\documentclass[10pt,a4paper]{report}

\usepackage{amsmath}
\usepackage{amssymb}
\usepackage{amsthm}
\usepackage{amsfonts}
\usepackage{microtype}
\usepackage{mathrsfs}
\usepackage{graphicx}
\usepackage{color}
\usepackage{latexsym}
\usepackage{rotating}
\usepackage{xspace}
\usepackage[all]{xy}
\usepackage{longtable}
\usepackage{leftidx}
\usepackage{mathtools}
\usepackage{titlesec}
\usepackage{setspace}
\usepackage{tikz}
\usetikzlibrary{calc}
\usetikzlibrary{arrows}
\usetikzlibrary{decorations.markings}
\usetikzlibrary{positioning}
\usepackage[top=1.2in, bottom=1.2in]{geometry}
\usepackage{geometry}
\usepackage[nottoc,notlot,notlof]{tocbibind}
\usepackage{makeidx}
\usepackage[final]{pdfpages}

\newtheorem{theorem}{Theorem}[section]
\numberwithin{equation}{theorem}
\newtheorem{lemma}[theorem]{Lemma}
\newtheorem{proposition}[theorem]{Proposition}
\newtheorem{corollary}[theorem]{Corollary}

\theoremstyle{definition}
\newtheorem{definition}[theorem]{Definition}
\newtheorem{example}[theorem]{Example}
\newtheorem{Notation}[theorem]{Notation}
\newtheorem{remark}[theorem]{Remark}
\newtheorem{construction}[theorem]{Construction}
\theoremstyle{conjecture}
\newtheorem{conjecture}[theorem]{Conjecture}
\newtheorem{question}[theorem]{Question}

\newcommand{\Ass}{\operatorname{Ass}}
\newcommand{\im}{\operatorname{im}}

\newcommand{\Min}{\operatorname{Min}}
\newcommand{\Assh}{\operatorname{Assh}}
\newcommand{\Spec}{\operatorname{Spec}}

\newcommand{\ara}{\operatorname{ara}}

\newcommand{\ann}{\operatorname{ann}}
\newcommand{\rank}{\operatorname{rank}}

\newcommand{\amp}{\operatorname{amp}}

\newcommand{\cd}{\operatorname{cd}}

\newcommand{\hd}{\operatorname{hd}}
\newcommand{\Ht}{\operatorname{ht}}
\newcommand{\id}{\operatorname{id}}
\newcommand{\fd}{\operatorname{fd}}
\newcommand{\Gid}{\operatorname{Gid}}
\newcommand{\pd}{\operatorname{pd}}
\newcommand{\Gpd}{\operatorname{Gpd}}
\newcommand{\V}{\operatorname{V}}

\newcommand{\Cone}{\operatorname{Cone}}

\newcommand{\Ext}{\operatorname{Ext}}
\newcommand{\Supp}{\operatorname{Supp}}

\newcommand{\Tor}{\operatorname{Tor}}
\newcommand{\Hom}{\operatorname{Hom}}
\newcommand{\Att}{\operatorname{Att}}

\newcommand{\depth}{\operatorname{depth}}
\newcommand{\width}{\operatorname{width}}
\newcommand{\Coass}{\operatorname{Coass}}
\newcommand{\coker}{\operatorname{coker}}

\newcommand{\Max}{\operatorname{Max}}

\newcommand{\fm}{\frak{m}}
\newcommand{\fp}{\frak{p}}

\newcommand{\fa}{\frak{a}}
\newcommand{\fb}{\frak{b}}

\newcommand{\suchthat}{\;\ifnum\currentgrouptype=16 \middle\fi|\;}

\newenvironment{prf}[1][Proof]{\begin{proof}[\bf #1]}{\end{proof}}

\makeatletter
\newcommand{\hocolim@}[2]{%
  \vtop{\m@th\ialign{##\cr
    \hfil$#1\operator@font holim$\hfil\cr
    \noalign{\nointerlineskip\kern1.5\ex@}#2\cr
    \noalign{\nointerlineskip\kern-\ex@}\cr}}%
}
\newcommand{\hocolim}{%
  \mathop{\mathpalette\hocolim@{\rightarrowfill@\textstyle}}\nmlimits@
}
\makeatother

\makeatletter
\newcommand{\holim@}[2]{%
  \vtop{\m@th\ialign{##\cr
    \hfil$#1\operator@font holim$\hfil\cr
    \noalign{\nointerlineskip\kern1.5\ex@}#2\cr
    \noalign{\nointerlineskip\kern-\ex@}\cr}}%
}
\newcommand{\holim}{%
  \mathop{\mathpalette\holim@{\leftarrowfill@\textstyle}}\nmlimits@
}
\makeatother

\makeindex

\begin{document}

\sloppy

\begin{titlepage}
\begin{center}

\begin{spacing}{3}
\vspace*{.01\textheight}\textbf{\Huge Gorenstein Homology and Finiteness Properties of Local (Co)homology}\\[1cm]
\end{spacing}

\begin{spacing}{2}
\textbf{\large by}\\
\textbf{\huge Hossein Faridian}\\[2cm]
\end{spacing}

\begin{spacing}{1.7}
\textnormal{\Large A Dissertation Submitted
in Partial Fulfillment of the Requirements for the
Degree of Doctor of Philosophy in Mathematics}\\[1.5cm]
\end{spacing}

\begin{spacing}{1.5}
\textnormal{\Large Department of Mathematics and Statistics \\ Shahid Beheshti University \\ Tehran, Iran \\}
\end{spacing}

\begin{spacing}{3.5}
\textnormal{\Large January 2018 \\}
\end{spacing}

\end{center}
\end{titlepage}

\pagenumbering{roman}

\begin{abstract}
This thesis is comprised of three chapters. The first chapter deals with bounded complexes of Gorenstein projective and Gorenstein injective modules. Deploying methods of relative homological algebra, we approximate such complexes with bounded complexes of projective and injective modules, respectively. As an application, we investigate the Gorenstein version of the New Intersection Theorem.

The second chapter studies the notion of cofiniteness for modules and complexes set forth by Hartshorne. Recruiting techniques of derived category, we study this notion thoroughly, obtain novel results, and extend some of the results to stable under specialization sets.

The third chapter delves into the Greenlees-May Duality Theorem which is widely thought of as a far-reaching generalization of the Grothendieck's Local Duality Theorem. This theorem is not addressed in the literature as it merits and its proof is indeed a tangled web in a series of scattered papers. By carefully scrutinizing the requisite tools, we present a clear-cut well-documented proof of this theorem.
\end{abstract}

\clearpage
\thispagestyle{plain}
\par\vspace*{.2\textheight}{\huge Mathematics, rightly viewed, possesses not only truth, but supreme beauty; a beauty cold and austere, like that of sculpture, without appeal to any part of our weaker nature, without the gorgeous trappings of painting or music, yet sublimely pure, and capable of a stern perfection such as only the greatest art can show.\par}

\par\vspace*{.02\textheight}{\Large \hfill Bertrand Russell\par}

\chapter*{Acknowledgements}
It is a pleasure to express my sincerest gratitude to my advisors professor Massoud Tousi and professor Kamran Divaani-Aazar, without their supportive guidance and cordial help this very dissertation would not have seen the light of day. I am also deeply grateful to professors Lars W. Christensen, Edgar E. Enochs, Paul C. Roberts, Charles A. Weibel, Joseph Lipman, John Greenlees, Peter Schenzel, and Ryo Takahashi from whom I learned a whole lot through email correspondence during my Ph.D. career.

\tableofcontents

\chapter*{Preface}
\addcontentsline{toc}{chapter}{Preface}
The main theme of this dissertation hinges upon the pivotal applications of homological algebra in commutative algebra. Homological algebra is by and large comprised of three major realms; namely, classical homological algebra, relative homological algebra, and the theory of derived category. Strictly speaking, classical homological algebra develops classical resolutions and the derived functors constructed upon them. Relative homological algebra deploys, in essence, the theory of covers and envelopes to proceed even further, so as to construct relative resolutions and relative derived functors. However, the theory of derived category is recognized as the ultimate formulation of homological algebra whose demanding language lays on a formidable machinery to carry out homological methods in an even more supple fashion and instates the warp and woof of a convenient context within its expanse many overwhelming spectral sequences turn into straightforward isomorphisms.

In the first chapter of this work, we focus on a special relative theory, the so-called Gorenstein homology. The main ingredients of this relative theory are Gorenstein projective, Gorenstein injective, and Gorenstein flat modules which were defined by M. Auslander, M. Bridger, E. E. Enochs, and O. Jenda in several stages to take on the current full-fledged formalism. Based on these classes of modules, we are privileged to have Gorenstein homological dimensions and Gorenstein derived functors. However, we are solely interested in considering a bounded complex of Gorenstein projective modules and effectively approximate it with a bounded complex of projective modules of the same length, in a way that their homology modules are isomorphic as far as it is possible. We further establish analogous results for bounded complexes of Gorenstein injective modules. As an application, we scrupulously investigate the conjectural Gorenstein version of the New Intersection Theorem and settle some special cases of this conjecture.

In the second chapter, we crucially study the notion of cofiniteness of modules and complexes; a notion that dates back to Grothendieck's epoch-making seminars on local cohomology delivered at Harvard university in 1961-2. We address Hartshorne's questions on cofiniteness, collect all the results thus far procured in the literature with fully worked out and demystified proofs in a coherent and cohesive manner, and establish plenty of new results in this direction. In the course of our scrutiny, we deal efficaciously with both module and complex cases; in the latter, we are irrevocably driven towards the recruitment of the state-of-the-art techniques of derived category. At the culminating apex of our work, we curiously inspect the untoward correlation between Hartshorne's questions that divulges an intricate link-up which in turn revealingly spotlights on the raison d'\^{e}tre of these questions. We further strive to extend these results to local cohomology modules with respect to stable under specialization sets; a notion that could be vehemently visualized as an all-in-one generalization of all the generalized local cohomology theories introduced and delved by various authors.

In the third chapter, we are interested in the celebrated Greenlees-May Duality which is deemed as a far-reaching generalization of the Grothendieck's Local Duality. Despite its undeniable impact on the theory of derived local homology and cohomology, we regretfully notice that there is no comprehensive and accessible treatment of this theorem in the literature. There is a handful of papers that touch on the subject, each from a different perspective, but none of them affords the presentation of a clear-cut and thoroughgoing proof that is fairly readable by non-experts. To remedy this defect and bridge the gap in the literature, we punctiliously commence on probing this theorem by providing the prerequisites from scratch and build upon a well-documented rigorous proof which is basically tendered in layman's terms. In the course of our proof, some arguments are familiar while some others are novel. We finally elaborate on the highly non-trivial fact that the Greenlees-May Duality generalizes the Local Duality in simple and traceable steps.

\chapter{Approximation of Complexes in Gorenstein Homology}

\pagenumbering{arabic}

\section{Introduction}

Throughout this chapter, all rings are assumed to be commutative with identity.

Let $R$ be a ring. Given a bounded $R$-complex
$$X: 0\rightarrow X_{s} \xrightarrow {\partial^{X}_{s}} X_{s-1}\rightarrow \cdots \rightarrow X_{1} \xrightarrow {\partial^{X}_{1}} X_{0} \rightarrow 0$$
of Gorenstein projective modules, we would like to construct a bounded $R$-complex
$$F: 0\rightarrow F_{s} \xrightarrow {\partial^{F}_{s}} F_{s-1}\rightarrow \cdots \rightarrow F_{1} \xrightarrow {\partial^{F}_{1}} F_{0} \rightarrow 0$$
of free modules, whose homology modules are as close as possible to the homology modules of $X$. The best-case scenario would be $H_{i}(X) \cong H_{i}(F)$ for every $0 \leq i \leq s$. Unfortunately, this cannot occur in general. It turns out that the best possible result is the following which is proved in Section 1.2.

\begin{theorem} \label{1.1.1}
Let
$$X: 0 \rightarrow X_{s} \rightarrow X_{s-1} \rightarrow \cdots \rightarrow X_{1} \rightarrow X_{0} \rightarrow 0$$
be an $R$-complex of Gorenstein projective modules. Then there exists an $R$-complex
$$F: 0 \rightarrow F_{s} \rightarrow F_{s-1} \rightarrow \cdots \rightarrow F_{1} \rightarrow F_{0} \rightarrow 0,$$
of free modules that enjoys the following properties:
\begin{enumerate}
\item[(i)] $H_{i}(X)\cong H_{i}(F)$ for every $1\leq i \leq s$.
\item[(ii)] There is a short exact sequence of $R$-modules
$$0\rightarrow H_{0}(X) \rightarrow H_{0}(F) \rightarrow C \rightarrow 0,$$
where $C$ is Gorenstein projective.
\item[(iii)] If in particular, $(R,\mathfrak{m})$ is noetherian local and the homology modules of $X$ have all finite length, then the $R$-complex $F$ can be chosen in a way that the modules in the short exact sequence of (ii) are locally free on the punctured spectrum of $R$ with constant rank.
\end{enumerate}
\end{theorem}

As an immediate corollary to the above result, we recover the following well-known fact; see \cite[Lemma 2.17]{CFrH}.

\begin{corollary} \label{1.1.2}
Let $M$ be an $R$-module with $\Gpd_{R}(M)<\infty$. Then there exists a short exact sequence of $R$-modules
$$0 \rightarrow M \rightarrow L \rightarrow C \rightarrow 0,$$
where $\pd_{R}(L)= \Gpd_{R}(M)$, and $C$ is Gorenstein projective.
\end{corollary}

In a similar fashion, the dual result is proved in Section 1.3.

\begin{theorem} \label{1.1.3}
Let
$$Y: 0 \rightarrow Y_{s} \rightarrow Y_{s-1} \rightarrow \cdots \rightarrow Y_{1} \rightarrow Y_{0} \rightarrow 0$$
be an $R$-complex of Gorenstein injective modules. Then there exists an $R$-complex
$$I: 0 \rightarrow I_{s} \rightarrow I_{s-1} \rightarrow \cdots \rightarrow I_{1} \rightarrow I_{0} \rightarrow 0,$$
of injective modules that enjoys the following properties:
\begin{enumerate}
\item[(i)] $H_{i}(Y)\cong H_{i}(I)$ for every $0\leq i \leq s-1$.
\item[(ii)] There is a short exact sequence of $R$-modules
$$0\rightarrow K \rightarrow H_{0}(I) \rightarrow H_{0}(Y) \rightarrow 0,$$
where $K$ is Gorenstein injective.
\item[(iii)] If in particular, $(R,\mathfrak{m})$ is noetherian local ring which is a homomorphic image a Gorenstein local ring, and the homology modules of $Y$ have all finite length, then the $R$-complex $I$ can be chosen in a way that the modules in the short exact sequence of (ii) are locally injective on the punctured spectrum of $R$.
\end{enumerate}
\end{theorem}

As an immediate corollary to the above result, we recover the following well-known fact; see \cite[Lemma 2.18]{CFrH}.

\begin{corollary} \label{1.1.4}
Let $M$ be an $R$-module with $\Gid_{R}(M)<\infty$. Then there exists a short exact sequence of $R$-modules
$$0 \rightarrow K \rightarrow T \rightarrow M \rightarrow 0,$$
where $\id_{R}(T)= \Gid_{R}(M)$, and $K$ is Gorenstein injective.
\end{corollary}

In section 1.4, we deploy these results to investigate some special cases of the Gorenstein version of the New Intersection Theorem.

\section{Bounded Complexes of Gorenstein Projective Modules}

In this section, we develop some constructions to study bounded complexes of Gorenstein projective modules. As a result, we obtain an approximation for complexes of Gorenstein projective modules by complexes of free modules. We first collect some definitions and facts for the convenience of the reader.

We begin with the definition of Gorenstein projective modules and Gorenstein projective dimension.

\begin{definition} \label{1.2.1}
Let $\mathcal{A}$ be a class of $R$-modules. An $R$-complex
$$X: \cdots \rightarrow X_{i+1} \xrightarrow{\partial^{X}_{i+1}} X_{i} \xrightarrow{\partial^{X}_{i}} X_{i-1} \rightarrow \cdots$$
is said to be \textit{$\Hom_{R}(-,\mathcal{A})$-exact} if the $R$-complex
$$\cdots \rightarrow \Hom_{R}(X_{i-1},A) \xrightarrow{\Hom_{R}\left(\partial^{X}_{i},A\right)} \Hom_{R}(X_{i},A) \xrightarrow{\Hom_{R}\left(\partial^{X}_{i+1},A\right)} \Hom_{R}(X_{i+1},A) \rightarrow \cdots$$
is exact for every $A\in \mathcal{A}$.
\end{definition}

In the sequel, $\mathcal{P}$ will denote the class of projective $R$-modules.

\begin{definition} \label{1.2.2}
We have the following definitions:
\begin{enumerate}
\item[(i)] An $R$-module $Q$ is called \textit{Gorenstein projective} \index{Gorenstein projective module} if there exists a $\Hom_{R}(-,\mathcal{P})$-exact exact $R$-complex
$$P: \cdots \rightarrow P_{2} \xrightarrow {\partial^{P}_{2}} P_{1} \xrightarrow {\partial^{P}_{1}} P_{0} \xrightarrow {\partial^{P}_{0}} P_{-1} \xrightarrow {\partial^{P}_{-1}} P_{-2} \rightarrow \cdots,$$
consisting of projective $R$-modules, such that $Q \cong \im \partial^{P}_{1}$.
\item[(ii)] A \textit{Gorenstein projective resolution} \index{Gorenstein projective resolution} of an $R$-module $M$ is defined to be an $R$-complex
$$\cdots \rightarrow Q_{2} \rightarrow Q_{1} \rightarrow Q_{0} \rightarrow 0,$$
where $Q_{i}$ is a Gorenstein projective $R$-module for every $i \geq 0$, and there is an $R$-homomorphism $\varepsilon: Q_{0}\rightarrow M$ such that the augmented $R$-complex
$$\cdots \rightarrow Q_{2} \rightarrow Q_{1} \rightarrow Q_{0} \xrightarrow{\varepsilon} M \rightarrow 0$$
is exact.
\item[(iii)] The \textit{Gorenstein projective dimension} \index{Gorenstein projective dimension} of a nonzero $R$-module $M$ is defined to be
$$\Gpd_{R}(M):= \inf \left\{n \geq 0 \suchthat \begin{tabular}{ccc}
  \text{\textit{There is a Gorenstein projective resolution}}\\
 $0 \rightarrow Q_{n} \rightarrow \cdots \rightarrow Q_{1} \rightarrow Q_{0} \rightarrow$ 0 \text{\textit{of}} $M$\\
  \end{tabular} \right\},$$
with the convention that $\inf \emptyset := \infty$. Further, we define $\Gpd_{R}(0):=-\infty$.
\end{enumerate}
\end{definition}

The following facts on Gorenstein projective modules and Gorenstein projective dimension are required in the rest of the work.

\begin{proposition} \label{1.2.3}
Let
$$0\rightarrow M^{\prime}\rightarrow M\rightarrow M^{\prime\prime}\rightarrow 0$$
be a short exact sequence of $R$-modules with $M^{\prime\prime}$ Gorenstein projective. Then $M^{\prime}$ is Gorenstein projective if and only if $M$ is Gorenstein projective.
\end{proposition}

\begin{prf}
See \cite[Theorem 2.5]{Ho}.
\end{prf}

\begin{proposition} \label{1.2.4}
Let $R$ be a noetherian ring with finite Krull dimension. If $Q$ is a Gorenstein projective $R$-module, then $Q_{\mathfrak{p}}$ is a Gorenstein projective $R_{\mathfrak{p}}$-module for every $\mathfrak{p}\in \Spec(R)$.
\end{proposition}

\begin{prf}
See \cite[Proposition 2.17]{CFH1}.
\end{prf}

\begin{proposition} \label{1.2.5}
Let $M$ be an $R$-module with $\Gpd_{R}(M)<\infty$. Then we have
$$\Gpd_{R}(M)= \sup \left\{i \geq 0 \suchthat \Ext_{R}^{i}(M,P)\neq 0 \text{ for some } P \in \mathcal{P}\right\}.$$
\end{proposition}

\begin{prf}
See \cite[Theorem 2.20]{Ho}.
\end{prf}

\begin{proposition} \label{1.2.6}
Let $M$ be an $R$-module. Then $\Gpd_{R}(M)\leq \pd_{R}(M)$ with equality if $\pd_{R}(M)<\infty$.
\end{proposition}

\begin{prf}
See \cite[Proposition 2.12]{CFH1}.
\end{prf}

We further need the following definitions.

\begin{definition} \label{1.2.7}
Let $\mathcal{A}$ be a class of $R$-modules. Define the \textit{left orthogonal class} \index{left orthogonal class} of $\mathcal{A}$, denoted by $\leftidx{^\perp}{\mathcal{A}}$, to be the class of all $R$-modules $M$ with $\Ext^{1}_{R}(M,A)=0$ for every $A\in \mathcal{A}$.
\end{definition}

\begin{definition} \label{1.2.8}
Let $\mathcal{A}$ be a class of $R$-modules, and $M$ an $R$-module. Then we have the following definitions:
\begin{enumerate}
\item[(i)] By an \textit{$\mathcal{A}$-preenvelope} \index{preenvelope} of $M$, we mean an $R$-homomorphism $f:M\rightarrow A$ where $A\in \mathcal{A}$, with the property that the $R$-homomorphism
$$\Hom_{R}(f,B):\Hom_{R}(A,B)\rightarrow \Hom_{R}(M,B)$$
is surjective for every $B\in \mathcal{A}$.
\item[(ii)] An $\mathcal{A}$-preenvelope $f:M\rightarrow A$ is said to be \textit{special} \index{special preenvelope} whenever $f$ is injective and $\coker f \in \leftidx{^\perp}{\mathcal{A}}$.
\end{enumerate}
\end{definition}

The next proposition is the main ingredient of Construction \ref{1.2.12}.

\begin{proposition} \label{1.2.9}
Every Gorenstein projective $R$-module admits a special free preenvelope.
\end{proposition}

\begin{prf}
Let $M$ be a Gorenstein projective $R$-module. By \cite[Proposition 2.4]{Ho}, there exists a $\Hom_{R}(-,\mathcal{P})$-exact exact $R$-complex
\begin{equation} \label{eq:1.2.9.1}
F: \cdots \rightarrow F_{2} \xrightarrow {\partial^{F}_{2}} F_{1} \xrightarrow {\partial^{F}_{1}} F_{0} \xrightarrow {\partial^{F}_{0}} F_{-1} \xrightarrow {\partial^{F}_{-1}} F_{-2} \rightarrow \cdots,
\end{equation}
consisting of free $R$-modules, such that $M \cong \im \partial^{F}_{1}$. Let $\theta:M \rightarrow \im \partial^{F}_{1}$ be an isomorphism, and let $\lambda:= \iota\theta$, where $\iota:\im \partial^{F}_{1}\rightarrow F_{0}$ is the inclusion map. Therefore, there is a short exact sequence
$$0\rightarrow M \xrightarrow {\lambda} F_{0} \rightarrow C\rightarrow 0,$$
where $$C=\coker \lambda=\coker \partial^{F}_{1} \cong\im \partial^{F}_{0}$$
is Gorenstein projective by the symmetry in Definition \ref{1.2.2} (i).
We show that the $R$-homomorphism $\lambda:M \rightarrow F_{0}$ is a special free preenvelope of $M$.
Given any free $R$-module $E$, we need to argue that the $R$-homomorphism
$$\Hom_{R}(\lambda,E):\Hom_{R}(F_{0},E)\rightarrow \Hom_{R}(M,E)$$
is surjective. There is an exact sequence
$$\Hom_{R}(F_{0},E)\xrightarrow{\Hom_{R}(\lambda,E)} \Hom_{R}(M,E) \rightarrow \Ext^{1}_{R}(C,E)=0,$$
where the vanishing is due to the fact that $C$ is Gorenstein projective and $E$ is free. This establishes the claim.
\end{prf}

We further need the pushout construction.

\begin{definition} \label{1.2.10}
Let $f:X\rightarrow M$ and $g:X\rightarrow N$ be two $R$-homomorphisms. Define the \textit{pushout} \index{pushout} of the pair $(f,g)$ to be the $R$-module
$$M\sqcup_{X}N :=  \frac{M \oplus N}{\Im},$$
where
\begin{align*}
\Im:= \left\{\left(-f(x),g(x)\right)\suchthat x\in X\right\}.
\end{align*}
Further, complete the pair $(f,g)$ to the commutative diagram
\[\begin{tikzpicture}[every node/.style={midway}]
  \matrix[column sep={3em}, row sep={3em}]
  {\node(X) {$X$}; & \node(N) {$N$}; \\
  \node(M) {$M$}; & \node (T) {$M\sqcup_{X}N$};\\};
  \draw[decoration={markings,mark=at position 1 with {\arrow[scale=1.5]{>}}},postaction={decorate},shorten >=0.5pt] (X) -- (M) node[anchor=west] {$f$};
  \draw[decoration={markings,mark=at position 1 with {\arrow[scale=1.5]{>}}},postaction={decorate},shorten >=0.5pt] (X) -- (N) node[anchor=south] {$g$};
  \draw[decoration={markings,mark=at position 1 with {\arrow[scale=1.5]{>}}},postaction={decorate},shorten >=0.5pt] (N) -- (T) node[anchor=west] {$f^{\prime}$};
  \draw[decoration={markings,mark=at position 1 with {\arrow[scale=1.5]{>}}},postaction={decorate},shorten >=0.5pt] (M) -- (T) node[anchor=south] {$g^{\prime}$};
\end{tikzpicture}\]
where $f^{\prime}(n):= (0,n)+\Im$ and $g^{\prime}(m):= (m,0)+\Im$.
\end{definition}

The following proposition describes the most important property of the pushout construction.

\begin{proposition} \label{1.2.11}
Let $f:X\rightarrow M$ and $g:X\rightarrow N$ be two $R$-homomorphisms. Then there exists a commutative diagram with exact rows and columns as follows:
\[\begin{tikzpicture}[every node/.style={midway}]
  \matrix[column sep={3em}, row sep={3em}]
  {\node(a) {$X$}; & \node(b) {$N$}; & \node(c) {$\coker g$}; & \node(d) {$0$};\\
  \node(e) {$M$}; & \node (f) {$M\sqcup_{X}N$}; & \node(g) {$\coker g^{\prime}$}; & \node(h) {$0$};\\
  \node(i) {$\coker f$}; & \node (j) {$\coker f^{\prime}$}; & \node (k) {}; & \node (l) {};\\
  \node(m) {$0$}; & \node (n) {$0$}; & \node (o) {}; & \node (p) {};\\};
  \draw[decoration={markings,mark=at position 1 with {\arrow[scale=1.5]{>}}},postaction={decorate},shorten >=0.5pt] (a) -- (b) node[anchor=south] {$g$};
  \draw[decoration={markings,mark=at position 1 with {\arrow[scale=1.5]{>}}},postaction={decorate},shorten >=0.5pt] (b) -- (c) node[anchor=south] {};
  \draw[decoration={markings,mark=at position 1 with {\arrow[scale=1.5]{>}}},postaction={decorate},shorten >=0.5pt] (c) -- (d) node[anchor=south] {};
  \draw[decoration={markings,mark=at position 1 with {\arrow[scale=1.5]{>}}},postaction={decorate},shorten >=0.5pt] (e) -- (f) node[anchor=south] {$g^{\prime}$};
  \draw[decoration={markings,mark=at position 1 with {\arrow[scale=1.5]{>}}},postaction={decorate},shorten >=0.5pt] (f) -- (g) node[anchor=south] {};
  \draw[decoration={markings,mark=at position 1 with {\arrow[scale=1.5]{>}}},postaction={decorate},shorten >=0.5pt] (g) -- (h) node[anchor=south] {};
  \draw[decoration={markings,mark=at position 1 with {\arrow[scale=1.5]{>}}},postaction={decorate},shorten >=0.5pt] (i) -- (j) node[anchor=south] {$\widetilde{g^{\prime}}$};
  \draw[decoration={markings,mark=at position 1 with {\arrow[scale=1.5]{>}}},postaction={decorate},shorten >=0.5pt] (a) -- (e) node[anchor=west] {$f$};
  \draw[decoration={markings,mark=at position 1 with {\arrow[scale=1.5]{>}}},postaction={decorate},shorten >=0.5pt] (e) -- (i) node[anchor=east] {};
  \draw[decoration={markings,mark=at position 1 with {\arrow[scale=1.5]{>}}},postaction={decorate},shorten >=0.5pt] (i) -- (m) node[anchor=east] {};
  \draw[decoration={markings,mark=at position 1 with {\arrow[scale=1.5]{>}}},postaction={decorate},shorten >=0.5pt] (b) -- (f) node[anchor=west] {$f^{\prime}$};
  \draw[decoration={markings,mark=at position 1 with {\arrow[scale=1.5]{>}}},postaction={decorate},shorten >=0.5pt] (f) -- (j) node[anchor=west] {};
  \draw[decoration={markings,mark=at position 1 with {\arrow[scale=1.5]{>}}},postaction={decorate},shorten >=0.5pt] (j) -- (n) node[anchor=west] {};
  \draw[decoration={markings,mark=at position 1 with {\arrow[scale=1.5]{>}}},postaction={decorate},shorten >=0.5pt] (c) -- (g) node[anchor=west] {$\widetilde{f^{\prime}}$};
\end{tikzpicture}\]
where $\widetilde{f^{\prime}}$ and $\widetilde{g^{\prime}}$ are induced by $f^{\prime}$ and $g^{\prime}$, respectively. Furthermore, the following assertions hold:
\begin{enumerate}
\item[(i)] If $f$ is injective, then so is $f^{\prime}$.
\item[(ii)] If $g$ is injective, then so is $g^{\prime}$.
\item[(iii)] Both $\widetilde{f^{\prime}}$ and $\widetilde{g^{\prime}}$ are isomorphisms.
\end{enumerate}
\end{proposition}

\begin{prf}
(i): Suppose that $f^{\prime}(n)= (0,n)+ \Im =0$, so that $(0,n) \in \Im$. This means that $(0,n) = \left(-f(x),g(x)\right)$ for some $x \in X$. It follows that $f(x)=0$, and thus $x=0$ as $f$ is injective. Therefore, $n=0$. Hence $f^{\prime}$ is injective.

(ii): Similar to (i).

(iii): Suppose that $\widetilde{f^{\prime}}(n + \im g)= f^{\prime}(n) + \im g^{\prime}=0$, so that $f^{\prime}(n) = (0,n)+ \Im \in \im g^{\prime}$. It follows that $(0,n)+ \Im = g^{\prime}(m) = (m,0)+ \Im$ for some $m \in M$. Hence $(-m,n) \in \Im$, so $(-m,n) = \left(-f(x),g(x)\right)$ for some $x \in X$. In particular, $n=g(x)$, so $n + \im g =0$, implying that $\widetilde{f^{\prime}}$ is injective. On the other hand, every element of $\coker g^{\prime}$ is of the form $\left((m,n) + \Im\right) + \im g^{\prime}$, and we have
\begin{equation*}
\begin{split}
\left((m,n) + \Im\right) + \im g^{\prime} & = \left((m,0) + \Im\right) + \left((0,n) + \Im\right) + \im g^{\prime} \\
 & = g^{\prime}(m) + f^{\prime}(n) + \im g^{\prime} \\
 & = f^{\prime}(n) + \im g^{\prime}.
\end{split}
\end{equation*}
Now it is clear that
$$\widetilde{f^{\prime}}(n + \im g)= f^{\prime}(n) + \im g^{\prime} = \left((m,n) + \Im\right) + \im g^{\prime},$$
so $\widetilde{f^{\prime}}$ is surjective. It follows that $\widetilde{f^{\prime}}$ is an isomorphism. Similarly, one can see that $\widetilde{g^{\prime}}$ is an isomorphism.
\end{prf}

We are now ready to present the following construction.

\begin{construction} \label{1.2.12}
Let
$$X: 0\rightarrow X_{s} \xrightarrow {\partial^{X}_{s}} X_{s-1}\rightarrow \cdots \rightarrow X_{1} \xrightarrow {\partial^{X}_{1}} X_{0} \rightarrow 0,$$
be a bounded $R$-complex of Gorenstein projective modules. The proof of Proposition \ref{1.2.9} gives a short exact sequence
$$0\rightarrow X_{s} \xrightarrow {\lambda_{s}} F_{s} \rightarrow \coker \lambda_{s} \rightarrow 0,$$
where $F_{s}$ is a free $R$-module and $\coker \lambda_{s}$ is a Gorenstein projective $R$-module. Construct the following pushout diagram as in Proposition \ref{1.2.11}:
\[\begin{tikzpicture}[every node/.style={midway}]
  \matrix[column sep={3em}, row sep={3em}]
  {\node(a) {$0$}; & \node(b) {$0$}; & \node(c) {}; & \node(d) {};\\
  \node(e) {$X_{s}$}; & \node(f) {$X_{s-1}$}; & \node(g) {$\coker {\partial^{X}_{s}}$}; & \node(h) {$0$};\\
  \node(i) {$F_{s}$}; & \node (j) {$L_{s-1}$}; & \node(k) {$\coker \alpha_{s}$}; & \node(l) {$0$};\\
  \node(m) {$\coker \lambda_{s}$}; & \node (n) {$\coker \beta_{s-1}$}; & \node(o) {}; & \node(p) {};\\
  \node(q) {$0$}; & \node (r) {$0$}; & \node(s) {}; & \node(t) {};\\};
  \draw[decoration={markings,mark=at position 1 with {\arrow[scale=1.5]{>}}},postaction={decorate},shorten >=0.5pt] (f) -- (g) node[anchor=south] {};
  \draw[decoration={markings,mark=at position 1 with {\arrow[scale=1.5]{>}}},postaction={decorate},shorten >=0.5pt] (g) -- (h) node[anchor=south] {};
  \draw[decoration={markings,mark=at position 1 with {\arrow[scale=1.5]{>}}},postaction={decorate},shorten >=0.5pt] (e) -- (f) node[anchor=south] {${\partial^{X}_{s}}$};
  \draw[decoration={markings,mark=at position 1 with {\arrow[scale=1.5]{>}}},postaction={decorate},shorten >=0.5pt] (k) -- (l) node[anchor=south] {};
  \draw[decoration={markings,mark=at position 1 with {\arrow[scale=1.5]{>}}},postaction={decorate},shorten >=0.5pt] (i) -- (j) node[anchor=south] {$\alpha_{s}$};
  \draw[decoration={markings,mark=at position 1 with {\arrow[scale=1.5]{>}}},postaction={decorate},shorten >=0.5pt] (j) -- (k) node[anchor=south] {$\pi_{s-1}$};
  \draw[decoration={markings,mark=at position 1 with {\arrow[scale=1.5]{>}}},postaction={decorate},shorten >=0.5pt] (g) -- (h) node[anchor=south] {};
  \draw[decoration={markings,mark=at position 1 with {\arrow[scale=1.5]{>}}},postaction={decorate},shorten >=0.5pt] (m) -- (n) node[anchor=south] {$\widetilde{\alpha_{s}}$};
  \draw[decoration={markings,mark=at position 1 with {\arrow[scale=1.5]{>}}},postaction={decorate},shorten >=0.5pt] (e) -- (i) node[anchor=west] {$\lambda_{s}$};
  \draw[decoration={markings,mark=at position 1 with {\arrow[scale=1.5]{>}}},postaction={decorate},shorten >=0.5pt] (a) -- (e) node[anchor=east] {};
  \draw[decoration={markings,mark=at position 1 with {\arrow[scale=1.5]{>}}},postaction={decorate},shorten >=0.5pt] (i) -- (m) node[anchor=east] {};
  \draw[decoration={markings,mark=at position 1 with {\arrow[scale=1.5]{>}}},postaction={decorate},shorten >=0.5pt] (m) -- (q) node[anchor=east] {};
  \draw[decoration={markings,mark=at position 1 with {\arrow[scale=1.5]{>}}},postaction={decorate},shorten >=0.5pt] (f) -- (j) node[anchor=west] {$\beta_{s-1}$};
  \draw[decoration={markings,mark=at position 1 with {\arrow[scale=1.5]{>}}},postaction={decorate},shorten >=0.5pt] (b) -- (f) node[anchor=east] {};
  \draw[decoration={markings,mark=at position 1 with {\arrow[scale=1.5]{>}}},postaction={decorate},shorten >=0.5pt] (j) -- (n) node[anchor=east] {};
  \draw[decoration={markings,mark=at position 1 with {\arrow[scale=1.5]{>}}},postaction={decorate},shorten >=0.5pt] (n) -- (r) node[anchor=east] {};
  \draw[decoration={markings,mark=at position 1 with {\arrow[scale=1.5]{>}}},postaction={decorate},shorten >=0.5pt] (g) -- (k) node[anchor=west] {$\widetilde{\beta_{s-1}}$};
\end{tikzpicture}\]
where $L_{s-1}:= F_{s} \sqcup_{X_{s}} X_{s-1}$. Note that the injectivity of $\beta_{s-1}$ follows from that of $\lambda_{s}$ in view of Proposition \ref{1.2.11}. Since $\widetilde{\alpha_{s}}$ is an isomorphism, $\coker \beta_{s-1}$ is Gorenstein projective, and Proposition \ref{1.2.3} implies that $L_{s-1}$ is Gorenstein projective.
Let $\varphi_{s-1}$ be the composition of the $R$-homomorphisms
$$L_{s-1} \xrightarrow {\pi_{s-1}} \coker \alpha_{s} \xrightarrow {\widetilde{\beta_{s-1}}^{-1}} \coker {\partial^{X}_{s}} \xrightarrow {{\overline{\partial^{X}_{s-1}}}} X_{s-2},$$
where $\overline{\partial^{X}_{s-1}}$ is induced by $\partial^{X}_{s-1}$ in the obvious way.
One can easily check that the following sequence is an $R$-complex:
$$0\rightarrow F_{s} \xrightarrow {\alpha_{s}} L_{s-1} \xrightarrow {\varphi_{s-1}} X_{s-2} \xrightarrow {\partial^{X}_{s-2}} X_{s-3} \rightarrow \cdots \rightarrow X_{1} \xrightarrow {\partial^{X}_{1}} X_{0} \rightarrow 0$$
Continuing this construction in the same fashion, we get the following commutative diagram:
\[\begin{tikzpicture}[every node/.style={midway}]
  \matrix[column sep={3em}, row sep={3em}]
  {\node(1) {$X_{s}$}; & \node(2) {$X_{s-1}$}; & \node(3) {}; & \node(4) {}; & \node(5) {};\\
  \node(6) {$F_{s}$}; & \node(7) {$L_{s-1}$}; & \node(8) {$X_{s-2}$}; & \node(9) {}; & \node(10) {};\\
  \node(11) {}; & \node(12) {$F_{s-1}$}; & \node(13) {$L_{s-2}$}; & \node(14) {}; & \node(15) {};\\
  \node(16) {}; & \node(17) {}; & \node(18) {}; & \node(19) {$L_{1}$}; & \node(20) {$X_{0}$};\\
  \node(21) {}; & \node(22) {}; & \node(23) {}; & \node(24) {$F_{1}$}; & \node(25) {$L_{0}$};\\
  \node(26) {}; & \node(27) {}; & \node(28) {}; & \node(29) {}; & \node(30) {$F_{0}$};\\};
  \draw[decoration={markings,mark=at position 1 with {\arrow[scale=1.5]{>}}},postaction={decorate},shorten >=0.5pt] (1) -- (6) node[anchor=west] {$\lambda_{s}$};
  \draw[decoration={markings,mark=at position 1 with {\arrow[scale=1.5]{>}}},postaction={decorate},shorten >=0.5pt] (2) -- (7) node[anchor=west] {$\beta_{s-1}$};
  \draw[decoration={markings,mark=at position 1 with {\arrow[scale=1.5]{>}}},postaction={decorate},shorten >=0.5pt] (7) -- (12) node[anchor=west] {$\lambda_{s-1}$};
  \draw[decoration={markings,mark=at position 1 with {\arrow[scale=1.5]{>}}},postaction={decorate},shorten >=0.5pt] (8) -- (13) node[anchor=west] {$\beta_{s-2}$};
  \draw[decoration={markings,mark=at position 1 with {\arrow[scale=1.5]{>}}},postaction={decorate},shorten >=0.5pt] (19) -- (24) node[anchor=west] {$\lambda_{1}$};
  \draw[decoration={markings,mark=at position 1 with {\arrow[scale=1.5]{>}}},postaction={decorate},shorten >=0.5pt] (20) -- (25) node[anchor=west] {$\beta_{0}$};
  \draw[decoration={markings,mark=at position 1 with {\arrow[scale=1.5]{>}}},postaction={decorate},shorten >=0.5pt] (25) -- (30) node[anchor=west] {$\lambda_{0}$};
  \draw[decoration={markings,mark=at position 1 with {\arrow[scale=1.5]{>}}},postaction={decorate},shorten >=0.5pt] (1) -- (2) node[anchor=south] {$\partial^{X}_{s}$};
  \draw[decoration={markings,mark=at position 1 with {\arrow[scale=1.5]{>}}},postaction={decorate},shorten >=0.5pt] (6) -- (7) node[anchor=south] {$\alpha_{s}$};
  \draw[decoration={markings,mark=at position 1 with {\arrow[scale=1.5]{>}}},postaction={decorate},shorten >=0.5pt] (7) -- (8) node[anchor=south] {$\varphi_{s-1}$};
  \draw[decoration={markings,mark=at position 1 with {\arrow[scale=1.5]{>}}},postaction={decorate},shorten >=0.5pt] (12) -- (13) node[anchor=south] {$\alpha_{s-1}$};
  \draw[decoration={markings,mark=at position 1 with {\arrow[scale=1.5]{>}}},postaction={decorate},shorten >=0.5pt] (19) -- (20) node[anchor=south] {$\varphi_{1}$};
  \draw[decoration={markings,mark=at position 1 with {\arrow[scale=1.5]{>}}},postaction={decorate},shorten >=0.5pt] (24) -- (25) node[anchor=south] {$\alpha_{1}$};
  \path (13) -- (19) node [midway, sloped] {$\dots$};
\end{tikzpicture}\]
where $L_{i}:= F_{i+1} \sqcup_{L_{i+1}} X_{i}$ for every $0 \leq i \leq s-2$, and the homomorphisms $\lambda_{i}$, $\beta_{i}$, $\alpha_{i}$, and $\varphi_{i}$ are defined analogously. One can see by inspection that $\partial^{X}_{i}=\varphi_{i}\beta_{i}$ for every $1 \leq i \leq s-1$. Moreover, we let $\partial^{F}_{i}:=\lambda_{i-1}\alpha_{i}$ for every $1 \leq i \leq s$. We have
$$\partial^{F}_{i-1}\partial^{F}_{i}=\lambda_{i-2}\alpha_{i-1}\lambda_{i-1}\alpha_{i}=\lambda_{i-2}\beta_{i-2}\varphi_{i-1}\alpha_{i}=0,$$
since $\varphi_{i-1}\alpha_{i}=0$.
Hence $F$ is a bounded $R$-complex of free modules. Besides, we let $\sigma_{s}:=\lambda_{s}$, and $\sigma_{i}:=\lambda_{i}\beta_{i}$ for every $0 \leq i \leq s-1$. Then
$$\sigma_{s-1}\partial^{X}_{s}=\lambda_{s-1}\beta_{s-1}\partial^{X}_{s}=
\lambda_{s-1}\alpha_{s}\lambda_{s}=\partial^{F}_{s}\sigma_{s}$$
and
$$\sigma_{i-1}\partial^{X}_{i}=\lambda_{i-1}\beta_{i-1}\varphi_{i}\beta_{i}=
\lambda_{i-1}\alpha_{i}\lambda_{i}\beta_{i}=\partial^{F}_{i}\sigma_{i}$$
for every $1 \leq i \leq s-1$. As observed before, all the vertical homomorphisms are injective. Thus $\sigma=(\sigma_{i})_{0 \leq i \leq s}: X \rightarrow F$ is an injective morphism of $R$-complexes as follows:
\[\begin{tikzpicture}[every node/.style={midway}]
  \matrix[column sep={3em}, row sep={3em}]
  {\node(1) {$0$}; & \node(2) {$X_{s}$}; & \node(3) {$X_{s-1}$}; & \node(4) {$\cdots$}; & \node(5) {$X_{1}$}; & \node(6) {$X_{0}$}; & \node(7) {$0$};\\
  \node(8) {$0$}; & \node(9) {$F_{s}$}; & \node(10) {$F_{s-1}$}; & \node(11) {$\cdots$}; & \node(12) {$F_{1}$}; & \node(13) {$F_{0}$}; & \node(14) {$0$};\\};
  \draw[decoration={markings,mark=at position 1 with {\arrow[scale=1.5]{>}}},postaction={decorate},shorten >=0.5pt] (2) -- (9) node[anchor=west] {$\sigma_{s}$};
  \draw[decoration={markings,mark=at position 1 with {\arrow[scale=1.5]{>}}},postaction={decorate},shorten >=0.5pt] (3) -- (10) node[anchor=west] {$\sigma_{s-1}$};
  \draw[decoration={markings,mark=at position 1 with {\arrow[scale=1.5]{>}}},postaction={decorate},shorten >=0.5pt] (5) -- (12) node[anchor=west] {$\sigma_{1}$};
  \draw[decoration={markings,mark=at position 1 with {\arrow[scale=1.5]{>}}},postaction={decorate},shorten >=0.5pt] (6) -- (13) node[anchor=west] {$\sigma_{0}$};
  \draw[decoration={markings,mark=at position 1 with {\arrow[scale=1.5]{>}}},postaction={decorate},shorten >=0.5pt] (1) -- (2) node[anchor=south] {};
  \draw[decoration={markings,mark=at position 1 with {\arrow[scale=1.5]{>}}},postaction={decorate},shorten >=0.5pt] (8) -- (9) node[anchor=south] {};
  \draw[decoration={markings,mark=at position 1 with {\arrow[scale=1.5]{>}}},postaction={decorate},shorten >=0.5pt] (2) -- (3) node[anchor=south] {$\partial^{X}_{s}$};
  \draw[decoration={markings,mark=at position 1 with {\arrow[scale=1.5]{>}}},postaction={decorate},shorten >=0.5pt] (9) -- (10) node[anchor=south] {$\partial^{F}_{s}$};
  \draw[decoration={markings,mark=at position 1 with {\arrow[scale=1.5]{>}}},postaction={decorate},shorten >=0.5pt] (3) -- (4) node[anchor=south] {};
  \draw[decoration={markings,mark=at position 1 with {\arrow[scale=1.5]{>}}},postaction={decorate},shorten >=0.5pt] (10) -- (11) node[anchor=south] {};
  \draw[decoration={markings,mark=at position 1 with {\arrow[scale=1.5]{>}}},postaction={decorate},shorten >=0.5pt] (4) -- (5) node[anchor=south] {};
  \draw[decoration={markings,mark=at position 1 with {\arrow[scale=1.5]{>}}},postaction={decorate},shorten >=0.5pt] (11) -- (12) node[anchor=south] {};
  \draw[decoration={markings,mark=at position 1 with {\arrow[scale=1.5]{>}}},postaction={decorate},shorten >=0.5pt] (5) -- (6) node[anchor=south] {$\partial^{X}_{1}$};
  \draw[decoration={markings,mark=at position 1 with {\arrow[scale=1.5]{>}}},postaction={decorate},shorten >=0.5pt] (12) -- (13) node[anchor=south] {$\partial^{F}_{1}$};
  \draw[decoration={markings,mark=at position 1 with {\arrow[scale=1.5]{>}}},postaction={decorate},shorten >=0.5pt] (6) -- (7) node[anchor=south] {};
  \draw[decoration={markings,mark=at position 1 with {\arrow[scale=1.5]{>}}},postaction={decorate},shorten >=0.5pt] (13) -- (14) node[anchor=south] {};
\end{tikzpicture}\]
Let $C:= \coker \sigma$. By the above discussion, $C_{s}= \coker \sigma_{s}= \coker \lambda_{s}$ is Gorenstein projective. Now fix $0 \leq i \leq s-1$. There exists a short exact sequence of $R$-modules
$$0\rightarrow \coker \beta_{i} \xrightarrow {\overline{\lambda_{i}}} \coker \sigma_{i} \rightarrow \coker \lambda_{i} \rightarrow 0,$$
where $\overline{\lambda_{i}}$ is given by $\overline{\lambda_{i}}(x+\im \beta_{i})= \lambda_{i}(x)+ \im \sigma_{i}$.
By the above construction, $\coker \lambda_{i}$ is Gorenstein projective, and thus $\coker \beta_{i}\cong \coker \lambda_{i+1}$ is Gorenstein projective. It follows from Proposition \ref{1.2.3} that $C_{i}= \coker \sigma_{i}$ is Gorenstein projective. We can then form the short exact sequence of $R$-complexes
$$0\rightarrow X \xrightarrow {\sigma} F \rightarrow C \rightarrow 0,$$
where $C$ is an $R$-complex of Gorenstein projective modules.
\end{construction}

We can even say more about homology modules.

\begin{theorem} \label{1.2.13}
Let
$$X: 0 \rightarrow X_{s} \rightarrow X_{s-1} \rightarrow \cdots \rightarrow X_{1} \rightarrow X_{0} \rightarrow 0$$
be a bounded $R$-complex of Gorenstein projective modules. Then there exist $R$-complexes
$$F: 0 \rightarrow F_{s} \rightarrow F_{s-1} \rightarrow \cdots \rightarrow F_{1} \rightarrow F_{0} \rightarrow 0,$$
consisting of free modules, and
$$C: 0 \rightarrow C_{s} \rightarrow C_{s-1} \rightarrow \cdots \rightarrow C_{1} \rightarrow C_{0} \rightarrow 0,$$
consisting of Gorenstein projective modules, that fit into a short exact sequence of $R$-complexes
$$0\rightarrow X \xrightarrow{\sigma} F \rightarrow C \rightarrow 0.$$
In addition, the following assertions hold:
\begin{enumerate}
\item[(i)] The morphism $\sigma:X \rightarrow F$ induces isomorphisms $H_{i}(X)\cong H_{i}(F)$ for every $1\leq i \leq s$. As a result, $H_{i}(C)=0$ for every $1\leq i \leq s$.
\item[(ii)] There is a short exact sequence of $R$-modules
$$0\rightarrow H_{0}(X) \xrightarrow {H_{0}(\sigma)} H_{0}(F) \rightarrow H_{0}(C) \rightarrow 0,$$
where $H_{0}(C)$ is Gorenstein projective.
\item[(iii)] If in particular, $X$ is an $R$-complex of projective modules, then $C$ is also an $R$-complex of projective modules. Furthermore, $H_{0}(C)$ is projective, so in particular, the short exact sequence in (ii) splits.
\item[(iv)] If $X_{i}$ happens to be projective for every $1 \leq i \leq s$, then $C_{i}$ is projective for every $1 \leq i \leq s$, and $H_{0}(C)$ is isomorphic to a direct summand of $C_{0}$ with a projective complement.
\item[(v)] If $H_{i}(X)=0$ for every $1 \leq i \leq s$, then $\pd_{R} \left(H_{0}(F)\right) \leq s$.
\end{enumerate}
\end{theorem}

\begin{prf}
(i) and (ii): We apply Construction \ref{1.2.12} to the $R$-complex $X$ to get the desired short exact sequence of $R$-complexes
$$0\rightarrow X \xrightarrow {\sigma} F \rightarrow C \rightarrow 0.$$
We now focus on homology modules, sticking to the notation of Construction \ref{1.2.12}.

First, we consider $H_{s}(\sigma):H_{s}(X)\rightarrow H_{s}(F)$. Suppose that
$$H_{s}(\sigma)(x)=\sigma_{s}(x)=\lambda_{s}(x)=0,$$
for some $x\in H_{s}(X)= \ker \partial^{X}_{s}$. Then $x=0$, so $H_{s}(\sigma)$ is injective. Now let $y\in H_{s}(F)= \ker \partial^{F}_{s}$.
Then $\partial^{F}_{s}(y)=\lambda_{s-1}\left(\alpha_{s}(y)\right)$, so $\alpha_{s}(y)=0$.
By the structure of the pushout, $\alpha_{s}(y)= (y,0)+\Im_{s}$, where
$$\Im_{s}:= \left\{\left(-\lambda_{s}(g),\partial^{X}_{s}(g)\right) \suchthat g\in X_{s}\right\}.$$
Therefore, there exists an element $g\in X_{s}$ such that $(y,0)= \left(-\lambda_{s}(g),\partial^{X}_{s}(g)\right)$. Thus $y= \lambda_{s}(-g)$ and $\partial^{X}_{s}(-g)= 0$. That is to say
$$H_{s}(\sigma)(-g)=\sigma_{s}(-g)= \lambda _{s}(-g)= y$$
and $-g\in H_{s}(X)$, showing that $H_{s}(\sigma)$ is surjective, and thus an isomorphism.

Next, we consider $H_{i}(\sigma):H_{i}(X)\rightarrow H_{i}(F)$ for any given $0 \leq i \leq s-1$. Suppose that
$$H_{i}(\sigma)\left(x+ \im \partial^{X}_{i+1}\right)= \sigma_{i}(x)+ \im \partial^{F}_{i+1}= \lambda_{i}\left(\beta_{i}(x)\right)+ \im \partial^{F}_{i+1}=0,$$
for some $x\in \ker \partial^{X}_{i}$.
Then $\lambda_{i}\left(\beta_{i}(x)\right)=\partial^{F}_{i+1}(t)$ for some $t\in F_{i+1}$. But then $\lambda_{i}\left(\beta_{i}(x)\right)=\lambda_{i}\left(\alpha_{i+1}(t)\right)$, so
$$\beta_{i}(x)=(0,x)+ \Im_{i+1} =(t,0)+\Im_{i+1}=\alpha_{i+1}(t),$$
where
$$\Im_{i+1}:= \left\{\left(-\lambda_{i+1}(g),\varphi_{i+1}(g)\right) \suchthat g\in L_{i+1}\right\}.$$
Therefore, there exists an element $g\in L_{i+1}$, such that $(-t,x)=\left(-\lambda_{i+1}(g),\varphi_{i+1}(g)\right)$,
so
$$x=\varphi_{i+1}(g)= \overline{\partial^{X}_{i+1}}\left(\widetilde{\beta_{i+1}}^{-1}\left(\pi_{i+1}(g)\right)\right) \in \im \partial^{X}_{i+1},$$
showing that $H_{i}(\sigma)$ is injective.
Now suppose that $1 \leq i \leq s$, and let $y+ \im \partial^{F}_{i+1}\in H_{i}(F)$. Thus $\partial^{F}_{i}(y)=\lambda_{i-1}\left(\alpha_{i}(y)\right)=0$,
so $\alpha_{i}(y)=0$. Just as before, it follows that there exists an element $g\in L_{i}$, such that $(y,0)=\left(-\lambda_{i}(g),\varphi_{i}(g)\right)$, so $y=\lambda_{i}(-g)$ and $\varphi_{i}(-g)=0$.
By the definition of $\varphi_{i}$, we have
$$\overline{\partial^{X}_{i}}\left(\widetilde{\beta_{i}}^{-1}\left(\pi_{i}(-g)\right)\right)= \overline{\partial^{X}_{i}}\left(\widetilde{\beta_{i}}^{-1}(-g+ \im \alpha_{i+1})\right)=0.$$
If $\widetilde{\beta_{i}}^{-1}\left(-g+ \im \alpha_{i+1}\right)= h+ \im \varphi_{i+1}$, for some $h\in X_{i}$, then
$$\overline{\partial^{X}_{i}}\left(\widetilde{\beta_{i}}^{-1}(-g+ \im \alpha_{i+1})\right)=\overline{\partial^{X}_{i}}\left(h+ \im \varphi_{i+1}\right)=\partial^{X}_{i}(h)=0,$$
i.e. $h\in \ker \partial^{X}_{i}$.
On the other hand, by the choice of $h$, we have
$$\widetilde{\beta_{i}}\left(h+ \im \varphi_{i+1}\right)= \beta_{i}(h)+ \im \alpha_{i+1} =-g + \im \alpha_{i+1}.$$
Hence $\beta_{i}(h)+g=\alpha_{i+1}(e)$ for some $e\in F_{i+1}$. It follows that
\begin{equation*}
\begin{split}
H_{i}(\sigma)\left(h+ \im \varphi_{i+1}\right) & = \sigma_{i}(h)+ \im \partial^{F}_{i+1} \\
 & = \lambda_{i}\left(\beta_{i}(h)\right)+ \im \partial^{F}_{i+1} \\
 & = \lambda_{i}\left(-g+ \alpha_{i+1}(e)\right)+ \im \partial^{F}_{i+1} \\
 & = \lambda_{i}(-g)+ \lambda_{i}\left(\alpha_{i+1}(e)\right)+ \im \partial^{F}_{i+1} \\
 & = y+ \partial^{F}_{i+1}(e)+ \im \partial^{F}_{i+1} \\
 & = y+ \im \partial^{F}_{i+1}. \\
\end{split}
\end{equation*}
This shows that $H_{i}(\sigma)$ is surjective, and thus an isomorphism.

Consider the short exact sequence
$$0\rightarrow X \xrightarrow {\sigma} F \rightarrow C \rightarrow 0$$
of $R$-complexes, and its associated long exact sequence
$$0\rightarrow H_{s}(X) \xrightarrow {H_{s}(\sigma)} H_{s}(F) \rightarrow H_{s}(C) \rightarrow \cdots \rightarrow H_{1}(X) \xrightarrow {H_{1}(\sigma)} H_{1}(F) \rightarrow $$$$ H_{1}(C) \rightarrow H_{0}(X) \xrightarrow {H_{0}(\sigma)} H_{0}(F) \rightarrow H_{0}(C) \rightarrow 0,$$
of homology modules. For any given $1 \leq i \leq s$, the exact sequence
$$H_{i}(X) \xrightarrow {H_{i}(\sigma)} H_{i}(F) \xrightarrow {\vartheta} H_{i}(C) \rightarrow H_{i-1}(X) \xrightarrow {H_{i-1}(\sigma)} H_{i-1}(F)$$
shows that
$$0= \coker H_{i}(\sigma)\cong \im \vartheta$$
and
$$\coker \vartheta \cong \ker H_{i-1}(\sigma)=0,$$
hence $H_{i}(C)=0$.
Thus we get a short exact sequence
$$0\rightarrow H_{0}(X) \xrightarrow {H_{0}(\sigma)} H_{0}(F) \rightarrow H_{0}(C) \rightarrow 0$$
of $R$-modules. We further notice that
\begin{equation*}
\begin{split}
H_{0}(C) & \cong \coker H_{0}(\sigma) \\
 & = \frac{H_{0}(F)}{H_{0}(\sigma)\left(H_{0}(X)\right)} \\
 & \cong \frac{F_{0}}{\sigma_{0}(X_{0})+ \im \partial^{F}_{1}} \\
 & = \frac{F_{0}}{\lambda_{0}\left(\beta_{0}(X_{0})+\alpha_{1}(F_{1})\right)} \\
 & = \frac{F_{0}}{\lambda_{0}(L_{0})} \\
 & = \coker \lambda_{0}. \\
\end{split}
\end{equation*}
But $\coker \lambda_{0}$ is Gorenstein projective by the construction.

(iii): We consider the following commutative diagram with exact columns:
\[\begin{tikzpicture}[every node/.style={midway}]
  \matrix[column sep={3em}, row sep={3em}]
  {\node(1) {}; & \node(2) {$0$}; & \node(3) {$0$}; & \node(4) {}; & \node(5) {$0$}; & \node(6) {$0$}; & \node(7) {};\\
  \node(8) {$0$}; & \node(9) {$X_{s}$}; & \node(10) {$X_{s-1}$}; & \node(11) {$\cdots$}; & \node(12) {$X_{1}$}; & \node(13) {$X_{0}$}; & \node(14) {$0$};\\
  \node(15) {$0$}; & \node(16) {$F_{s}$}; & \node(17) {$F_{s-1}$}; & \node(18) {$\cdots$}; & \node(19) {$F_{1}$}; & \node(20) {$F_{0}$}; & \node(21) {$0$};\\
  \node(22) {$0$}; & \node(23) {$C_{s}$}; & \node(24) {$C_{s-1}$}; & \node(25) {$\cdots$}; & \node(26) {$C_{1}$}; & \node(27) {$C_{0}$}; & \node(28) {$0$};\\
  \node(29) {}; & \node(30) {$0$}; & \node(31) {$0$}; & \node(32) {}; & \node(33) {$0$}; & \node(34) {$0$}; & \node(35) {};\\};
  \draw[decoration={markings,mark=at position 1 with {\arrow[scale=1.5]{>}}},postaction={decorate},shorten >=0.5pt] (2) -- (9) node[anchor=west] {};
  \draw[decoration={markings,mark=at position 1 with {\arrow[scale=1.5]{>}}},postaction={decorate},shorten >=0.5pt] (3) -- (10) node[anchor=west] {};
  \draw[decoration={markings,mark=at position 1 with {\arrow[scale=1.5]{>}}},postaction={decorate},shorten >=0.5pt] (5) -- (12) node[anchor=west] {};
  \draw[decoration={markings,mark=at position 1 with {\arrow[scale=1.5]{>}}},postaction={decorate},shorten >=0.5pt] (6) -- (13) node[anchor=west] {};
  \draw[decoration={markings,mark=at position 1 with {\arrow[scale=1.5]{>}}},postaction={decorate},shorten >=0.5pt] (9) -- (16) node[anchor=west] {};
  \draw[decoration={markings,mark=at position 1 with {\arrow[scale=1.5]{>}}},postaction={decorate},shorten >=0.5pt] (10) -- (17) node[anchor=west] {};
  \draw[decoration={markings,mark=at position 1 with {\arrow[scale=1.5]{>}}},postaction={decorate},shorten >=0.5pt] (12) -- (19) node[anchor=west] {};
  \draw[decoration={markings,mark=at position 1 with {\arrow[scale=1.5]{>}}},postaction={decorate},shorten >=0.5pt] (13) -- (20) node[anchor=west] {};
  \draw[decoration={markings,mark=at position 1 with {\arrow[scale=1.5]{>}}},postaction={decorate},shorten >=0.5pt] (6) -- (13) node[anchor=west] {};
  \draw[decoration={markings,mark=at position 1 with {\arrow[scale=1.5]{>}}},postaction={decorate},shorten >=0.5pt] (16) -- (23) node[anchor=west] {};
  \draw[decoration={markings,mark=at position 1 with {\arrow[scale=1.5]{>}}},postaction={decorate},shorten >=0.5pt] (17) -- (24) node[anchor=west] {};
  \draw[decoration={markings,mark=at position 1 with {\arrow[scale=1.5]{>}}},postaction={decorate},shorten >=0.5pt] (19) -- (26) node[anchor=west] {};
  \draw[decoration={markings,mark=at position 1 with {\arrow[scale=1.5]{>}}},postaction={decorate},shorten >=0.5pt] (20) -- (27) node[anchor=west] {};
  \draw[decoration={markings,mark=at position 1 with {\arrow[scale=1.5]{>}}},postaction={decorate},shorten >=0.5pt] (23) -- (30) node[anchor=west] {};
  \draw[decoration={markings,mark=at position 1 with {\arrow[scale=1.5]{>}}},postaction={decorate},shorten >=0.5pt] (24) -- (31) node[anchor=west] {};
  \draw[decoration={markings,mark=at position 1 with {\arrow[scale=1.5]{>}}},postaction={decorate},shorten >=0.5pt] (26) -- (33) node[anchor=west] {};
  \draw[decoration={markings,mark=at position 1 with {\arrow[scale=1.5]{>}}},postaction={decorate},shorten >=0.5pt] (27) -- (34) node[anchor=west] {};
  \draw[decoration={markings,mark=at position 1 with {\arrow[scale=1.5]{>}}},postaction={decorate},shorten >=0.5pt] (8) -- (9) node[anchor=south] {};
  \draw[decoration={markings,mark=at position 1 with {\arrow[scale=1.5]{>}}},postaction={decorate},shorten >=0.5pt] (9) -- (10) node[anchor=south] {$\partial^{X}_{s}$};
  \draw[decoration={markings,mark=at position 1 with {\arrow[scale=1.5]{>}}},postaction={decorate},shorten >=0.5pt] (10) -- (11) node[anchor=south] {};
  \draw[decoration={markings,mark=at position 1 with {\arrow[scale=1.5]{>}}},postaction={decorate},shorten >=0.5pt] (11) -- (12) node[anchor=south] {};
  \draw[decoration={markings,mark=at position 1 with {\arrow[scale=1.5]{>}}},postaction={decorate},shorten >=0.5pt] (12) -- (13) node[anchor=south] {$\partial^{X}_{1}$};
  \draw[decoration={markings,mark=at position 1 with {\arrow[scale=1.5]{>}}},postaction={decorate},shorten >=0.5pt] (13) -- (14) node[anchor=south] {};
  \draw[decoration={markings,mark=at position 1 with {\arrow[scale=1.5]{>}}},postaction={decorate},shorten >=0.5pt] (15) -- (16) node[anchor=south] {};
  \draw[decoration={markings,mark=at position 1 with {\arrow[scale=1.5]{>}}},postaction={decorate},shorten >=0.5pt] (16) -- (17) node[anchor=south] {$\partial^{F}_{s}$};
  \draw[decoration={markings,mark=at position 1 with {\arrow[scale=1.5]{>}}},postaction={decorate},shorten >=0.5pt] (17) -- (18) node[anchor=south] {};
  \draw[decoration={markings,mark=at position 1 with {\arrow[scale=1.5]{>}}},postaction={decorate},shorten >=0.5pt] (18) -- (19) node[anchor=south] {};
  \draw[decoration={markings,mark=at position 1 with {\arrow[scale=1.5]{>}}},postaction={decorate},shorten >=0.5pt] (19) -- (20) node[anchor=south] {$\partial^{F}_{1}$};
  \draw[decoration={markings,mark=at position 1 with {\arrow[scale=1.5]{>}}},postaction={decorate},shorten >=0.5pt] (20) -- (21) node[anchor=south] {};
  \draw[decoration={markings,mark=at position 1 with {\arrow[scale=1.5]{>}}},postaction={decorate},shorten >=0.5pt] (22) -- (23) node[anchor=south] {};
  \draw[decoration={markings,mark=at position 1 with {\arrow[scale=1.5]{>}}},postaction={decorate},shorten >=0.5pt] (23) -- (24) node[anchor=south] {$\partial^{C}_{s}$};
  \draw[decoration={markings,mark=at position 1 with {\arrow[scale=1.5]{>}}},postaction={decorate},shorten >=0.5pt] (24) -- (25) node[anchor=south] {};
  \draw[decoration={markings,mark=at position 1 with {\arrow[scale=1.5]{>}}},postaction={decorate},shorten >=0.5pt] (25) -- (26) node[anchor=south] {};
  \draw[decoration={markings,mark=at position 1 with {\arrow[scale=1.5]{>}}},postaction={decorate},shorten >=0.5pt] (26) -- (27) node[anchor=south] {$\partial^{C}_{1}$};
  \draw[decoration={markings,mark=at position 1 with {\arrow[scale=1.5]{>}}},postaction={decorate},shorten >=0.5pt] (27) -- (28) node[anchor=south] {};
\end{tikzpicture}\]
For any given $0 \leq i \leq s$, the exact column
$$0 \rightarrow X_{i} \rightarrow F_{i} \rightarrow C_{i} \rightarrow 0$$
in the above diagram shows that $\pd_{R}(C_{i})<\infty$, so Proposition \ref{1.2.6} implies that $\pd_{R}(C_{i})= \Gpd_{R}(C_{i})=0$, i.e. $C_{i}$ is projective. Consider the exact augmented $R$-complex
$$C^{+}: 0 \rightarrow C_{s} \xrightarrow {\partial^{C}_{s}} C_{s-1} \rightarrow \cdots \rightarrow C_{1} \xrightarrow {\partial^{C}_{1}} C_{0} \rightarrow H_{0}(C) \rightarrow 0.$$
It shows that $\pd_{R}\left(H_{0}(C)\right)<\infty$, so Proposition \ref{1.2.6} implies that $$\pd_{R}\left(H_{0}(C)\right)= \Gpd_{R}\left(H_{0}(C)\right)=0,$$ i.e. $H_{0}(C)$ is projective.

(iv): We first note that if $s=0$, then the condition is void, and the assertion trivially holds. Now let $s\geq 1$. Since $X_{i}$ is projective for every $1 \leq i \leq s$, as in the proof of part (iii), it follows that $C_{i}$ is projective for every $1 \leq i \leq s$. On the other hand, $C_{0}$ and $H_{0}(C)$ are Gorenstein projective, so using Proposition \ref{1.2.3}, the short exact sequence
\begin{equation} \label{eq:1.2.13.1}
0 \rightarrow \im \partial^{C}_{1} \rightarrow C_{0} \rightarrow H_{0}(C) \rightarrow 0
\end{equation}
implies that $\im \partial^{C}_{1}$ is Gorenstein projective. On the other hand, the exact sequence
$$0 \rightarrow C_{s} \xrightarrow {\partial^{C}_{s}} C_{s-1} \rightarrow \cdots \rightarrow C_{2} \xrightarrow {\partial^{C}_{2}} C_{1} \rightarrow \im \partial^{C}_{1} \rightarrow 0$$
shows that $\pd_{R}\left(\im \partial^{C}_{1}\right)<\infty$, so Proposition \ref{1.2.6} implies that
$$\pd_{R}\left(\im \partial^{C}_{1}\right)= \Gpd_{R}\left(\im \partial^{C}_{1}\right)=0,$$
i.e. $\im \partial^{C}_{1}$ is projective. Therefore, $\Ext^{1}_{R}\left(H_{0}(C),\im \partial^{C}_{1}\right)=0$, so the short exact sequence \eqref{eq:1.2.13.1} splits and yields the assertions.

(v): If $s=0$, then $H_{0}(F)$ is indeed free. Now let $s \geq 1$. Then $H_{i}(F) \cong H_{i}(X)=0$ for every $1 \leq i \leq s$. Thus the exact sequence
$$0 \rightarrow F_{s} \xrightarrow {\partial^{F}_{s}} F_{s-1} \rightarrow \cdots \rightarrow F_{1} \xrightarrow {\partial^{F}_{1}} F_{0} \rightarrow H_{0}(F) \rightarrow 0,$$
shows that $\pd_{R} \left(H_{0}(F)\right) \leq s$.
\end{prf}

We immediately draw the following conclusion.

\begin{corollary} \label{1.2.14}
Let $M$ be an $R$-module with $\Gpd_{R}(M)< \infty$. Then there exists a short exact sequence
$$0 \rightarrow M \rightarrow L \rightarrow T \rightarrow 0$$
of $R$-modules, where $\pd_{R}(L)= \Gpd_{R}(M)$, and $T$ is Gorenstein projective.
\end{corollary}

\begin{prf}
Let $\Gpd_{R}(M)=s$. The case $s=0$ is established in the proof of Proposition \ref{1.2.9}, so we may assume that $s\geq 1$. There exists an exact sequence
$$0 \rightarrow X_{s} \rightarrow X_{s-1} \rightarrow \cdots \rightarrow X_{1} \rightarrow X_{0} \rightarrow M \rightarrow 0,$$
where $X_{i}$ is a Gorenstein projective $R$-module for every $0 \leq i \leq s$. Now consider the $R$-complex
$$X: 0 \rightarrow X_{s} \rightarrow X_{s-1} \rightarrow \cdots \rightarrow X_{1} \rightarrow X_{0} \rightarrow 0,$$
and apply Theorem \ref{1.2.13} to get the $R$-complexes $F$ and $C$, and the short exact sequence of $R$-modules
$$0 \rightarrow H_{0}(X) \rightarrow H_{0}(F) \rightarrow H_{0}(C) \rightarrow 0.$$
Obviously, $H_{0}(X)\cong M$. Let $L= H_{0}(F)$ and $T= H_{0}(C)$.
Theorem \ref{1.2.13} implies that $\pd_{R}(L) \leq s$ and $T$ is Gorenstein projective. By Proposition \ref{1.2.5}, we may choose a projective $R$-module $P$ for which $\Ext_{R}^{s}(M,P)\neq 0$.
If $\pd_{R}(L) < s$, then in view of Theorem \ref{1.2.13}, we get the exact sequence
$$0= \Ext_{R}^{s}(L,P) \rightarrow \Ext_{R}^{s}(M,P) \rightarrow \Ext_{R}^{s+1}(T,P)=0,$$
which gives $\Ext_{R}^{s}(M,P)=0$, a contradiction. Hence $\pd_{R}(L) = s$.
\end{prf}

\begin{remark} \label{1.2.15}
As mentioned in the Introduction, given a bounded $R$-complex
$$X: 0 \rightarrow X_{s} \rightarrow X_{s-1} \rightarrow  \cdots \rightarrow X_{1} \rightarrow X_{0} \rightarrow 0$$
of Gorenstein projective modules, we are interested in finding an $R$-complex
$$F: 0 \rightarrow F_{s} \rightarrow F_{s-1} \rightarrow \cdots \rightarrow F_{1} \rightarrow F_{0} \rightarrow 0$$
of free $R$-modules, with the property that the homology modules of $F$ are as close as possible to the homology modules of $X$. The best-case scenario is that $H_{i}(X) \cong H_{i}(F)$ for every $0 \leq i \leq s$. But as the example below suggests, this cannot occur in general. Therefore Theorem \ref{1.2.13} gives, in some sense, the best approximation one could hope for.
\end{remark}

\begin{example} \label{1.2.16}
Let $(R,\mathfrak{m},k)$ be a non-regular Gorenstein local ring. We know that $\Gpd_{R}(k) < \infty$. Set $s:=\Gpd_{R}(k)$, and take a Gorenstein projective resolution
$$X: 0 \rightarrow X_{s} \rightarrow X_{s-1} \rightarrow \cdots \rightarrow X_{1} \rightarrow X_{0} \rightarrow 0$$
for $k$. If there were an $R$-complex
$$F: 0 \rightarrow F_{s} \rightarrow F_{s-1} \rightarrow \cdots \rightarrow F_{1} \rightarrow F_{0} \rightarrow 0$$
of free $R$-modules with the property that $H_{i}(X)\cong H_{i}(F)$ for every $0 \leq i \leq s$, then $F$ would be a free resolution for $k$, implying that $\pd_{R}(k)<\infty$. This contradicts the non-regularity of $R$.
\end{example}

Our next goal is to apply Construction \ref{1.2.12} to bounded $R$-complexes of Gorenstein projective modules whose homology modules have finite length.

We first need the following two propositions.

\begin{proposition} \label{1.2.17}
For every $R$-complex $X$, there is a an $R$-complex $F$ of free modules and a quasi-isomorphism $F \xrightarrow {\simeq} X$, such that $F_{i}=0$ for every $i< \inf X$.
\end{proposition}

\begin{prf}
See \cite[Theorem 5.1.12]{CFH2}.
\end{prf}

\begin{proposition} \label{1.2.18}
Let $X$ be a homologically right bounded $R$-complex, and $P$ a right-bounded $R$-complex of Gorenstein projective modules such that $X \simeq P$ in the derived category $\mathcal{D}(R)$. If $\Gpd_{R}(X)\leq s$, then $\coker \partial^{P}_{s+1}$ is Gorenstein projective.
\end{proposition}

\begin{prf}
See \cite[Theorem 3.1]{CFrH}.
\end{prf}

We can now present the following construction.

\begin{construction} \label{1.2.19}
Let
$$X: 0\rightarrow X_{s} \xrightarrow {\partial^{X}_{s}} X_{s-1}\rightarrow \cdots \rightarrow X_{1} \xrightarrow {\partial^{X}_{1}} X_{0} \rightarrow 0$$
be an $R$-complex of Gorenstein projective modules. By Proposition \ref{1.2.17}, there exists an $R$-complex $E$ of free modules and a quasi-isomorphism $f:E \xrightarrow {\simeq} X$. Furthermore, since $\inf X \geq 0$, we may choose $E_{i}=0$ for every $i<0$. We then have the following commutative diagram:
\[\begin{tikzpicture}[every node/.style={midway}]
\matrix[column sep={3em}, row sep={3em}]
  {\node(1) {$\cdots$}; & \node(2) {$E_{s+1}$}; & \node(3) {$E_{s}$}; & \node(4) {$E_{s-1}$}; & \node(5) {$\cdots$}; & \node(6) {$E_{1}$}; & \node(7) {$E_{0}$}; & \node(8) {$0$};\\
  \node(9) {}; & \node(10) {$0$}; & \node(11) {$X_{s}$}; & \node(12) {$X_{s-1}$}; & \node(13) {$\cdots$}; & \node(14) {$X_{1}$}; & \node(15) {$X_{0}$}; & \node(16) {$0$};\\};
  \draw[decoration={markings,mark=at position 1 with {\arrow[scale=1.5]{>}}},postaction={decorate},shorten >=0.5pt] (1) -- (2) node[anchor=south] {};
  \draw[decoration={markings,mark=at position 1 with {\arrow[scale=1.5]{>}}},postaction={decorate},shorten >=0.5pt] (2) -- (3) node[anchor=south] {$\partial^{E}_{s+1}$};
  \draw[decoration={markings,mark=at position 1 with {\arrow[scale=1.5]{>}}},postaction={decorate},shorten >=0.5pt] (3) -- (4) node[anchor=south] {$\partial^{E}_{s}$};
  \draw[decoration={markings,mark=at position 1 with {\arrow[scale=1.5]{>}}},postaction={decorate},shorten >=0.5pt] (4) -- (5) node[anchor=south] {};
  \draw[decoration={markings,mark=at position 1 with {\arrow[scale=1.5]{>}}},postaction={decorate},shorten >=0.5pt] (5) -- (6) node[anchor=south] {};
  \draw[decoration={markings,mark=at position 1 with {\arrow[scale=1.5]{>}}},postaction={decorate},shorten >=0.5pt] (6) -- (7) node[anchor=south] {$\partial^{E}_{1}$};
  \draw[decoration={markings,mark=at position 1 with {\arrow[scale=1.5]{>}}},postaction={decorate},shorten >=0.5pt] (7) -- (8) node[anchor=south] {};
  \draw[decoration={markings,mark=at position 1 with {\arrow[scale=1.5]{>}}},postaction={decorate},shorten >=0.5pt] (10) -- (11) node[anchor=south] {};
  \draw[decoration={markings,mark=at position 1 with {\arrow[scale=1.5]{>}}},postaction={decorate},shorten >=0.5pt] (11) -- (12) node[anchor=south] {$\partial^{X}_{s}$};
  \draw[decoration={markings,mark=at position 1 with {\arrow[scale=1.5]{>}}},postaction={decorate},shorten >=0.5pt] (12) -- (13) node[anchor=south] {};
  \draw[decoration={markings,mark=at position 1 with {\arrow[scale=1.5]{>}}},postaction={decorate},shorten >=0.5pt] (13) -- (14) node[anchor=south] {};
  \draw[decoration={markings,mark=at position 1 with {\arrow[scale=1.5]{>}}},postaction={decorate},shorten >=0.5pt] (14) -- (15) node[anchor=south] {$\partial^{X}_{1}$};
  \draw[decoration={markings,mark=at position 1 with {\arrow[scale=1.5]{>}}},postaction={decorate},shorten >=0.5pt] (15) -- (16) node[anchor=south] {};
  \draw[decoration={markings,mark=at position 1 with {\arrow[scale=1.5]{>}}},postaction={decorate},shorten >=0.5pt] (7) -- (15) node[anchor=west] {$f_{0}$};
  \draw[decoration={markings,mark=at position 1 with {\arrow[scale=1.5]{>}}},postaction={decorate},shorten >=0.5pt] (6) -- (14) node[anchor=west] {$f_{1}$};
  \draw[decoration={markings,mark=at position 1 with {\arrow[scale=1.5]{>}}},postaction={decorate},shorten >=0.5pt] (4) -- (12) node[anchor=west] {$f_{s-1}$};
  \draw[decoration={markings,mark=at position 1 with {\arrow[scale=1.5]{>}}},postaction={decorate},shorten >=0.5pt] (3) -- (11) node[anchor=west] {$f_{s}$};
  \draw[decoration={markings,mark=at position 1 with {\arrow[scale=1.5]{>}}},postaction={decorate},shorten >=0.5pt] (2) -- (10) node[anchor=west] {$f_{s+1}$};
\end{tikzpicture}\]
The above diagram shows that $\im \partial^{E}_{s+1} \subseteq \ker f_{s}$. Let $K_{s}:= \coker \partial^{E}_{s+1}$, and consider the softly truncated $R$-complex $E_{s \subset}$. For convenience, we denote $E_{s \subset}$ by $E$ in what follows. Then we obtain the following commutative diagram:
\[\begin{tikzpicture}[every node/.style={midway}]
  \matrix[column sep={3em}, row sep={3em}]
  {\node(1) {$0$}; & \node(2) {$K_{s}$}; & \node(3) {$E_{s-1}$}; & \node(4) {$\cdots$}; & \node(5) {$E_{1}$}; & \node(6) {$E_{0}$}; & \node(7) {$0$};\\
  \node(8) {$0$}; & \node(9) {$X_{s}$}; & \node(10) {$X_{s-1}$}; & \node(11) {$\cdots$}; & \node(12) {$X_{1}$}; & \node(13) {$X_{0}$}; & \node(14) {$0$};\\};
  \draw[decoration={markings,mark=at position 1 with {\arrow[scale=1.5]{>}}},postaction={decorate},shorten >=0.5pt] (2) -- (9) node[anchor=west] {$\overline{f_{s}}$};
  \draw[decoration={markings,mark=at position 1 with {\arrow[scale=1.5]{>}}},postaction={decorate},shorten >=0.5pt] (3) -- (10) node[anchor=west] {$f_{s-1}$};
  \draw[decoration={markings,mark=at position 1 with {\arrow[scale=1.5]{>}}},postaction={decorate},shorten >=0.5pt] (5) -- (12) node[anchor=west] {$f_{1}$};
  \draw[decoration={markings,mark=at position 1 with {\arrow[scale=1.5]{>}}},postaction={decorate},shorten >=0.5pt] (6) -- (13) node[anchor=west] {$f_{0}$};
  \draw[decoration={markings,mark=at position 1 with {\arrow[scale=1.5]{>}}},postaction={decorate},shorten >=0.5pt] (1) -- (2) node[anchor=south] {};
  \draw[decoration={markings,mark=at position 1 with {\arrow[scale=1.5]{>}}},postaction={decorate},shorten >=0.5pt] (8) -- (9) node[anchor=south] {};
  \draw[decoration={markings,mark=at position 1 with {\arrow[scale=1.5]{>}}},postaction={decorate},shorten >=0.5pt] (2) -- (3) node[anchor=south] {$\overline{\partial^{E}_{s}}$};
  \draw[decoration={markings,mark=at position 1 with {\arrow[scale=1.5]{>}}},postaction={decorate},shorten >=0.5pt] (9) -- (10) node[anchor=south] {$\partial^{X}_{s}$};
  \draw[decoration={markings,mark=at position 1 with {\arrow[scale=1.5]{>}}},postaction={decorate},shorten >=0.5pt] (3) -- (4) node[anchor=south] {};
  \draw[decoration={markings,mark=at position 1 with {\arrow[scale=1.5]{>}}},postaction={decorate},shorten >=0.5pt] (10) -- (11) node[anchor=south] {};
  \draw[decoration={markings,mark=at position 1 with {\arrow[scale=1.5]{>}}},postaction={decorate},shorten >=0.5pt] (4) -- (5) node[anchor=south] {};
  \draw[decoration={markings,mark=at position 1 with {\arrow[scale=1.5]{>}}},postaction={decorate},shorten >=0.5pt] (11) -- (12) node[anchor=south] {};
  \draw[decoration={markings,mark=at position 1 with {\arrow[scale=1.5]{>}}},postaction={decorate},shorten >=0.5pt] (5) -- (6) node[anchor=south] {$\partial^{E}_{1}$};
  \draw[decoration={markings,mark=at position 1 with {\arrow[scale=1.5]{>}}},postaction={decorate},shorten >=0.5pt] (12) -- (13) node[anchor=south] {$\partial^{X}_{1}$};
  \draw[decoration={markings,mark=at position 1 with {\arrow[scale=1.5]{>}}},postaction={decorate},shorten >=0.5pt] (6) -- (7) node[anchor=south] {};
  \draw[decoration={markings,mark=at position 1 with {\arrow[scale=1.5]{>}}},postaction={decorate},shorten >=0.5pt] (13) -- (14) node[anchor=south] {};
\end{tikzpicture}\]
where $\overline{\partial^{E}_{s}}:K_{s} \rightarrow E_{s-1}$ is given by $\overline{\partial^{E}_{s}}\left(x+ \im \partial^{E}_{s+1}\right)= \partial^{E}_{s}(x)$ for every $x\in E_{s}$, and $\overline{f_{s}}:K_{s}\rightarrow X_{s}$ is given by $\overline{f_{s}}\left(x+ \im \partial^{E}_{s+1}\right)= f_{s}(x)$ for every $x\in E_{s}$. Note that these homomorphisms are well-defined since $\im \partial^{E}_{s+1} \subseteq \ker \partial^{E}_{s} \cap \ker f_{s}$. It is easily seen that the new $f:E \rightarrow X$, obtained in the above diagram, is also a quasi-isomorphism, so that $H_{i}(E)\cong H_{i}(X)$ for every $0 \leq i \leq s$. On the other hand, as $X$ is an $R$-complex of Gorenstein projective modules, it is clear that $\Gpd_{R} (X) \leq s$, so Proposition \ref{1.2.18} implies that $K_{s}$ is Gorenstein projective.
\end{construction}

To study complexes of Gorenstein projective modules whose homology modules have finite length, we need the following notion.

\begin{definition} \label{1.2.20}
Let $(R,\mathfrak{m})$ be a local ring. An $R$-module $M$ is said to be \textit{locally free on the punctured spectrum of} $R$ \index{locally free on the punctured spectrum}, if $M_{\mathfrak{p}}$ is a free $R_{\mathfrak{p}}$-module for every $\mathfrak{p}\in \Spec(R)\backslash \{\mathfrak{m}\}$. Further, given an integer $r \geq 0$, $M$ is said to be \textit{locally free on the punctured spectrum of} $R$ \textit{with constant rank} $r$, if $M_{\mathfrak{p}}$ is a free $R_{\mathfrak{p}}$-module of rank $r$ for every $\mathfrak{p}\in \Spec(R)\backslash \{\mathfrak{m}\}$.
\end{definition}

The following lemma may be of independent interest.

\begin{lemma} \label{1.2.21}
Let $(R,\mathfrak{m})$ be a noetherian local ring, and let
$$X: 0\rightarrow X_{s} \rightarrow X_{s-1} \rightarrow \cdots \rightarrow X_{1} \rightarrow X_{0} \rightarrow 0$$
be an $R$-complex whose homology modules have all finite length.
Then the following assertions hold:
\begin{enumerate}
\item[(i)] If $X_{i}$ is locally free on the punctured spectrum of $R$ with constant rank for every $0 \leq i \leq s-1$, then $X_{s}$ is locally free on the punctured spectrum of $R$ with constant rank.
\item[(ii)] If $X_{i}$ is locally free on the punctured spectrum of $R$ with constant rank for every $1 \leq i \leq s$, and $X_{0}$ is Gorenstein projective, then $X_{0}$ is locally free on the punctured spectrum of $R$ with constant rank.
\end{enumerate}
\end{lemma}

\begin{prf}
Let $\mathfrak{p}\in \Spec(R)\backslash \{\mathfrak{m}\}$. Since $\ell_{R}\left(H_{i}(X)\right)<\infty$ for every $0 \leq i \leq s$, we have $H_{i}\left(X_{\mathfrak{p}}\right)= H_{i}(X)_{\mathfrak{p}}=0$, so we get the exact sequence
\begin{equation} \label{eq:1.2.21.1}
X_{\mathfrak{p}}: 0 \rightarrow (X_{s})_{\mathfrak{p}} \rightarrow (X_{s-1})_{\mathfrak{p}} \rightarrow \cdots \rightarrow (X_{1})_{\mathfrak{p}} \rightarrow (X_{0})_{\mathfrak{p}} \rightarrow 0.
\end{equation}

In the situation of part (i), $(X_{i})_{\mathfrak{p}}$ is a free $R_{\mathfrak{p}}$-module for every $0 \leq i \leq s-1$.
Extracting a short exact sequence
\begin{equation} \label{eq:1.2.21.2}
0\rightarrow K \rightarrow (X_{1})_{\mathfrak{p}} \rightarrow (X_{0})_{\mathfrak{p}} \rightarrow 0
\end{equation}
from $X_{\mathfrak{p}}$, we are left with an exact sequence
\begin{equation} \label{eq:1.2.21.3}
0 \rightarrow (X_{s})_{\mathfrak{p}} \rightarrow (X_{s-1})_{\mathfrak{p}} \rightarrow \cdots \rightarrow (X_{2})_{\mathfrak{p}} \rightarrow K \rightarrow 0.
\end{equation}
Since $(X_{0})_{\mathfrak{p}}$ is free, the short exact sequence \eqref{eq:1.2.21.2} splits, showing that $K$ is projective and hence free.
Continuing this way with \eqref{eq:1.2.21.3}, we conclude that $(X_{s})_{\mathfrak{p}}$ is free. Moreover, the number
$$\rank_{R} \left((X_{s})_{\mathfrak{p}}\right)=\sum_{i=0}^{s-1} (-1)^{i}\rank_{R} \left((X_{s-i-1})_{\mathcal{\mathfrak{p}}}\right)$$
is independent of the choice of $\mathfrak{p}$, since each $\rank_{R} \left((X_{s-i-1})_{\mathcal{\mathfrak{p}}}\right)$ is so.

In the situation of part (ii), $(X_{i})_{\mathfrak{p}}$ is a free $R_{\mathfrak{p}}$-module for every $1 \leq i \leq s$.
Now the sequence \eqref{eq:1.2.21.1} shows that $\pd_{R_{\mathfrak{p}}}\left((X_{0})_{\mathfrak{p}}\right)<\infty$.
On the other hand, it follows from Proposition \ref{1.2.4} that $(X_{0})_{\mathfrak{p}}$ is a Gorenstein projective $R_{\mathfrak{p}}$-module. Now Proposition \ref{1.2.6} implies that $\pd_{R_{\mathfrak{p}}}\left((X_{0})_{\mathfrak{p}}\right)= \Gpd_{R_{\mathfrak{p}}}\left((X_{0})_{\mathfrak{p}}\right)=0$, so $(X_{0})_{\mathfrak{p}}$ is projective and thus free.
Besides, the number
$$\rank_{R} \left((X_{0})_{\mathfrak{p}}\right)=\sum_{i=1}^{s} (-1)^{i+1}\rank_{R} \left((X_{i})_{\mathcal{\mathfrak{p}}}\right)$$
is independent of the choice of $\mathfrak{p}$, since each $\rank_{R} \left((X_{i})_{\mathcal{\mathfrak{p}}}\right)$ is so.
\end{prf}

We now present the anticipated theorem.

\begin{theorem} \label{1.2.22}
Let $(R,\mathfrak{m})$ be a noetherian local ring, and let
$$X: 0\rightarrow X_{s} \rightarrow X_{s-1} \rightarrow \cdots \rightarrow X_{1} \rightarrow X_{0} \rightarrow 0$$
be an $R$-complex of Gorenstein projective modules whose homology modules have all finite length. Then there exists an $R$-complex
$$F: 0\rightarrow F_{s} \rightarrow F_{s-1} \rightarrow \cdots \rightarrow F_{1} \rightarrow F_{0} \rightarrow 0,$$
consisting of free modules, with the following properties:
\begin{enumerate}
\item[(i)] $H_{i}(X)\cong H_{i}(F)$ for every $1 \leq i \leq s$.
\item[(ii)] There is a short exact sequence
$$0 \rightarrow H_{0}(X) \rightarrow H_{0}(F) \rightarrow C \rightarrow 0$$
of locally free modules on the punctured spectrum of $R$, such that $C$ is Gorenstein projective.
\end{enumerate}
\end{theorem}

\begin{prf}
We apply Construction \ref{1.2.19} to the $R$-complex $X$ to get an $R$-complex $E$ with $E_{i}$ free for every $0 \leq i \leq s-1$ and $K_{s}$ Gorenstein projective, and a quasi-isomorphism $f:E \rightarrow X$ as in the following diagram:
\[\begin{tikzpicture}[every node/.style={midway}]
  \matrix[column sep={3em}, row sep={3em}]
  {\node(1) {$0$}; & \node(2) {$K_{s}$}; & \node(3) {$E_{s-1}$}; & \node(4) {$\cdots$}; & \node(5) {$E_{1}$}; & \node(6) {$E_{0}$}; & \node(7) {$0$};\\
  \node(8) {$0$}; & \node(9) {$X_{s}$}; & \node(10) {$X_{s-1}$}; & \node(11) {$\cdots$}; & \node(12) {$X_{1}$}; & \node(13) {$X_{0}$}; & \node(14) {$0$};\\};
  \draw[decoration={markings,mark=at position 1 with {\arrow[scale=1.5]{>}}},postaction={decorate},shorten >=0.5pt] (2) -- (9) node[anchor=west] {$\overline{f_{s}}$};
  \draw[decoration={markings,mark=at position 1 with {\arrow[scale=1.5]{>}}},postaction={decorate},shorten >=0.5pt] (3) -- (10) node[anchor=west] {$f_{s-1}$};
  \draw[decoration={markings,mark=at position 1 with {\arrow[scale=1.5]{>}}},postaction={decorate},shorten >=0.5pt] (5) -- (12) node[anchor=west] {$f_{1}$};
  \draw[decoration={markings,mark=at position 1 with {\arrow[scale=1.5]{>}}},postaction={decorate},shorten >=0.5pt] (6) -- (13) node[anchor=west] {$f_{0}$};
  \draw[decoration={markings,mark=at position 1 with {\arrow[scale=1.5]{>}}},postaction={decorate},shorten >=0.5pt] (1) -- (2) node[anchor=south] {};
  \draw[decoration={markings,mark=at position 1 with {\arrow[scale=1.5]{>}}},postaction={decorate},shorten >=0.5pt] (8) -- (9) node[anchor=south] {};
  \draw[decoration={markings,mark=at position 1 with {\arrow[scale=1.5]{>}}},postaction={decorate},shorten >=0.5pt] (2) -- (3) node[anchor=south] {$\overline{\partial^{E}_{s}}$};
  \draw[decoration={markings,mark=at position 1 with {\arrow[scale=1.5]{>}}},postaction={decorate},shorten >=0.5pt] (9) -- (10) node[anchor=south] {$\partial^{X}_{s}$};
  \draw[decoration={markings,mark=at position 1 with {\arrow[scale=1.5]{>}}},postaction={decorate},shorten >=0.5pt] (3) -- (4) node[anchor=south] {};
  \draw[decoration={markings,mark=at position 1 with {\arrow[scale=1.5]{>}}},postaction={decorate},shorten >=0.5pt] (10) -- (11) node[anchor=south] {};
  \draw[decoration={markings,mark=at position 1 with {\arrow[scale=1.5]{>}}},postaction={decorate},shorten >=0.5pt] (4) -- (5) node[anchor=south] {};
  \draw[decoration={markings,mark=at position 1 with {\arrow[scale=1.5]{>}}},postaction={decorate},shorten >=0.5pt] (11) -- (12) node[anchor=south] {};
  \draw[decoration={markings,mark=at position 1 with {\arrow[scale=1.5]{>}}},postaction={decorate},shorten >=0.5pt] (5) -- (6) node[anchor=south] {$\partial^{E}_{1}$};
  \draw[decoration={markings,mark=at position 1 with {\arrow[scale=1.5]{>}}},postaction={decorate},shorten >=0.5pt] (12) -- (13) node[anchor=south] {$\partial^{X}_{1}$};
  \draw[decoration={markings,mark=at position 1 with {\arrow[scale=1.5]{>}}},postaction={decorate},shorten >=0.5pt] (6) -- (7) node[anchor=south] {};
  \draw[decoration={markings,mark=at position 1 with {\arrow[scale=1.5]{>}}},postaction={decorate},shorten >=0.5pt] (13) -- (14) node[anchor=south] {};
\end{tikzpicture}\]
Since $H_{i}(E)\cong H_{i}(X)$ for every $0 \leq i \leq s$, Lemma \ref{1.2.21} implies that $K_{s}$ is locally free on the punctured spectrum of $R$ with constant rank. Now apply Construction \ref{1.2.12} to the $R$-complex $E$ to get the following commutative diagram with exact columns:
\[\begin{tikzpicture}[every node/.style={midway}]
  \matrix[column sep={3em}, row sep={3em}]
  {\node(1) {}; & \node(2) {$0$}; & \node(3) {$0$}; & \node(4) {}; & \node(5) {$0$}; & \node(6) {$0$}; & \node(7) {};\\
  \node(8) {$0$}; & \node(9) {$K_{s}$}; & \node(10) {$E_{s-1}$}; & \node(11) {$\cdots$}; & \node(12) {$E_{1}$}; & \node(13) {$E_{0}$}; & \node(14) {$0$};\\
  \node(15) {$0$}; & \node(16) {$F_{s}$}; & \node(17) {$F_{s-1}$}; & \node(18) {$\cdots$}; & \node(19) {$F_{1}$}; & \node(20) {$F_{0}$}; & \node(21) {$0$};\\
  \node(22) {$0$}; & \node(23) {$C_{s}$}; & \node(24) {$C_{s-1}$}; & \node(25) {$\cdots$}; & \node(26) {$C_{1}$}; & \node(27) {$C_{0}$}; & \node(28) {$0$};\\
  \node(29) {}; & \node(30) {$0$}; & \node(31) {$0$}; & \node(32) {}; & \node(33) {$0$}; & \node(34) {$0$}; & \node(35) {};\\};
  \draw[decoration={markings,mark=at position 1 with {\arrow[scale=1.5]{>}}},postaction={decorate},shorten >=0.5pt] (2) -- (9) node[anchor=west] {};
  \draw[decoration={markings,mark=at position 1 with {\arrow[scale=1.5]{>}}},postaction={decorate},shorten >=0.5pt] (3) -- (10) node[anchor=west] {};
  \draw[decoration={markings,mark=at position 1 with {\arrow[scale=1.5]{>}}},postaction={decorate},shorten >=0.5pt] (5) -- (12) node[anchor=west] {};
  \draw[decoration={markings,mark=at position 1 with {\arrow[scale=1.5]{>}}},postaction={decorate},shorten >=0.5pt] (6) -- (13) node[anchor=west] {};
  \draw[decoration={markings,mark=at position 1 with {\arrow[scale=1.5]{>}}},postaction={decorate},shorten >=0.5pt] (9) -- (16) node[anchor=west] {};
  \draw[decoration={markings,mark=at position 1 with {\arrow[scale=1.5]{>}}},postaction={decorate},shorten >=0.5pt] (10) -- (17) node[anchor=west] {};
  \draw[decoration={markings,mark=at position 1 with {\arrow[scale=1.5]{>}}},postaction={decorate},shorten >=0.5pt] (12) -- (19) node[anchor=west] {};
  \draw[decoration={markings,mark=at position 1 with {\arrow[scale=1.5]{>}}},postaction={decorate},shorten >=0.5pt] (13) -- (20) node[anchor=west] {};
  \draw[decoration={markings,mark=at position 1 with {\arrow[scale=1.5]{>}}},postaction={decorate},shorten >=0.5pt] (6) -- (13) node[anchor=west] {};
  \draw[decoration={markings,mark=at position 1 with {\arrow[scale=1.5]{>}}},postaction={decorate},shorten >=0.5pt] (16) -- (23) node[anchor=west] {};
  \draw[decoration={markings,mark=at position 1 with {\arrow[scale=1.5]{>}}},postaction={decorate},shorten >=0.5pt] (17) -- (24) node[anchor=west] {};
  \draw[decoration={markings,mark=at position 1 with {\arrow[scale=1.5]{>}}},postaction={decorate},shorten >=0.5pt] (19) -- (26) node[anchor=west] {};
  \draw[decoration={markings,mark=at position 1 with {\arrow[scale=1.5]{>}}},postaction={decorate},shorten >=0.5pt] (20) -- (27) node[anchor=west] {};
  \draw[decoration={markings,mark=at position 1 with {\arrow[scale=1.5]{>}}},postaction={decorate},shorten >=0.5pt] (23) -- (30) node[anchor=west] {};
  \draw[decoration={markings,mark=at position 1 with {\arrow[scale=1.5]{>}}},postaction={decorate},shorten >=0.5pt] (24) -- (31) node[anchor=west] {};
  \draw[decoration={markings,mark=at position 1 with {\arrow[scale=1.5]{>}}},postaction={decorate},shorten >=0.5pt] (26) -- (33) node[anchor=west] {};
  \draw[decoration={markings,mark=at position 1 with {\arrow[scale=1.5]{>}}},postaction={decorate},shorten >=0.5pt] (27) -- (34) node[anchor=west] {};
  \draw[decoration={markings,mark=at position 1 with {\arrow[scale=1.5]{>}}},postaction={decorate},shorten >=0.5pt] (8) -- (9) node[anchor=south] {};
  \draw[decoration={markings,mark=at position 1 with {\arrow[scale=1.5]{>}}},postaction={decorate},shorten >=0.5pt] (9) -- (10) node[anchor=south] {$\overline{\partial^{E}_{s}}$};
  \draw[decoration={markings,mark=at position 1 with {\arrow[scale=1.5]{>}}},postaction={decorate},shorten >=0.5pt] (10) -- (11) node[anchor=south] {};
  \draw[decoration={markings,mark=at position 1 with {\arrow[scale=1.5]{>}}},postaction={decorate},shorten >=0.5pt] (11) -- (12) node[anchor=south] {};
  \draw[decoration={markings,mark=at position 1 with {\arrow[scale=1.5]{>}}},postaction={decorate},shorten >=0.5pt] (12) -- (13) node[anchor=south] {$\partial^{E}_{1}$};
  \draw[decoration={markings,mark=at position 1 with {\arrow[scale=1.5]{>}}},postaction={decorate},shorten >=0.5pt] (13) -- (14) node[anchor=south] {};
  \draw[decoration={markings,mark=at position 1 with {\arrow[scale=1.5]{>}}},postaction={decorate},shorten >=0.5pt] (15) -- (16) node[anchor=south] {};
  \draw[decoration={markings,mark=at position 1 with {\arrow[scale=1.5]{>}}},postaction={decorate},shorten >=0.5pt] (16) -- (17) node[anchor=south] {$\partial^{F}_{s}$};
  \draw[decoration={markings,mark=at position 1 with {\arrow[scale=1.5]{>}}},postaction={decorate},shorten >=0.5pt] (17) -- (18) node[anchor=south] {};
  \draw[decoration={markings,mark=at position 1 with {\arrow[scale=1.5]{>}}},postaction={decorate},shorten >=0.5pt] (18) -- (19) node[anchor=south] {};
  \draw[decoration={markings,mark=at position 1 with {\arrow[scale=1.5]{>}}},postaction={decorate},shorten >=0.5pt] (19) -- (20) node[anchor=south] {$\partial^{F}_{1}$};
  \draw[decoration={markings,mark=at position 1 with {\arrow[scale=1.5]{>}}},postaction={decorate},shorten >=0.5pt] (20) -- (21) node[anchor=south] {};
  \draw[decoration={markings,mark=at position 1 with {\arrow[scale=1.5]{>}}},postaction={decorate},shorten >=0.5pt] (22) -- (23) node[anchor=south] {};
  \draw[decoration={markings,mark=at position 1 with {\arrow[scale=1.5]{>}}},postaction={decorate},shorten >=0.5pt] (23) -- (24) node[anchor=south] {$\partial^{C}_{s}$};
  \draw[decoration={markings,mark=at position 1 with {\arrow[scale=1.5]{>}}},postaction={decorate},shorten >=0.5pt] (24) -- (25) node[anchor=south] {};
  \draw[decoration={markings,mark=at position 1 with {\arrow[scale=1.5]{>}}},postaction={decorate},shorten >=0.5pt] (25) -- (26) node[anchor=south] {};
  \draw[decoration={markings,mark=at position 1 with {\arrow[scale=1.5]{>}}},postaction={decorate},shorten >=0.5pt] (26) -- (27) node[anchor=south] {$\partial^{C}_{1}$};
  \draw[decoration={markings,mark=at position 1 with {\arrow[scale=1.5]{>}}},postaction={decorate},shorten >=0.5pt] (27) -- (28) node[anchor=south] {};
\end{tikzpicture}\]
where $F$ is an $R$-complex of free modules, and $C$ is an $R$-complex of Gorenstein projective modules. It follows from Theorem \ref{1.2.13} that $H_{i}(E)\cong H_{i}(F)$ and $H_{i}(C)=0$ for every $1 \leq i \leq s$, and there is a short exact sequence of $R$-modules
\begin{equation} \label{eq:1.2.22.1}
0 \rightarrow H_{0}(E) \rightarrow H_{0}(F) \rightarrow H_{0}(C) \rightarrow 0,
\end{equation}
where $H_{0}(C)$ is Gorenstein projective.
For any given $0 \leq i \leq s-1$, the exact column
$$0 \rightarrow E_{i} \rightarrow F_{i} \rightarrow C_{i} \rightarrow 0$$
in the above diagram shows that $\pd_{R}(C_{i}) < \infty$, so by Proposition \ref{1.2.6}, $\pd_{R}(C_{i})= \Gpd_{R}(C_{i})= 0$, i.e. $C_{i}$ is projective, hence free. Considering the short exact sequence
$$0 \rightarrow K_{s} \rightarrow F_{s} \rightarrow C_{s} \rightarrow 0,$$
Lemma \ref{1.2.21} implies that $C_{s}$ is locally free on the punctured spectrum of $R$ with constant rank. Therefore, considering the exact sequence
$$0 \rightarrow C_{s} \rightarrow C_{s-1} \rightarrow \cdots \rightarrow C_{1} \rightarrow C_{0} \rightarrow H_{0}(C) \rightarrow 0,$$
another application of Lemma \ref{1.2.21}, shows that $H_{0}(C)$ is locally free on the punctured spectrum of $R$ with constant rank.
Since $H_{0}(E) \cong H_{0}(X)$ has finite length, localizing the short exact sequence \eqref{eq:1.2.22.1} at any given $\mathfrak{p}\in \Spec (R)\backslash \{\mathfrak{m}\}$, yields $\left(H_{0}(F)\right)_{\mathfrak{p}} \cong \left(H_{0}(C)\right)_{\mathfrak{p}}$. Therefore, $H_{0}(F)$ is locally free on the punctured spectrum of $R$ with constant rank. Letting $C := H_{0}(C)$, we get from \eqref{eq:1.2.22.1}, the short exact sequence
$$0 \rightarrow H_{0}(X) \rightarrow H_{0}(F) \rightarrow C \rightarrow 0.$$
\end{prf}

Recall that an $R$-complex $X$ is said to be minimal if every homotopy equivalence $f:X \rightarrow X$ is an isomorphism. Given a bounded $R$-complex $X$ of of finitely generated Gorenstein projective modules, we can make the $R$-complex $F$ of Theorem \ref{1.2.13} minimal, at the cost of losing the injectivity of the morphism $\sigma: X \rightarrow F$.

We first need the following proposition.

\begin{proposition} \label{1.2.23}
Let $(R,\mathfrak{m})$ be a noetherian local ring, and let
$$F: 0 \rightarrow F_{s} \rightarrow F_{s-1} \rightarrow \cdots \rightarrow F_{1} \rightarrow F_{0} \rightarrow 0$$
be an $R$-complex of finitely generated free modules. Then there exists a decomposition $F=E\oplus X$ of $R$-complexes of finitely generated free $R$-modules, in which $E$ is minimal and $X$ is contractible. In particular, the projection $\eta: F \rightarrow E$ is a surjective quasi-isomorphism.
\end{proposition}

\begin{prf}
See \cite[Theorem 5.4.50]{CFH2}, noting that local rings are semi-perfect.
\end{prf}

\begin{theorem} \label{1.2.24}
Let $(R,\mathfrak{m})$ be a noetherian local ring, and let
$$X: 0 \rightarrow X_{s} \rightarrow X_{s-1} \rightarrow \cdots \rightarrow X_{1} \rightarrow X_{0} \rightarrow 0$$
be an $R$-complex of finitely generated Gorenstein projective modules. Then there exists a minimal $R$-complex
$$E: 0 \rightarrow E_{s} \rightarrow E_{s-1} \rightarrow \cdots \rightarrow E_{1} \rightarrow E_{0} \rightarrow 0$$
of finitely generated free modules, together with a morphism of complexes $\xi: X \rightarrow E$, with the property that $\xi$ induces isomorphisms $H_{i}(X)\cong H_{i}(E)$ for every $1 \leq i \leq s$, and a short exact sequence of $R$-modules
$$0 \rightarrow H_{0}(X) \xrightarrow {H_{0}(\xi)} H_{0}(E) \rightarrow D \rightarrow 0,$$
where $D$ is finitely generated Gorenstein projective.
\end{theorem}

\begin{prf}
For the $R$-complex $X$, choose an $R$-complex $F$ as in Theorem \ref{1.2.13}. On can easily observe that $F$ is degreewise finitely generated. We then have an injective morphism of $R$-complexes $\sigma: X \rightarrow F$, with the property that it induces isomorphisms $H_{i}(X)\cong H_{i}(F)$ for every $1 \leq i \leq s$, and a short exact sequence of $R$-modules
$$0 \rightarrow H_{0}(X) \xrightarrow {H_{0}(\sigma)} H_{0}(F) \rightarrow C \rightarrow 0,$$
where $C$ is finitely generated Gorenstein projective. By Proposition \ref{1.2.23}, there is a minimal $R$-complex
$$E: 0 \rightarrow E_{s} \rightarrow E_{s-1} \rightarrow \cdots \rightarrow E_{1} \rightarrow E_{0} \rightarrow 0$$
of finitely generated free modules, and a surjective quasi-isomorphism $\eta: F \rightarrow E$. Let $\xi=\eta\sigma: X \rightarrow E$. Then it is clear that $\xi$ induces isomorphisms $H_{i}(X)\cong H_{i}(E)$ for every $1 \leq i \leq s$. Further, $H_{0}(\xi)$ is injective. Let $D:= \coker H_{0}(\xi)$, and consider the following completed commutative diagram with exact rows:
\[\begin{tikzpicture}[every node/.style={midway},]
  \matrix[column sep={3em}, row sep={3em}]
  {\node(1) {$0$}; & \node(2) {$H_{0}(X)$}; & \node(3) {$H_{0}(F)$}; & \node(4) {$C$}; & \node(5) {$0$};\\
  \node(6) {$0$}; & \node(7) {$H_{0}(X)$}; & \node(8) {$H_{0}(E)$}; & \node(9) {$D$}; & \node(10) {$0$};\\};
  \draw[double distance=1.5pt] (2) -- (7) node[anchor=west] {};
  \draw[decoration={markings,mark=at position 1 with {\arrow[scale=1.5]{>}}},postaction={decorate},shorten >=0.5pt] (3) -- (8) node[anchor=east] {$\cong$} node[anchor=west] {$H_{0}(\eta)$};
  \draw[dashed,decoration={markings,mark=at position 1 with {\arrow[scale=1.5]{>}}},postaction={decorate},shorten >=0.5pt] (4) -- (9) node[anchor=west] {$\exists \psi$};
  \draw[decoration={markings,mark=at position 1 with {\arrow[scale=1.5]{>}}},postaction={decorate},shorten >=0.5pt] (1) -- (2) node[anchor=south] {};
  \draw[decoration={markings,mark=at position 1 with {\arrow[scale=1.5]{>}}},postaction={decorate},shorten >=0.5pt] (2) -- (3) node[anchor=south] {$H_{0}(\sigma)$};
  \draw[decoration={markings,mark=at position 1 with {\arrow[scale=1.5]{>}}},postaction={decorate},shorten >=0.5pt] (3) -- (4) node[anchor=south] {};
  \draw[decoration={markings,mark=at position 1 with {\arrow[scale=1.5]{>}}},postaction={decorate},shorten >=0.5pt] (4) -- (5) node[anchor=south] {};
  \draw[decoration={markings,mark=at position 1 with {\arrow[scale=1.5]{>}}},postaction={decorate},shorten >=0.5pt] (6) -- (7) node[anchor=south] {};
  \draw[decoration={markings,mark=at position 1 with {\arrow[scale=1.5]{>}}},postaction={decorate},shorten >=0.5pt] (7) -- (8) node[anchor=south] {$H_{0}(\xi)$};
  \draw[decoration={markings,mark=at position 1 with {\arrow[scale=1.5]{>}}},postaction={decorate},shorten >=0.5pt] (8) -- (9) node[anchor=south] {};
  \draw[decoration={markings,mark=at position 1 with {\arrow[scale=1.5]{>}}},postaction={decorate},shorten >=0.5pt](9) -- (10) node[anchor=south] {};
\end{tikzpicture}\]
The Short Five Lemma implies that $\psi$ is an isomorphism. Therefore $D$ is a finitely generated Gorenstein projective $R$-module.
\end{prf}

\section{Bounded Complexes of Gorenstein Injective Modules}

In this section, we develop the dual constructions to the previous section, to study bounded complexes of Gorenstein injective modules. As a result, we obtain an approximation for complexes of Gorenstein injective modules by complexes of injective modules. As in the previous section, we collect some definitions and facts for the convenience of the reader.

We begin with the definition of Gorenstein injective modules and Gorenstein injective dimension.

\begin{definition} \label{1.3.1}
Let $\mathcal{A}$ be a class of $R$-modules. An $R$-complex
$$X: \cdots \rightarrow X_{i+1} \xrightarrow{\partial^{X}_{i+1}} X_{i} \xrightarrow{\partial^{X}_{i}} X_{i-1} \rightarrow \cdots$$
is said to be \textit{$\Hom_{R}(\mathcal{A},-)$-exact} if the $R$-complex
$$\cdots \rightarrow \Hom_{R}(A,X_{i+1}) \xrightarrow{\Hom_{R}\left(A,\partial^{X}_{i+1}\right)} \Hom_{R}(A,X_{i}) \xrightarrow{\Hom_{R}\left(A,\partial^{X}_{i}\right)} \Hom_{R}(A,X_{i-1}) \rightarrow \cdots$$
is exact for every $A\in \mathcal{A}$.
\end{definition}

In the sequel, $\mathcal{I}$ will denote the class of injective $R$-modules.

\begin{definition} \label{1.3.2}
We have the following definitions:
\begin{enumerate}
\item[(i)] An $R$-module $J$ is called \textit{Gorenstein injective} \index{Gorenstein injective module} if there exists a $\Hom_{R}(\mathcal{I},-)$-exact exact $R$-complex
$$I: \cdots \rightarrow I_{2} \xrightarrow{\partial^{I}_{2}} I_{1} \xrightarrow{\partial^{I}_{1}} I_{0} \xrightarrow{\partial^{I}_{0}} I_{-1} \xrightarrow{\partial^{I}_{-1}} I_{-2} \rightarrow \cdots,$$
consisting of injective $R$-modules, such that $J\cong \im \partial^{I}_{1}$.
\item[(ii)] A \textit{Gorenstein injective resolution} \index{Gorenstein injective resolution} of an $R$-module $M$ is defined to be an $R$-complex
$$0 \rightarrow J_{0} \rightarrow J_{-1} \rightarrow J_{-2} \rightarrow \cdots,$$
where $J_{-i}$ is a Gorenstein injective $R$-module for every $i \geq 0$, and there is a homomorphism $\varepsilon: M\rightarrow J_{0}$ such that the augmented $R$-complex
$$0 \rightarrow M \xrightarrow{\varepsilon} J_{0} \rightarrow J_{-1} \rightarrow J_{-2} \rightarrow \cdots$$
is exact.
\item[(iii)]The \textit{Gorenstein injective dimension} \index{Gorenstein injective dimension} of a nonzero $R$-module $M$ is defined to be
$$\Gid_{R}(M):= \inf \left\{n \geq 0 \suchthat \begin{tabular}{ccc}
  \text{\textit{There is a Gorenstein injective resolution}}\\
 $0 \rightarrow J_{0} \rightarrow J_{-1} \rightarrow \cdots \rightarrow J_{-n} \rightarrow$ 0 \text{\textit{of}} $M$\\
  \end{tabular} \right\},$$
with the convention that $\inf \emptyset := \infty$. Further, we define $\Gid_{R}(0):=-\infty$.
\end{enumerate}
\end{definition}

The following facts on Gorenstein injective modules and Gorenstein injective dimension are required in the rest of the work.

\begin{proposition} \label{1.3.3}
Let $0\rightarrow M^{\prime}\rightarrow M\rightarrow M^{\prime\prime}\rightarrow 0$ be a short exact sequence of $R$-modules with $M^{\prime}$ Gorenstein injective. Then $M$ is Gorenstein injective if and only if $M^{\prime\prime}$ is Gorenstein injective.
\end{proposition}

\begin{prf}
\cite[Theorem 2.6]{Ho}.
\end{prf}

\begin{proposition} \label{1.3.4}
Let $R$ be a noetherian ring which is a homomorphic image of a Gorenstein ring with finite Krull dimension. If $J$ is a Gorenstein injective $R$-module, then $J_{\mathfrak{p}}$ is a Gorenstein injective $R_{\mathfrak{p}}$-module for every $\mathfrak{p}\in \Spec(R)$.
\end{proposition}

\begin{prf}
See \cite[Proposition 3.20]{CFH1}.
\end{prf}

\begin{proposition} \label{1.3.5}
Let $M$ be an $R$-module with $\Gid_{R}(M)<\infty$. Then we have
$$\Gid_{R}(M)= \sup \left\{i \geq 0 \suchthat \Ext_{R}^{i}(I,M)\neq 0 \text{ for some } I \in \mathcal{I}\right\}.$$
\end{proposition}

\begin{prf}
See \cite[Theorem 2.22]{Ho}.
\end{prf}

\begin{proposition} \label{1.3.6}
Let $M$ be an $R$-module. Then $\Gid_{R}(M)\leq \id_{R}(M)$ with equality if $\id_{R}(M)<\infty$.
\end{proposition}

\begin{prf}
See \cite[Proposition 3.10]{CFH1}.
\end{prf}

We further need the following definition.

\begin{definition} \label{1.3.7}
Let $\mathcal{A}$ be a class of $R$-modules. Define the \textit{right orthogonal class} \index{right orthogonal class} of $\mathcal{A}$, denoted by $\mathcal{A}^{\perp}$, to be the class of all $R$-modules $M$ with $\Ext^{1}_{R}(A,M)=0$ for every $A\in \mathcal{A}$.
\end{definition}

\begin{definition} \label{1.3.8}
Let $\mathcal{A}$ be a class of $R$-modules, and $M$ an $R$-module. Then we have the following definitions:
\begin{enumerate}
\item[(i)] By an \textit{$\mathcal{A}$-precover} \index{precover} of $M$, we mean a homomorphism $f:A\rightarrow M$ where $A\in \mathcal{A}$, with the property that the homomorphism
$$\Hom_{R}(B,f):\Hom_{R}(B,A)\rightarrow \Hom_{R}(B,M)$$
is surjective for every $B\in \mathcal{A}$.
\item[(ii)] An $\mathcal{A}$-precover $f:A\rightarrow M$ is said to be \textit{special} \index{special precover} whenever $f$ is surjective and $\ker f \in \mathcal{A}^{\perp}$.
\end{enumerate}
\end{definition}

The next proposition is the main ingredient of Construction \ref{1.3.12}.

\begin{proposition} \label{1.3.9}
Every Gorenstein injective $R$-module admits a special injective precover.
\end{proposition}

\begin{prf}
Let $M$ be a Gorenstein injective $R$-module. By definition, there exists a $\Hom_{R}(\mathcal{I},-)$-exact exact $R$-complex
\begin{equation} \label{eq:1.3.9.1}
I: \cdots \rightarrow I_{2} \xrightarrow{\partial^{I}_{2}} I_{1} \xrightarrow{\partial^{I}_{1}} I_{0} \xrightarrow{\partial^{I}_{0}} I_{-1} \xrightarrow{\partial^{I}_{-1}} I_{-2} \rightarrow \cdots,
\end{equation}
consisting of injective modules, such that $M\cong \im \partial^{I}_{1}$. Let $\theta:M\rightarrow \im \partial^{I}_{1}$ be any isomorphism, and let $\pi=\theta^{-1}\overline{\partial^{I}_{1}}$ where $\overline{\partial^{I}_{1}}:I_{1}\rightarrow \im \partial^{I}_{1}$ is given by $\overline{\partial^{I}_{1}}(x)=\partial^{I}_{1}(x)$ for every $x\in I_{1}$.
Therefore, there is a short exact sequence
$$0\rightarrow K \rightarrow I_{1} \xrightarrow {\pi} M \rightarrow 0,$$
where
$$K=\ker \pi=\ker \partial^{I}_{1} \cong \im \partial^{I}_{2}$$
is Gorenstein injective by the symmetry Definition \ref{1.3.2}.
We show that the $R$-homomorphism $\pi:I_{1}\rightarrow M$ is a special injective precover of $M$.
Given any injective $R$-module $I$, we need to argue that the $R$-homomorphism
$$\Hom_{R}(I,\pi):\Hom_{R}(I,I_{1})\rightarrow \Hom_{R}(I,M)$$
is surjective.
There is an exact sequence
$$\Hom_{R}(I,I_{1})\xrightarrow{\Hom_{R}(I,\pi)} \Hom_{R}(I,M) \rightarrow \Ext^{1}_{R}(I,K)=0,$$
where the vanishing is due to the fact that $K$ is Gorenstein injective and $I$ is injective. This established the claim.
\end{prf}

We further need the pullback construction.

\begin{definition} \label{1.3.10}
Let $f:M\rightarrow Y$ and $g:N\rightarrow Y$ be two $R$-homomorphisms. Define the \textit{pullback} \index{pullback} of the pair $(f,g)$ to be the $R$-module
$$M\sqcap_{Y}N := \left\{(m,n)\in M \oplus N \suchthat f(m)=g(n)\right\}.$$
Further, complete the pair $(f,g)$ to the commutative diagram
\[\begin{tikzpicture}[every node/.style={midway}]
  \matrix[column sep={2.5em}, row sep={2.5em}]
  {\node(1) {$M\sqcap_{Y}N$}; & \node(2) {$N$}; \\
  \node(3) {$M$}; & \node (4) {$Y$};\\};
  \draw[decoration={markings,mark=at position 1 with {\arrow[scale=1.5]{>}}},postaction={decorate},shorten >=0.5pt] (1) -- (2) node[anchor=south] {$f^{\prime}$};
  \draw[decoration={markings,mark=at position 1 with {\arrow[scale=1.5]{>}}},postaction={decorate},shorten >=0.5pt] (3) -- (4) node[anchor=south] {$f$};
  \draw[decoration={markings,mark=at position 1 with {\arrow[scale=1.5]{>}}},postaction={decorate},shorten >=0.5pt] (1) -- (3) node[anchor=west] {$g^{\prime}$};
  \draw[decoration={markings,mark=at position 1 with {\arrow[scale=1.5]{>}}},postaction={decorate},shorten >=0.5pt] (2) -- (4) node[anchor=west] {$g$};
\end{tikzpicture}\]
where $f^{\prime}(m,n):= n$ and $g^{\prime}(m,n):= m$.
\end{definition}

The following proposition describes the most important property of the pullback construction.

\begin{proposition} \label{1.3.11}
Let $f:M\rightarrow Y$ and $g:N\rightarrow Y$ be two $R$-homomorphisms. Then there exists a commutative diagram with exact rows and columns as follows:
\[\begin{tikzpicture}[every node/.style={midway}]
  \matrix[column sep={2.5em}, row sep={2.5em}]
  {\node(a) {}; & \node(b) {}; & \node(c) {$0$}; & \node(d) {$0$};\\
  \node(e) {}; & \node (f) {}; & \node(g) {$\ker g^{\prime}$}; & \node(h) {$\ker g$};\\
  \node(i) {$0$}; & \node (j) {$\ker f^{\prime}$}; & \node(k) {$M\sqcap_{Y}N$}; & \node (l) {$N$};\\
  \node(m) {$0$}; & \node (n) {$\ker f$}; & \node(o) {$M$}; & \node (p) {$Y$};\\};
  \draw[decoration={markings,mark=at position 1 with {\arrow[scale=1.5]{>}}},postaction={decorate},shorten >=0.5pt] (g) -- (h) node[anchor=south] {$\widetilde{f^{\prime}}$};
  \draw[decoration={markings,mark=at position 1 with {\arrow[scale=1.5]{>}}},postaction={decorate},shorten >=0.5pt] (i) -- (j) node[anchor=south] {};
  \draw[decoration={markings,mark=at position 1 with {\arrow[scale=1.5]{>}}},postaction={decorate},shorten >=0.5pt] (j) -- (k) node[anchor=south] {};
  \draw[decoration={markings,mark=at position 1 with {\arrow[scale=1.5]{>}}},postaction={decorate},shorten >=0.5pt] (k) -- (l) node[anchor=south] {$f^{\prime}$};
  \draw[decoration={markings,mark=at position 1 with {\arrow[scale=1.5]{>}}},postaction={decorate},shorten >=0.5pt] (m) -- (n) node[anchor=south] {};
  \draw[decoration={markings,mark=at position 1 with {\arrow[scale=1.5]{>}}},postaction={decorate},shorten >=0.5pt] (n) -- (o) node[anchor=south] {};
  \draw[decoration={markings,mark=at position 1 with {\arrow[scale=1.5]{>}}},postaction={decorate},shorten >=0.5pt] (o) -- (p) node[anchor=south] {$f$};
  \draw[decoration={markings,mark=at position 1 with {\arrow[scale=1.5]{>}}},postaction={decorate},shorten >=0.5pt] (c) -- (g) node[anchor=west] {};
  \draw[decoration={markings,mark=at position 1 with {\arrow[scale=1.5]{>}}},postaction={decorate},shorten >=0.5pt] (d) -- (h) node[anchor=east] {};
  \draw[decoration={markings,mark=at position 1 with {\arrow[scale=1.5]{>}}},postaction={decorate},shorten >=0.5pt] (g) -- (k) node[anchor=east] {};
  \draw[decoration={markings,mark=at position 1 with {\arrow[scale=1.5]{>}}},postaction={decorate},shorten >=0.5pt] (h) -- (l) node[anchor=west] {};
  \draw[decoration={markings,mark=at position 1 with {\arrow[scale=1.5]{>}}},postaction={decorate},shorten >=0.5pt] (j) -- (n) node[anchor=west] {$\widetilde{g^{\prime}}$};
  \draw[decoration={markings,mark=at position 1 with {\arrow[scale=1.5]{>}}},postaction={decorate},shorten >=0.5pt] (k) -- (o) node[anchor=west] {$g^{\prime}$};
  \draw[decoration={markings,mark=at position 1 with {\arrow[scale=1.5]{>}}},postaction={decorate},shorten >=0.5pt] (l) -- (p) node[anchor=west] {$g$};
\end{tikzpicture}\]
where $\widetilde{f^{\prime}}$ and $\widetilde{g^{\prime}}$ are induced by $f^{\prime}$ and $g^{\prime}$, respectively. Furthermore, the following assertions hold:
\begin{enumerate}
\item[(i)] If $f$ is surjective, then so is $f^{\prime}$.
\item[(ii)] If $g$ is surjective, then so is $g^{\prime}$.
\item[(iii)] Both $\widetilde{f^{\prime}}$ and $\widetilde{g^{\prime}}$ are isomorphisms.
\end{enumerate}
\end{proposition}

\begin{prf}
(i): Let $n\in N$. Since $f$ is surjective, there is an $m\in M$ such that $f(m)=g(n)$. Hence $(m,n)\in M\sqcap_{Y}N$, and $f^{\prime}(m,n)=n$.

(ii): similar to (i).

(iii): Let $(m,n)\in \ker g^{\prime}$, and suppose that $\widetilde{f^{\prime}}(m,n)=f^{\prime}(m,n)=n=0$. Further, $g^{\prime}(m,n)=m=0$. It follows that $(m,n)=0$, so $\widetilde{f^{\prime}}$ is injective. On the other hand, let $n\in \ker g$. Then $g(n)=0=f(0)$, so $(0,n)\in M\sqcap_{Y}N$. Moreover, $(0,n)\in \ker g^{\prime}$ and $\widetilde{f^{\prime}}(0,n)=f^{\prime}(0,n)=n$. Hence $\widetilde{f^{\prime}}$ is surjective. It follows that $\widetilde{f^{\prime}}$ is an isomorphism. Similarly, one can see that $\widetilde{g^{\prime}}$ is an isomorphism.
\end{prf}

The following construction is dual to Construction \ref{1.2.12}.

\begin{construction} \label{1.3.12}
Let
$$Y: 0\rightarrow Y_{s} \xrightarrow {\partial^{Y}_{s}} Y_{s-1}\rightarrow \cdots \rightarrow Y_{1} \xrightarrow {\partial^{Y}_{1}} Y_{0} \rightarrow 0,$$
be an $R$-complex of Gorenstein injective modules. The proof of Proposition \ref{1.3.9} gives a short exact sequence
$$0 \rightarrow \ker \pi_{0} \rightarrow I_{0} \xrightarrow {\pi_{0}} Y_{0} \rightarrow 0,$$
where $I_{0}$ is an injective $R$-module, and $\ker \pi_{0}$ is a Gorenstein injective $R$-module.
Construct the following pullback diagram as in Proposition \ref{1.3.11}:
\[\begin{tikzpicture}[every node/.style={midway}]
  \matrix[column sep={2.5em}, row sep={2.5em}]
  {\node(a) {}; & \node(b) {}; & \node(c) {$0$}; & \node(d) {$0$};\\
  \node(e) {}; & \node (f) {}; & \node(g) {$\ker \theta_{1}$}; & \node(h) {$\ker \pi_{0}$};\\
  \node(i) {$0$}; & \node (j) {$\ker \eta_{1}$}; & \node(k) {$L_{1}$}; & \node (l) {$I_{0}$};\\
  \node(m) {$0$}; & \node (n) {$\ker \partial^{Y}_{1}$}; & \node(o) {$Y_{1}$}; & \node (p) {$Y_{0}$};\\
  \node(q) {}; & \node (r) {}; & \node(s) {$0$}; & \node (t) {$0$};\\};
  \draw[decoration={markings,mark=at position 1 with {\arrow[scale=1.5]{>}}},postaction={decorate},shorten >=0.5pt] (g) -- (h) node[anchor=south] {$\widetilde{\eta_{1}}$};
  \draw[decoration={markings,mark=at position 1 with {\arrow[scale=1.5]{>}}},postaction={decorate},shorten >=0.5pt] (i) -- (j) node[anchor=south] {};
  \draw[decoration={markings,mark=at position 1 with {\arrow[scale=1.5]{>}}},postaction={decorate},shorten >=0.5pt] (j) -- (k) node[anchor=south] {$\iota_{1}$};
  \draw[decoration={markings,mark=at position 1 with {\arrow[scale=1.5]{>}}},postaction={decorate},shorten >=0.5pt] (k) -- (l) node[anchor=south] {$\eta_{1}$};
  \draw[decoration={markings,mark=at position 1 with {\arrow[scale=1.5]{>}}},postaction={decorate},shorten >=0.5pt] (m) -- (n) node[anchor=south] {};
  \draw[decoration={markings,mark=at position 1 with {\arrow[scale=1.5]{>}}},postaction={decorate},shorten >=0.5pt] (n) -- (o) node[anchor=south] {};
  \draw[decoration={markings,mark=at position 1 with {\arrow[scale=1.5]{>}}},postaction={decorate},shorten >=0.5pt] (o) -- (p) node[anchor=south] {$\partial^{Y}_{1}$};
  \draw[decoration={markings,mark=at position 1 with {\arrow[scale=1.5]{>}}},postaction={decorate},shorten >=0.5pt] (c) -- (g) node[anchor=west] {};
  \draw[decoration={markings,mark=at position 1 with {\arrow[scale=1.5]{>}}},postaction={decorate},shorten >=0.5pt] (d) -- (h) node[anchor=east] {};
  \draw[decoration={markings,mark=at position 1 with {\arrow[scale=1.5]{>}}},postaction={decorate},shorten >=0.5pt] (g) -- (k) node[anchor=east] {};
  \draw[decoration={markings,mark=at position 1 with {\arrow[scale=1.5]{>}}},postaction={decorate},shorten >=0.5pt] (h) -- (l) node[anchor=west] {};
  \draw[decoration={markings,mark=at position 1 with {\arrow[scale=1.5]{>}}},postaction={decorate},shorten >=0.5pt] (j) -- (n) node[anchor=west] {$\widetilde{\theta_{1}}$};
  \draw[decoration={markings,mark=at position 1 with {\arrow[scale=1.5]{>}}},postaction={decorate},shorten >=0.5pt] (k) -- (o) node[anchor=west] {$\theta_{1}$};
  \draw[decoration={markings,mark=at position 1 with {\arrow[scale=1.5]{>}}},postaction={decorate},shorten >=0.5pt] (l) -- (p) node[anchor=west] {$\pi_{0}$};
  \draw[decoration={markings,mark=at position 1 with {\arrow[scale=1.5]{>}}},postaction={decorate},shorten >=0.5pt] (o) -- (s) node[anchor=west] {};
  \draw[decoration={markings,mark=at position 1 with {\arrow[scale=1.5]{>}}},postaction={decorate},shorten >=0.5pt] (p) -- (t) node[anchor=west] {};
\end{tikzpicture}\]
where $L_{1}:= Y_{1} \sqcap_{Y_{0}} I_{0}$. Note that the surjectivity of $\theta_{1}$ follows from that of $\pi_{0}$ in view of Proposition \ref{1.3.11}. Since $\widetilde{\eta_{1}}$ is an isomorphism, $\ker \theta_{1}$ is Gorenstein injective, and Proposition \ref{1.3.3} implies that $L_{1}$ is Gorenstein injective.
Let $\psi_{2}$ be the composition
$$Y_{2} \xrightarrow {\overline{\partial^{Y}_{2}}} \ker \partial^{Y}_{1} \xrightarrow {\widetilde{\theta_{1}}^{-1}} \ker {\eta}_{1} \xrightarrow {\iota_{1}} L_{1},$$
where $\overline{\partial^{Y}_{2}}$ is induced by $\partial^{Y}_{2}$ in the obvious way.
One can easily check that the following sequence is an $R$-complex:
$$0\rightarrow Y_{s} \xrightarrow {\partial^{Y}_{s}} Y_{s-1} \rightarrow \cdots \rightarrow Y_{3} \xrightarrow {\partial^{Y}_{3}} Y_{2} \xrightarrow {\psi_{2}} L_{1} \xrightarrow {\eta_{1}} I_{0} \rightarrow 0$$
Continuing this construction, we get the following commutative diagram:
\[\begin{tikzpicture}[every node/.style={midway}]
  \matrix[column sep={2.5em}, row sep={2.5em}]
  {\node(1) {$I_{s}$}; & \node(2) {}; & \node(3) {}; & \node(4) {}; & \node(5) {};\\
  \node(6) {$L_{s}$}; & \node(7) {$I_{s-1}$}; & \node(8) {}; & \node(9) {}; & \node(10) {};\\
  \node(11) {$Y_{s}$}; & \node(12) {$L_{s-1}$}; & \node(13) {}; & \node(14) {}; & \node(15) {};\\
  \node(16) {}; & \node(17) {}; & \node(18) {$L_{2}$}; & \node(19) {$I_{1}$}; & \node(20) {};\\
  \node(21) {}; & \node(22) {}; & \node(23) {$Y_{2}$}; & \node(24) {$L_{1}$}; & \node(25) {$I_{0}$};\\
  \node(26) {}; & \node(27) {}; & \node(28) {}; & \node(29) {$Y_{1}$}; & \node(30) {$Y_{0}$};\\};
  \draw[decoration={markings,mark=at position 1 with {\arrow[scale=1.5]{>}}},postaction={decorate},shorten >=0.5pt] (1) -- (6) node[anchor=west] {$\pi_{s}$};
  \draw[decoration={markings,mark=at position 1 with {\arrow[scale=1.5]{>}}},postaction={decorate},shorten >=0.5pt] (6) -- (11) node[anchor=west] {$\theta_{s}$};
  \draw[decoration={markings,mark=at position 1 with {\arrow[scale=1.5]{>}}},postaction={decorate},shorten >=0.5pt] (7) -- (12) node[anchor=west] {$\pi_{s-1}$};
  \draw[decoration={markings,mark=at position 1 with {\arrow[scale=1.5]{>}}},postaction={decorate},shorten >=0.5pt] (18) -- (23) node[anchor=west] {$\theta_{2}$};
  \draw[decoration={markings,mark=at position 1 with {\arrow[scale=1.5]{>}}},postaction={decorate},shorten >=0.5pt] (19) -- (24) node[anchor=west] {$\pi_{1}$};
  \draw[decoration={markings,mark=at position 1 with {\arrow[scale=1.5]{>}}},postaction={decorate},shorten >=0.5pt] (24) -- (29) node[anchor=west] {$\theta_{1}$};
  \draw[decoration={markings,mark=at position 1 with {\arrow[scale=1.5]{>}}},postaction={decorate},shorten >=0.5pt] (25) -- (30) node[anchor=west] {$\pi_{0}$};
  \draw[decoration={markings,mark=at position 1 with {\arrow[scale=1.5]{>}}},postaction={decorate},shorten >=0.5pt] (6) -- (7) node[anchor=south] {$\eta_{s}$};
  \draw[decoration={markings,mark=at position 1 with {\arrow[scale=1.5]{>}}},postaction={decorate},shorten >=0.5pt] (11) -- (12) node[anchor=south] {$\psi_{s}$};
  \draw[decoration={markings,mark=at position 1 with {\arrow[scale=1.5]{>}}},postaction={decorate},shorten >=0.5pt] (18) -- (19) node[anchor=south] {$\eta_{2}$};
  \draw[decoration={markings,mark=at position 1 with {\arrow[scale=1.5]{>}}},postaction={decorate},shorten >=0.5pt] (23) -- (24) node[anchor=south] {$\psi_{2}$};
  \draw[decoration={markings,mark=at position 1 with {\arrow[scale=1.5]{>}}},postaction={decorate},shorten >=0.5pt] (24) -- (25) node[anchor=south] {$\eta_{1}$};
  \draw[decoration={markings,mark=at position 1 with {\arrow[scale=1.5]{>}}},postaction={decorate},shorten >=0.5pt] (29) -- (30) node[anchor=south] {$\partial^{Y}_{1}$};
  \path (12) -- (18) node [midway, sloped] {$\dots$};
\end{tikzpicture}\]
where $L_{i}:= Y_{i} \sqcap_{L_{i-1}} I_{i-1}$ for every $2 \leq i \leq s$, and the homomorphisms $\theta_{i}$, $\pi_{i}$, $\psi_{i}$, and $\eta_{i}$ are defined analogously. One can see by inspection that $\partial^{Y}_{i}=\theta_{i-1}\psi_{i}$ for every $2 \leq i \leq s$. Moreover, we let $\partial^{I}_{i}:=\eta_{i}\pi_{i}$ for every $1 \leq i \leq s$. We have
$$\partial^{I}_{i-1}\partial^{I}_{i}=\eta_{i-1}\pi_{i-1}\eta_{i}\pi_{i}=\eta_{i-1}\psi_{i}\theta_{i}\pi_{i}=0,$$
since $\eta_{i-1}\psi_{i}=0$.
Hence $I$ is an $R$-complex of injective modules. Besides, we let $\sigma_{0}:=\pi_{0}$, and $\sigma_{i}:=\theta_{i}\pi_{i}$ for every $1 \leq i \leq s$. Then
$$\sigma_{0}\partial^{I}_{1}=\pi_{0}\eta_{1}\pi_{1}= \partial^{Y}_{1}\theta_{1}\pi_{1}=\partial^{Y}_{1}\sigma_{1}$$
and
$$\sigma_{i-1}\partial^{I}_{i}=\theta_{i-1}\pi_{i-1}\eta_{i}\pi_{i}=\theta_{i-1}\psi_{i}\theta_{i}\pi_{i}=\partial^{Y}_{i}\sigma_{i}$$
for every $2 \leq i \leq s$.
As observed before, all the vertical homomorphisms are surjective. Thus $\sigma=(\sigma_{i})_{0 \leq i \leq s}: I \rightarrow Y$ is a surjective morphism of complexes as follows:
\[\begin{tikzpicture}[every node/.style={midway}]
  \matrix[column sep={2.5em}, row sep={2.5em}]
  {\node(1) {$0$}; & \node(2) {$I_{s}$}; & \node(3) {$I_{s-1}$}; & \node(4) {$\cdots$}; & \node(5) {$I_{1}$}; & \node(6) {$I_{0}$}; & \node(7) {$0$};\\
  \node(8) {$0$}; & \node(9) {$Y_{s}$}; & \node(10) {$Y_{s-1}$}; & \node(11) {$\cdots$}; & \node(12) {$Y_{1}$}; & \node(13) {$Y_{0}$}; & \node(14) {$0$};\\};
  \draw[decoration={markings,mark=at position 1 with {\arrow[scale=1.5]{>}}},postaction={decorate},shorten >=0.5pt] (2) -- (9) node[anchor=west] {$\sigma_{s}$};
  \draw[decoration={markings,mark=at position 1 with {\arrow[scale=1.5]{>}}},postaction={decorate},shorten >=0.5pt] (3) -- (10) node[anchor=west] {$\sigma_{s-1}$};
  \draw[decoration={markings,mark=at position 1 with {\arrow[scale=1.5]{>}}},postaction={decorate},shorten >=0.5pt] (5) -- (12) node[anchor=west] {$\sigma_{1}$};
  \draw[decoration={markings,mark=at position 1 with {\arrow[scale=1.5]{>}}},postaction={decorate},shorten >=0.5pt] (6) -- (13) node[anchor=west] {$\sigma_{0}$};
  \draw[decoration={markings,mark=at position 1 with {\arrow[scale=1.5]{>}}},postaction={decorate},shorten >=0.5pt] (1) -- (2) node[anchor=south] {};
  \draw[decoration={markings,mark=at position 1 with {\arrow[scale=1.5]{>}}},postaction={decorate},shorten >=0.5pt] (8) -- (9) node[anchor=south] {};
  \draw[decoration={markings,mark=at position 1 with {\arrow[scale=1.5]{>}}},postaction={decorate},shorten >=0.5pt] (2) -- (3) node[anchor=south] {$\partial^{I}_{s}$};
  \draw[decoration={markings,mark=at position 1 with {\arrow[scale=1.5]{>}}},postaction={decorate},shorten >=0.5pt] (9) -- (10) node[anchor=south] {$\partial^{Y}_{s}$};
  \draw[decoration={markings,mark=at position 1 with {\arrow[scale=1.5]{>}}},postaction={decorate},shorten >=0.5pt] (3) -- (4) node[anchor=south] {};
  \draw[decoration={markings,mark=at position 1 with {\arrow[scale=1.5]{>}}},postaction={decorate},shorten >=0.5pt] (10) -- (11) node[anchor=south] {};
  \draw[decoration={markings,mark=at position 1 with {\arrow[scale=1.5]{>}}},postaction={decorate},shorten >=0.5pt] (4) -- (5) node[anchor=south] {};
  \draw[decoration={markings,mark=at position 1 with {\arrow[scale=1.5]{>}}},postaction={decorate},shorten >=0.5pt] (11) -- (12) node[anchor=south] {};
  \draw[decoration={markings,mark=at position 1 with {\arrow[scale=1.5]{>}}},postaction={decorate},shorten >=0.5pt] (5) -- (6) node[anchor=south] {$\partial^{I}_{1}$};
  \draw[decoration={markings,mark=at position 1 with {\arrow[scale=1.5]{>}}},postaction={decorate},shorten >=0.5pt] (12) -- (13) node[anchor=south] {$\partial^{Y}_{1}$};
  \draw[decoration={markings,mark=at position 1 with {\arrow[scale=1.5]{>}}},postaction={decorate},shorten >=0.5pt] (6) -- (7) node[anchor=south] {};
  \draw[decoration={markings,mark=at position 1 with {\arrow[scale=1.5]{>}}},postaction={decorate},shorten >=0.5pt] (13) -- (14) node[anchor=south] {};
\end{tikzpicture}\]
Let $K:= \ker \sigma$. By the above discussion, $K_{0}= \ker \sigma_{0}= \ker \pi_{0}$ is Gorenstein injective. Now fix $1 \leq i \leq s$. There exists a short exact sequence of $R$-modules
$$0\rightarrow \ker \pi_{i} \rightarrow \ker \sigma_{i} \xrightarrow {\pi_{i}|_{\ker \sigma_{i}}} \ker \theta_{i} \rightarrow 0.$$
By the above construction, $\ker \pi_{i}$ is Gorenstein injective, and thus $\ker \theta_{i}\cong \ker \pi_{i-1}$ is Gorenstein injective. It follows from Proposition \ref{1.3.3} that $K_{i}:= \ker \sigma_{i}$ is Gorenstein injective. We can then form the short exact sequence of complexes
$$0\rightarrow K \rightarrow I \xrightarrow {\sigma} Y \rightarrow 0,$$
where $K$ is an $R$-complex of Gorenstein injective modules.
\end{construction}

We can even say more about homology modules.

\begin{theorem} \label{1.3.13}
Let
$$Y: 0 \rightarrow Y_{s} \rightarrow Y_{s-1} \rightarrow \cdots \rightarrow Y_{1} \rightarrow Y_{0} \rightarrow 0$$
be an $R$-complex of Gorenstein injective modules. Then there exist $R$-complexes
$$I: 0 \rightarrow I_{s} \rightarrow I_{s-1} \rightarrow \cdots \rightarrow I_{1} \rightarrow I_{0} \rightarrow 0,$$
consisting of injective modules, and
$$K: 0 \rightarrow K_{s} \rightarrow K_{s-1} \rightarrow \cdots \rightarrow K_{1} \rightarrow K_{0} \rightarrow 0,$$
consisting of Gorenstein injective modules, that fit into a short exact sequence of $R$-complexes
$$0\rightarrow K \rightarrow I \xrightarrow{\sigma} Y \rightarrow 0.$$
In addition, the following assertions hold:
\begin{enumerate}
\item[(i)] The morphism $\sigma:I \rightarrow Y$ induces isomorphisms $H_{i}(I)\cong H_{i}(Y)$ for every $0\leq i \leq s-1$. As a result, $H_{i}(K)=0$ for every $0\leq i \leq s-1$.
\item[(ii)] There is a short exact sequence of $R$-modules
$$0\rightarrow H_{s}(K) \rightarrow H_{s}(I) \xrightarrow {H_{s}(\sigma)} H_{s}(Y) \rightarrow 0,$$
where $H_{s}(K)$ is Gorenstein injective.
\item[(iii)] If in particular, $Y$ is an $R$-complex of injective modules, then $K$ is also an $R$-complex of injective modules. Furthermore, $H_{s}(K)$ is injective, so in particular, the short exact sequence in (ii) splits.
\item[(iv)] If $Y_{i}$ happens to be injective for every $0 \leq i \leq s-1$, then $K_{i}$ is injective for every $0 \leq i \leq s-1$, and $H_{s}(K)$ is isomorphic to a direct summand of $K_{s}$ with an injective complement.
\item[(v)] If $H_{i}(Y)=0$ for every $0 \leq i \leq s-1$, then $\id_{R} \left(H_{s}(I)\right) \leq s$.
\end{enumerate}
\end{theorem}

\begin{prf}
(i) and (ii): We apply Construction \ref{1.3.12} to the $R$-complex $Y$ to get the desired short exact sequence of $R$-complexes
$$0\rightarrow K \rightarrow I \xrightarrow {\sigma} Y \rightarrow 0.$$
We now focus on homology modules, sticking to the notations of Construction \ref{1.3.12}.

First, we consider $H_{0}(\sigma):H_{0}(I)\rightarrow H_{0}(Y)$. Suppose that
$$H_{0}(\sigma)\left(x+ \im \partial^{I}_{1}\right)=\sigma_{0}(x)+ \im \partial^{Y}_{1} = \pi_{0}(x) + \im \partial^{Y}_{1}=0,$$
for some $x\in I_{0}$. Then there is an element $y\in Y_{1}$ such that $\pi_{0}(x)= \partial^{Y}_{1}(y)$. By the construction of the pullback, we have $(y,x)\in L_{1}$. Since $\pi_{1}$ is surjective, there is an element $z\in I_{1}$ such that $(y,x)=\pi_{1}(z)$. Hence $$x=\eta_{1}(y,x)=\eta_{1}\left(\pi_{1}(z)\right)=\partial^{I}_{1}(z),$$
i.e. $x\in \im \partial^{I}_{1}$. Therefore, $H_{0}(\sigma)$ is injective. Now let $y+ \im \partial^{Y}_{1} \in H_{0}(Y)$. Since $\sigma_{0}$ is surjective, there is an element $x\in I_{0}$ such that $y= \sigma_{0}(x)$. Thus
$$H_{0}(\sigma)\left(x+ \im \partial^{I}_{1}\right)= \sigma_{0}(x)+\im \partial^{Y}_{1}=y+ \im \partial^{Y}_{1},$$
showing that $H_{0}(\sigma)$ is surjective, and so an isomorphism.

Next, we consider $H_{i}(\sigma):H_{i}(I)\rightarrow H_{i}(Y)$ for any given $1 \leq i \leq s$. Let $y+ \im \partial^{Y}_{i+1}\in H_{i}(Y)$. Thus $\partial^{Y}_{i}(y)=0$, so by definition of $\psi_{i}$, $$\psi_{i}(y)=\iota_{i-1}\left(\widetilde{\theta_{i-1}}^{-1}\left(\overline{\partial^{Y}_{i}}(y)\right)\right)=0.$$
It follows that $(y,0)\in L_{i}$. Now by the surjectivity of $\pi_{i}$, there is an element $x\in I_{i}$ such that $(y,0)=\pi_{i}(x)$. Therefore, $$\partial^{I}_{i}(x)= \eta_{i}\left(\pi_{i}(x)\right)= \eta_{i}(y,0)=0,$$
and
\begin{equation*}
\begin{split}
H_{i}(\sigma)\left(x+ \im \partial^{I}_{i+1} \right) & = \sigma_{i}(x)+ \im \partial^{Y}_{i+1} \\
 & = \theta_{i}\left(\pi_{i}(x)\right)+ \im \partial^{Y}_{i+1} \\
 & = \theta_{i}(y,0)+ \im \partial^{Y}_{i+1} \\
 & = y+ \im \partial^{Y}_{i+1}. \\
\end{split}
\end{equation*}
This shows that $H_{i}(\sigma)$ is surjective. Now suppose that $1 \leq i \leq s-1$, and
$$H_{i}(\sigma)\left(x+ \im \partial^{I}_{i+1}\right)= \sigma_{i}(x)+ \im \partial^{Y}_{i+1} = \theta_{i}\left(\pi_{i}(x)\right)+ \im \partial^{Y}_{i+1}=0,$$
for some $x\in \ker \partial^{I}_{i}$.
Then
$$\theta_{i}\left(\pi_{i}(x)\right)=\partial^{Y}_{i+1}(t)=\theta_{i}\left(\psi_{i+1}(t)\right)$$
for some $t\in Y_{i+1}$. Let $\pi_{i}(x)=(m,n)\in L_{i}$. Now the hypothesis $x\in \ker \partial^{I}_{i}$ implies that $$0=\partial^{I}_{i}(x)=\eta_{i}\left(\pi_{i}(x)\right)=\eta_{i}(m,n)=n,$$
so $\pi_{i}(x)=(m,0)$. On the other hand, letting $\psi_{i+1}(t)=(m^{\prime},n^{\prime})\in L_{i}$, we get $$0=\eta_{i}\left(\psi_{i+1}(t)\right)=\eta_{i}(m^{\prime},n^{\prime})=n^{\prime},$$
so $\psi_{i+1}(t)=(m^{\prime},0)$. Further,
$$m=\theta_{i}(m,n)=\theta_{i}\left(\pi_{i}(x)\right)=\theta_{i}\left(\psi_{i+1}(t)\right)=\theta_{i}(m^{\prime},n^{\prime})=m^{\prime}.$$
Therefore,
$$\pi_{i}(x)=(m,0)=(m^{\prime},0)= \psi_{i+1}(t).$$
It follows that $(t,x)\in L_{i+1}$. But by the surjectivity of $\pi_{i+1}$, there is an element $z\in I_{i+1}$ such that $(t,x)=\pi_{i+1}(z)$. It then follows that
$$\partial^{I}_{i+1}(z)=\eta_{i+1}\left(\pi_{i+1}(z)\right)=\eta_{i+1}(t,x)=x,$$
i.e. $x\in \im \partial^{I}_{i+1}$. Therefore, $H_{i}(\sigma)$ is injective, and thus an isomorphism.

Consider the short exact sequence
$$0\rightarrow K \rightarrow I \xrightarrow {\sigma} Y \rightarrow 0$$
of $R$-complexes, and its associated long exact sequence
$$0\rightarrow H_{s}(K) \rightarrow H_{s}(I) \xrightarrow {H_{s}(\sigma)} H_{s}(Y) \rightarrow H_{s-1}(K) \rightarrow H_{s-1}(I) \xrightarrow {H_{s-1}(\sigma)} H_{s-1}(Y) \rightarrow \cdots \rightarrow $$$$ H_{0}(K) \rightarrow H_{0}(I) \xrightarrow {H_{0}(\sigma)} H_{0}(Y) \rightarrow 0,$$
of homology modules. For any given $0 \leq i \leq s-1$, the exact sequence
$$H_{i+1}(I) \xrightarrow {H_{i+1}(\sigma)} H_{i+1}(Y) \xrightarrow {\vartheta} H_{i}(K) \rightarrow H_{i}(I) \xrightarrow {H_{i}(\sigma)} H_{i}(Y)$$
shows that
$$0= \coker H_{i+1}(\sigma)\cong \im \vartheta$$
and
$$\coker \vartheta \cong \ker H_{i}(\sigma)=0,$$
hence $H_{i}(K)=0$.
Thus we get a short exact sequence
$$0\rightarrow H_{s}(K) \rightarrow H_{s}(I) \xrightarrow {H_{s}(\sigma)} H_{s}(Y) \rightarrow 0.$$
We further notice that
$$H_{s}(K)\cong \ker H_{s}(\sigma) = \ker \pi_{s}.$$
In fact, let $x\in \ker \pi_{s}$, then $\partial^{I}_{s}(x)=\eta_{s}\left(\pi_{s}(x)\right)=0$ and $$H_{s}(\sigma)(x)=\sigma_{s}(x)=\theta_{s}\left(\pi_{s}(x)\right)=0,$$
so $x\in \ker H_{s}(\sigma)$.
Conversely, let $x\in \ker H_{s}(\sigma)$ and $\pi_{s}(x)=(m,n)$. Then
$$0=\partial^{I}_{s}(x)=\eta_{s}\left(\pi_{s}(x)\right)=\eta_{s}(m,n)=n,$$
and
$$0=H_{s}(\sigma)(x)=\sigma_{s}(x)=\theta_{s}\left(\pi_{s}(x)\right)=\theta_{s}(m,n)=m.$$
Thus $x\in \ker \pi_{s}$. It follows that $H_{s}(K)$ is Gorenstein injective by the construction.

(iii): We consider the following commutative diagram with exact columns:
\[\begin{tikzpicture}[every node/.style={midway}]
  \matrix[column sep={2.5em}, row sep={2.5em}]
  {\node(1) {}; & \node(2) {$0$}; & \node(3) {$0$}; & \node(4) {}; & \node(5) {$0$}; & \node(6) {$0$}; & \node(7) {};\\
  \node(8) {$0$}; & \node(9) {$K_{s}$}; & \node(10) {$K_{s-1}$}; & \node(11) {$\cdots$}; & \node(12) {$K_{1}$}; & \node(13) {$K_{0}$}; & \node(14) {$0$};\\
  \node(15) {$0$}; & \node(16) {$I_{s}$}; & \node(17) {$I_{s-1}$}; & \node(18) {$\cdots$}; & \node(19) {$I_{1}$}; & \node(20) {$I_{0}$}; & \node(21) {$0$};\\
  \node(22) {$0$}; & \node(23) {$Y_{s}$}; & \node(24) {$Y_{s-1}$}; & \node(25) {$\cdots$}; & \node(26) {$Y_{1}$}; & \node(27) {$Y_{0}$}; & \node(28) {$0$};\\
  \node(29) {}; & \node(30) {$0$}; & \node(31) {$0$}; & \node(32) {}; & \node(33) {$0$}; & \node(34) {$0$}; & \node(35) {};\\};
  \draw[decoration={markings,mark=at position 1 with {\arrow[scale=1.5]{>}}},postaction={decorate},shorten >=0.5pt] (2) -- (9) node[anchor=west] {};
  \draw[decoration={markings,mark=at position 1 with {\arrow[scale=1.5]{>}}},postaction={decorate},shorten >=0.5pt] (3) -- (10) node[anchor=west] {};
  \draw[decoration={markings,mark=at position 1 with {\arrow[scale=1.5]{>}}},postaction={decorate},shorten >=0.5pt] (5) -- (12) node[anchor=west] {};
  \draw[decoration={markings,mark=at position 1 with {\arrow[scale=1.5]{>}}},postaction={decorate},shorten >=0.5pt] (6) -- (13) node[anchor=west] {};
  \draw[decoration={markings,mark=at position 1 with {\arrow[scale=1.5]{>}}},postaction={decorate},shorten >=0.5pt] (9) -- (16) node[anchor=west] {};
  \draw[decoration={markings,mark=at position 1 with {\arrow[scale=1.5]{>}}},postaction={decorate},shorten >=0.5pt] (10) -- (17) node[anchor=west] {};
  \draw[decoration={markings,mark=at position 1 with {\arrow[scale=1.5]{>}}},postaction={decorate},shorten >=0.5pt] (12) -- (19) node[anchor=west] {};
  \draw[decoration={markings,mark=at position 1 with {\arrow[scale=1.5]{>}}},postaction={decorate},shorten >=0.5pt] (13) -- (20) node[anchor=west] {};
  \draw[decoration={markings,mark=at position 1 with {\arrow[scale=1.5]{>}}},postaction={decorate},shorten >=0.5pt] (6) -- (13) node[anchor=west] {};
  \draw[decoration={markings,mark=at position 1 with {\arrow[scale=1.5]{>}}},postaction={decorate},shorten >=0.5pt] (16) -- (23) node[anchor=west] {};
  \draw[decoration={markings,mark=at position 1 with {\arrow[scale=1.5]{>}}},postaction={decorate},shorten >=0.5pt] (17) -- (24) node[anchor=west] {};
  \draw[decoration={markings,mark=at position 1 with {\arrow[scale=1.5]{>}}},postaction={decorate},shorten >=0.5pt] (19) -- (26) node[anchor=west] {};
  \draw[decoration={markings,mark=at position 1 with {\arrow[scale=1.5]{>}}},postaction={decorate},shorten >=0.5pt] (20) -- (27) node[anchor=west] {};
  \draw[decoration={markings,mark=at position 1 with {\arrow[scale=1.5]{>}}},postaction={decorate},shorten >=0.5pt] (23) -- (30) node[anchor=west] {};
  \draw[decoration={markings,mark=at position 1 with {\arrow[scale=1.5]{>}}},postaction={decorate},shorten >=0.5pt] (24) -- (31) node[anchor=west] {};
  \draw[decoration={markings,mark=at position 1 with {\arrow[scale=1.5]{>}}},postaction={decorate},shorten >=0.5pt] (26) -- (33) node[anchor=west] {};
  \draw[decoration={markings,mark=at position 1 with {\arrow[scale=1.5]{>}}},postaction={decorate},shorten >=0.5pt] (27) -- (34) node[anchor=west] {};
  \draw[decoration={markings,mark=at position 1 with {\arrow[scale=1.5]{>}}},postaction={decorate},shorten >=0.5pt] (8) -- (9) node[anchor=south] {};
  \draw[decoration={markings,mark=at position 1 with {\arrow[scale=1.5]{>}}},postaction={decorate},shorten >=0.5pt] (9) -- (10) node[anchor=south] {$\partial^{K}_{s}$};
  \draw[decoration={markings,mark=at position 1 with {\arrow[scale=1.5]{>}}},postaction={decorate},shorten >=0.5pt] (10) -- (11) node[anchor=south] {};
  \draw[decoration={markings,mark=at position 1 with {\arrow[scale=1.5]{>}}},postaction={decorate},shorten >=0.5pt] (11) -- (12) node[anchor=south] {};
  \draw[decoration={markings,mark=at position 1 with {\arrow[scale=1.5]{>}}},postaction={decorate},shorten >=0.5pt] (12) -- (13) node[anchor=south] {$\partial^{K}_{1}$};
  \draw[decoration={markings,mark=at position 1 with {\arrow[scale=1.5]{>}}},postaction={decorate},shorten >=0.5pt] (13) -- (14) node[anchor=south] {};
  \draw[decoration={markings,mark=at position 1 with {\arrow[scale=1.5]{>}}},postaction={decorate},shorten >=0.5pt] (15) -- (16) node[anchor=south] {};
  \draw[decoration={markings,mark=at position 1 with {\arrow[scale=1.5]{>}}},postaction={decorate},shorten >=0.5pt] (16) -- (17) node[anchor=south] {$\partial^{I}_{s}$};
  \draw[decoration={markings,mark=at position 1 with {\arrow[scale=1.5]{>}}},postaction={decorate},shorten >=0.5pt] (17) -- (18) node[anchor=south] {};
  \draw[decoration={markings,mark=at position 1 with {\arrow[scale=1.5]{>}}},postaction={decorate},shorten >=0.5pt] (18) -- (19) node[anchor=south] {};
  \draw[decoration={markings,mark=at position 1 with {\arrow[scale=1.5]{>}}},postaction={decorate},shorten >=0.5pt] (19) -- (20) node[anchor=south] {$\partial^{I}_{1}$};
  \draw[decoration={markings,mark=at position 1 with {\arrow[scale=1.5]{>}}},postaction={decorate},shorten >=0.5pt] (20) -- (21) node[anchor=south] {};
  \draw[decoration={markings,mark=at position 1 with {\arrow[scale=1.5]{>}}},postaction={decorate},shorten >=0.5pt] (22) -- (23) node[anchor=south] {};
  \draw[decoration={markings,mark=at position 1 with {\arrow[scale=1.5]{>}}},postaction={decorate},shorten >=0.5pt] (23) -- (24) node[anchor=south] {$\partial^{Y}_{s}$};
  \draw[decoration={markings,mark=at position 1 with {\arrow[scale=1.5]{>}}},postaction={decorate},shorten >=0.5pt] (24) -- (25) node[anchor=south] {};
  \draw[decoration={markings,mark=at position 1 with {\arrow[scale=1.5]{>}}},postaction={decorate},shorten >=0.5pt] (25) -- (26) node[anchor=south] {};
  \draw[decoration={markings,mark=at position 1 with {\arrow[scale=1.5]{>}}},postaction={decorate},shorten >=0.5pt] (26) -- (27) node[anchor=south] {$\partial^{Y}_{1}$};
  \draw[decoration={markings,mark=at position 1 with {\arrow[scale=1.5]{>}}},postaction={decorate},shorten >=0.5pt] (27) -- (28) node[anchor=south] {};
\end{tikzpicture}\]
For any given $0 \leq i \leq s$, the exact column
$$0 \rightarrow K_{i} \rightarrow I_{i} \rightarrow Y_{i} \rightarrow 0$$
in the above diagram shows that $\id_{R}(K_{i})<\infty$, so Proposition \ref{1.3.6} implies that $\id_{R}(K_{i})= \Gid_{R}(K_{i})=0$, i.e. $K_{i}$ is injective. Consider the exact augmented $R$-complex
$$K^{+}: 0\rightarrow H_{s}(K) \rightarrow K_{s} \xrightarrow {\partial^{K}_{s}} K_{s-1} \rightarrow \cdots \rightarrow K_{1} \xrightarrow {\partial^{K}_{1}} K_{0} \rightarrow 0.$$
It shows that $\id_{R}\left(H_{s}(K)\right)<\infty$, so Proposition \ref{1.3.6} implies that
$$\id_{R}\left(H_{s}(K)\right)= \Gid_{R}\left(H_{s}(K)\right)=0,$$
i.e. $H_{s}(K)$ is injective.

(iv): We first note that if $s=0$, then the condition is void, and the assertion trivially holds. Now let $s\geq 1$. Since $Y_{i}$ is injective for every $0 \leq i \leq s-1$, as in the proof of part (iii), it follows that $K_{i}$ is injective for every $0 \leq i \leq s-1$. On the other hand, $K_{s}$ and $H_{s}(K)$ are Gorenstein injective, so using Proposition \ref{1.3.3}, the short exact sequence
\begin{equation} \label{eq:1.3.13.1}
0 \rightarrow H_{s}(K) \rightarrow K_{s} \rightarrow \im \partial^{K}_{s} \rightarrow 0
\end{equation}
implies that $\im \partial^{K}_{s}$ is Gorenstein injective. On the other hand, the exact sequence
$$0 \rightarrow \im \partial^{K}_{s} \rightarrow K_{s-1} \xrightarrow {\partial^{K}_{s-1}} K_{s-2} \rightarrow \cdots \rightarrow K_{1} \xrightarrow {\partial^{K}_{1}} K_{0} \rightarrow 0$$
shows that $\id_{R}\left(\im \partial^{K}_{s}\right)<\infty$, so Proposition \ref{1.3.6} implies that
$$\id_{R}\left(\im \partial^{K}_{s}\right)= \Gid_{R}\left(\im \partial^{K}_{s}\right)=0,$$
i.e. $\im \partial^{K}_{s}$ is injective. Therefore, $\Ext^{1}_{R}\left(\im \partial^{K}_{s},H_{s}(K)\right)=0$, so the short exact sequence \eqref{eq:1.3.13.1} splits and yields the assertions.

(v): If $s=0$, then $H_{s}(I)$ is indeed injective. Now let $s \geq 1$. Then $H_{i}(I) \cong H_{i}(Y)=0$ for every $0 \leq i \leq s-1$. Thus the exact sequence
$$0\rightarrow H_{s}(I) \rightarrow I_{s} \xrightarrow {\partial^{I}_{s}} I_{s-1} \rightarrow \cdots \rightarrow I_{1} \xrightarrow {\partial^{I}_{1}} I_{0} \rightarrow 0$$
shows that $\id_{R} \left(H_{s}(I)\right) \leq s$.
\end{prf}

We immediately draw the following conclusion.

\begin{corollary} \label{1.3.14}
Let $M$ be an $R$-module with $\Gid_{R}(M)< \infty$. Then there exists a short exact sequence of $R$-modules
$$0 \rightarrow L \rightarrow T \rightarrow M \rightarrow 0,$$
where $\id_{R}(T)= \Gid_{R}(M)$, and $L$ is Gorenstein injective.
\end{corollary}

\begin{prf}
Let $\Gid_{R}(M)=s$. The case $s=0$ is established in the proof of Proposition \ref{1.3.9}, so we may assume that $s\geq 1$. There exists an exact sequence
$$0 \rightarrow M\rightarrow Y_{s} \rightarrow Y_{s-1} \rightarrow \cdots \rightarrow Y_{1} \rightarrow Y_{0} \rightarrow 0,$$
where $Y_{i}$ is a Gorenstein injective $R$-module for every $0 \leq i \leq s$. Now consider the $R$-complex
$$Y: 0 \rightarrow Y_{s} \rightarrow Y_{s-1} \rightarrow \cdots \rightarrow Y_{1} \rightarrow Y_{0} \rightarrow 0,$$
and apply Theorem \ref{1.3.13} to get the complexes $I$ and $K$, and the short exact sequence of $R$-modules
$$0 \rightarrow H_{s}(K) \rightarrow H_{s}(I) \rightarrow H_{s}(Y) \rightarrow 0.$$
Obviously, $H_{s}(Y)\cong M$. Let $L= H_{s}(K)$ and $T= H_{s}(I)$.
Theorem \ref{1.3.13} implies that $\id_{R}(T) \leq s$ and $L$ is Gorenstein injective. By Proposition \ref{1.3.5}, we may choose an injective $R$-module $I$ for which $\Ext_{R}^{s}(I,M)\neq 0$.
If $\id_{R}(T) < s$, then in view of Theorem \ref{1.3.13}, we get the exact sequence
$$0= \Ext_{R}^{s}(I,T) \rightarrow \Ext_{R}^{s}(I,M) \rightarrow \Ext_{R}^{s+1}(I,L)=0,$$
which gives $\Ext_{R}^{s}(I,M)=0$, a contradiction. Hence $\id_{R}(T) = s$.
\end{prf}

\begin{remark} \label{1.3.15}
As mentioned in the Introduction, given an $R$-complex
$$Y: 0 \rightarrow Y_{s} \rightarrow Y_{s-1} \rightarrow  \cdots \rightarrow Y_{1} \rightarrow Y_{0} \rightarrow 0$$
of Gorenstein injective modules, we are interested in finding an $R$-complex
$$I: 0 \rightarrow I_{s} \rightarrow I_{s-1} \rightarrow \cdots \rightarrow I_{1} \rightarrow I_{0} \rightarrow 0$$
of injective modules, with the property that the homology modules of $I$ are as close as possible to the homology modules of $Y$. The best-case scenario is that $H_{i}(Y) \cong H_{i}(I)$ for every $0 \leq i \leq s$. But as the example below suggests, this cannot occur in general. Therefore Theorem \ref{1.3.13} gives, in a sense, the best approximation one could hope for.
\end{remark}

\begin{example} \label{1.3.16}
Let $(R,\mathfrak{m},k)$ be a non-regular Gorenstein local ring. We know that $\Gid_{R}(k) < \infty$. Let $s=\Gid_{R}(k)$, and take a Gorenstein injective resolution
$$Y:0 \rightarrow Y_{0} \rightarrow Y_{-1} \rightarrow \cdots \rightarrow Y_{-s} \rightarrow 0$$
for $k$. If there were an $R$-complex
$$I: 0 \rightarrow I_{0} \rightarrow I_{-1} \rightarrow \cdots \rightarrow I_{-s} \rightarrow 0$$
of injective modules with the property that $H_{i}(Y)\cong H_{i}(I)$ for every $0 \leq i \leq s$, then $I$
would be an injective resolution for $k$, implying that $\id_{R}(k)<\infty$. This contradicts the non-regularity of $R$.
\end{example}

Our next goal is to apply Construction \ref{1.3.12} to bounded complexes of Gorenstein injective modules whose homology modules have finite length.

We first need the following two propositions.

\begin{proposition} \label{1.3.17}
For every $R$-complex $Y$, there is a an $R$-complex of injective modules $I$, and a quasi-isomorphism $Y\xrightarrow {\simeq} I$, such that $I_{i}=0$ for every $i> \sup Y$.
\end{proposition}

\begin{prf}
See \cite[Theorem 5.3.26]{CFH2}.
\end{prf}

\begin{proposition} \label{1.3.18}
Let $Y$ be a homologically left-bounded $R$-complex, and $I$ be a left-bounded $R$-complex of Gorenstein injective modules such that $Y\simeq I$ in the derived category $\mathcal{D}(R)$. If $\Gid_{R}(Y)\leq s$, then $\ker \partial^{I}_{-s}$ is Gorenstein injective.
\end{proposition}

\begin{prf}
See \cite[Theorem 3.3]{CFrH}.
\end{prf}

We can now present the following construction.

\begin{construction} \label{1.3.19}
Let
$$Y: 0\rightarrow Y_{s} \xrightarrow {\partial^{Y}_{s}} Y_{s-1}\rightarrow \cdots \rightarrow Y_{1} \xrightarrow {\partial^{Y}_{1}} Y_{0} \rightarrow 0$$
be an $R$-complex of Gorenstein injective modules. By Proposition \ref{1.3.17}, there exists an $R$-complex $J$ of injective modules and a quasi-isomorphism $g:Y\xrightarrow {\simeq} J$. Furthermore, since $\sup Y \leq s$, we may choose $J_{i}=0$ for every $i>s$. We then have the following commutative diagram:
\[\begin{tikzpicture}[every node/.style={midway}]
\matrix[column sep={2.5em}, row sep={2.5em}]
  {\node(1) {$0$}; & \node(2) {$Y_{s}$}; & \node(3) {$Y_{s-1}$}; & \node(4) {$\cdots$}; & \node(5) {$Y_{1}$}; & \node(6) {$Y_{0}$}; & \node(7) {$0$}; & \node(8) {};\\
  \node(9) {$0$}; & \node(10) {$J_{s}$}; & \node(11) {$J_{s-1}$}; & \node(12) {$\cdots$}; & \node(13) {$J_{1}$}; & \node(14) {$J_{0}$}; & \node(15) {$J_{-1}$}; & \node(16) {$\cdots$};\\};
  \draw[decoration={markings,mark=at position 1 with {\arrow[scale=1.5]{>}}},postaction={decorate},shorten >=0.5pt] (1) -- (2) node[anchor=south] {};
  \draw[decoration={markings,mark=at position 1 with {\arrow[scale=1.5]{>}}},postaction={decorate},shorten >=0.5pt] (2) -- (3) node[anchor=south] {$\partial^{Y}_{s}$};
  \draw[decoration={markings,mark=at position 1 with {\arrow[scale=1.5]{>}}},postaction={decorate},shorten >=0.5pt] (3) -- (4) node[anchor=south] {};
  \draw[decoration={markings,mark=at position 1 with {\arrow[scale=1.5]{>}}},postaction={decorate},shorten >=0.5pt] (4) -- (5) node[anchor=south] {};
  \draw[decoration={markings,mark=at position 1 with {\arrow[scale=1.5]{>}}},postaction={decorate},shorten >=0.5pt] (5) -- (6) node[anchor=south] {$\partial^{Y}_{1}$};
  \draw[decoration={markings,mark=at position 1 with {\arrow[scale=1.5]{>}}},postaction={decorate},shorten >=0.5pt] (6) -- (7) node[anchor=south] {};
  \draw[decoration={markings,mark=at position 1 with {\arrow[scale=1.5]{>}}},postaction={decorate},shorten >=0.5pt] (9) -- (10) node[anchor=south] {};
  \draw[decoration={markings,mark=at position 1 with {\arrow[scale=1.5]{>}}},postaction={decorate},shorten >=0.5pt] (10) -- (11) node[anchor=south] {$\partial^{J}_{s}$};
  \draw[decoration={markings,mark=at position 1 with {\arrow[scale=1.5]{>}}},postaction={decorate},shorten >=0.5pt] (11) -- (12) node[anchor=south] {};
  \draw[decoration={markings,mark=at position 1 with {\arrow[scale=1.5]{>}}},postaction={decorate},shorten >=0.5pt] (12) -- (13) node[anchor=south] {};
  \draw[decoration={markings,mark=at position 1 with {\arrow[scale=1.5]{>}}},postaction={decorate},shorten >=0.5pt] (13) -- (14) node[anchor=south] {$\partial^{J}_{1}$};
  \draw[decoration={markings,mark=at position 1 with {\arrow[scale=1.5]{>}}},postaction={decorate},shorten >=0.5pt] (14) -- (15) node[anchor=south] {$\partial^{J}_{0}$};
  \draw[decoration={markings,mark=at position 1 with {\arrow[scale=1.5]{>}}},postaction={decorate},shorten >=0.5pt] (15) -- (16) node[anchor=south] {};
  \draw[decoration={markings,mark=at position 1 with {\arrow[scale=1.5]{>}}},postaction={decorate},shorten >=0.5pt] (7) -- (15) node[anchor=west] {$g_{-1}$};
  \draw[decoration={markings,mark=at position 1 with {\arrow[scale=1.5]{>}}},postaction={decorate},shorten >=0.5pt] (6) -- (14) node[anchor=west] {$g_{0}$};
  \draw[decoration={markings,mark=at position 1 with {\arrow[scale=1.5]{>}}},postaction={decorate},shorten >=0.5pt] (5) -- (13) node[anchor=west] {$g_{1}$};
  \draw[decoration={markings,mark=at position 1 with {\arrow[scale=1.5]{>}}},postaction={decorate},shorten >=0.5pt] (3) -- (11) node[anchor=west] {$g_{s-1}$};
  \draw[decoration={markings,mark=at position 1 with {\arrow[scale=1.5]{>}}},postaction={decorate},shorten >=0.5pt] (2) -- (10) node[anchor=west] {$g_{s}$};
\end{tikzpicture}\]
The above diagram shows that $\im g_{0} \subseteq \ker \partial^{J}_{0}$. Let $C_{0}:= \ker \partial^{J}_{0}$, and consider the softly truncated $R$-complex $J_{\supset 0}$. For convenience, we denote $J_{\supset 0}$ by $J$ in what follows. Then we obtain the following commutative diagram:
\[\begin{tikzpicture}[every node/.style={midway}]
  \matrix[column sep={2.5em}, row sep={2.5em}]
  {\node(1) {$0$}; & \node(2) {$Y_{s}$}; & \node(3) {$Y_{s-1}$}; & \node(4) {$\cdots$}; & \node(5) {$Y_{1}$}; & \node(6) {$Y_{0}$}; & \node(7) {$0$};\\
  \node(8) {$0$}; & \node(9) {$J_{s}$}; & \node(10) {$J_{s-1}$}; & \node(11) {$\cdots$}; & \node(12) {$J_{1}$}; & \node(13) {$C_{0}$}; & \node(14) {$0$};\\};
  \draw[decoration={markings,mark=at position 1 with {\arrow[scale=1.5]{>}}},postaction={decorate},shorten >=0.5pt] (2) -- (9) node[anchor=west] {$g_{s}$};
  \draw[decoration={markings,mark=at position 1 with {\arrow[scale=1.5]{>}}},postaction={decorate},shorten >=0.5pt] (3) -- (10) node[anchor=west] {$g_{s-1}$};
  \draw[decoration={markings,mark=at position 1 with {\arrow[scale=1.5]{>}}},postaction={decorate},shorten >=0.5pt] (5) -- (12) node[anchor=west] {$g_{1}$};
  \draw[decoration={markings,mark=at position 1 with {\arrow[scale=1.5]{>}}},postaction={decorate},shorten >=0.5pt] (6) -- (13) node[anchor=west] {$\overline{g_{0}}$};
  \draw[decoration={markings,mark=at position 1 with {\arrow[scale=1.5]{>}}},postaction={decorate},shorten >=0.5pt] (1) -- (2) node[anchor=south] {};
  \draw[decoration={markings,mark=at position 1 with {\arrow[scale=1.5]{>}}},postaction={decorate},shorten >=0.5pt] (8) -- (9) node[anchor=south] {};
  \draw[decoration={markings,mark=at position 1 with {\arrow[scale=1.5]{>}}},postaction={decorate},shorten >=0.5pt] (2) -- (3) node[anchor=south] {$\partial^{Y}_{s}$};
  \draw[decoration={markings,mark=at position 1 with {\arrow[scale=1.5]{>}}},postaction={decorate},shorten >=0.5pt] (9) -- (10) node[anchor=south] {$\partial^{J}_{s}$};
  \draw[decoration={markings,mark=at position 1 with {\arrow[scale=1.5]{>}}},postaction={decorate},shorten >=0.5pt] (3) -- (4) node[anchor=south] {};
  \draw[decoration={markings,mark=at position 1 with {\arrow[scale=1.5]{>}}},postaction={decorate},shorten >=0.5pt] (10) -- (11) node[anchor=south] {};
  \draw[decoration={markings,mark=at position 1 with {\arrow[scale=1.5]{>}}},postaction={decorate},shorten >=0.5pt] (4) -- (5) node[anchor=south] {};
  \draw[decoration={markings,mark=at position 1 with {\arrow[scale=1.5]{>}}},postaction={decorate},shorten >=0.5pt] (11) -- (12) node[anchor=south] {};
  \draw[decoration={markings,mark=at position 1 with {\arrow[scale=1.5]{>}}},postaction={decorate},shorten >=0.5pt] (5) -- (6) node[anchor=south] {$\partial^{Y}_{1}$};
  \draw[decoration={markings,mark=at position 1 with {\arrow[scale=1.5]{>}}},postaction={decorate},shorten >=0.5pt] (12) -- (13) node[anchor=south] {$\overline{\partial^{J}_{1}}$};
  \draw[decoration={markings,mark=at position 1 with {\arrow[scale=1.5]{>}}},postaction={decorate},shorten >=0.5pt] (6) -- (7) node[anchor=south] {};
  \draw[decoration={markings,mark=at position 1 with {\arrow[scale=1.5]{>}}},postaction={decorate},shorten >=0.5pt] (13) -- (14) node[anchor=south] {};
\end{tikzpicture}\]
where $\overline{\partial^{J}_{1}}:J_{1} \rightarrow C_{0}$ is given by $\overline{\partial^{J}_{1}}(x)= \partial^{J}_{1}(x)$ for every $x\in J_{1}$, and $\overline{g_{0}}:Y_{0}\rightarrow C_{0}$ is given by $\overline{g_{0}}(x)= g_{0}(x)$ for every $x\in Y_{0}$. Note that these homomorphisms are well-defined since $\im g_{0} \cup \im \partial^{J}_{1} \subseteq \ker \partial^{J}_{0}$. It is easily seen that the new $g:Y \rightarrow J$, obtained in the above diagram, is also a quasi-isomorphism, so that $H_{i}(Y)\cong H_{i}(J)$ for every $0 \leq i \leq s$. On the other hand, as $Y$ is an $R$-complex of Gorenstein injective modules, it is clear that $\Gid_{R} (Y) \leq s$, so Proposition \ref{1.3.18} implies that $C_{0}$ is Gorenstein injective.
\end{construction}

To study complexes of Gorenstein injective modules whose homology modules have finite length, we need the following notion.

\begin{definition} \label{1.3.20}
Let $(R,\mathfrak{m})$ be a local ring. An $R$-module $M$ is said to be \textit{locally injective on the punctured spectrum of} \index{locally injective on the punctured spectrum} $R$, if $M_{\mathfrak{p}}$ is an injective $R_{\mathfrak{p}}$-module for every $\mathfrak{p}\in \Spec(R)\backslash \{\mathfrak{m}\}$.
\end{definition}

The following lemma may be of independent interest.

\begin{lemma} \label{1.3.21}
Let $(R,\mathfrak{m})$ be a noetherian local ring which is a homomorphic image of a Gorenstein local ring, and let
$$Y: 0\rightarrow Y_{s} \rightarrow Y_{s-1} \rightarrow \cdots \rightarrow Y_{1} \rightarrow Y_{0} \rightarrow 0$$
be an $R$-complex whose homology modules have all finite length.
Then the following assertions hold:
\begin{enumerate}
\item[(i)] If $Y_{i}$ is locally injective on the punctured spectrum of $R$ for every $1 \leq i \leq s$, then $Y_{0}$ is locally injective on the punctured spectrum of $R$.
\item[(ii)] If $Y_{i}$ is locally injective on the punctured spectrum of $R$ for every $0 \leq i \leq s-1$, and $Y_{s}$ is Gorenstein injective, then $Y_{s}$ is locally injective on the punctured spectrum of $R$.
\end{enumerate}
\end{lemma}

\begin{prf}
Let $\mathfrak{p}\in \Spec(R)\backslash \{\mathfrak{m}\}$. As in the proof of Lemma \ref{1.2.21}, we get the exact sequence
\begin{equation} \label{eq:1.3.21.1}
Y_{\mathfrak{p}}: 0 \rightarrow (Y_{s})_{\mathfrak{p}} \rightarrow (Y_{s-1})_{\mathfrak{p}} \rightarrow \cdots \rightarrow (Y_{1})_{\mathfrak{p}} \rightarrow (Y_{0})_{\mathfrak{p}} \rightarrow 0.
\end{equation}

In the situation of part (i), $(Y_{i})_{\mathfrak{p}}$ is an injective $R_{\mathfrak{p}}$-module for every $1 \leq i \leq s$.
Extracting a short exact sequence
\begin{equation} \label{eq:1.3.21.2}
0\rightarrow (Y_{s})_{\mathfrak{p}} \rightarrow (Y_{s-1})_{\mathfrak{p}} \rightarrow C \rightarrow 0
\end{equation}
from $Y_{\mathfrak{p}}$, we are left with an exact sequence
\begin{equation} \label{eq:1.3.21.3}
0 \rightarrow C \rightarrow (Y_{s-2})_{\mathfrak{p}} \rightarrow \cdots \rightarrow (Y_{1})_{\mathfrak{p}} \rightarrow (Y_{0})_{\mathfrak{p}} \rightarrow 0.
\end{equation}
Since $(Y_{s})_{\mathfrak{p}}$ is injective, the short exact sequence \eqref{eq:1.3.21.2}, showing that $C$ is injective.
Continuing this way with \eqref{eq:1.3.21.3}, we conclude that $(Y_{0})_{\mathfrak{p}}$ is injective.

In the situation of part (ii), $(Y_{i})_{\mathfrak{p}}$ is an injective $R_{\mathfrak{p}}$-module for every $0 \leq i \leq s-1$.
Now the sequence \eqref{eq:1.3.21.1} shows that $\id_{R_{\mathfrak{p}}}\left((Y_{s})_{\mathfrak{p}}\right)<\infty$.
On the other hand, it follows from Proposition \ref{1.3.4} that $(Y_{s})_{\mathfrak{p}}$ is a Gorenstein injective $R_{\mathfrak{p}}$-module. Now Proposition \ref{1.3.6} implies that $\id_{R_{\mathfrak{p}}}\left((Y_{s})_{\mathfrak{p}}\right)= \Gid_{R_{\mathfrak{p}}}\left((Y_{s})_{\mathfrak{p}}\right)=0$, so $(Y_{s})_{\mathfrak{p}}$ is injective.
\end{prf}

We now present the anticipated theorem.

\begin{theorem} \label{1.3.22}
Let $(R,\mathfrak{m})$ be a noetherian local ring which is a homomorphic image of a Gorenstein local ring, and let
$$Y: 0\rightarrow Y_{s} \rightarrow Y_{s-1} \rightarrow \cdots \rightarrow Y_{1} \rightarrow Y_{0} \rightarrow 0$$
be an $R$-complex of Gorenstein injective modules whose homology modules have all finite length. Then there exists an $R$-complex
$$I: 0\rightarrow I_{s} \rightarrow I_{s-1} \rightarrow \cdots \rightarrow I_{1} \rightarrow I_{0} \rightarrow 0,$$
consisting of injective modules, with the following properties:
\begin{enumerate}
\item[(i)] $H_{i}(Y)\cong H_{i}(I)$ for every $0 \leq i \leq s-1$.
\item[(ii)] There is a short exact sequence
$$0 \rightarrow K \rightarrow H_{s}(I) \rightarrow H_{s}(Y) \rightarrow 0$$
of locally injective modules on the punctured spectrum of $R$, such that $K$ is Gorenstein injective.
\end{enumerate}
\end{theorem}

\begin{prf}
We apply Construction \ref{1.3.19} to the complex $Y$, to get a $R$-complex $J$ with $J_{i}$ injective for every $1 \leq i \leq s$ and $C_{0}$ Gorenstein injective, and a quasi-isomorphism $g:Y \rightarrow J$ as in the following diagram:
\[\begin{tikzpicture}[every node/.style={midway}]
  \matrix[column sep={2.5em}, row sep={2.5em}]
  {\node(1) {$0$}; & \node(2) {$Y_{s}$}; & \node(3) {$Y_{s-1}$}; & \node(4) {$\cdots$}; & \node(5) {$Y_{1}$}; & \node(6) {$Y_{0}$}; & \node(7) {$0$};\\
  \node(8) {$0$}; & \node(9) {$J_{s}$}; & \node(10) {$J_{s-1}$}; & \node(11) {$\cdots$}; & \node(12) {$J_{1}$}; & \node(13) {$C_{0}$}; & \node(14) {$0$};\\};
  \draw[decoration={markings,mark=at position 1 with {\arrow[scale=1.5]{>}}},postaction={decorate},shorten >=0.5pt] (2) -- (9) node[anchor=west] {$g_{s}$};
  \draw[decoration={markings,mark=at position 1 with {\arrow[scale=1.5]{>}}},postaction={decorate},shorten >=0.5pt] (3) -- (10) node[anchor=west] {$g_{s-1}$};
  \draw[decoration={markings,mark=at position 1 with {\arrow[scale=1.5]{>}}},postaction={decorate},shorten >=0.5pt] (5) -- (12) node[anchor=west] {$g_{1}$};
  \draw[decoration={markings,mark=at position 1 with {\arrow[scale=1.5]{>}}},postaction={decorate},shorten >=0.5pt] (6) -- (13) node[anchor=west] {$\overline{g_{0}}$};
  \draw[decoration={markings,mark=at position 1 with {\arrow[scale=1.5]{>}}},postaction={decorate},shorten >=0.5pt] (1) -- (2) node[anchor=south] {};
  \draw[decoration={markings,mark=at position 1 with {\arrow[scale=1.5]{>}}},postaction={decorate},shorten >=0.5pt] (8) -- (9) node[anchor=south] {};
  \draw[decoration={markings,mark=at position 1 with {\arrow[scale=1.5]{>}}},postaction={decorate},shorten >=0.5pt] (2) -- (3) node[anchor=south] {$\partial^{Y}_{s}$};
  \draw[decoration={markings,mark=at position 1 with {\arrow[scale=1.5]{>}}},postaction={decorate},shorten >=0.5pt] (9) -- (10) node[anchor=south] {$\partial^{J}_{s}$};
  \draw[decoration={markings,mark=at position 1 with {\arrow[scale=1.5]{>}}},postaction={decorate},shorten >=0.5pt] (3) -- (4) node[anchor=south] {};
  \draw[decoration={markings,mark=at position 1 with {\arrow[scale=1.5]{>}}},postaction={decorate},shorten >=0.5pt] (10) -- (11) node[anchor=south] {};
  \draw[decoration={markings,mark=at position 1 with {\arrow[scale=1.5]{>}}},postaction={decorate},shorten >=0.5pt] (4) -- (5) node[anchor=south] {};
  \draw[decoration={markings,mark=at position 1 with {\arrow[scale=1.5]{>}}},postaction={decorate},shorten >=0.5pt] (11) -- (12) node[anchor=south] {};
  \draw[decoration={markings,mark=at position 1 with {\arrow[scale=1.5]{>}}},postaction={decorate},shorten >=0.5pt] (5) -- (6) node[anchor=south] {$\partial^{Y}_{1}$};
  \draw[decoration={markings,mark=at position 1 with {\arrow[scale=1.5]{>}}},postaction={decorate},shorten >=0.5pt] (12) -- (13) node[anchor=south] {$\overline{\partial^{J}_{1}}$};
  \draw[decoration={markings,mark=at position 1 with {\arrow[scale=1.5]{>}}},postaction={decorate},shorten >=0.5pt] (6) -- (7) node[anchor=south] {};
  \draw[decoration={markings,mark=at position 1 with {\arrow[scale=1.5]{>}}},postaction={decorate},shorten >=0.5pt] (13) -- (14) node[anchor=south] {};
\end{tikzpicture}\]
Since $H_{i}(Y)\cong H_{i}(J)$ for every $0 \leq i \leq s$, Lemma \ref{1.3.21} implies that $C_{0}$ is locally injective on the punctured spectrum of $R$. Now apply Construction \ref{1.3.12} to the complex $J$ to get the following commutative diagram with exact columns:
\[\begin{tikzpicture}[every node/.style={midway}]
  \matrix[column sep={2.5em}, row sep={2.5em}]
  {\node(1) {}; & \node(2) {$0$}; & \node(3) {$0$}; & \node(4) {}; & \node(5) {$0$}; & \node(6) {$0$}; & \node(7) {};\\
  \node(8) {$0$}; & \node(9) {$K_{s}$}; & \node(10) {$K_{s-1}$}; & \node(11) {$\cdots$}; & \node(12) {$K_{1}$}; & \node(13) {$K_{0}$}; & \node(14) {$0$};\\
  \node(15) {$0$}; & \node(16) {$I_{s}$}; & \node(17) {$I_{s-1}$}; & \node(18) {$\cdots$}; & \node(19) {$I_{1}$}; & \node(20) {$I_{0}$}; & \node(21) {$0$};\\
  \node(22) {$0$}; & \node(23) {$J_{s}$}; & \node(24) {$J_{s-1}$}; & \node(25) {$\cdots$}; & \node(26) {$J_{1}$}; & \node(27) {$C_{0}$}; & \node(28) {$0$};\\
  \node(29) {}; & \node(30) {$0$}; & \node(31) {$0$}; & \node(32) {}; & \node(33) {$0$}; & \node(34) {$0$}; & \node(35) {};\\};
  \draw[decoration={markings,mark=at position 1 with {\arrow[scale=1.5]{>}}},postaction={decorate},shorten >=0.5pt] (2) -- (9) node[anchor=west] {};
  \draw[decoration={markings,mark=at position 1 with {\arrow[scale=1.5]{>}}},postaction={decorate},shorten >=0.5pt] (3) -- (10) node[anchor=west] {};
  \draw[decoration={markings,mark=at position 1 with {\arrow[scale=1.5]{>}}},postaction={decorate},shorten >=0.5pt] (5) -- (12) node[anchor=west] {};
  \draw[decoration={markings,mark=at position 1 with {\arrow[scale=1.5]{>}}},postaction={decorate},shorten >=0.5pt] (6) -- (13) node[anchor=west] {};
  \draw[decoration={markings,mark=at position 1 with {\arrow[scale=1.5]{>}}},postaction={decorate},shorten >=0.5pt] (9) -- (16) node[anchor=west] {};
  \draw[decoration={markings,mark=at position 1 with {\arrow[scale=1.5]{>}}},postaction={decorate},shorten >=0.5pt] (10) -- (17) node[anchor=west] {};
  \draw[decoration={markings,mark=at position 1 with {\arrow[scale=1.5]{>}}},postaction={decorate},shorten >=0.5pt] (12) -- (19) node[anchor=west] {};
  \draw[decoration={markings,mark=at position 1 with {\arrow[scale=1.5]{>}}},postaction={decorate},shorten >=0.5pt] (13) -- (20) node[anchor=west] {};
  \draw[decoration={markings,mark=at position 1 with {\arrow[scale=1.5]{>}}},postaction={decorate},shorten >=0.5pt] (6) -- (13) node[anchor=west] {};
  \draw[decoration={markings,mark=at position 1 with {\arrow[scale=1.5]{>}}},postaction={decorate},shorten >=0.5pt] (16) -- (23) node[anchor=west] {};
  \draw[decoration={markings,mark=at position 1 with {\arrow[scale=1.5]{>}}},postaction={decorate},shorten >=0.5pt] (17) -- (24) node[anchor=west] {};
  \draw[decoration={markings,mark=at position 1 with {\arrow[scale=1.5]{>}}},postaction={decorate},shorten >=0.5pt] (19) -- (26) node[anchor=west] {};
  \draw[decoration={markings,mark=at position 1 with {\arrow[scale=1.5]{>}}},postaction={decorate},shorten >=0.5pt] (20) -- (27) node[anchor=west] {};
  \draw[decoration={markings,mark=at position 1 with {\arrow[scale=1.5]{>}}},postaction={decorate},shorten >=0.5pt] (23) -- (30) node[anchor=west] {};
  \draw[decoration={markings,mark=at position 1 with {\arrow[scale=1.5]{>}}},postaction={decorate},shorten >=0.5pt] (24) -- (31) node[anchor=west] {};
  \draw[decoration={markings,mark=at position 1 with {\arrow[scale=1.5]{>}}},postaction={decorate},shorten >=0.5pt] (26) -- (33) node[anchor=west] {};
  \draw[decoration={markings,mark=at position 1 with {\arrow[scale=1.5]{>}}},postaction={decorate},shorten >=0.5pt] (27) -- (34) node[anchor=west] {};
  \draw[decoration={markings,mark=at position 1 with {\arrow[scale=1.5]{>}}},postaction={decorate},shorten >=0.5pt] (8) -- (9) node[anchor=south] {};
  \draw[decoration={markings,mark=at position 1 with {\arrow[scale=1.5]{>}}},postaction={decorate},shorten >=0.5pt] (9) -- (10) node[anchor=south] {$\partial^{K}_{s}$};
  \draw[decoration={markings,mark=at position 1 with {\arrow[scale=1.5]{>}}},postaction={decorate},shorten >=0.5pt] (10) -- (11) node[anchor=south] {};
  \draw[decoration={markings,mark=at position 1 with {\arrow[scale=1.5]{>}}},postaction={decorate},shorten >=0.5pt] (11) -- (12) node[anchor=south] {};
  \draw[decoration={markings,mark=at position 1 with {\arrow[scale=1.5]{>}}},postaction={decorate},shorten >=0.5pt] (12) -- (13) node[anchor=south] {$\partial^{K}_{1}$};
  \draw[decoration={markings,mark=at position 1 with {\arrow[scale=1.5]{>}}},postaction={decorate},shorten >=0.5pt] (13) -- (14) node[anchor=south] {};
  \draw[decoration={markings,mark=at position 1 with {\arrow[scale=1.5]{>}}},postaction={decorate},shorten >=0.5pt] (15) -- (16) node[anchor=south] {};
  \draw[decoration={markings,mark=at position 1 with {\arrow[scale=1.5]{>}}},postaction={decorate},shorten >=0.5pt] (16) -- (17) node[anchor=south] {$\partial^{I}_{s}$};
  \draw[decoration={markings,mark=at position 1 with {\arrow[scale=1.5]{>}}},postaction={decorate},shorten >=0.5pt] (17) -- (18) node[anchor=south] {};
  \draw[decoration={markings,mark=at position 1 with {\arrow[scale=1.5]{>}}},postaction={decorate},shorten >=0.5pt] (18) -- (19) node[anchor=south] {};
  \draw[decoration={markings,mark=at position 1 with {\arrow[scale=1.5]{>}}},postaction={decorate},shorten >=0.5pt] (19) -- (20) node[anchor=south] {$\partial^{I}_{1}$};
  \draw[decoration={markings,mark=at position 1 with {\arrow[scale=1.5]{>}}},postaction={decorate},shorten >=0.5pt] (20) -- (21) node[anchor=south] {};
  \draw[decoration={markings,mark=at position 1 with {\arrow[scale=1.5]{>}}},postaction={decorate},shorten >=0.5pt] (22) -- (23) node[anchor=south] {};
  \draw[decoration={markings,mark=at position 1 with {\arrow[scale=1.5]{>}}},postaction={decorate},shorten >=0.5pt] (23) -- (24) node[anchor=south] {$\partial^{J}_{s}$};
  \draw[decoration={markings,mark=at position 1 with {\arrow[scale=1.5]{>}}},postaction={decorate},shorten >=0.5pt] (24) -- (25) node[anchor=south] {};
  \draw[decoration={markings,mark=at position 1 with {\arrow[scale=1.5]{>}}},postaction={decorate},shorten >=0.5pt] (25) -- (26) node[anchor=south] {};
  \draw[decoration={markings,mark=at position 1 with {\arrow[scale=1.5]{>}}},postaction={decorate},shorten >=0.5pt] (26) -- (27) node[anchor=south] {$\partial^{J}_{1}$};
  \draw[decoration={markings,mark=at position 1 with {\arrow[scale=1.5]{>}}},postaction={decorate},shorten >=0.5pt] (27) -- (28) node[anchor=south] {};
\end{tikzpicture}\]
where $I$ is an $R$-complex of injective $R$-modules, and $K$ is an $R$-complex of Gorenstein injective modules. It follows from Theorem \ref{1.3.13} that $H_{i}(I)\cong H_{i}(J)$ and $H_{i}(K)=0$ for every $0 \leq i \leq s-1$, and there is a short exact sequence of $R$-modules
\begin{equation} \label{eq:1.3.22.1}
0 \rightarrow H_{s}(K) \rightarrow H_{s}(I) \rightarrow H_{s}(J) \rightarrow 0,
\end{equation}
where $H_{s}(K)$ is Gorenstein injective.
For any given $1 \leq i \leq s$, the exact column
$$0 \rightarrow K_{i} \rightarrow I_{i} \rightarrow J_{i} \rightarrow 0$$
in the above diagram shows that $\id_{R}(K_{i})<\infty$, so by Proposition \ref{1.3.6}, $\id_{R}(K_{i})= \Gid_{R}(K_{i})= 0$, i.e. $K_{i}$ is injective. Considering the short exact sequence
$$0 \rightarrow K_{0} \rightarrow I_{0} \rightarrow C_{0} \rightarrow 0,$$
Lemma \ref{1.3.21} implies that $K_{0}$ is locally injective on the punctured spectrum of $R$.
Therefore, considering the exact sequence
$$0 \rightarrow H_{s}(K) \rightarrow K_{s} \rightarrow K_{s-1} \rightarrow \cdots \rightarrow K_{1} \rightarrow K_{0} \rightarrow 0,$$
another application of Lemma \ref{1.3.21}, shows that $H_{s}(K)$ is locally injective on the punctured spectrum of $R$.
Since $H_{s}(J) \cong H_{s}(Y)$ has finite length, localizing the short exact sequence \eqref{eq:1.3.22.1} at any given $\mathfrak{p}\in \Spec (R)\backslash \{\mathfrak{m}\}$, yields $H_{s}(K)_{\mathfrak{p}} \cong H_{s}(I)_{\mathfrak{p}}$. Therefore, $H_{s}(I)$ is locally injective on the punctured spectrum of $R$. Letting $K=H_{s}(K)$, we get from \eqref{eq:1.3.22.1}, the short exact sequence
$$0 \rightarrow K \rightarrow H_{s}(I) \rightarrow H_{s}(Y) \rightarrow 0.$$
\end{prf}

Given a bounded $R$-complex $Y$ of of Gorenstein injective modules, we can make the complex $I$ of Theorem \ref{1.3.13} minimal, at the cost of losing the surjectivity of the morphism $\sigma: I \rightarrow Y$.

We first need the following proposition.

\begin{proposition} \label{1.3.23}
Let
$$I: 0 \rightarrow I_{s} \rightarrow I_{s-1} \rightarrow \cdots \rightarrow I_{1} \rightarrow I_{0} \rightarrow 0$$
be an $R$-complex of injective modules. Then there exists a decomposition $I=J \oplus X$ of $R$-complexes of injective modules, in which $J$ is minimal and $X$ is contractible. In particular, the injection $\eta: J \rightarrow I$ is an injective quasi-isomorphism.
\end{proposition}

\begin{prf}
See \cite[Theorem 5.4.17]{CFH2}.
\end{prf}

\begin{theorem} \label{1.3.24}
Let
$$Y: 0 \rightarrow Y_{s} \rightarrow Y_{s-1} \rightarrow \cdots \rightarrow Y_{1} \rightarrow Y_{0} \rightarrow 0$$
be an $R$-complex of Gorenstein injective modules. Then there exists a minimal $R$-complex
$$J: 0 \rightarrow J_{s} \rightarrow J_{s-1} \rightarrow \cdots \rightarrow J_{1} \rightarrow J_{0} \rightarrow 0$$
of injective modules, together with a morphism of $R$-complexes $\xi: J \rightarrow Y$, with the property that $\xi$ induces isomorphisms $H_{i}(Y)\cong H_{i}(J)$ for every $0 \leq i \leq s-1$, and a short exact sequence of $R$-modules
$$0 \rightarrow L \rightarrow H_{s}(J) \xrightarrow {H_{s}(\xi)} H_{s}(Y) \rightarrow 0,$$
where $L$ is Gorenstein injective.
\end{theorem}

\begin{prf}
For the $R$-complex $Y$, choose an $R$-complex $I$ as in Theorem \ref{1.3.13}. We then have a surjective morphism of $R$-complexes $\sigma: I \rightarrow Y$, with the property that it induces isomorphisms $H_{i}(Y)\cong H_{i}(I)$ for every $0 \leq i \leq s-1$, and a short exact sequence of $R$-modules
$$0 \rightarrow K \rightarrow H_{s}(I) \xrightarrow {H_{s}(\sigma)} H_{s}(Y) \rightarrow 0,$$
where $K$ is Gorenstein injective. By Proposition \ref{1.3.23}, there is a minimal $R$-complex
$$J: 0 \rightarrow J_{s} \rightarrow J_{s-1} \rightarrow \cdots \rightarrow J_{1} \rightarrow J_{0} \rightarrow 0$$
of injective $R$-modules, and an injective quasi-isomorphism $\eta: J \rightarrow I$. Let $\xi=\sigma\eta: J \rightarrow Y$.
Then it is clear that $\xi$ induces isomorphisms $H_{i}(J)\cong H_{i}(Y)$ for every $0 \leq i \leq s-1$. Further, $H_{s}(\xi)$ is surjective. Let $L= \ker H_{s}(\xi)$, and consider the following completed commutative diagram with exact rows:
\[\begin{tikzpicture}[every node/.style={midway},]
  \matrix[column sep={2.5em}, row sep={2.5em}]
  {\node(1) {$0$}; & \node(2) {$L$}; & \node(3) {$H_{s}(J)$}; & \node(4) {$H_{s}(Y)$}; & \node(5) {$0$};\\
  \node(6) {$0$}; & \node(7) {$K$}; & \node(8) {$H_{s}(I)$}; & \node(9) {$H_{s}(Y)$}; & \node(10) {$0$};\\};
  \draw[double distance=1.5pt] (4) -- (9) node[anchor=west] {};
  \draw[decoration={markings,mark=at position 1 with {\arrow[scale=1.5]{>}}},postaction={decorate},shorten >=0.5pt] (3) -- (8) node[anchor=east] {$\cong$} node[anchor=west] {$H_{s}(\eta)$};
  \draw[dashed,decoration={markings,mark=at position 1 with {\arrow[scale=1.5]{>}}},postaction={decorate},shorten >=0.5pt] (2) -- (7) node[anchor=west] {$\exists \psi$};
  \draw[decoration={markings,mark=at position 1 with {\arrow[scale=1.5]{>}}},postaction={decorate},shorten >=0.5pt] (1) -- (2) node[anchor=south] {};
  \draw[decoration={markings,mark=at position 1 with {\arrow[scale=1.5]{>}}},postaction={decorate},shorten >=0.5pt] (2) -- (3) node[anchor=south] {};
  \draw[decoration={markings,mark=at position 1 with {\arrow[scale=1.5]{>}}},postaction={decorate},shorten >=0.5pt] (3) -- (4) node[anchor=south] {$H_{s}(\xi)$};
  \draw[decoration={markings,mark=at position 1 with {\arrow[scale=1.5]{>}}},postaction={decorate},shorten >=0.5pt] (4) -- (5) node[anchor=south] {};
  \draw[decoration={markings,mark=at position 1 with {\arrow[scale=1.5]{>}}},postaction={decorate},shorten >=0.5pt] (6) -- (7) node[anchor=south] {};
  \draw[decoration={markings,mark=at position 1 with {\arrow[scale=1.5]{>}}},postaction={decorate},shorten >=0.5pt] (7) -- (8) node[anchor=south] {};
  \draw[decoration={markings,mark=at position 1 with {\arrow[scale=1.5]{>}}},postaction={decorate},shorten >=0.5pt] (8) -- (9) node[anchor=south] {$H_{s}(\sigma)$};
  \draw[decoration={markings,mark=at position 1 with {\arrow[scale=1.5]{>}}},postaction={decorate},shorten >=0.5pt] (9) -- (10) node[anchor=south] {};
\end{tikzpicture}\]
The Short Five Lemma implies that $\psi$ is an isomorphism. Therefore $L$ is a Gorenstein injective $R$-module.
\end{prf}

\section{Gorenstein New Intersection Conjecture}

The celebrated New Intersection Theorem is a remarkable theorem at the crossroad of homological algebra and commutative algebra. We state this theorem as follows:

\begin{theorem} \label{1.4.1}
Let $(R,\mathfrak{m})$ be a noetherian local ring. Let
$$F: 0 \rightarrow F_{s}\rightarrow F_{s-1}\rightarrow \cdots \rightarrow F_{1}\rightarrow F_{0}\rightarrow 0$$
be a non-exact $R$-complex of finitely generated free modules such that $\ell_{R}\left(H_{i}(F)\right)< \infty$ for every $0 \leq i \leq s$. Then $\dim(R)\leq s$.
\end{theorem}

In 1973, Peskine and Szpiro proved the New Intersection Theorem in prime characteristic, using the Frobenius functor. Their work ushered prime characteristic methods to the forefront of commutative algebra. By 1975, Hochster's work established a reduction to prime characteristic from equicharacteristic zero, to give a proof of the New Intersection Theorem in all equicharacteristic local rings. In 1987, Roberts proved this theorem for mixed characteristic rings using the machinery of Chern classes, settling the theorem completely in all cases. As a generalization of the New Intersection Theorem, we have the following Improved New Intersection Theorem:

\begin{theorem} \label{1.4.2}
Let $(R,\mathfrak{m})$ be a noetherian local ring. Let
$$F: 0\rightarrow F_{s}\rightarrow F_{s-1}\rightarrow \cdots \rightarrow F_{1}\rightarrow F_{0}\rightarrow 0$$
be an $R$-complex of finitely generated free modules such that $\ell_{R}\left(H_{i}(F)\right)< \infty$ for every $1 \leq i \leq s$ and there exists an element $e\in H_{0}(F)\backslash \mathfrak{m}H_{0}(F)$ with $\ell_{R}(Re)<\infty$. Then $\dim(R)\leq s$.
\end{theorem}

This theorem follows from the existence of a big Cohen-Macaulay module which has been recently established in all cases.

It is widely believed that most of the results in classical homological algebra enjoy their Gorenstein counterparts. Accordingly, it is natural to ask if there would possibly exist a Gorenstein analogue of the New Intersection Theorem, the so-called Gorenstein New Intersection Theorem. In this direction, we formulate a conjecture.

We first need a definition.

\begin {definition} \label{1.4.3}
Let $R$ be noetherian ring. A finitely generated $R$-module $M$ is called \textit{totally reflexive} \index{totally reflexive module} if the following conditions are satisfied:
\begin{enumerate}
\item [(i)] The biduality map $\mu^{M}:M\rightarrow M^{\ast\ast}$, given by $\mu^{M}(x)(f):=f(x)$ for every $f\in M^{\ast}$ and every $x\in M$, is an isomorphism.
\item [(ii)] $\Ext^{i}_{R}(M,R)=0=\Ext^{i}_{R}(M^{\ast},R)$ for every $i \geq 1$.
\end{enumerate}
\end{definition}

One can see that over a noetherian ring, a finitely generated module is totally reflexive if and only if it is Gorenstein projective.

We now pose the following conjecture:

\begin{conjecture} \label{1.4.4}
Let $(R,\mathfrak{m})$ be a noetherian local ring. Let
$$G: 0\rightarrow G_{s}\rightarrow G_{s-1}\rightarrow \cdots \rightarrow G_{1}\rightarrow G_{0}\rightarrow 0$$
be a non-exact $R$-complex of totally reflexive modules such that $\ell_{R}\left(H_{i}(G)\right)< \infty$ for every $0 \leq i \leq s$. Then $\dim(R)\leq s$.
\end{conjecture}

We may observe that since every finitely generated free module is totally reflexive, the above conjecture generalizes the classical New Intersection Theorem.

In this section, we plan to investigate this conjecture by deploying the constructions developed in Section 1.2.

\begin{Notation} \label{1.4.5}
Let $(R,\mathfrak{m})$ be a noetherian local ring. Denote by $\mathcal{LG}(R)$ the class of all totally reflexive $R$-modules which are also locally free on the punctured spectrum of $R$.
\end{Notation}

The class $\mathcal{LG}(R)$ trivially contains all finitely generated free $R$-modules. However, the following example indicates that the class $\mathcal{LG}(R)$ can be quite bigger than the class of finitely generated free $R$-modules even over a non-Cohen-Macaulay local ring.

\begin{example}  \label{1.4.6}
Let $k$ be a field, $R=k[[X,Y,Z]]/(X^{2},XY,Z^{2})$, and $M=R/(Z)$. Then $R$ is a non-Cohen-Macaulay local ring, and $M$ is a non-free totally reflexive $R$-module. Now it follows from \cite[Corollary 5.6]{Ta} that there exists a nonfree module $G\in \mathcal{LG}(R)$.
\end{example}

The next proposition sheds some light on the set of associated primes and support of modules in the class $\mathcal{LG}(R)$ with constant rank.

\begin{proposition}  \label{1.4.7}
Let $(R,\mathfrak{m})$ be a noetherian local ring, and $G\in \mathcal{LG}(R)$ with constant $r$. Then the following assertions hold:
\begin{enumerate}
\item [(i)] If $r=0$, then $\ell_{R}(G)< \infty$. In particular, $\Ass_{R}(G)= \Supp_{R}(G) \subseteq \{ \mathfrak{m} \}$.
\item [(ii)] If $r>0$, then $\Ass_{R}(G)=\Ass(R)$ and $\Supp_{R}(G)=\Spec(R)$.
\end{enumerate}
\end{proposition}

\begin{prf}
(i): Clear.

(ii): Take any $\mathfrak{p}\in \Spec(R)\backslash \{\mathfrak{m}\}$, so that $G_{\mathfrak{p}} \cong R^{r}_{\mathfrak{p}}$. Now the equality $\Supp_{R}(G)=\Spec(R)$ is clear. Further, we have:
$$\mathfrak{p}\in \Ass_{R}(G) \Longleftrightarrow \mathfrak{p}R_{\mathfrak{p}}\in \Ass_{R_{\mathfrak{p}}}\left(G_{\mathfrak{p}}\right) \Longleftrightarrow \mathfrak{p}R_{\mathfrak{p}}\in \Ass\left(R_{\mathfrak{p}}\right) \Longleftrightarrow \mathfrak{p}\in \Ass(R)$$
It follows that $\Ass_{R}(G)\backslash \{\mathfrak{m}\} = \Ass(R)\backslash \{\mathfrak{m}\}$. On the other hand, the Auslander-Bridger formula implies that $\depth_{R}(G) = \depth(R)$. Thus we have:
$$\mathfrak{m}\in \Ass_{R}(G) \Longleftrightarrow \depth_{R}(G) =0 \Longleftrightarrow \depth(R)=0 \Longleftrightarrow \mathfrak{m}\in \Ass(R)$$
Putting everything together, we get $\Ass_{R}(G)=\Ass(R)$.
\end{prf}

We next establish a variant of Theorem \ref{1.2.22}. We need a lemma.

\begin{lemma} \label{1.4.8}
Let $R$ be a noetherian ring, and $M$ a finitely generated $R$-module such that $M_{\mathfrak{p}}$ is a free $R_{\mathfrak{p}}$-module of constant rank $r$ for every $\mathfrak{p}\in \Ass(R)$. Then $M$ contains a free submodule $L$ of rank $r$ such that $M/L$ is a torsion $R$-module.
\end{lemma}

\begin{prf}
See \cite[Proposition 1.4.3]{BH}.
\end{prf}

\begin{proposition} \label{1.4.9}
Let $(R,\mathfrak{m})$ be a noetherian local ring with $\depth(R)>0$. Let
$$G: 0\rightarrow G_{s} \rightarrow G_{s-1}\rightarrow \cdots \rightarrow G_{1} \rightarrow G_{0} \rightarrow 0$$
be an $R$-complex of totally reflexive modules such that $\ell_{R}\left(H_{i}(G)\right)< \infty$ for every $0 \leq i \leq s$. Then there exists an $R$-complex
$$F: 0\rightarrow F_{s} \rightarrow F_{s-1}\rightarrow \cdots \rightarrow F_{1} \rightarrow F_{0} \rightarrow 0,$$
consisting of finitely generated free modules, and a nonzerodivisor $a\in R$, with the following properties:
\begin{enumerate}
\item [(i)] $\ell_{R}\left(H_{i}(F)\right)< \infty$ for every $1 \leq i \leq s$.
\item [(ii)] $H_{0}(F)$ is a finitely generated faithful $R$-module which is locally free on the punctured spectrum of $R$ with constant rank.
\item [(iii)] $\Gamma_{\mathfrak{m}}\left(H_{0}(F)\right)=(0:_{H_{0}(F)}a)$.
\item [(iv)] $aH_{0}(F)\in \mathcal{LG}(R)$.
\end{enumerate}
\end{proposition}

\begin{prf}
(i) and (ii): We apply Theorem \ref{1.2.23} to the $R$-complex $G$ to obtain an $R$-complex $F$, consisting of finitely generated free modules, with $H_{i}(G)\cong H_{i}(F)$ for every $1 \leq i \leq s$, and a short exact sequence of $R$-modules
\begin{equation} \label{eq:1.4.9.1}
0\rightarrow H_{0}(G) \xrightarrow {\iota} H_{0}(F) \rightarrow C \rightarrow 0,
\end{equation}
where $H_{0}(F)$ is locally free on the punctured spectrum of $R$ with constant rank and $C\in \mathcal{LG}(R)$. Since $\ell_{R}\left(H_{0}(G)\right)<\infty$, localizing the short exact sequence \eqref{eq:1.4.9.1} at any given $\mathfrak{p}\in \Spec(R)\backslash \{\mathfrak{m}\}$, implies that $H_{0}(F)$ and $C$ have the same constant rank, say $r$. But $\depth_{R}(C)=\depth(R)>0$, so by Proposition \ref{1.4.7}, we must have $r>0$. Therefore, We have
$$\ann_{R}\left(H_{0}(F)\right)_{\mathfrak{p}}= \ann_{R_{\mathfrak{p}}}\left(H_{0}(F)_{\mathfrak{p}}\right)=0,$$
so $\ann_{R}\left(H_{0}(F)\right)$ is an artinian submodule of $R$.
If $\ann_{R}\left(H_{0}(F)\right)\neq 0$, then
$$\{\mathfrak{m}\}=\Ass\left(\ann_{R}\left(H_{0}(F)\right)\right)\subseteq \Ass(R),$$
which contradicts the hypothesis $\depth(R)> 0$. Therefore, $\ann_{R}\left(H_{0}(F)\right)=0$, i.e. $H_{0}(F)$ is faithful.

(iii) and (iv): We note that since $\mathfrak{m}\notin \Ass(R)$, $H_{0}(F)_{\mathfrak{p}}$ is a free $R_{\mathfrak{p}}$-module of constant rank $r$ for every $\mathfrak{p}\in \Ass(R)$. By Lemma \ref{1.4.8}, $H_{0}(F)$ has a free submodule $L$ of rank $r$ such that $H_{0}(F)/L$ is a torsion $R$-module. Since $H_{0}(F)/L$ is finitely generated, it follows that there exists a nonzerodivisor $a\in R$, such that $aH_{0}(F)\subseteq L$.
We show that
$$(0:_{H_{0}(F)}a)=\iota\left(H_{0}(G)\right).$$
Indeed, if $x\in H_{0}(F)$ be such that $ax=0$, then $a\left(x+\iota \left(H_{0}(G)\right)\right)=0$ in $H_{0}(F)/\iota \left(H_{0}(G)\right)$.
As $r>0$, Proposition \ref{1.4.7} implies that $\Ass_{R}(C)=\Ass(R)$, so $a$ is also a nonzerodivisor on $C\cong H_{0}(F)/\iota \left(H_{0}(G)\right)$, implying that $x\in \iota\left(H_{0}(G)\right)$. Conversely, if $x\in \iota\left(H_{0}(G)\right)$, then $ax\in \iota\left(H_{0}(G)\right) \cap L$. But
$$\Ass_{R}\left(\iota\left(H_{0}(G)\right) \cap L \right)\subseteq \Ass_{R}\left(\iota\left(H_{0}(G)\right) \right)\cap \Ass_{R}(L)=\{\mathfrak{m}\}\cap \Ass(R)=\emptyset,$$
so $\iota\left(H_{0}(G)\right) \cap L =0$, and in particular $ax=0$, i.e. $x\in (0:_{H_{0}(F)}a)$. It follows that
$$(0:_{H_{0}(F)}a)= \ker \left(H_{0}(F)\xrightarrow {a} H_{0}(F)\right)=\iota\left(H_{0}(G)\right),$$
and thus there is a short exact sequence
$$0\rightarrow H_{0}(G) \xrightarrow {\iota} H_{0}(F) \xrightarrow {a} aH_{0}(F) \rightarrow 0,$$
which in turn yields
$$aH_{0}(F)\cong H_{0}(F)/\iota\left(H_{0}(G)\right)\cong C \in \mathcal{LG}(R).$$
On the other hand, applying the functor $\Gamma_{\mathfrak{m}}(-)$ to the short exact sequence \eqref{eq:1.4.9.1}, gives the exact sequence
$$0\rightarrow H_{0}(G) \xrightarrow {\iota} \Gamma_{\mathfrak{m}}\left(H_{0}(F)\right) \rightarrow \Gamma_{\mathfrak{m}}(C) = 0,$$
where the vanishing comes from the fact that $\depth_{R}(C)>0$. Therefore,
$$(0:_{H_{0}(F)}a)= \iota\left(H_{0}(G)\right) = \Gamma_{\mathfrak{m}}\left(H_{0}(F)\right).$$
\end{prf}

Next, we deal with some special cases of Conjecture \ref{1.4.4}. We need the Acyclicity Lemma.

\begin{lemma} \label{1.4.10}
Let $(R,\mathfrak{m})$ be a noetherian local ring, and let
$$L: 0\rightarrow L_{s}\rightarrow L_{s-1}\rightarrow \cdots \rightarrow L_{1}\rightarrow L_{0}\rightarrow 0$$
be an $R$-complex. Suppose that for every $1 \leq i \leq s$, the following conditions are satisfied:
\begin{enumerate}
\item [(a)] $\depth_{R}(L_{i})\geq i$.
\item [(b)] $H_{i}(L)=0$ or $\depth_{R}\left(H_{i}(L)\right)= 0$.
\end{enumerate}
Then $H_{i}(L)=0$ for every $1 \leq i \leq s$.
\end{lemma}

\begin{prf}
See \cite[Exercise 1.4.24]{BH}.
\end{prf}

\begin{theorem} \label{1.4.11}
Let $(R,\mathfrak{m})$ be a noetherian local ring. Let
$$G: 0\rightarrow G_{s}\rightarrow G_{s-1}\rightarrow \cdots \rightarrow G_{1}\rightarrow G_{0}\rightarrow 0$$
be a non-exact $R$-complex of totally reflexive modules such that $\ell_{R}\left(H_{i}(G)\right)< \infty$ for every $0 \leq i \leq s$. Assume that one of the following conditions holds:
\begin{enumerate}
\item [(a)] $R$ is Cohen-Macaulay.
\item [(b)] $R$ is reduced and $s=0$.
\item [(c)] $H_{t}(G)\in \mathcal{LG}(R)^{\perp}$, where $t:= \min \left\{0 \leq i \leq s \suchthat H_{i}(G)\neq 0\right\}$.
\end{enumerate}
Then $\dim(R)\leq s$.
\end{theorem}

\begin{prf}
First suppose that (a) holds. Consider the following augmented $R$-complex:
$$G^{+}: 0\rightarrow G_{s}\rightarrow G_{s-1}\rightarrow \cdots \rightarrow G_{1}\rightarrow G_{0}\rightarrow 0 \rightarrow 0$$
If $H_{i}(G^{+})\neq 0$ for some $0 \leq i \leq s$, then $\depth_{R}\left(H_{i}(G^{+})\right)=0$, since $\ell_{R}\left(H_{i}(G^{+})\right)< \infty$. As the $R$-complex $G$ is non-exact, Lemma \ref{1.4.10} implies that $\depth_{R}(G_{i})< i+1$ for some $0 \leq i \leq s$. But then
$$\dim(R)=\depth(R)= \depth_{R}(G_{i})< i+1 \leq s+1,$$
i.e. $\dim(R) \leq s$.

Next suppose that (b) holds. Since $s=0$, we have a nonzero totally reflexive $R$-module $G_{0}$ with finite length. Therefore, $\depth(R)=\depth_{R}(G_{0})=0$. It follows that $\Gamma_{\mathfrak{m}}(R)\neq 0$. If $\Gamma_{\mathfrak{m}}(R) \subseteq \mathfrak{m}$, then choose any $0 \neq a\in \Gamma_{\mathfrak{m}}(R)$. There is some integer $t \geq 1$ such that $\mathfrak{m}^{t}a=0$. But $a\in \mathfrak{m}$, so $a^{t+1}=0$. Since $R$ is reduced, we must have $a=0$, which is a contradiction. Therefore, $\Gamma_{\mathfrak{m}}(R)=R$. It follows that there is an integer $u \geq 1$ such that $\mathfrak{m}^{u}1=\mathfrak{m}^{u}=0$. In particular, every element of $\mathfrak{m}$ is nilpotent. As $R$ is reduced, $\mathfrak{m}=0$, i.e. $R$ is a field, so $\dim(R)=0 \leq s$.

Finally suppose that (c) holds. Without loss of generality, we may assume that $t=0$, i.e. $H_{0}(G)\neq 0$, since otherwise we can extract a short exact sequence
$$0 \rightarrow K \rightarrow G_{1} \rightarrow G_{0} \rightarrow 0$$
from the $R$-complex $G$, so that $K$ is totally reflexive. We may then work with the shorter $R$-complex
$$G: 0\rightarrow G_{s}\rightarrow G_{s-1}\rightarrow \cdots \rightarrow G_{2}\rightarrow K\rightarrow 0$$
instead. We now apply Theorem \ref{1.2.22} to the $R$-complex $G$. We obtain an $R$-complex
$$F: 0\rightarrow F_{s}\rightarrow F_{s-1}\rightarrow \cdots \rightarrow F_{1}\rightarrow F_{0}\rightarrow 0,$$
consisting of finitely generated free $R$-modules, such that $H_{i}(F)\cong H_{i}(G)$ for every $1 \leq i \leq s$, and a short exact sequence of $R$-modules
\begin{equation} \label{eq:1.4.11.1}
0\rightarrow H_{0}(G) \xrightarrow {\iota} H_{0}(F) \rightarrow C \rightarrow 0
\end{equation}
where $H_{0}(F)$ is a locally free module on the punctured spectrum of $R$ with constant rank, and $C\in \mathcal{LG}(R)$. The hypothesis implies that $\Ext^{1}_{R}\left(C,H_{0}(G)\right)=0$, so the short exact sequence \eqref{eq:1.4.11.1} splits. Applying the functor $R/\mathfrak{m}\otimes_{R}-$ to \eqref{eq:1.4.11.1}, gives the short exact sequence
$$0\rightarrow H_{0}(G)/\mathfrak{m}H_{0}(G) \xrightarrow {\overline{\iota}} H_{0}(F)/\mathfrak{m}H_{0}(F) \rightarrow C/\mathfrak{m}C \rightarrow 0.$$
If $\iota \left(H_{0}(G)\right)\subseteq \mathfrak{m}H_{0}(F)$, then $\overline{\iota}=0$, implying that $H_{0}(G)/\mathfrak{m}H_{0}(G)=0$, and that $H_{0}(G)=0$ by Nakayama's Lemma, which contradicts our assumption. Therefore, $\iota \left(H_{0}(G)\right)\nsubseteq \mathfrak{m}H_{0}(F)$. Choose any $e\in \iota\left(H_{0}(G)\right) \backslash \mathfrak{m}H_{0}(F)$. Clearly, $e\in H_{0}(F) \backslash \mathfrak{m}H_{0}(F)$ and $\ell_{R}(Re)< \infty$. Now Theorem \ref{1.4.2} implies that $\dim(R) \leq s$.
\end{prf}

In the sequel, we are going to establish another special case of Conjecture \ref{1.4.4}. To this end, we need the notion of a big Cohen-Macaulay module.

\begin{definition} \label{1.4.12}
Let $(R,\mathfrak{m})$ be a noetherian local ring. An $R$-module $B$ is called a \textit{big Cohen-Macaulay} \index{big Cohen-Macaulay module} module if there exists a system of parameters for $R$ which is $B$-regular. Further, $B$ is called a \textit{balanced big Cohen-Macaulay} \index{balanced big Cohen-Macaulay module} module if every system of parameters for $R$ is $B$-regular.
\end{definition}

The following proposition collects some important facts about big Cohen-Macaulay modules which will be used in Theorem \ref{1.4.16}.

\begin{proposition} \label{1.4.13}
Let $(R,\mathfrak{m},k)$ be a noetherian local ring, and $B$ an $R$-module. Then the following assertions hold:
\begin{enumerate}
\item [(i)] $R$ possesses a big Cohen-Macaulay module.
\item [(ii)] If $B$ is a big Cohen-Macaulay $R$-module, then $\widehat{B}^{\mathfrak{m}}$ is a balanced big Cohen-Macaulay $R$-module.
\item [(iii)] $\widehat{B}^{\mathfrak{m}}$ is a balanced big Cohen-Macaulay $R$-module if and only if $\depth_{R}(B)=\depth(R)$.
\item [(iv)] If $B$ is a big Cohen-Macaulay $R$-module and $M$ a nonzero finitely generated $R$-module, then $M\otimes _{R} B\neq 0$.
\end{enumerate}
\end{proposition}

\begin{prf}
(i): See \cite[Corollary 8.5.3]{BH} and \cite[Therem 0.7.1]{An}.

(ii) and (iii): See \cite[Corollary 8.5.3, and Exercise 9.1.12]{BH}.

(iv): Let $\underline{a}=a_{1},...,a_{n}$ be a system of parameters for $R$ which is $B$-regular. In particular, $(\underline{a})B \neq B$. As $\sqrt{(\underline{a})} = \mathfrak{m}$, there is an integer $t\geq 1$ such that $\mathfrak{m}^{t} \subseteq (\underline{a})$. In particular, $\mathfrak{m}^{t}B\neq B$, whence $\mathfrak{m}B\neq B$. On the other hand, the Nakayama's Lemma implies that $\mathfrak{m}M \neq M$. We thus have
$$(M\otimes_{R}B)\otimes_{R} k \cong (M/\mathfrak{m}M)\otimes_{k}(B/\mathfrak{m}B)\neq 0$$
as it is the tensor product of two nonzero vector spaces over the field $k$. It follows that $M\otimes_{R}B \neq 0$.
\end{prf}

We need two lemmas.

\begin{lemma} \label{1.4.14}
Let $(R,\mathfrak{m})$ be a noetherian local ring, and $B$ an $R$-module. Then the following conditions are equivalent:
\begin{enumerate}
\item [(i)] $\Ext^{1}_{R}(G,B)=0$ for every $G\in \mathcal{LG}(R)$.
\item [(ii)] $\Tor^{R}_{1}(G,B)=0$ for every $G\in \mathcal{LG}(R)$.
\end{enumerate}
\end{lemma}

\begin{prf}
In this proof, for any $i\geq 1$, we denote an $i$th syzygy of an $R$-module $M$ in any projective resolution by $M_{i}$. It is clear that $(M_{i})_{j}=M_{i+j}$. There is a short exact sequence
\begin{equation} \label{eq:1.4.14.1}
0 \rightarrow B_{1} \rightarrow F \rightarrow B \rightarrow 0,
\end{equation}
where $F$ is a free $R$-module.
For any $R$-module $B$ and $G\in \mathcal{LG}(R)$, define the homomorphism
$$\theta_{GB}: G \otimes_{R} B \rightarrow \Hom_{R}(G^{\ast},B),$$
given by $\theta_{GB}(g \otimes b)(\varphi):= \varphi(g)b$. As $G^{\ast}$ is totally reflexive, the short exact sequence \eqref{eq:1.4.14.1}, gives the following diagram with exact rows:
\[\begin{tikzpicture}[every node/.style={midway},]
  \matrix[column sep={1em}, row sep={2.5em}]
  {\node(1) {$0$}; & \node(2) {$\Tor^{R}_{1}(G,B)$}; & \node(3) {$G\otimes_{R}B_{1}$}; & \node(4) {$G\otimes_{R}F$}; & \node(5) {$G\otimes_{R}B$}; & \node(6) {$0$}; & \node(7) {};\\
  \node(8) {}; & \node(9) {$0$}; & \node(10) {$\Hom_{R}(G^{\ast},B_{1})$}; & \node(11) {$\Hom_{R}(G^{\ast},F)$}; & \node(12) {$\Hom_{R}(G^{\ast},B)$}; & \node(13) {$\Ext^{1}_{R}(G^{\ast},B_{1})$}; & \node(14) {$0$}; \\};
  \draw[decoration={markings,mark=at position 1 with {\arrow[scale=1.5]{>}}},postaction={decorate},shorten >=0.5pt] (1) -- (2) node[anchor=west] {};
  \draw[decoration={markings,mark=at position 1 with {\arrow[scale=1.5]{>}}},postaction={decorate},shorten >=0.5pt] (2) -- (3) node[anchor=west] {};
  \draw[decoration={markings,mark=at position 1 with {\arrow[scale=1.5]{>}}},postaction={decorate},shorten >=0.5pt] (3) -- (4) node[anchor=south] {$\zeta$};
  \draw[decoration={markings,mark=at position 1 with {\arrow[scale=1.5]{>}}},postaction={decorate},shorten >=0.5pt] (4) -- (5) node[anchor=south] {$\pi$};
  \draw[decoration={markings,mark=at position 1 with {\arrow[scale=1.5]{>}}},postaction={decorate},shorten >=0.5pt] (5) -- (6) node[anchor=west] {};
  \draw[decoration={markings,mark=at position 1 with {\arrow[scale=1.5]{>}}},postaction={decorate},shorten >=0.5pt] (9) -- (10) node[anchor=west] {};
  \draw[decoration={markings,mark=at position 1 with {\arrow[scale=1.5]{>}}},postaction={decorate},shorten >=0.5pt] (10) -- (11) node[anchor=south] {$\iota$};
  \draw[decoration={markings,mark=at position 1 with {\arrow[scale=1.5]{>}}},postaction={decorate},shorten >=0.5pt] (11) -- (12) node[anchor=south] {$\eta$};
  \draw[decoration={markings,mark=at position 1 with {\arrow[scale=1.5]{>}}},postaction={decorate},shorten >=0.5pt] (12) -- (13) node[anchor=west] {};
  \draw[decoration={markings,mark=at position 1 with {\arrow[scale=1.5]{>}}},postaction={decorate},shorten >=0.5pt] (13) -- (14) node[anchor=west] {};
  \draw[decoration={markings,mark=at position 1 with {\arrow[scale=1.5]{>}}},postaction={decorate},shorten >=0.5pt] (3) -- (10) node[anchor=west] {$\theta_{GB_{1}}$};
  \draw[decoration={markings,mark=at position 1 with {\arrow[scale=1.5]{>}}},postaction={decorate},shorten >=0.5pt] (4) -- (11) node[anchor=east] {$\cong$} node[anchor=west] {$\theta_{GF}$};
  \draw[decoration={markings,mark=at position 1 with {\arrow[scale=1.5]{>}}},postaction={decorate},shorten >=0.5pt] (5) -- (12) node[anchor=west] {$\theta_{GB}$};
\end{tikzpicture}\]
It follows from the above diagram that
$$\im \theta_{GB} = \im (\theta_{GB}\pi) = \im (\eta \theta_{GF})= \im \eta,$$
so
\begin{equation} \label{eq:1.4.14.2}
\coker \theta_{GB}= \coker \eta \cong \Ext^{1}_{R}(G^{\ast},B_{1}).
\end{equation}
Now consider the short exact sequence
$$0 \rightarrow G_{1} \rightarrow E \rightarrow G \rightarrow 0,$$
where $E$ is a finitely generated free $R$-module, so that $G_{1}$ is totally reflexive by Proposition \ref{1.2.3}. As $G$ is totally reflexive, we get the short exact sequence
$$0 \rightarrow G^{\ast} \rightarrow E \rightarrow (G_{1})^{\ast} \rightarrow 0,$$
which in turn yields the exact sequence
$$0=\Ext^{1}_{R}(E,B_{1}) \rightarrow \Ext^{1}_{R}(G^{\ast},B_{1}) \rightarrow \Ext^{2}_{R}\left((G_{1})^{\ast},B_{1}\right) \rightarrow \Ext^{2}_{R}(E,B_{1})=0,$$
implying the isomorphism
\begin{equation} \label{eq:1.4.14.3}
\Ext^{1}_{R}(G^{\ast},B_{1}) \cong \Ext^{2}_{R}\left((G_{1})^{\ast},B_{1}\right).
\end{equation}
Since $(G_{1})^{\ast}$ is totally reflexive, the short exact sequence \eqref{eq:1.4.14.1} yields the exact sequence
$$0=\Ext^{1}_{R}\left((G_{1})^{\ast},F\right) \rightarrow \Ext^{1}_{R}\left((G_{1})^{\ast},B\right) \rightarrow \Ext^{2}_{R}\left((G_{1})^{\ast},B_{1}\right) \rightarrow \Ext^{2}_{R}\left((G_{1})^{\ast},F\right)=0,$$
which gives the isomorphism
\begin{equation} \label{eq:1.4.14.4}
\Ext^{1}_{R}\left((G_{1})^{\ast},B\right) \cong \Ext^{2}_{R}\left((G_{1})^{\ast},B_{1}\right).
\end{equation}
From \eqref{eq:1.4.14.2}, \eqref{eq:1.4.14.3}, and \eqref{eq:1.4.14.4} we get
\begin{equation} \label{eq:1.4.14.5}
\coker \theta_{GB}\cong \Ext^{1}_{R}(G^{\ast},B_{1}) \cong \Ext^{1}_{R}\left((G_{1})^{\ast},B\right).
\end{equation}
From the above diagram, using the Snake Lemma and \eqref{eq:1.4.14.5}, we get
\begin{equation} \label{eq:1.4.14.6}
\begin{split}
\ker \theta_{GB} & \cong \coker \theta_{GB_{1}} \\
 & \cong \Ext^{1}_{R}(G^{\ast},B_{2}) \\
 & \cong \Ext^{1}_{R}\left((G_{1})^{\ast},B_{1}\right) \\
 & \cong \Ext^{1}_{R}\left((G_{2})^{\ast},B\right). \\
\end{split}
\end{equation}
From \eqref{eq:1.4.14.5} and \eqref{eq:1.4.14.6}, we get
\begin{equation} \label{eq:1.4.14.7}
\ker \theta_{GB_{1}} \cong \Ext^{1}_{R}\left((G_{2})^{\ast},B_{1}\right) \cong \Ext^{1}_{R}\left((G_{3})^{\ast},B\right). \\
\end{equation}
On the other hand, from the above diagram we have
\begin{equation} \label{eq:1.4.14.8}
\ker \theta_{GB_{1}} = \ker (\iota\theta_{GB_{1}}) = \ker (\theta_{GF}\zeta) \\ = \ker \zeta \\ \cong \Tor_{1}^{R}(G,B). \\
\end{equation}
From \eqref{eq:1.4.14.7} and \eqref{eq:1.4.14.8}, we get
\begin{equation} \label{eq:1.4.14.9}
\Ext^{1}_{R}\left((G_{3})^{\ast},B\right) \cong \Tor_{1}^{R}(G,B). \\
\end{equation}
Next we note that as $G\in \mathcal{LG}(R)$, Proposition \ref{1.2.3} together with Lemma \ref{1.2.22} imply that $G_{i}^{\ast}\in \mathcal{LG}(R)$ for every $i\geq 1$. On the other hand, we have $G^{\ast}\in \mathcal{LG}(R)$, so there is a $\Hom_{R}(-,\mathcal{P})$-exact exact sequence
$$F: \cdots \rightarrow F_{i+1} \xrightarrow {\partial^{F}_{i+1}} F_{i} \xrightarrow {\partial^{F}_{i}} F_{i-1} \rightarrow \cdots,$$
consisting of finitely generated free $R$-modules, such that $G^{\ast}\cong \im \partial^{F}_{i}$. If $L= \im \partial^{F}_{0}$, then $L$ is totally reflexive, and we get the exact sequence
$$0 \rightarrow G^{\ast} \rightarrow F_{i-1} \rightarrow F_{i-2} \rightarrow \cdots \rightarrow F_{0}\rightarrow L \rightarrow 0,$$
showing that $G^{\ast} = L_{i}$.
Further, Lemma \ref{1.2.22} implies that $L\in \mathcal{LG}(R)$.
Therefore, $G \cong G^{\ast\ast} = (L_{i})^{\ast}$, i.e. $G$ can be isomorphic to the dual of an $i$th syzygy of an element of $\mathcal{LG}(R)$ for any $i \geq 1$. From this observation and \eqref{eq:1.4.14.9}, the assertion is immediate.
\end{prf}

\begin{lemma} \label{1.4.15}
Let $(R,\mathfrak{m})$ be a noetherian local ring and $B$ an $R$-module. Let
$$F: 0\rightarrow F_{s}\xrightarrow {\partial^{F}_{s}} F_{s-1}\rightarrow \cdots \rightarrow F_{1} \xrightarrow {\partial^{F}_{1}} F_{0}\rightarrow 0$$
be an $R$-complex of finitely generated free modules such that $H_{i}(F\otimes_{R}B)=0$ for every $1 \leq i \leq s$. Then $$\depth_{R}\left(H_{0}(F)\otimes_{R}B\right) \geq \depth_{R}(B)-s.$$
\end{lemma}

\begin{prf}
Consider the $R$-complex
$$F\otimes_{R}B: 0 \rightarrow F_{s}\otimes_{R}B \xrightarrow {\partial^{F}_{s}\otimes_{R} B} F_{s-1}\otimes_{R}B \rightarrow \cdots \rightarrow F_{1}\otimes_{R} B \xrightarrow {\partial^{F}_{1}\otimes_{R}B} F_{0}\otimes_{R}B \rightarrow 0.$$
From the short exact sequence
$$0 \rightarrow \im \left(\partial^{F}_{1}\otimes_{R}B\right) \rightarrow F_{0}\otimes_{R}B \rightarrow H_{0}(F\otimes_{R}B) \rightarrow 0,$$
we get
\begin{equation} \label{eq:1.4.15.1}
\begin{split}
\depth_{R}\left(H_{0}(F)\otimes_{R}B\right) & = \depth_{R}\left(\left(\coker \partial^{F}_{1}\right) \otimes_{R}B\right) \\
 & = \depth_{R}\left(\coker \left(\partial^{F}_{1}\otimes_{R}B\right)\right) \\
 & = \depth_{R}\left(H_{0}(F\otimes_{R}B)\right) \\
 & \geq \min \left\{\depth_{R}\left(\im \left(\partial^{F}_{1}\otimes_{R}B\right)\right)-1, \depth_{R}(F_{0}\otimes_{R}B) \right\} \\
 & = \min \left\{\depth_{R}\left(\im \left(\partial^{F}_{1}\otimes_{R}B\right)\right)-1, \depth_{R}(B) \right\}. \\
\end{split}
\end{equation}
Now consider the short exact sequence
$$0 \rightarrow \im \left(\partial^{F}_{2}\otimes_{R}B\right) \rightarrow F_{1}\otimes_{R}B \rightarrow \im \left(\partial^{F}_{1}\otimes_{R}B\right) \rightarrow 0.$$
Similarly, we get
\begin{equation} \label{eq:1.4.15.2}
\depth_{R}\left(\im \left(\partial^{F}_{1}\otimes_{R}B\right)\right) \geq \min \left\{\depth_{R}\left(\im \left(\partial^{F}_{2}\otimes_{R}B\right)\right)-1, \depth_{R}(B) \right\}.
\end{equation}
From \eqref{eq:1.4.15.1} and \eqref{eq:1.4.15.2}, we get
$$\depth_{R}\left(H_{0}(F)\otimes_{R}B\right)\geq \min \left\{\depth_{R}\left(\im \left(\partial^{F}_{2}\otimes_{R}B\right)\right)-2, \depth_{R}(B) \right\}.$$
We continue this way until we reach at
\begin{equation} \label{eq:1.4.15.3}
\depth_{R}\left(H_{0}(F)\otimes_{R}B\right)\geq \min \left\{\depth_{R}\left(\im \left(\partial^{F}_{s-1}\otimes_{R}B\right)\right)-(s-1), \depth_{R}(B) \right\}.
\end{equation}
Now the short exact sequence
$$0 \rightarrow F_{s}\otimes_{R}B \rightarrow F_{s-1}\otimes_{R}B \rightarrow \im \left(\partial^{F}_{s-1}\otimes_{R}B\right) \rightarrow 0$$
gives
\begin{equation} \label{eq:1.4.15.4}
\depth_{R}\left(\im \left(\partial^{F}_{s-1}\otimes_{R}B\right)\right) \geq \min \left\{\depth_{R}(B)-1, \depth_{R}(B) \right\}=\depth_{R}(B)-1.
\end{equation}
From \eqref{eq:1.4.15.3} and \eqref{eq:1.4.15.4}, we get
\begin{equation*}
\begin{split}
\depth_{R}\left(H_{0}(F)\otimes_{R}B\right) & \geq \min \left\{\depth_{R}\left(\im \left(\partial^{F}_{s-1}\otimes_{R}B\right)\right)-(s-1), \depth_{R}(B) \right\} \\
 & \geq \min \left\{\depth_{R}(B)-1-(s-1)), \depth_{R}(B) \right\} \\
 & = \depth_{R}(B)-s. \\
\end{split}
\end{equation*}
\end{prf}

We are now ready to present the main theorem of this section.

\begin{theorem} \label{1.4.16}
Let $(R,\mathfrak{m})$ be a noetherian local ring. Let
$$G: 0\rightarrow G_{s}\rightarrow G_{s-1}\rightarrow \cdots \rightarrow G_{1}\rightarrow G_{0}\rightarrow 0$$
be a non-exact $R$-complex of totally reflexive modules whose homology modules have all finite length. If there exists a balanced big Cohen-Macaulay module $B\in \mathcal{LG}(R)^{\perp}$, then $\dim(R) \leq s$.
\end{theorem}

\begin{prf}
As observed in the proof of Theorem \ref{1.4.11}, we may assume that $H_{0}(G)\neq 0$. We apply Theorem \ref{1.2.22} to the $R$-complex $G$ to get an $R$-complex
$$F: 0\rightarrow F_{s} \rightarrow F_{s-1} \rightarrow \cdots \rightarrow F_{1} \rightarrow F_{0} \rightarrow 0,$$
consisting of finitely generated free modules, with $H_{i}(G)\cong H_{i}(F)$ for every $1 \leq i \leq s$, and a short exact sequence of $R$-modules
\begin{equation} \label{eq:1.4.16.1}
0\rightarrow H_{0}(G) \rightarrow H_{0}(F) \rightarrow C \rightarrow 0,
\end{equation}
where $H_{0}(F)$ is locally free on the punctured spectrum of $R$ with constant rank and $C\in \mathcal{LG}(R)$. We now assume to the contrary that $\dim(R)>s$, and aim at a contradiction.

We claim that
\begin{equation} \label{eq:1.4.16.2}
\dim(R)-\dim\left(\frac{R}{I_{r_{i}}\left(\partial^{F}_{i}\right)}\right) \geq i
\end{equation}
for every $1 \leq i \leq s$, where $r_{i}$ is the expected rank of $\partial^{F}_{i}$ i.e. $r_{i}= \sum_{j=i}^{s}(-1)^{j-i}\rank_{R}(F_{j})$, and $I_{r_{i}}\left(\partial^{F}_{i}\right)$ is the $r_{i}$th Fitting ideal of $\partial^{F}_{i}$. Indeed, if for some $1 \leq j \leq s$, we have
$$\dim(R)-\dim\left(\frac{R}{I_{r_{j}}\left(\partial^{F}_{j}\right)}\right) < j,$$
then the hypothesis $\dim(R) > s$ implies that $\dim\left(R/I_{r_{j}}\left(\partial^{F}_{j}\right)\right) \neq 0$. Take any $\mathfrak{p} \in \Assh_{R}\left(R/I_{r_{j}}\left(\partial^{F}_{j}\right)\right)$. It is clear that $\mathfrak{p} \neq \mathfrak{m}$. The hypothesis $\ell_{R}\left(H_{i}(F)\right)<\infty$ implies that $H_{i}(R_{\mathfrak{p}}\otimes_{R}F) \cong H_{i}(F)_{\mathfrak{p}}=0$ for every $1 \leq i \leq s$. Therefore, \cite[Lemma 9.1.6]{BH} implies that
$$\depth_{R_{\mathfrak{p}}}\left(I_{r_{i}}\left(R_{\mathfrak{p}}\otimes_{R}\partial_{i}^{F}\right),R_{\mathfrak{p}}\right) \geq i$$
for every $1 \leq i \leq s$. We thus have
\begin{equation*}
\begin{split}
j & \leq \depth_{R_{\mathfrak{p}}}\left(I_{r_{j}}\left(R_{\mathfrak{p}}\otimes_{R}\partial_{j}^{F}\right),R_{\mathfrak{p}}\right) \\
 & = \depth_{R_{\mathfrak{p}}}\left(I_{r_{j}}\left(\partial_{j}^{F}\right)R_{\mathfrak{p}},R_{\mathfrak{p}}\right) \\
 & \leq \Ht\left(I_{r_{j}}\left(\partial_{j}^{F}\right)R_{\mathfrak{p}}\right) \\
 & \leq \Ht\left(\mathfrak{p}R_{\mathfrak{p}}\right) \\
 & = \Ht(\mathfrak{p}) \\
 & \leq \dim(R)- \dim(R/\mathfrak{p}) \\
 & = \dim(R)- \dim \left(\frac{R}{I_{r_{i}}\left(\partial^{F}_{i}\right)}\right) \\
 & <j. \\
\end{split}
\end{equation*}
This contradiction establishes the claim.

Now fix $1 \leq i \leq s$. If
$$I_{r_{i}}\left(\partial^{F}_{i}\right) \subseteq \bigcup_{\mathfrak{p}\in\Assh(R)}\mathfrak{p},$$
then there is a $\mathfrak{p} \in \Assh(R)$ with $I_{r_{i}}\left(\partial^{F}_{i}\right) \subseteq \mathfrak{p}$. Hence by \eqref{eq:1.4.16.2}, we have
$$\dim(R)=\dim(R/\mathfrak{p}) \leq \dim\left(\frac{R}{I_{r_{i}}\left(\partial^{F}_{i}\right)}\right) \leq \dim(R)-i$$
which is a contradiction. Thus there is an element
$$a_{1} \in I_{r_{i}}\left(\partial^{F}_{i}\right) \setminus \bigcup_{\mathfrak{p}\in\Assh(R)}\mathfrak{p},$$
so that $a_{1}$ is part of a system of parameters for $R$. If $i \geq 2$ and
$$I_{r_{i}}\left(\partial^{F}_{i}\right) \subseteq \bigcup_{\mathfrak{p}\in\Assh\left(R/(a_{1})\right)}\mathfrak{p},$$
then there is a $\mathfrak{p} \in \Assh\left(R/(a_{1})\right)$ with $I_{r_{i}}\left(\partial^{F}_{i}\right) \subseteq \mathfrak{p}$. Hence by \eqref{eq:1.4.16.2}, we have
$$\dim(R)-1=\dim\left(R/(a_{1})\right)=\dim(R/\mathfrak{p}) \leq \dim\left(\frac{R}{I_{r_{i}}\left(\partial^{F}_{i}\right)}\right) \leq \dim(R)-i$$
which is a contradiction. Thus there is an element
$$a_{2} \in I_{r_{i}}\left(\partial^{F}_{i}\right) \setminus \bigcup_{\mathfrak{p}\in\Assh\left(R/(a_{1})\right)}\mathfrak{p},$$
so that $a_{1},a_{2}$ is part of a system of parameters for $R$. Continuing in this manner, we see that there is a sequence $a_{1},...,a_{i} \in I_{r_{i}}(\partial^{F}_{i})$, which is part of a system of parameters for $R$. Since $B$ is a balanced big Cohen-Macaulay $R$-module, $a_{1},...,a_{i}$ is a regular sequence on $B$, implying that $\depth_{R}\left(I_{r_{i}}(\partial^{F}_{i}),B\right) \geq i$. Now \cite[Lemma 9.1.6]{BH} implies that
\begin{equation} \label{eq:1.4.16.3}
H_{i}(F \otimes _{R}B)=0
\end{equation}
for every $1 \leq i \leq s$.
Therefore, Lemma \ref{1.4.15} yields that
$$\depth_{R}\left(H_{0}(F)\otimes_{R}B\right) \geq \depth_{R}(B)-s = \dim(R)-s \geq 1,$$
and since
$$\depth_{R}\left(H_{0}(F)\otimes_{R}B\right) = \inf \left\{i \geq 0 \suchthat H^{i}_{\mathfrak{m}}\left(H_{0}(F)\otimes_{R}B\right)\neq 0 \right\},$$
we conclude that
\begin{equation} \label{eq:1.4.16.4}
\Gamma_{\mathfrak{m}}\left(H_{0}(F)\otimes_{R}B\right)=0.
\end{equation}
Since $\inf F \geq 0$, we have an isomorphism $F \simeq F_{\supset 0}$ in the derived category $\mathcal{D}(R)$, and thus we get the distinguished triangle
$$F_{\supset 1} \rightarrow F \rightarrow H_{0}(F) \rightarrow.$$
Applying the triangulated functor $-\otimes_{R}^{\bf L}B$ to this distinguished triangle, we get the distinguished triangle
\begin{equation} \label{eq:1.4.16.5}
F_{\supset 1} \otimes_{R}^{\bf L} B \rightarrow F \otimes_{R}^{\bf L} B \rightarrow H_{0}(F) \otimes_{R}^{\bf L} B \rightarrow.
\end{equation}
As $F$ is an $R$-complex of free modules, we have
$$H_{1}\left(F \otimes_{R}^{\bf L} B\right) = H_{1}(F \otimes_{R} B) = 0$$
by \eqref{eq:1.4.16.3}.
On the other hand, we have
$$\inf \left(F_{\supset 1} \otimes_{R}^{\bf L} B\right) \geq \inf \left(F_{\supset 1}\right) + \inf B \geq 1.$$
Therefore, the long exact sequence of homology modules associated with the distinguished triangle \eqref{eq:1.4.16.5}, gives the exact sequence
$$0=H_{1}\left(F \otimes_{R}^{\bf L} B\right) \rightarrow H_{1}\left(H_{0}(F) \otimes_{R}^{\bf L} B\right) \rightarrow H_{0}\left(F_{\supset 1} \otimes_{R}^{\bf L} B\right)=0,$$
so $$\Tor_{1}^{R}\left(H_{0}(F),B\right)= H_{1}\left(H_{0}(F) \otimes_{R}^{\bf L} B\right)=0.$$
Thus, from the short exact sequence \eqref{eq:1.4.16.1}, we get the exact sequence
\begin{equation} \label{eq:1.4.16.6}
0 \rightarrow \Tor_{1}^{R}(C,B) \rightarrow H_{0}(G)\otimes_{R}B \rightarrow H_{0}(F)\otimes_{R}B \rightarrow C\otimes_{R}B \rightarrow 0.
\end{equation}
Since $C$ is locally free on the punctured spectrum of $R$, we observe that
$$\Tor_{1}^{R}(C,B)_{\mathfrak{p}}\cong \Tor_{1}^{R_{\mathfrak{p}}}(C_{\mathfrak{p}},B_{\mathfrak{p}})=0,$$
for every $\mathfrak{p}\in \Spec(R)\backslash \{\mathfrak{m}\}$, so $\Tor_{1}^{R}(C,B)$ is $\mathfrak{m}$-torsion. Moreover, since $\ell _{R}\left(H_{0}(G)\right) < \infty$, we see that
$$\left(H_{0}(G)\otimes_{R}B\right)_{\mathfrak{p}}\cong H_{0}(G)_{\mathfrak{p}}\otimes_{R_{\mathfrak{p}}}B_{\mathfrak{p}}=0,$$
for every $\mathfrak{p}\in \Spec(R)\backslash \{\mathfrak{m}\}$, so $H_{0}(G)\otimes_{R}B$ is $\mathfrak{m}$-torsion.
Therefore, applying the functor $\Gamma_{\mathfrak{m}}(-)$ to \eqref{eq:1.4.16.6} yields the exact sequence
$$0 \rightarrow \Tor_{1}^{R}(C,B) \rightarrow H_{0}(G)\otimes_{R}B \rightarrow \Gamma_{\mathfrak{m}}\left(H_{0}(F)\otimes_{R}B\right)=0,$$
where the vanishing comes from \eqref{eq:1.4.16.4}. Hence
\begin{equation} \label{eq:1.4.16.7}
\Tor_{1}^{R}(C,B) \cong H_{0}(G)\otimes_{R}B.
\end{equation}
By Proposition \ref{1.4.13}, $H_{0}(G)\otimes_{R}B \neq 0$. However, $B\in \mathcal{LG}(R)^{\perp}$ so $\Ext^{1}_{R}(G,B)=0$ for every $G\in \mathcal{LG}(R)$. Lemma \ref{1.4.14} implies that $\Tor^{R}_{1}(G,B)=0$ for every $G\in \mathcal{LG}(R)$. As $C\in \mathcal{LG}(R)$, we have $\Tor^{R}_{1}(C,B)=0$. Now the left hand side of \eqref{eq:1.4.16.7} is zero while the right hand side is nonzero. This contradiction establishes the theorem.
\end{prf}

We note that if $R$ is a Cohen-Macaulay local ring, then $R$ is a balanced big Cohen-Macaulay $R$-module such that $R\in \mathcal{LG}(R)^{\perp}$, i.e. the condition in Theorem \ref{1.4.16} is satisfied.

\chapter{Hartshorne's Notion of Cofiniteness}

\section{Introduction}

Throughout this chapter, all rings are assumed to be commutative noetherian with identity.

In his algebraic geometry seminars of 1961-2, Grothendieck founded the theory of local cohomology and raised, along the way, a few questions on the finiteness properties of local cohomology modules; see \cite[Conjectures 1.1 and 1.2]{Gr}. He specifically asked whether the $R$-modules $\Hom_{R}\left(R/\mathfrak{a},H^{i}_{\mathfrak{a}}(M)\right)$ were finitely generated for every ideal $\mathfrak{a}$ of $R$ and every finitely generated $R$-module $M$, which had been answered affirmatively in the same seminar when $(R,\mathfrak{m})$ is local and $\mathfrak{a}=\mathfrak{m}$.

In 1969, Hartshorne provided a counterexample in \cite[Section 3]{Ha1}, to show that this question does not have an affirmative answer in general. He then defined an $R$-module $M$ to be $\mathfrak{a}$-cofinite whenever $\Supp_{R}(M) \subseteq \V(\mathfrak{a})$ and $\Ext^{i}_{R}(R/\mathfrak{a},M)$ is a finitely generated $R$-module for every $i \geq 0$. Moreover, in the case where $R$ is an $\mathfrak{a}$-adically complete regular ring of finite Krull dimension, he defined an $R$-complex $X$ to be $\mathfrak{a}$-cofinite if $X \simeq {\bf R}\Hom_{R}\left(Y,{\bf R} \Gamma_{\mathfrak{a}}(R)\right)$ for some $R$-complex $Y$ with finitely generated homology modules. Hartshorne proceeded to pose three questions in this respect. We paraphrase the main themes of his questions as follows.

\begin{question} \label{2.1.1}
Are the local cohomology modules $H^{i}_{\mathfrak{a}}(M)$, $\mathfrak{a}$-cofinite for every finitely generated $R$-module $M$
and every $i \geq 0$?
\end{question}

\begin{question} \label{2.1.2}
Is the category $\mathcal{M}(R,\mathfrak{a})_{cof}$ of $\mathfrak{a}$-cofinite $R$-modules an abelian subcategory of the category $\mathcal{M}(R)$ of $R$-modules?
\end{question}

\begin{question} \label{2.1.3}
Is it true that an $R$-complex $X$ is $\mathfrak{a}$-cofinite if and only if its homology module $H_{i}(X)$ is $\mathfrak{a}$-cofinite for
every $i\in \mathbb{Z}$?
\end{question}

By providing the following counterexample, he showed that the answers to these questions are negative in general.

\begin{example} \label{2.1.4}
Let $k$ be a field, $R=k[X,Y][[U,V]]$, $\mathfrak{a}=(U,V)$, $\mathfrak{p}=(XV+YU)$, and $T=R/\mathfrak{p}$. Then $R$ is a regular domain
of dimension $4$ and $T$ is a non-regular domain. It is shown that $\Hom_{R}\left(R/\mathfrak{a},H^{2}_{\mathfrak{a}}(T)\right)$ is not
finitely generated, so in particular, $H^{2}_{\mathfrak{a}}(T)$ is not $\mathfrak{a}$-cofinite. This takes care of Question \ref{2.1.1}.
Furthermore, there is an exact sequence $$0 \rightarrow H^{1}_{\mathfrak{a}}(T) \rightarrow H^{2}_{\mathfrak{a}}(R) \rightarrow H^{2}_{\mathfrak{a}}(R)
\rightarrow H^{2}_{\mathfrak{a}}(T)\rightarrow 0.$$
The local cohomology module $H^{2}_{\mathfrak{a}}(R)$ turns out to be $\mathfrak{a}$-cofinite, whereas $H^{2}_{\mathfrak{a}}(T)$ is not
$\mathfrak{a}$-cofinite, answering Question \ref{2.1.2}. Finally, it is shown that ${\bf R}\Gamma_{\mathfrak{a}}(T)$ is an $\mathfrak{a}$-cofinite
$R$-complex, while $H^{2}_{\mathfrak{a}}(T)= H_{-2}\left({\bf R}\Gamma_{\mathfrak{a}}(T)\right)$ is not an $\mathfrak{a}$-cofinite $R$-module,
dealing with Question \ref{2.1.3}.
\end{example}

Hartshorne further established affirmative answers to these questions in the case where $\mathfrak{a}$ is a principal ideal generated by a nonzerodivisor and $R$
is an $\mathfrak{a}$-adically complete regular ring of finite Krull dimension, and also in the case where $\mathfrak{a}$ is a prime ideal with $\dim(R/\mathfrak{a})=1$ and $R$ is a complete regular local ring; see \cite[Propositions 6.1 and 6.2, Corollary 6.3, Theorem 7.5, Proposition 7.6 and Corollary 7.7]{Ha1}.

In the following years, Hartshorne's results on Questions \ref{2.1.1} and \ref{2.1.2}, were systematically extended and polished by commutative algebra practitioners in several stages to take the following culminating form.

\begin{theorem} \label{2.1.5}
Let $\mathfrak{a}$ be an ideal of $R$, and $M$ a finitely generated $R$-module. Then the following assertions hold for every $i \geq 0$:
\begin{enumerate}
\item[(i)] If $\cd(\mathfrak{a},R) \leq 1$, then $H^{i}_{\mathfrak{a}}(M)$ is $\mathfrak{a}$-cofinite, and if furthermore $R$ is local, then $\mathcal{M}(R,\mathfrak{a})_{cof}$ is abelian.
\item[(ii)] If $\dim(R) \leq 2$, then $H^{i}_{\mathfrak{a}}(M)$ is $\mathfrak{a}$-cofinite and $\mathcal{M}(R,\mathfrak{a})_{cof}$ is abelian.
\item[(iii)] If $\dim(R/\mathfrak{a})\leq 1$, then $H^{i}_{\mathfrak{a}}(M)$ is $\mathfrak{a}$-cofinite and $\mathcal{M}(R,\mathfrak{a})_{cof}$ is abelian.
\end{enumerate}
\end{theorem}

\begin{prf}
(i): See \cite[Corollary 3.14]{Me2} and \cite[Theorem 2.4]{PAB}.

(ii): See \cite[Theorem 7.10]{Me2} and \cite[Theorem 7.4]{Me2}.

(iii): See \cite[Theorem 2.6 and Corollary 2.12]{Me1}, \cite[Corollary 2.8]{BNS1}, and \cite[Corollary 2.7]{BN}.
\end{prf}

In this chapter, we relax the "local" assumption in Theorem \ref{2.1.5} (i) by deploying the theory of local homology. Subsequently, we will be privileged with the following full-fledged form of the above theorem.

\begin{theorem} \label{2.1.6}
Let $\mathfrak{a}$ be an ideal of $R$ such that either $\cd(\mathfrak{a},R) \leq 1$, or $\dim(R) \leq 2$, or $\dim(R/\mathfrak{a})\leq 1$. Then $\mathcal{M}(R,\mathfrak{a})_{cof}$ is abelian, and $H^{i}_{\mathfrak{a}}(M)$ is $\mathfrak{a}$-cofinite for every finitely generated $R$-module $M$ and every $i \geq 0$.
\end{theorem}

In the year 1990, Huneke posed some questions in \cite{Hu} on local cohomology modules. In particular, he asked whether the set of associated primes $\Ass_{R}\left(H^{i}_{\mathfrak{a}}(M)\right)$ is finite for every finitely generated $R$-module $M$ and every $i \geq 0$, or whether the Bass number $\mu^{j}_{R}\left(\mathfrak{p},H^{i}_{\mathfrak{a}}(M)\right)$ is finite for every $\mathfrak{p}\in \Spec(R)$ and every $i,j \geq 0$. These questions and the subsequent research motivated by them took the subject at hand to a whole new level. In the year 2000, Singh presented a counterexample in \cite{Sn} to set up a negative answer to Huneke's first question. As a consequence, his question was narrowed down to the following question of Lyubeznik.

\begin{question} \label{2.1.7}
Let $R$ be a regular ring and $\mathfrak{a}$ an ideal of $R$. Is $\Ass_{R}\left(H^{i}_{\mathfrak{a}}(R)\right)$ finite for every $i \geq 0$?
\end{question}

It is worth noting that Lyubeznik's question has been answered affirmatively thus far in the cases where $R$ is equicharacteristic, and remains open in the mixed characteristic case.

One can easily observe that if the local cohomology module $H^{i}_{\mathfrak{a}}(M)$ is $\mathfrak{a}$-cofinite, then the set of associated primes $\Ass_{R}\left(H^{i}_{\mathfrak{a}}(M)\right)$ and the Bass numbers $\mu^{j}_{R}\left(\mathfrak{p},H^{i}_{\mathfrak{a}}(M)\right)$ are finite for every $\mathfrak{p} \in \Spec(R)$ and every $i,j \geq 0$. This observation reveals the utmost significance of Hartshorne's Question \ref{2.1.1}, as an affirmative answer to this question in any case would readily set fourth affirmative answers to Huneke's questions, to say the least.

Proceeding further to the world of complexes, it turns out that to establish the cofiniteness of $H^{i}_{\mathfrak{a}}(X)$ for any $R$-complex $X$ with finitely generated homology modules, all we need to know is the cofiniteness of $H^{i}_{\mathfrak{a}}(M)$ for any finitely generated $R$-module $M$ and the abelianness of $\mathcal{M}(R,\mathfrak{a})_{cof}$. This in turn unmasks the importance of Hartshorne's Question \ref{2.1.2}. The crucial step to achieve this is to recruit the technique of way-out functors which is to be implemented in this chapter.

Questions \ref{2.1.1} and \ref{2.1.2} have been high-profile among researchers, whereas not much attention has been brought to Question \ref{2.1.3}. The most striking result on this question is the following; see \cite[Theorem 1]{EK}.

\begin{theorem} \label{2.1.8}
Let $R$ be a complete Gorenstein local domain and $\mathfrak{a}$ an ideal of $R$ with $\dim(R/\mathfrak{a})=1$. Then an $R$-complex $X\in \mathcal{D}_{\sqsubset}(R)$ is $\mathfrak{a}$-cofinite if and only if $H_{i}(X)$ is $\mathfrak{a}$-cofinite for every $i \in \mathbb{Z}$.
\end{theorem}

To be consistent and coherent in both module and complex cases, we define an $R$-complex $X$ to be $\mathfrak{a}$-cofinite if $\Supp_{R}(X)\subseteq
\V(\mathfrak{a})$ and ${\bf R}\Hom_{R}(R/\mathfrak{a},X)$ has finitely generated homology modules. We show that this definition coincides with that of Hartshorne.

At first glance, Hartshorne's questions seem to be irrelevant to each other and rather independent of one another. But on second thought, a curious idea that may strike one's mind swirls around the possible existence of profound connections between them way beyond what meets the eye. These questions sound to be inextricably bound up in an elegant way whose interrelations are yet to be unraveled. An ample evidence that motivates and supports this speculation relies on the very fact that these questions stand true all together in one case and fail to hold collectively in another case. As to shed some light on this revelation, we unearth a connection between Hartshorne's questions and subsequently answer Question \ref{2.1.3} completely in the cases $\cd(\mathfrak{a},R) \leq 1$, $\dim(R) \leq 2$, and $\dim(R/\mathfrak{a})\leq 1$.

Several authors have strived to extend the results of Theorem \ref{2.1.6} to generalized local cohomology modules. However, it is folklore that all the generalizations $H_{\varphi}^{i}(M)$, $H_{\mathcal{Z}}^{i}(M)$, $H_{\mathcal{\mathfrak{a},\mathfrak{b}}}^{i}(M)$, $H_{\mathcal{\mathfrak{a}}}^{i}(M,N)$, $H_{\varphi}^{i}(M,N)$, and $H_{\mathcal{\mathfrak{a},\mathfrak{b}}}^{i}(M,N)$ of the local cohomology module $H_{\mathcal{\mathfrak{a}}}^{i}(M)$ of an $R$-module $M$, are special cases of the local cohomology module $H_{\mathcal{Z}}^{i}(X)$ of an $R$-complex $X$ with support in a stable under specialization subset $\mathcal{Z}$ of $\Spec(R)$. Therefore, any established result on $H_{\mathcal{Z}}^{i}(X)$ encompasses all the previously known results on each of these generalized local cohomology modules. In this direction, we aspire to define the general notion of $\mathcal{Z}$-cofiniteness and establish the analogous results for the local cohomology modules $H_{\mathcal{Z}}^{i}(X)$.

\section{Containment in Serre Classes}

An excursion among the results \cite[Theorem 2.10]{AM1}, \cite[Theorem 2.9]{AM2}, \cite[Corollary 3.1]{BA}, \cite[Lemma 2.1]{BKN}, \cite[Proposition 1 and Corollary 1]{DM}, \cite[Lemma 4.2]{HK}, \cite[Theorem 2.1]{Me2}, \cite[Propositions 7.1, 7.2 and 7.4]{WW}, and \cite[Lemma 1.2]{Yo} reveals an in-depth connection between local cohomology, Ext modules, Tor modules, and Koszul homology in terms of their containment in a Serre class of modules. The purpose of this section is to bring the local homology into play and uncover the true connection between all these homology and cohomology theories, and consequently, to illuminate and enhance the aforementioned results.

To shed some light on this revelation, we observe that given elements $\underline{a}= a_{1},\ldots,a_{n} \in R$ and an $R$-module $M$, the Koszul homology $H_{i}(\underline{a};M)=0$ for every $i < 0$ or $i > n$, and further
$$H_{n}(\underline{a};M)\cong \left(0:_{M}(\underline{a})\right) \cong \Ext_{R}^{0}\left(R/(\underline{a}),M\right)$$
and
$$H_{0}(\underline{a};M)\cong M/(\underline{a})M \cong \Tor^{R}_{0}\left(R/(\underline{a}),M\right).$$
These isomorphisms suggest the existence of an intimate connection between the Koszul homology on the one hand, and Ext and Tor modules on the other hand. The connection seems to manifest itself as in the following casual diagram for any given integer $s \geq 0$:\\

\begin{center}
\begin{tikzpicture}
  \node(1) at (-2.23,1.5) {$\underbrace{\Ext_{R}^{0}\left(R/(\underline{a}),M\right),\ldots,\Ext_{R}^{s}\left(R/(\underline{a}),M\right)}$};
  \node(2) at (-2.23,0.78) {$\updownarrow$};
  \node(3) at (0,0) {$\overbrace{H_{n}(\underline{a};M),\ldots,H_{n-s}(\underline{a};M)},\ldots,\underbrace{H_{s}(\underline{a};M),\ldots,H_{0}(\underline{a};M)}$};
  \node(4) at (2.43,-0.78) {$\updownarrow$};
  \node(5) at (2.43,-1.5) {$\overbrace{\Tor_{s}^{R}\left(R/(\underline{a}),M\right),\ldots,\Tor_{0}^{R}\left(R/(\underline{a}),M\right)}$};
\end{tikzpicture}
\end{center}

As the above diagram depicts, the Koszul homology acts as a connecting device in attaching the Ext modules to Tor modules in an elegant way. One such connection, one may surmise, could be the containment in some Serre class of modules that we are about to probe in this section.

First and foremost, we recall the notion of a Serre class of modules.

\begin{definition} \label{2.2.1}
A class $\mathcal{S}$ of $R$-modules is said to be a \textit{Serre class} \index{Serre class of modules} if, given any short exact sequence
$$0 \rightarrow M^{\prime} \rightarrow M \rightarrow M^{\prime\prime} \rightarrow 0$$
of $R$-modules, we have $M \in \mathcal{S}$ if and only if $M^{\prime}, M^{\prime\prime} \in \mathcal{S}$.
\end{definition}

The following example showcases the stereotypical instances of Serre classes.

\begin{example} \label{2.2.2}
Let $\mathfrak{a}$ be an ideal of $R$. Then the following classes of modules are Serre classes:
\begin{enumerate}
\item[(i)] The zero class.
\item[(ii)] The class of all noetherian $R$-modules.
\item[(iii)] The class of all artinian $R$-modules.
\item[(iv)] The class of all minimax $R$-modules.
\item[(v)] The class of all minimax and $\mathfrak{a}$-cofinite $R$-modules; see \cite[Corollary 4.4]{Me1}.
\item[(vi)] The class of all weakly Laskerian $R$-modules.
\item[(v)] The class of all Matlis reflexive $R$-modules.
\end{enumerate}
\end{example}

We will be concerned with a change of rings when dealing with Serre classes. Therefore, we adopt the following notion of Serre property to exclude any possible ambiguity in the statement of the results.

\begin{definition} \label{2.2.3}
A property $\mathcal{P}$ concerning modules is said to be a \textit{Serre property} if
$$\mathcal{S}_{\mathcal{P}}(R):= \left\{M\in \mathcal{M}(R) \suchthat M \text{ satisfies the property } \mathcal{P} \right\}$$
is a Serre class for every ring $R$.
\end{definition}

Given a Serre property, there naturally arise the corresponding notions of depth and width.

\begin{definition} \label{2.2.4}
Let $\mathfrak{a}$ be an ideal of $R$, and $M$ an $R$-module. Let $\mathcal{P}$ be a Serre property. Then:
\begin{enumerate}
\item[(i)] Define the $\mathcal{P}$-\textit{depth} \index{depth} of $M$ with respect to $\mathfrak{a}$ to be
\begin{center}
$\mathcal{P}$-$\depth_{R}(\mathfrak{a},M):= \inf \left\{i\geq 0 \suchthat \Ext^{i}_{R}(R/\mathfrak{a},M) \notin \mathcal{S}_{\mathcal{P}}(R) \right\}.$
\end{center}
\item[(ii)] Define the $\mathcal{P}$-\textit{width} \index{width} of $M$ with respect to $\mathfrak{a}$ to be
\begin{center}
$\mathcal{P}$-$\width_{R}(\mathfrak{a},M):= \inf \left\{i\geq 0 \suchthat \Tor^{R}_{i}(R/\mathfrak{a},M) \notin \mathcal{S}_{\mathcal{P}}(R) \right\}.$
\end{center}
\end{enumerate}
\end{definition}

It is clear that upon letting $\mathcal{P}$ be the Serre property of being zero, we recover the classical notions of depth and width.

We next remind the definition of local cohomology and local homology functors. Let $\mathfrak{a}$ be an ideal of $R$. We let
$$\Gamma_{\mathfrak{a}}(M):=\left\{x\in M \suchthat \mathfrak{a}^{t}x=0 \text{ for some } t\geq 1 \right\}$$
for an $R$-module $M$, and $\Gamma_{\mathfrak{a}}(f):=f|_{\Gamma_{\mathfrak{a}}(M)}$ for an $R$-homomorphism $f:M\rightarrow N$. This provides us with the so-called $\mathfrak{a}$-torsion functor $\Gamma_{\mathfrak{a}}(-)$ on the category of $R$-modules. The $i$th local cohomology functor with respect to $\mathfrak{a}$ is defined to be
$$H^{i}_{\mathfrak{a}}(-):= R^{i}\Gamma_{\mathfrak{a}}(-)$$
for every $i \geq 0$. In addition, the cohomological dimension \index{cohomological dimension of an ideal} of $\mathfrak{a}$ with respect to $M$ is
$$\cd(\mathfrak{a},M) := \sup \left\{i \geq 0 \suchthat H^{i}_{\mathfrak{a}}(M)\neq 0 \right\}.$$

Recall that an $R$-module $M$ is said to be $\mathfrak{a}$-torsion if $M=\Gamma_{\mathfrak{a}}(M)$. It is well-known that any $\mathfrak{a}$-torsion $R$-module $M$ enjoys a natural $\widehat{R}^{\mathfrak{a}}$-module structure in such a way that its lattices of $R$-submodules and $\widehat{R}^{\mathfrak{a}}$-submodules coincide. Besides, we have an isomorphism $\widehat{R}^{\mathfrak{a}}\otimes_{R}M \cong M$ both as $R$-modules and $\widehat{R}^{\mathfrak{a}}$-modules.

Likewise, we let
$$\Lambda^{\mathfrak{a}}(M):=\widehat{M}^{\mathfrak{a}}= \underset{n}\varprojlim (M/\mathfrak{a}^{n}M)$$
for an $R$-module $M$, and $\Lambda^{\mathfrak{a}}(f):=\widehat{f}$ for an $R$-homomorphism $f:M\rightarrow N$. This provides us with the so-called $\mathfrak{a}$-adic completion functor $\Lambda^{\mathfrak{a}}(-)$ on the category of $R$-modules. The $i$th local homology functor with respect to $\mathfrak{a}$ is defined to be
$$H^{\mathfrak{a}}_{i}(-):= L_{i}\Lambda^{\mathfrak{a}}(-)$$
for every $i \geq 0$. Besides, the homological dimension \index{homological dimension of an ideal} of $\mathfrak{a}$ with respect to $M$ is
$$\hd(\mathfrak{a},M) := \sup \left\{i \geq 0 \suchthat H^{\mathfrak{a}}_{i}(M)\neq 0 \right\}.$$

The next conditions are also required to be imposed on the Serre classes when we intend to bring the local homology and local cohomology into the picture.

\begin{definition} \label{2.2.5}
Let $\mathcal{P}$ be a Serre property, and $\mathfrak{a}$ an ideal of $R$. We say that $\mathcal{P}$ satisfies the condition $\mathfrak{D}_{\mathfrak{a}}$ if the following statements hold:
\begin{enumerate}
\item[(i)] If $R$ is $\mathfrak{a}$-adically complete and $M/\mathfrak{a}M \in \mathcal{S}_{\mathcal{P}}(R)$ for some $R$-module $M$, then $H^{\mathfrak{a}}_{0}(M) \in \mathcal{S}_{\mathcal{P}}(R)$.
\item[(ii)] For any $\mathfrak{a}$-torsion $R$-module $M$, we have $M\in \mathcal{S}_{\mathcal{P}}(R)$ if and only if $M\in \mathcal{S}_{\mathcal{P}}\left(\widehat{R}^{\mathfrak{a}}\right)$.
\end{enumerate}
\end{definition}

\begin{example} \label{2.2.6}
Let $\mathfrak{a}$ be an ideal of $R$. Then we have:
\begin{enumerate}
\item[(i)] The Serre property of being zero satisfies the condition $\mathfrak{D}_{\mathfrak{a}}$. Use \cite[Lemma 5.1 (iii)]{Si1}.
\item[(ii)] The Serre property of being noetherian satisfies the condition $\mathfrak{D}_{\mathfrak{a}}$. Refer to the next section.
\end{enumerate}
\end{example}

Melkersson \cite[Definition 2.1]{AM3} defines the condition  $\mathfrak{C}_{\mathfrak{a}}$  for Serre classes. Adopting his definition, we have:

\begin{definition} \label{2.2.7}
Let $\mathcal{P}$ be a Serre property, and $\mathfrak{a}$ an ideal of $R$. We say that $\mathcal{P}$ satisfies the condition $\mathfrak{C}_{\mathfrak{a}}$ if the containment $(0:_{M}\mathfrak{a}) \in \mathcal{S}_{\mathcal{P}}(R)$ for some $R$-module $M$, implies that $\Gamma_{\mathfrak{a}}(M) \in \mathcal{S}_{\mathcal{P}}(R)$.
\end{definition}

\begin{example} \label{2.2.8}
Let $\mathfrak{a}$ be an ideal of $R$. Then we have:
\begin{enumerate}
\item[(i)] The Serre property of being zero satisfies the condition $\mathfrak{C}_{\mathfrak{a}}$. Indeed, if $M$ is an $R$-module such that $(0:_{M}\mathfrak{a})=0$, then it can be easily seen that $\Gamma_{\mathfrak{a}}(M)=0$.
\item[(ii)] The Serre property of being artinian satisfies the condition $\mathfrak{C}_{\mathfrak{a}}$. This follows from the Melkersson's Criterion \cite[Theorem 1.3]{Me2}.
\item[(iii)] A Serre class which is closed under taking direct limits is called a torsion theory. There exists a whole bunch of torsion theories, e.g. given any $R$-module $L$, the class
    $$\mathfrak{T}_{L}:= \left\{M\in \mathcal{M}(R) \suchthat \Supp_{R}(M)\subseteq \Supp_{R}(L) \right\}$$
    is a torsion theory. One can easily see that any torsion theory satisfies the condition $\mathfrak{C}_{\mathfrak{a}}$.
\end{enumerate}
\end{example}

We finally describe the Koszul homology briefly. The Koszul complex $K^{R}(a)$ on an element $a \in R$ is the $R$-complex
$$K^{R}(a):=\Cone(R \xrightarrow{a} R),$$
and the Koszul complex $K^{R}(\underline{a})$ on a sequence of elements $\underline{a} = a_{1},\ldots,a_{n} \in R$ is the $R$-complex
$$K^{R}(\underline{a}):= K^{R}(a_{1}) \otimes_{R} \cdots \otimes_{R} K^{R}(a_{n}).$$
It is easy to see that $K^{R}(\underline{a})$ is a complex of finitely generated free $R$-modules concentrated in degrees $n,\ldots,0$.
Given any $R$-module $M$, there is an isomorphism of $R$-complexes
$$K^{R}(\underline{a})\otimes_{R}M \cong \Sigma^{n} \Hom_{R}\left(K^{R}(\underline{a}),M\right),$$
which is sometimes referred to as the self-duality property of the Koszul complex. Accordingly, we feel free to define the Koszul homology of the sequence $\underline{a}$ with coefficients in $M$, by setting
$$H_{i}(\underline{a};M):= H_{i}\left(K^{R}(\underline{a})\otimes_{R}M\right) \cong H_{i-n}\left(\Hom_{R}\left(K^{R}(\underline{a}),M\right)\right)$$
for every $i\geq 0$.

Now we are ready to study the containment of Koszul homology, Ext modules, Tor modules, local homology, and local cohomology in Serre classes.

We need a lemma.

\begin{lemma} \label{2.2.9}
Let $\mathfrak{a}$ be an ideal of $R$ and $M$ an $R$-module. Then the following assertions hold:
\begin{enumerate}
\item[(i)] There is a first quadrant spectral sequence \index{spectral sequence}
$$E^{2}_{p,q}= \Tor^{R}_{p}\left(R/\mathfrak{a},H^{\mathfrak{a}}_{q}(M)\right) \underset {p} \Rightarrow \Tor^{R}_{p+q}(R/\mathfrak{a},M).$$
\item[(ii)] There is a third quadrant spectral sequence
$$E^{2}_{p,q}= \Ext^{-p}_{R}\left(R/\mathfrak{a},H^{-q}_{\mathfrak{a}}(M)\right) \underset {p} \Rightarrow \Ext^{-p-q}_{\mathfrak{a}}(R/\mathfrak{a},M).$$
\end{enumerate}
\end{lemma}

\begin{prf}
(i): Let $\mathcal{F}=(R/\mathfrak{a})\otimes_{R}-$, and $\mathcal{G}=\Lambda^{\mathfrak{a}}(-)$. Then $\mathcal{F}$ is right exact, and
$\mathcal{G}(P)$ is left $\mathcal{F}$-acyclic for every projective $R$-module $P,$ since the $\mathfrak{a}$-adic completion of a flat
$R$-module is flat by \cite[1.4.7]{B}. Therefore, by \cite[Theorem 10.48]{Ro}, there is a Grothendieck spectral sequence
$$E^{2}_{p,q}= L_{p}\mathcal{F}\left(L_{q}\mathcal{G}(M)\right) \underset {p} \Rightarrow L_{p+q}(\mathcal{FG})(M).$$
Let $F$ be a free resolution of $M$. By \cite[Theorem 1.3.1]{B} or \cite[Theorem 15]{Mat1}, we have
$$(R/\mathfrak{a})\otimes_{R}\Lambda^{\mathfrak{a}}(F)\cong (R/\mathfrak{a})\otimes_{R}F,$$
whence
$$L_{p+q}(\mathcal{FG})(M) = H_{p+q}\left((\mathcal{FG})(F)\right) \cong \Tor^{R}_{p+q}(R/\mathfrak{a},M).$$

(ii): Similar to (i).
\end{prf}

In the proof of the following theorems, we use the straightforward observation that given elements $\underline{a}= a_{1},\dots,a_{n} \in R$, a finitely generated $R$-module $N$, and an $R$-module $M$, if $M$ belongs to a Serre class $\mathcal{S}$, then $H_{i}(\underline{a};M)\in \mathcal{S}$, $\Ext^{i}_{R}(N,M) \in \mathcal{S}$, and $\Tor^{R}_{i}(N,M) \in \mathcal{S}$ for every $i \geq 0$. For a proof refer to \cite[Lemma 2.1]{AT}.

\begin{theorem} \label{2.2.10}
Let $\mathfrak{a}=(a_{1},\ldots,a_{n})$ be an ideal of $R$, $\underline{a}=a_{1},\ldots,a_{n}$, and $M$ an $R$-module. Let $\mathcal{P}$ be a Serre property. Then the following three conditions are equivalent for any given $s \geq 0$:
\begin{enumerate}
\item[(i)] $H_{i}(\underline{a};M) \in \mathcal{S}_{\mathcal{P}}(R)$ for every $0 \leq i \leq s$.
\item[(ii)] $\Tor^{R}_{i}(N,M) \in \mathcal{S}_{\mathcal{P}}(R)$ for every finitely generated $R$-module $N$ with $\Supp_{R}(N)\subseteq V(\mathfrak{a})$, and for every $0 \leq i \leq s$.
\item[(iii)] $\Tor^{R}_{i}(N,M) \in \mathcal{S}_{\mathcal{P}}(R)$ for some finitely generated $R$-module $N$ with $\Supp_{R}(N)= V(\mathfrak{a})$, and for every $0 \leq i \leq s$.
\end{enumerate}
If in addition, $\mathcal{P}$ satisfies the condition $\mathfrak{D}_{\mathfrak{a}}$, then the above three conditions are equivalent to the following condition:
\begin{enumerate}
\item[(iv)] $H^{\mathfrak{a}}_{i}(M)\in \mathcal{S}_{\mathcal{P}}\left(\widehat{R}^{\mathfrak{a}}\right)$ for every $0 \leq i \leq s$.
\end{enumerate}
\end{theorem}

\begin{prf}
(i) $\Rightarrow$ (iii): Set $N=R/\mathfrak{a}$. Let $F$ be a free resolution of $R/\mathfrak{a}$ consisting of finitely generated $R$-modules. Then the $R$-complex $F\otimes_{R}M$ is isomorphic to an $R$-complex of the form
$$\cdots \rightarrow M^{s_{2}} \xrightarrow{\partial_{2}} M^{s_{1}} \xrightarrow{\partial_{1}} M^{s_{0}} \rightarrow 0.$$
We note that
\begin{equation} \label{eq:2.2.10.1}
\Tor_{0}^{R}(R/\mathfrak{a},M) \cong \coker \partial_{1} \cong H_{0}(\underline{a};M) \in \mathcal{S}_{\mathcal{P}}(R),
\end{equation}
by the assumption. If $s=0$, then we are done. Suppose that $s \geq 1$. The short exact sequence
$$0 \rightarrow \im \partial_{1} \rightarrow M^{s_{0}} \rightarrow \coker \partial_{1} \rightarrow 0,$$
induces the exact sequence
$$H_{i+1}\left(\underline{a};\coker \partial_{1}\right) \rightarrow H_{i}\left(\underline{a};\im \partial_{1}\right) \rightarrow H_{i}(\underline{a};M^{s_{0}}).$$
The assumption together with the display \eqref{eq:2.2.10.1} imply that the two lateral terms of the above exact sequence belong to $\mathcal{S}_{\mathcal{P}}(R)$, so $H_{i}\left(\underline{a};\im \partial_{1}\right) \in \mathcal{S}_{\mathcal{P}}(R)$ for every $0 \leq i \leq s$. The short exact sequence
$$0 \rightarrow \ker \partial_{1} \rightarrow M^{s_{1}} \rightarrow \im \partial_{1} \rightarrow 0,$$
yields the exact sequence
$$H_{i+1}\left(\underline{a};\im \partial_{1}\right) \rightarrow H_{i}\left(\underline{a};\ker \partial_{1}\right) \rightarrow H_{i}(\underline{a};M^{s_{1}}).$$
Therefore, $H_{i}\left(\underline{a};\ker \partial_{1}\right) \in \mathcal{S}_{\mathcal{P}}(R)$ for every $0 \leq i \leq s-1$. As $s \geq 1$, we see that $H_{0}\left(\underline{a};\ker \partial_{1}\right) \in \mathcal{S}_{\mathcal{P}}(R)$. The short exact sequence
\begin{equation} \label{eq:2.2.10.2}
0 \rightarrow \im \partial_{2} \rightarrow \ker \partial_{1} \rightarrow \Tor_{1}^{R}(R/\mathfrak{a},M) \rightarrow 0,
\end{equation}
implies the exact sequence
$$H_{0}\left(\underline{a};\ker \partial_{1}\right) \rightarrow H_{0}\left(\underline{a};\Tor_{1}^{R}(R/\mathfrak{a},M)\right) \rightarrow 0.$$
Therefore, $H_{0}\left(\underline{a};\Tor_{1}^{R}(R/\mathfrak{a},M)\right) \in \mathcal{S}_{\mathcal{P}}(R)$.
But $\mathfrak{a} \Tor^{R}_{1}(R/\mathfrak{a},M)=0$, so
$$\Tor^{R}_{1}(R/\mathfrak{a},M) \cong H_{0}\left(\underline{a};\Tor_{1}^{R}(R/\mathfrak{a},M)\right) \in \mathcal{S}_{\mathcal{P}}(R).$$
If $s=1$, then we are done. Suppose that $s \geq 2$. The short exact sequence \eqref{eq:2.2.10.2} induces the exact sequence
$$H_{i+1}\left(\underline{a};\Tor_{1}^{R}(R/\mathfrak{a},M)\right) \rightarrow H_{i}\left(\underline{a};\im \partial_{2}\right) \rightarrow H_{i}\left(\underline{a};\ker \partial_{1}\right).$$
It follows that $H_{i}\left(\underline{a};\im \partial_{2}\right) \in \mathcal{S}_{\mathcal{P}}(R)$ for every $0 \leq i \leq s-1$. The short exact sequence
$$0 \rightarrow \ker \partial_{2} \rightarrow M^{s_{2}} \rightarrow \im \partial_{2} \rightarrow 0,$$
yields the exact sequence
$$H_{i+1}\left(\underline{a};\im \partial_{2}\right) \rightarrow H_{i}\left(\underline{a};\ker \partial_{2}\right) \rightarrow H_{i}(\underline{a};M^{s_{2}}).$$
Thus $H_{i}\left(\underline{a};\ker \partial_{2}\right) \in \mathcal{S}_{\mathcal{P}}(R)$ for every $0 \leq i \leq s-2$.
As $s \geq 2$, we see that $H_{0}\left(\underline{a};\ker \partial_{2}\right) \in \mathcal{S}_{\mathcal{P}}(R)$. The short exact sequence
$$0 \rightarrow \im \partial_{3} \rightarrow \ker \partial_{2} \rightarrow \Tor_{2}^{R}\left(R/\mathfrak{a},M\right) \rightarrow 0,$$
yields the exact sequence
$$H_{0}\left(\underline{a};\ker \partial_{2}\right) \rightarrow H_{0}\left(\underline{a};\Tor_{2}^{R}(R/\mathfrak{a},M)\right) \rightarrow 0.$$
As before, we conclude that
$$\Tor^{R}_{2}(R/\mathfrak{a},M) \cong H_{0}\left(\underline{a};\Tor_{2}^{R}(R/\mathfrak{a},M)\right) \in \mathcal{S}_{\mathcal{P}}(R).$$
If $s=2$, then we are done. Proceeding in this manner, we see that $\Tor^{R}_{i}(R/\mathfrak{a},M) \in \mathcal{S}_{\mathcal{P}}(R)$ for every $0 \leq i \leq s$.

(iii) $\Rightarrow$ (ii): Let $L$ be a finitely generated $R$-module with $\Supp_{R}(L) \subseteq V(\mathfrak{a})$. By induction on $s$, we show that $\Tor_{i}^{R}(L,M)\in \mathcal{S}_{\mathcal{P}}(R)$ for every $0 \leq i \leq s$. Let $s=0$. Then $N \otimes_{R} M \in \mathcal{S}(R)$. Using Gruson's Filtration Lemma \cite[Theorem 4.1]{Va}, there is a finite filtration
$$0=L_{0} \subseteq L_{1} \subseteq \cdots \subseteq L_{t-1} \subseteq L_{t} = L,$$
such that $L_{j}/L_{j-1}$ is isomorphic to a quotient of a finite direct sum of copies of $N$ for every $1 \leq j \leq t$. Given any $1 \leq j \leq t$, the exact sequence
\begin{equation} \label{eq:2.2.10.3}
0 \rightarrow K_{j} \rightarrow N^{r_{j}} \rightarrow L_{j}/L_{j-1} \rightarrow 0,
\end{equation}
induces the exact sequence
$$N^{r_{j}} \otimes_{R} M \rightarrow (L_{j}/L_{j-1}) \otimes_{R} M \rightarrow 0.$$
It follows that $(L_{j}/L_{j-1}) \otimes_{R} M \in \mathcal{S}_{\mathcal{P}}(R)$.
Now the short exact sequence
\begin{equation} \label{eq:2.2.10.4}
0 \rightarrow L_{j-1} \rightarrow L_{j} \rightarrow L_{j}/L_{j-1} \rightarrow 0,
\end{equation}
yields the exact sequence
$$L_{j-1} \otimes_{R} M \rightarrow L_{j} \otimes_{R} M \rightarrow (L_{j}/L_{j-1}) \otimes_{R} M.$$
A successive use of the above exact sequence, letting $j=1,\ldots,t$, implies that $L\otimes_{R}M \in \mathcal{S}_{\mathcal{P}}(R)$.

Now let $s \geq 1$, and suppose that the results holds true for $s-1$. The induction hypothesis implies that $\Tor^{R}_{i}(L,M) \in \mathcal{S}_{\mathcal{P}}(R)$ for every $0 \leq i \leq s-1$. Hence, it suffices to show that $\Tor^{R}_{s}(L,M) \in \mathcal{S}_{\mathcal{P}}(R)$. Given any $1 \leq j \leq t$, the short exact sequence \eqref{eq:2.2.10.3} induces the exact sequence $$\Tor^{R}_{s}(N^{r_{j}},M) \rightarrow \Tor^{R}_{s}\left(L_{j}/L_{j-1},M\right) \rightarrow \Tor^{R}_{s-1}(K_{j},M).$$
The induction hypothesis shows that $\Tor^{R}_{s-1}(K_{j},M) \in \mathcal{S}_{\mathcal{P}}(R)$, so from the above exact sequence we get that $\Tor^{R}_{s}\left(L_{j}/L_{j-1},M\right) \in \mathcal{S}_{\mathcal{P}}(R)$. Now, the short exact sequence \eqref{eq:2.2.10.4} yields the exact sequence
$$\Tor^{R}_{s}(L_{j-1},M) \rightarrow \Tor^{R}_{s}(L_{j},M) \rightarrow \Tor^{R}_{s}\left(L_{j}/L_{j-1},M\right).$$
A successive use of the above exact sequence, letting $j=1,\ldots,t$, implies that $\Tor^{R}_{s}(L,M) \in \mathcal{S}_{\mathcal{P}}(R)$.

(ii) $\Rightarrow$ (i): The hypothesis implies that $\Tor^{R}_{i}(R/\mathfrak{a},M) \in \mathcal{S}_{\mathcal{P}}(R)$ for every $0 \leq i \leq s$. The $R$-complex $K^{R}(\underline{a})\otimes_{R}M$ is isomorphic to an $R$-complex of the form
$$0 \rightarrow M \xrightarrow{\partial_{n}} M^{t_{n-1}} \rightarrow \cdots \rightarrow M^{t_{2}} \xrightarrow{\partial_{2}} M^{t_{1}} \xrightarrow{\partial_{1}} M \rightarrow 0.$$
We have
\begin{equation} \label{eq:2.2.10.5}
H_{0}(\underline{a};M)\cong \coker \partial_{1} \cong \Tor^{R}_{0}(R/\mathfrak{a},M) \in \mathcal{S}_{\mathcal{P}}(R),
\end{equation}
by the assumption. If $s=0$, then we are done. Suppose that $s \geq 1$. The short exact sequence
$$0 \rightarrow \im \partial_{1} \rightarrow M \rightarrow \coker \partial_{1} \rightarrow 0,$$
induces the exact sequence
$$\Tor^{R}_{i+1}(R/\mathfrak{a},\coker \partial_{1}) \rightarrow \Tor^{R}_{i}(R/\mathfrak{a},\im \partial_{1}) \rightarrow \Tor^{R}_{i}(R/\mathfrak{a},M).$$
The assumption along with the display \eqref{eq:2.2.10.5} imply that the two lateral terms of the above exact sequence belong to $\mathcal{S}_{\mathcal{P}}(R)$, so $\Tor^{R}_{i}(R/\mathfrak{a},\im \partial_{1}) \in \mathcal{S}_{\mathcal{P}}(R)$ for every $0 \leq i \leq s$. The short exact sequence
$$0 \rightarrow \ker \partial_{1} \rightarrow M^{t_{1}} \rightarrow \im \partial_{1} \rightarrow 0,$$
yields the exact sequence
$$\Tor^{R}_{i+1}(R/\mathfrak{a},\im \partial_{1}) \rightarrow \Tor^{R}_{i}(R/\mathfrak{a},\ker \partial_{1}) \rightarrow \Tor^{R}_{i}(R/\mathfrak{a},M^{t_{1}}).$$
Therefore, $\Tor^{R}_{i}(R/\mathfrak{a},\ker \partial_{1}) \in \mathcal{S}_{\mathcal{P}}(R)$ for every $0 \leq i \leq s-1$.
As $s \geq 1$, we see that $(R/\mathfrak{a})\otimes_{R} \ker \partial_{1} \in \mathcal{S}_{\mathcal{P}}(R)$. The short exact sequence
\begin{equation} \label{eq:2.2.10.6}
0 \rightarrow \im \partial_{2} \rightarrow \ker \partial_{1} \rightarrow H_{1}(\underline{a};M) \rightarrow 0,
\end{equation}
implies the exact sequence
$$(R/\mathfrak{a})\otimes_{R} \ker \partial_{1} \rightarrow (R/\mathfrak{a})\otimes_{R} H_{1}(\underline{a};M) \rightarrow 0.$$
Therefore, $(R/\mathfrak{a})\otimes_{R} H_{1}(\underline{a};M) \in \mathcal{S}_{\mathcal{P}}(R)$.
But $\mathfrak{a} H_{1}(\underline{a};M)=0$, so
$$H_{1}(\underline{a};M) \cong (R/\mathfrak{a})\otimes_{R} H_{1}(\underline{a};M) \in \mathcal{S}_{\mathcal{P}}(R).$$
If $s=1$, then we are done. Suppose that $s \geq 2$. The short exact sequence \eqref{eq:2.2.10.6} induces the exact sequence
$$\Tor^{R}_{i+1}\left(R/\mathfrak{a},H_{1}(\underline{a};M)\right) \rightarrow \Tor^{R}_{i}(R/\mathfrak{a},\im \partial_{2}) \rightarrow \Tor^{R}_{i}(R/\mathfrak{a},\ker \partial_{1}).$$
It follows that $\Tor^{R}_{i}(R/\mathfrak{a},\im \partial_{2}) \in \mathcal{S}_{\mathcal{P}}(R)$ for every $0 \leq i \leq s-1$. The short exact sequence
$$0 \rightarrow \ker \partial_{2} \rightarrow M^{t_{2}} \rightarrow \im \partial_{2} \rightarrow 0,$$
yields the exact sequence
$$\Tor^{R}_{i+1}\left(R/\mathfrak{a},\im \partial_{2}\right) \rightarrow \Tor^{R}_{i}(R/\mathfrak{a},\ker \partial_{2}) \rightarrow \Tor^{R}_{i}(R/\mathfrak{a},M^{t_{2}}).$$
Thus $\Tor^{R}_{i}(R/\mathfrak{a},\ker \partial_{2}) \in \mathcal{S}_{\mathcal{P}}(R)$ for every $0 \leq i \leq s-2$. As $s \geq 2$, we see that $(R/\mathfrak{a}) \otimes_{R} \ker \partial_{2} \in \mathcal{S}_{\mathcal{P}}(R)$. The short exact sequence
$$0 \rightarrow \im \partial_{3} \rightarrow \ker \partial_{2} \rightarrow H_{2}(\underline{a};M) \rightarrow 0,$$
implies the exact sequence
$$(R/\mathfrak{a}) \otimes_{R} \ker \partial_{2} \rightarrow (R/\mathfrak{a}) \otimes_{R} H_{2}(\underline{a};M) \rightarrow 0.$$
As before, we conclude that
$$H_{2}(\underline{a};M) \cong (R/\mathfrak{a}) \otimes_{R} H_{2}(\underline{a};M) \in \mathcal{S}_{\mathcal{P}}(R).$$
If $s=2$, then we are done. Proceeding in this manner, we infer that $H_{i}(\underline{a};M) \in \mathcal{S}_{\mathcal{P}}(R)$ for every $0 \leq i \leq s$.

Now, assume that the Serre property $\mathcal{P}$ satisfies the condition $\mathfrak{D}_{\mathfrak{a}}$. We first note that since $\Tor^{R}_{i}(R/\mathfrak{a},M)$ is an $\mathfrak{a}$-torsion $R$-module, it has an $\widehat{R}^{\mathfrak{a}}$-module structure such that
$$\Tor^{R}_{i}(R/\mathfrak{a},M) \cong \widehat{R}^{\mathfrak{a}} \otimes_{R} \Tor^{R}_{i}(R/\mathfrak{a},M) \cong \Tor^{\widehat{R}^{\mathfrak{a}}}_{i}\left(\widehat{R}^{\mathfrak{a}}/\mathfrak{a}\widehat{R}^{\mathfrak{a}},\widehat{R}^{\mathfrak{a}} \otimes_{R} M\right),$$
for every $i \geq 0$ both as $R$-modules and $\widehat{R}^{\mathfrak{a}}$-modules. Moreover, by \cite[Lemma 2.3]{Si2}, we have
$$H^{\mathfrak{a}}_{i}(M)\cong H^{\mathfrak{a}\widehat{R}^{\mathfrak{a}}}_{i}\left(\widehat{R}^{\mathfrak{a}}\otimes_{R}M\right),$$
for every $i \geq 0$ both as $R$-modules and $\widehat{R}^{\mathfrak{a}}$-modules. With this preparation we prove:

(iv) $\Rightarrow$ (iii): We let $N=R/\mathfrak{a}$, and show that $\Tor_{i}^{R}(R/\mathfrak{a},M)\in \mathcal{S}_{\mathcal{P}}(R)$ for every $0 \leq i \leq s$. First suppose that $R$ is $\mathfrak{a}$-adically complete. By Lemma \ref{2.2.9}, there is a first quadrant spectral sequence
\begin{equation} \label{eq:2.2.10.7}
E^{2}_{p,q}= \Tor^{R}_{p}\left(R/\mathfrak{a},H^{\mathfrak{a}}_{q}(M)\right) \underset {p} \Rightarrow \Tor^{R}_{p+q}(R/\mathfrak{a},M).
\end{equation}
The hypothesis implies that $E^{2}_{p,q}\in \mathcal{S}_{\mathcal{P}}(R)$ for every $p \geq 0$ and $0 \leq q \leq s$. Let $0 \leq i \leq s$. There is a finite filtration
$$0=U^{-1} \subseteq U^{0} \subseteq \cdots \subseteq U^{i} = \Tor^{R}_{i}(R/\mathfrak{a},M),$$
such that $U^{p}/U^{p-1}\cong E^{\infty}_{p,i-p}$ for every $0 \leq p \leq i$. Since $E^{\infty}_{p,i-p}$ is a subquotient of $E^{2}_{p,i-p}$ and $0 \leq i-p \leq s$, we infer that
$$U^{p}/U^{p-1}\cong E^{\infty}_{p,i-p}\in \mathcal{S}_{\mathcal{P}}(R)$$
for every $0 \leq p \leq i$.
A successive use of the short exact sequence
$$0 \rightarrow U^{p-1} \rightarrow U^{p} \rightarrow U^{p}/U^{p-1} \rightarrow 0,$$
by letting $p=0,\ldots,i$, implies that $\Tor^{R}_{i}(R/\mathfrak{a},M)\in \mathcal{S}_{\mathcal{P}}(R)$.

Now, consider the general case. Since
$$H^{\mathfrak{a}\widehat{R}^{\mathfrak{a}}}_{i}\left(\widehat{R}^{\mathfrak{a}}\otimes_{R}M\right) \cong H^{\mathfrak{a}}_{i}(M) \in \mathcal{S}_{\mathcal{P}}\left(\widehat{R}^{\mathfrak{a}}\right),$$
for every $0 \leq i \leq s$, the special case implies that
$$\Tor^{R}_{i}(R/\mathfrak{a},M) \cong \Tor^{\widehat{R}^{\mathfrak{a}}}_{i}\left(\widehat{R}^{\mathfrak{a}}/\mathfrak{a}\widehat{R}^{\mathfrak{a}},\widehat{R}^{\mathfrak{a}} \otimes_{R} M\right)\in \mathcal{S}_{\mathcal{P}}\left(\widehat{R}^{\mathfrak{a}}\right),$$
for every $0 \leq i \leq s$. However, $\Tor^{R}_{i}(R/\mathfrak{a},M)$ is $\mathfrak{a}$-torsion, so the condition $\mathfrak{D}_{\mathfrak{a}}$ implies that $\Tor^{R}_{i}(R/\mathfrak{a},M)\in \mathcal{S}_{\mathcal{P}}(R)$ for every $0 \leq i \leq s$.

(ii) $\Rightarrow$ (iv): It follows from the hypothesis that $\Tor^{R}_{i}(R/\mathfrak{a},M)\in \mathcal{S}_{\mathcal{P}}(R)$ for every $0 \leq i \leq s$. First suppose that $R$ is $\mathfrak{a}$-adically complete. We argue by induction on $s$. Let $s=0$. Then $M/\mathfrak{a}M \in \mathcal{S}_{\mathcal{P}}(R)$, whence $H^{\mathfrak{a}}_{0}(M)\in \mathcal{S}_{\mathcal{P}}(R)$ by the condition $\mathfrak{D}_{\mathfrak{a}}$.

Now let $s \geq 1$, and suppose that the result holds true for $s-1$. The induction hypothesis implies that $H^{\mathfrak{a}}_{i}(M)\in \mathcal{S}_{\mathcal{P}}(R)$ for every $0 \leq i \leq s-1$. Hence, it suffices to show that $H^{\mathfrak{a}}_{s}(M)\in \mathcal{S}_{\mathcal{P}}(R)$.
Consider the spectral sequence \eqref{eq:2.2.10.7}. By the hypothesis, $E^{2}_{p,q}\in \mathcal{S}_{\mathcal{P}}(R)$ for every $p \in \mathbb{Z}$ and $0 \leq q \leq s-1$.
There is a finite filtration
$$0=U^{-1} \subseteq U^{0} \subseteq \cdots \subseteq U^{s} = \Tor^{R}_{s}(R/\mathfrak{a},M),$$
such that $U^{p}/U^{p-1}\cong E^{\infty}_{p,s-p}$ for every $0 \leq p \leq s$. As $\Tor^{R}_{s}(R/\mathfrak{a},M)\in \mathcal{S}_{\mathcal{P}}(R)$,
we conclude that
$$E^{\infty}_{0,s}\cong U^{0}/U^{-1} \cong U^{0}\in \mathcal{S}_{\mathcal{P}}(R).$$
Let $r \geq 2$, and consider the differentials
$$E^{r}_{r,s-r+1}\xrightarrow {d^{r}_{r,s-r+1}} E^{r}_{0,s}\xrightarrow {d^{r}_{0,s}} E^{r}_{-r,s+r-1}=0.$$
Since $s-r+1 \leq s-1$ and $E^{r}_{r,s-r+1}$ is a subquotient of $E^{2}_{r,s-r+1}$, the hypothesis implies that $E^{r}_{r,s-r+1}\in \mathcal{S}_{\mathcal{P}}(R)$, and consequently $\im d^{r}_{r,s-r+1}\in \mathcal{S}_{\mathcal{P}}(R)$ for every $r\geq 2$.
We thus obtain
$$E^{r+1}_{0,s} \cong \ker d^{r}_{0,s}/ \im d^{r}_{r,s-r+1} = E^{r}_{0,s}/ \im d^{r}_{r,s-r+1},$$
and consequently a short exact sequence
\begin{equation} \label{eq:2.2.10.8}
0 \rightarrow \im d^{r}_{r,s-r+1} \rightarrow E^{r}_{0,s} \rightarrow E^{r+1}_{0,s} \rightarrow 0.
\end{equation}
There is an integer $r_{0}\geq 2$, such that $E^{\infty}_{0,s} \cong E^{r+1}_{0,s}$ for every $r\geq r_{0}$. It follows that $E^{r_{0}+1}_{0,s}\in \mathcal{S}_{\mathcal{P}}(R)$.
Now, the short exact sequence \eqref{eq:2.2.10.8} implies that $E^{r_{0}}_{0,s}\in \mathcal{S}_{\mathcal{P}}(R)$. Using the short exact sequence
\eqref{eq:2.2.10.8} inductively, we conclude that
$$H^{\mathfrak{a}}_{s}(M)/\mathfrak{a}H^{\mathfrak{a}}_{s}(M)\cong E^{2}_{0,s}\in \mathcal{S}_{\mathcal{P}}(R).$$
Therefore, by Lemma \ref{2.3.7} and the condition $\mathfrak{D}_{\mathfrak{a}}$, we get
$$H^{\mathfrak{a}}_{s}(M) \cong H^{\mathfrak{a}}_{0}\left(H^{\mathfrak{a}}_{s}(M)\right)\in \mathcal{S}_{\mathcal{P}}(R).$$

Now, consider the general case. Since $\Tor^{R}_{i}(R/\mathfrak{a},M)$ is an $\mathfrak{a}$-torsion $R$-module such that $\Tor^{R}_{i}(R/\mathfrak{a},M)\in \mathcal{S}_{\mathcal{P}}(R)$ for every $0 \leq i \leq s$, we deduce that
$$\Tor^{\widehat{R}^{\mathfrak{a}}}_{i}\left(\widehat{R}^{\mathfrak{a}}/\mathfrak{a}\widehat{R}^{\mathfrak{a}},\widehat{R}^{\mathfrak{a}} \otimes_{R} M\right) \cong \Tor^{R}_{i}(R/\mathfrak{a},M) \in \mathcal{S}_{\mathcal{P}}\left(\widehat{R}^{\mathfrak{a}}\right),$$
for every $0 \leq i \leq s$. But the special case yields that
$$H^{\mathfrak{a}}_{i}(M) \cong H^{\mathfrak{a}\widehat{R}^{\mathfrak{a}}}_{i}\left(\widehat{R}^{\mathfrak{a}}\otimes_{R}M\right) \in \mathcal{S}_{\mathcal{P}}\left(\widehat{R}^{\mathfrak{a}}\right),$$
for every $0 \leq i \leq s$.
\end{prf}

The following special case may be of independent interest.

\begin{corollary} \label{2.2.11}
Let $\mathfrak{a}=(a_{1},\ldots,a_{n})$ be an ideal of $R$, $\underline{a}=a_{1},\ldots,a_{n}$, and $M$ an $R$-module. Let $\mathcal{P}$ be a Serre property satisfying the condition $\mathfrak{D}_{\mathfrak{a}}$. Then the following conditions are equivalent:
\begin{enumerate}
\item[(i)] $M/\mathfrak{a}M \in \mathcal{S}_{\mathcal{P}}(R)$.
\item[(ii)] $H^{\mathfrak{a}}_{0}(M) \in \mathcal{S}_{\mathcal{P}}\left(\widehat{R}^{\mathfrak{a}}\right)$.
\item[(iii)] $\widehat{M}^{\mathfrak{a}}\in \mathcal{S}_{\mathcal{P}}\left(\widehat{R}^{\mathfrak{a}}\right)$.
\end{enumerate}
\end{corollary}

\begin{prf}
(i) $\Leftrightarrow$ (ii): Follows from Theorem \ref{2.2.10} upon letting $s=0$.

(ii) $\Rightarrow$ (iii): By \cite[Lemma 5.1 (i)]{Si1}, the natural homomorphism $H^{\mathfrak{a}}_{0}(M) \rightarrow \widehat{M}^{\mathfrak{a}}$ is surjective. Thus the result follows.

(iii) $\Rightarrow$ (i): Since $\widehat{M}^{\mathfrak{a}}\in \mathcal{S}_{\mathcal{P}}\left(\widehat{R}^{\mathfrak{a}}\right)$, we see that $\widehat{M}^{\mathfrak{a}}/\mathfrak{a}\widehat{M}^{\mathfrak{a}} \in \mathcal{S}_{\mathcal{P}}\left(\widehat{R}^{\mathfrak{a}}\right)$. However, by \cite[Theorem 1.1]{Si1}, we have $\widehat{M}^{\mathfrak{a}}/\mathfrak{a}\widehat{M}^{\mathfrak{a}} \cong M/\mathfrak{a}M$. It follows that $M/\mathfrak{a}M \in \mathcal{S}_{\mathcal{P}}\left(\widehat{R}^{\mathfrak{a}}\right)$. But $M/\mathfrak{a}M$ is $\mathfrak{a}$-torsion, so by the condition $\mathfrak{D}_{\mathfrak{a}}$, we have $M/\mathfrak{a}M \in \mathcal{S}_{\mathcal{P}}(R)$.
\end{prf}

Using similar arguments, we can prove the dual result to Theorem \ref{2.2.10}. It is worth noting that the following result is proved in \cite[Theorem 2.9]{AM3}. However, the condition $\mathfrak{C}_{\mathfrak{a}}$ is assumed to be satisfied in all four statements, but here we only require that the condition $\mathfrak{C}_{\mathfrak{a}}$ is satisfied for the last statement. Moreover, the techniques used there are quite different than those used here. We include a proof for the sake of completeness.

\begin{theorem} \label{2.2.12}
Let $\mathfrak{a}=(a_{1},\ldots,a_{n})$ be an ideal of $R$, $\underline{a}=a_{1},\ldots,a_{n}$, and $M$ an $R$-module. Let $\mathcal{P}$ be a Serre property. Then the following three conditions are equivalent for any given $s \geq 0$:
\begin{enumerate}
\item[(i)] $H_{n-i}(\underline{a};M) \in \mathcal{S}_{\mathcal{P}}(R)$ for every $0 \leq i \leq s$.
\item[(ii)] $\Ext^{i}_{R}(N,M) \in \mathcal{S}_{\mathcal{P}}(R)$ for every finitely generated $R$-module $N$ with $\Supp_{R}(N)\subseteq V(\mathfrak{a})$, and for every $0 \leq i \leq s$.
\item[(iii)] $\Ext^{i}_{R}(N,M) \in \mathcal{S}_{\mathcal{P}}(R)$ for some finitely generated $R$-module $N$ with $\Supp_{R}(N)= V(\mathfrak{a})$, and for every $0 \leq i \leq s$.
\end{enumerate}
If in addition, $\mathcal{P}$ satisfies the condition $\mathfrak{C}_{\mathfrak{a}}$, then the above three conditions are equivalent to the following condition:
\begin{enumerate}
\item[(iv)] $H^{i}_{\mathfrak{a}}(M)\in \mathcal{S}_{\mathcal{P}}(R)$ for every $0 \leq i \leq s$.
\end{enumerate}
\end{theorem}

\begin{prf}
(i) $\Rightarrow$ (iii): Set $N=R/\mathfrak{a}$. Let $F$ be a free resolution of $R/\mathfrak{a}$ consisting of finitely generated $R$-modules. Then the $R$-complex $\Hom_{R}(F,M)$ is isomorphic to an $R$-complex of the form
$$0 \rightarrow M^{s_{0}} \xrightarrow{\partial_{0}} M^{s_{1}} \xrightarrow{\partial_{-1}} M^{s_{2}} \rightarrow \cdots.$$
We note that
\begin{equation} \label{eq:2.2.12.1}
\Ext_{R}^{0}(R/\mathfrak{a},M) \cong \ker \partial_{0} \cong H_{n}(\underline{a};M) \in \mathcal{S}_{\mathcal{P}}(R),
\end{equation}
by the assumption. If $s=0$, then we are done. Suppose that $s \geq 1$. The short exact sequence
$$0 \rightarrow \ker \partial_{0} \rightarrow M^{s_{0}} \rightarrow \im \partial_{0} \rightarrow 0,$$
induces the exact sequence
$$H_{n-i}\left(\underline{a};M^{s_{0}}\right) \rightarrow H_{n-i}\left(\underline{a};\im \partial_{0}\right) \rightarrow H_{n-i-1}(\underline{a};\ker \partial_{0}).$$
The assumption together with the display \eqref{eq:2.2.12.1} imply that the two lateral terms of the above exact sequence belong to $\mathcal{S}_{\mathcal{P}}(R)$, so $H_{n-i}\left(\underline{a};\im \partial_{0}\right) \in \mathcal{S}_{\mathcal{P}}(R)$ for every $0 \leq i \leq s$. As $s \geq 1$, we have $H_{n-1}\left(\underline{a};\im \partial_{0}\right) \in \mathcal{S}_{\mathcal{P}}(R)$. Moreover, since $\ker \partial_{-1} \subseteq M^{s_{1}}$, we get the exact sequence
$$0 \rightarrow H_{n}(\underline{a};\ker \partial_{-1}) \rightarrow H_{n}(\underline{a};M^{s_{1}}).$$
But $H_{n}(\underline{a};M^{s_{1}}) \in \mathcal{S}_{\mathcal{P}}(R)$, so $H_{n}(\underline{a};\ker \partial_{-1}) \in \mathcal{S}_{\mathcal{P}}(R)$.
The short exact sequence
\begin{equation} \label{eq:2.2.12.2}
0 \rightarrow \im \partial_{0} \rightarrow \ker \partial_{-1} \rightarrow \Ext^{1}_{R}(R/\mathfrak{a},M) \rightarrow 0,
\end{equation}
yields the exact sequence
$$H_{n}\left(\underline{a};\ker \partial_{-1}\right) \rightarrow H_{n}\left(\underline{a};\Ext^{1}_{R}(R/\mathfrak{a},M)\right) \rightarrow H_{n-1}(\underline{a};\im \partial_{0}).$$
Therefore, $H_{n}\left(\underline{a};\Ext^{1}_{R}(R/\mathfrak{a},M)\right) \in \mathcal{S}_{\mathcal{P}}(R)$. But $\mathfrak{a} \Ext^{1}_{R}(R/\mathfrak{a},M)=0$, so
$$\Ext^{1}_{R}(R/\mathfrak{a},M) \cong H_{n}\left(\underline{a};\Ext^{1}_{R}(R/\mathfrak{a},M)\right) \in \mathcal{S}_{\mathcal{P}}(R).$$
If $s=1$, then we are done. Suppose that $s \geq 2$.
The short exact sequence \eqref{eq:2.2.12.2} induces the exact sequence
$$H_{n-i}\left(\underline{a};\im \partial_{0}\right) \rightarrow H_{n-i}\left(\underline{a};\ker \partial_{-1}\right) \rightarrow H_{i}\left(\underline{a};\Ext^{1}_{R}(R/\mathfrak{a},M)\right).$$
It follows that $H_{n-i}\left(\underline{a};\ker \partial_{-1}\right) \in \mathcal{S}_{\mathcal{P}}(R)$ for every $0 \leq i \leq s$.
The short exact sequence
$$0 \rightarrow \ker \partial_{-1} \rightarrow M^{s_{1}} \rightarrow \im \partial_{-1} \rightarrow 0,$$
yields the exact sequence
$$H_{n-i}\left(\underline{a};M^{s_{1}}\right) \rightarrow H_{n-i}\left(\underline{a};\im \partial_{-1}\right) \rightarrow H_{n-i-1}\left(\underline{a};\ker \partial_{-1}\right).$$
Thus $H_{n-i}\left(\underline{a};\im \partial_{-1}\right) \in \mathcal{S}_{\mathcal{P}}(R)$ for every $0 \leq i \leq s-1$.
As $s \geq 2$, we see that $H_{n-1}\left(\underline{a};\im \partial_{-1}\right) \in \mathcal{S}_{\mathcal{P}}(R)$. Moreover, since $\ker \partial_{-2} \subseteq M^{s_{2}}$, we get the exact sequence
$$0 \rightarrow H_{n}\left(\underline{a};\ker \partial_{-2}\right) \rightarrow H_{n}\left(\underline{a};M^{s_{2}}\right).$$
But $H_{n}\left(\underline{a};M^{s_{2}}\right) \in \mathcal{S}_{\mathcal{P}}(R)$, so $H_{n}\left(\underline{a};\ker \partial_{-2}\right) \in \mathcal{S}_{\mathcal{P}}(R)$.
The short exact sequence
$$0 \rightarrow \im \partial_{-1} \rightarrow \ker \partial_{-2} \rightarrow \Ext_{R}^{2}\left(R/\mathfrak{a},M\right) \rightarrow 0,$$
yields the exact sequence
$$H_{n}\left(\underline{a};\ker \partial_{-2}\right) \rightarrow H_{n}\left(\underline{a};\Ext_{R}^{2}(R/\mathfrak{a},M)\right) \rightarrow H_{n-1}\left(\underline{a};\im \partial_{-1}\right).$$
Therefore, $H_{n}\left(\underline{a};\Ext_{R}^{2}(R/\mathfrak{a},M)\right) \in \mathcal{S}_{\mathcal{P}}(R)$.
As before, we conclude that
$$\Ext^{2}_{R}(R/\mathfrak{a},M) \cong H_{n}\left(\underline{a};\Ext_{R}^{2}(R/\mathfrak{a},M)\right) \in \mathcal{S}_{\mathcal{P}}(R).$$
If $s=2$, then we are done. Proceeding in this manner, we see that $\Ext^{i}_{R}(R/\mathfrak{a},M) \in \mathcal{S}_{\mathcal{P}}(R)$ for every $0 \leq i \leq s$.

(iii) $\Rightarrow$ (ii): Let $L$ be a finitely generated $R$-module with $\Supp_{R}(L) \subseteq V(\mathfrak{a})$. By induction on $s$, we show that $\Ext_{R}^{i}(L,M)\in \mathcal{S}_{\mathcal{P}}(R)$ for every $0 \leq i \leq s$. Let $s=0$. Then $\Hom_{R}(N,M) \in \mathcal{S}(R)$. Using Gruson's Filtration Lemma \cite[Theorem 4.1]{Va}, there is a finite filtration
$$0=L_{0} \subseteq L_{1} \subseteq \cdots \subseteq L_{t-1} \subseteq L_{t} = L,$$
such that $L_{j}/L_{j-1}$ is isomorphic to a quotient of a finite direct sum of copies of $N$ for every $1 \leq j \leq t$. Given any $1 \leq j \leq t$, the exact sequence
\begin{equation} \label{eq:2.2.12.3}
0 \rightarrow K_{j} \rightarrow N^{r_{j}} \rightarrow L_{j}/L_{j-1} \rightarrow 0,
\end{equation}
induces the exact sequence
$$0 \rightarrow \Hom_{R}(L_{j}/L_{j-1},M) \rightarrow \Hom_{R}(N^{r_{j}},M).$$
It follows that $\Hom_{R}(L_{j}/L_{j-1},M) \in \mathcal{S}_{\mathcal{P}}(R)$.
Now the short exact sequence
\begin{equation} \label{eq:2.2.12.4}
0 \rightarrow L_{j-1} \rightarrow L_{j} \rightarrow L_{j}/L_{j-1} \rightarrow 0,
\end{equation}
yields the exact sequence
$$\Hom_{R}(L_{j}/L_{j-1},M) \rightarrow \Hom_{R}(L_{j},M) \rightarrow \Hom_{R}(L_{j-1},M).$$
A successive use of the above exact sequence, letting $j=1,\ldots,t$, implies that $\Hom_{R}(L,M) \in \mathcal{S}_{\mathcal{P}}(R)$.

Now let $s \geq 1$, and suppose that the results holds true for $s-1$. The induction hypothesis implies that $\Ext^{i}_{R}(L,M) \in \mathcal{S}_{\mathcal{P}}(R)$ for every $0 \leq i \leq s-1$. Hence, it suffices to show that $\Ext^{s}_{R}(L,M) \in \mathcal{S}_{\mathcal{P}}(R)$. Given any $1 \leq j \leq t$, the short exact sequence \eqref{eq:2.2.12.3} induces the exact sequence
$$\Ext^{s-1}_{R}(K_{j},M) \rightarrow \Ext^{s}_{R}\left(L_{j}/L_{j-1},M\right) \rightarrow \Ext^{s}_{R}(N^{r_{j}},M).$$
The induction hypothesis shows that $\Ext^{s-1}_{R}(K_{j},M) \in \mathcal{S}_{\mathcal{P}}(R)$, so from the above exact sequence we get that $\Ext^{s}_{R}\left(L_{j}/L_{j-1},M\right) \in \mathcal{S}_{\mathcal{P}}(R)$. Now, the short exact sequence \eqref{eq:2.2.12.4} yields the exact sequence
$$\Ext^{s}_{R}(L_{j}/L_{j-1},M) \rightarrow \Ext^{s}_{R}(L_{j},M) \rightarrow \Ext^{s}_{R}\left(L_{j-1},M\right).$$
A successive use of the above exact sequence, letting $j=1,\ldots,t$, implies that $\Ext^{s}_{R}(L,M) \in \mathcal{S}_{\mathcal{P}}(R)$.

(ii) $\Rightarrow$ (i): The hypothesis implies that $\Ext^{i}_{R}(R/\mathfrak{a},M) \in \mathcal{S}_{\mathcal{P}}(R)$ for every $0 \leq i \leq s$. The $R$-complex $K^{R}(\underline{a})\otimes_{R}M$ is isomorphic to an $R$-complex of the form
$$0 \rightarrow M \xrightarrow{\partial_{n}} M^{t_{n-1}} \rightarrow \cdots \rightarrow M^{t_{2}} \xrightarrow{\partial_{2}} M^{t_{1}} \xrightarrow{\partial_{1}} M \rightarrow 0.$$
We have
\begin{equation} \label{eq:2.2.12.5}
H_{n}(\underline{a};M)\cong \ker \partial_{n} \cong \Ext^{0}_{R}(R/\mathfrak{a},M) \in \mathcal{S}_{\mathcal{P}}(R),
\end{equation}
by the assumption. If $s=0$, then we are done. Suppose that $s \geq 1$. The short exact sequence
$$0 \rightarrow \ker \partial_{n} \rightarrow M \rightarrow \im \partial_{n} \rightarrow 0,$$
induces the exact sequence
$$\Ext^{i}_{R}(R/\mathfrak{a},M) \rightarrow \Ext^{i}_{R}(R/\mathfrak{a},\im \partial_{n}) \rightarrow \Ext^{i+1}_{R}(R/\mathfrak{a},\ker \partial_{n}).$$
The assumption along with the display \eqref{eq:2.2.12.5} imply that the two lateral terms of the above exact sequence belong to $\mathcal{S}_{\mathcal{P}}(R)$, so $\Ext^{i}_{R}(R/\mathfrak{a},\im \partial_{n}) \in \mathcal{S}_{\mathcal{P}}(R)$ for every $0 \leq i \leq s$. As $s \geq 1$, we have $\Ext^{1}_{R}(R/\mathfrak{a},\im \partial_{n}) \in \mathcal{S}_{\mathcal{P}}(R)$. Moreover, since $\ker \partial_{n-1} \subseteq M^{t_{n-1}}$, we get the exact sequence
$$0 \rightarrow \Hom_{R}(R/\mathfrak{a},\ker \partial_{n-1}) \rightarrow \Hom_{R}(R/\mathfrak{a},M^{t_{n-1}}).$$
But $\Hom_{R}(R/\mathfrak{a},M^{t_{n-1}}) \in \mathcal{S}_{\mathcal{P}}(R)$, so $\Hom_{R}(R/\mathfrak{a},\ker \partial_{n-1}) \in \mathcal{S}_{\mathcal{P}}(R)$.
The short exact sequence
\begin{equation} \label{eq:2.2.12.6}
0 \rightarrow \im \partial_{n} \rightarrow \ker \partial_{n-1} \rightarrow H_{n-1}(\underline{a};M) \rightarrow 0,
\end{equation}
yields the exact sequence
$$\Hom_{R}(R/\mathfrak{a},\ker \partial_{n-1}) \rightarrow \Hom_{R}\left(R/\mathfrak{a},H_{n-1}(\underline{a};M)\right) \rightarrow \Ext^{1}_{R}(R/\mathfrak{a},\im \partial_{n}).$$
Therefore, $\Hom_{R}\left(R/\mathfrak{a},H_{n-1}(\underline{a};M)\right) \in \mathcal{S}_{\mathcal{P}}(R)$.
But $\mathfrak{a} H_{n-1}(\underline{a};M)=0$, so
$$H_{n-1}(\underline{a};M) \cong \Hom_{R}\left(R/\mathfrak{a},H_{n-1}(\underline{a};M)\right) \in \mathcal{S}_{\mathcal{P}}(R).$$
If $s=1$, then we are done.
Suppose that $s \geq 2$. The short exact sequence \eqref{eq:2.2.12.6} induces the exact sequence
$$\Ext^{i}_{R}\left(R/\mathfrak{a},\im \partial_{n}\right) \rightarrow \Ext^{i}_{R}(R/\mathfrak{a},\ker \partial_{n-1}) \rightarrow \Ext^{i}_{R}\left(R/\mathfrak{a},H_{n-1}(\underline{a};M)\right).$$
It follows that $\Ext^{i}_{R}(R/\mathfrak{a},\ker \partial_{n-1}) \in \mathcal{S}_{\mathcal{P}}(R)$ for every $0 \leq i \leq s$.
The short exact sequence
$$0 \rightarrow \ker \partial_{n-1} \rightarrow M^{t_{n-1}} \rightarrow \im \partial_{n-1} \rightarrow 0,$$
yields the exact sequence
$$\Ext^{i}_{R}\left(R/\mathfrak{a},M^{t_{n-1}}\right) \rightarrow \Ext^{i}_{R}\left(R/\mathfrak{a},\im \partial_{n-1}\right) \rightarrow \Ext^{i+1}_{R}\left(R/\mathfrak{a},\ker \partial_{n-1}\right).$$
Thus $\Ext^{i}_{R}\left(R/\mathfrak{a},\im \partial_{n-1}\right) \in \mathcal{S}_{\mathcal{P}}(R)$ for every $0 \leq i \leq s-1$.
As $s \geq 2$, we see that $\Ext^{1}_{R}\left(R/\mathfrak{a},\im \partial_{n-1}\right) \in \mathcal{S}_{\mathcal{P}}(R)$. Moreover, since $\ker \partial_{n-2} \subseteq M^{t_{n-2}}$, we get the exact sequence
$$0 \rightarrow \Hom_{R}(R/\mathfrak{a},\ker \partial_{n-2}) \rightarrow \Hom_{R}(R/\mathfrak{a},M^{t_{n-2}}).$$
But $\Hom_{R}(R/\mathfrak{a},M^{t_{n-2}}) \in \mathcal{S}_{\mathcal{P}}(R)$, so $\Hom_{R}(R/\mathfrak{a},\ker \partial_{n-2}) \in \mathcal{S}_{\mathcal{P}}(R)$. The short exact sequence
$$0 \rightarrow \im \partial_{n-1} \rightarrow \ker \partial_{n-2} \rightarrow H_{n-2}(\underline{a};M) \rightarrow 0,$$
implies the exact sequence
$$\Hom_{R}(R/\mathfrak{a},\ker \partial_{n-2}) \rightarrow \Hom_{R}\left(R/\mathfrak{a},H_{n-2}(\underline{a};M)\right) \rightarrow \Ext_{R}^{1}(R/\mathfrak{a},\im \partial_{n-1}).$$
Therefore, $\Hom_{R}\left(R/\mathfrak{a},H_{n-2}(\underline{a};M)\right) \in \mathcal{S}_{\mathcal{P}}(R)$.
As before, we conclude that
$$H_{n-2}(\underline{a};M) \cong \Hom_{R}\left(R/\mathfrak{a},H_{n-2}(\underline{a};M)\right) \in \mathcal{S}_{\mathcal{P}}(R).$$
If $s=2$, then we are done. Proceeding in this manner, we infer that $H_{n-i}(\underline{a};M) \in \mathcal{S}_{\mathcal{P}}(R)$ for every $0 \leq i \leq s$.

Now, assume that the Serre property $\mathcal{P}$ satisfies the condition $\mathfrak{C}_{\mathfrak{a}}$.

(iv) $\Rightarrow$ (iii): We let $N=R/\mathfrak{a}$, and show that $\Ext_{R}^{i}(R/\mathfrak{a},M)\in \mathcal{S}_{\mathcal{P}}(R)$ for every $0 \leq i \leq s$. By Lemma \ref{2.2.9}, there is a third quadrant spectral sequence
\begin{equation} \label{eq:2.2.12.7}
E^{2}_{p,q}= \Ext^{-p}_{R}\left(R/\mathfrak{a},H^{-q}_{\mathfrak{a}}(M)\right) \underset {p} \Rightarrow \Ext^{-p-q}_{R}(R/\mathfrak{a},M).
\end{equation}
The hypothesis implies that $E^{2}_{p,q}\in \mathcal{S}_{\mathcal{P}}(R)$ for every $p \in \mathbb{Z}$ and $-s \leq q \leq 0$. Let $-s \leq i \leq 0$. There is a finite filtration
$$0=U^{i-1} \subseteq U^{i} \subseteq \cdots \subseteq U^{0} = \Ext^{-i}_{R}(R/\mathfrak{a},M),$$
such that $U^{p}/U^{p-1}\cong E^{\infty}_{p,i-p}$ for every $i \leq p \leq 0$. Since $E^{\infty}_{p,i-p}$ is a subquotient of $E^{2}_{p,i-p}$ and $-s \leq i-p \leq 0$, we infer that
$$U^{p}/U^{p-1}\cong E^{\infty}_{p,i-p}\in \mathcal{S}_{\mathcal{P}}(R)$$
for every $i \leq p \leq 0$.
A successive use of the short exact sequence
$$0 \rightarrow U^{p-1} \rightarrow U^{p} \rightarrow U^{p}/U^{p-1} \rightarrow 0,$$
by letting $p=i,\ldots,0$, implies that $\Ext^{-i}_{R}(R/\mathfrak{a},M)\in \mathcal{S}_{\mathcal{P}}(R)$.

(ii) $\Rightarrow$ (iv): It follows from the hypothesis that $\Ext^{i}_{R}(R/\mathfrak{a},M)\in \mathcal{S}_{\mathcal{P}}(R)$ for every $0 \leq i \leq s$. We argue by induction on $s$. Let $s=0$. Then
$$(0:_{M}\mathfrak{a}) \cong \Hom_{R}(R/\mathfrak{a},M) \in \mathcal{S}_{\mathcal{P}}(R),$$
whence
$$H^{0}_{s}(M) \cong \Gamma_{\mathfrak{a}}(M) \in \mathcal{S}_{\mathcal{P}}(R)$$
by the condition $\mathfrak{C}_{\mathfrak{a}}$.

Now let $s \geq 1$, and suppose that the result holds true for $s-1$. The induction hypothesis implies that $H^{i}_{\mathfrak{a}}(M)\in \mathcal{S}_{\mathcal{P}}(R)$ for every $0 \leq i \leq s-1$. Hence, it suffices to show that $H^{s}_{\mathfrak{a}}(M)\in \mathcal{S}_{\mathcal{P}}(R)$.
Consider the spectral sequence \eqref{eq:2.2.12.7}. By the hypothesis, $E^{2}_{p,q}\in \mathcal{S}_{\mathcal{P}}(R)$ for every $p \in \mathbb{Z}$ and $1-s \leq q \leq 0$.
There is a finite filtration
$$0=U^{-s-1} \subseteq U^{-s} \subseteq \cdots \subseteq U^{0} = \Ext^{s}_{R}(R/\mathfrak{a},M),$$
such that $U^{p}/U^{p-1}\cong E^{\infty}_{p,-s-p}$ for every $-s \leq p \leq 0$. As $\Ext^{s}_{R}(R/\mathfrak{a},M)\in \mathcal{S}_{\mathcal{P}}(R)$,
we conclude that
$$E^{\infty}_{0,-s}\cong U^{0}/U^{-1} \in \mathcal{S}_{\mathcal{P}}(R).$$
Let $r \geq 2$, and consider the differentials
$$0=E^{r}_{r,-s-r+1}\xrightarrow {d^{r}_{r,-s-r+1}} E^{r}_{0,-s}\xrightarrow {d^{r}_{0,-s}} E^{r}_{-r,-s+r-1}.$$
Since $1-s \leq -s+r-1$ and $E^{r}_{-r,-s+r-1}$ is a subquotient of $E^{2}_{-r,-s+r-1}$, the hypothesis implies that $E^{r}_{-r,-s+r-1}\in \mathcal{S}_{\mathcal{P}}(R)$, and consequently $\im d^{r}_{0,-s}\in \mathcal{S}_{\mathcal{P}}(R)$ for every $r\geq 2$.
We thus obtain
$$E^{r+1}_{0,-s} \cong \ker d^{r}_{0,-s}/ \im d^{r}_{r,-s-r+1} = \ker d^{r}_{0,-s},$$
and consequently a short exact sequence
\begin{equation} \label{eq:2.2.12.8}
0 \rightarrow E^{r+1}_{0,-s} \rightarrow E^{r}_{0,-s} \rightarrow \im d^{r}_{0,-s} \rightarrow 0.
\end{equation}
There is an integer $r_{0}\geq 2$, such that $E^{\infty}_{0,-s} \cong E^{r+1}_{0,-s}$ for every $r\geq r_{0}$. It follows that $E^{r_{0}+1}_{0,-s}\in \mathcal{S}_{\mathcal{P}}(R)$.
Now, the short exact sequence \eqref{eq:2.2.12.8} implies that $E^{r_{0}}_{0,-s}\in \mathcal{S}_{\mathcal{P}}(R)$. Using the short exact sequence
\eqref{eq:2.2.12.8} inductively, we conclude that
$$(0:_{H^{s}_{\mathfrak{a}}(M)}\mathfrak{a}) \cong \Hom_{R}\left(R/\mathfrak{a},H^{s}_{\mathfrak{a}}(M)\right)\cong E^{2}_{0,-s}\in \mathcal{S}_{\mathcal{P}}(R).$$
Therefore, by the condition $\mathfrak{C}_{\mathfrak{a}}$, we get
$$H^{s}_{\mathfrak{a}}(M) \cong \Gamma_{\mathfrak{a}}\left(H^{s}_{\mathfrak{a}}(M)\right)\in \mathcal{S}_{\mathcal{P}}(R).$$
\end{prf}

If we let the integer $s$ exhaust the whole nonzero range of Koszul homology, i.e. $s=n$, then we can effectively combine Theorems \ref{2.2.10} and \ref{2.2.12} to obtain the following result which in turn generalizes \cite[Corollary 3.1]{BA}.

\begin{corollary} \label{2.2.13}
Let $\mathfrak{a}=(a_{1},\ldots,a_{n})$ be an ideal of $R$, $\underline{a}=a_{1},\ldots,a_{n}$, and $M$ an $R$-module. Let $\mathcal{P}$ be a Serre property. Then the following conditions are equivalent:
\begin{enumerate}
\item[(i)] $H_{i}(\underline{a};M)\in \mathcal{S}_{\mathcal{P}}(R)$ for every $0 \leq i \leq n$.
\item[(ii)] $\Tor^{R}_{i}(N,M) \in \mathcal{S}_{\mathcal{P}}(R)$ for every finitely generated $R$-module $N$ with $\Supp_{R}(N)\subseteq V(\mathfrak{a})$, and for every $i\geq 0$.
\item[(iii)] $\Tor^{R}_{i}(N,M) \in \mathcal{S}_{\mathcal{P}}(R)$ for some finitely generated $R$-module $N$ with $\Supp_{R}(N)= V(\mathfrak{a})$, and for every $0 \leq i \leq n$.
\item[(iv)] $\Ext^{i}_{R}(N,M) \in \mathcal{S}_{\mathcal{P}}(R)$ for every finitely generated $R$-module $N$ with $\Supp_{R}(N)\subseteq V(\mathfrak{a})$, and for every $i\geq 0$.
\item[(v)] $\Ext^{i}_{R}(N,M) \in \mathcal{S}_{\mathcal{P}}(R)$ for some finitely generated $R$-module $N$ with $\Supp_{R}(N)= V(\mathfrak{a})$, and for every $0 \leq i \leq n$.
\end{enumerate}
\end{corollary}

\begin{prf}
(i) $\Leftrightarrow$ (ii) and (i) $\Leftrightarrow$ (iv): Since $H_{i}(\underline{a};M)=0$ for every $i > n$, these equivalences follow from Theorems \ref{2.2.10} and \ref{2.2.12}, respectively.

(i) $\Leftrightarrow$ (iii): Follows from Theorem \ref{2.2.10}.

(i) $\Leftrightarrow$ (v): Follows from Theorem \ref{2.2.12}.
\end{prf}

The following corollaries describe the numerical invariants $\mathcal{P}$-depth and $\mathcal{P}$-width in terms of Koszul homology, local homology, and local cohomology.

\begin{corollary} \label{2.2.14}
Let $\mathfrak{a}=(a_{1},\ldots,a_{n})$ be an ideal of $R$, $\underline{a}=a_{1},\ldots,a_{n}$, and $M$ an $R$-module. Let $\mathcal{P}$ be a Serre property. Then the following assertions hold:
\begin{enumerate}
\item[(i)] $\mathcal{P}$-$\depth_{R}(\mathfrak{a},M) = \inf \left\{i \geq 0 \suchthat H_{n-i}(\underline{a};M) \notin \mathcal{S}_{\mathcal{P}}(R) \right\}.$
\item[(ii)] $\mathcal{P}$-$\width_{R}(\mathfrak{a},M) = \inf \left\{i \geq 0 \suchthat H_{i}(\underline{a};M) \notin \mathcal{S}_{\mathcal{P}}(R) \right\}.$
\item[(iii)] We have $\mathcal{P}$-$\depth_{R}(\mathfrak{a},M) < \infty$ if and only if $\mathcal{P}$-$\width_{R}(\mathfrak{a},M) < \infty$. Moreover in this case, we have
\begin{center}
$\sup\left\{i \geq 0 \suchthat H_{i}(\underline{a};M) \notin \mathcal{S}_{\mathcal{P}}(R) \right\} +$ $\mathcal{P}$-$\depth_{R}(\mathfrak{a},M)=n.$
\end{center}
\end{enumerate}
\end{corollary}

\begin{prf}
(i) and (ii): Follows from Theorems \ref{2.2.10} and \ref{2.2.12}.

(iii): The first assertion follows from Corollary \ref{2.2.11}. For the second assertion, note that
$$\inf \left\{i \geq 0 \suchthat H_{n-i}(\underline{a};M) \notin \mathcal{S}_{\mathcal{P}}(R) \right\} = n- \sup \left\{i \geq 0 \suchthat H_{i}(\underline{a};M) \notin \mathcal{S}_{\mathcal{P}}(R) \right\}.$$
\end{prf}

\begin{corollary} \label{2.2.15}
Let $\mathfrak{a}$ be an ideal of $R$, and $M$ an $R$-module. Let $\mathcal{P}$ be a Serre property. Then the following assertions hold:
\begin{enumerate}
\item[(i)] If $\mathcal{P}$ satisfies the condition $\mathfrak{C}_{\mathfrak{a}}$, then
\begin{center}
$\mathcal{P}$-$\depth_{R}(\mathfrak{a},M) = \inf \left\{i \geq 0 \suchthat H^{i}_{\mathfrak{a}}(M) \notin \mathcal{S}_{\mathcal{P}}(R) \right\}.$
\end{center}
\item[(ii)] If $\mathcal{P}$ satisfies the condition $\mathfrak{D}_{\mathfrak{a}}$, then
\begin{center}
$\mathcal{P}$-$\width_{R}(\mathfrak{a},M) = \inf \left\{i \geq 0 \suchthat H^{\mathfrak{a}}_{i}(M) \notin \mathcal{S}_{\mathcal{P}}\left(\widehat{R}^{\mathfrak{a}}\right) \right\}.$
\end{center}
\item[(iii)] If $\mathcal{P}$ satisfies both conditions $\mathfrak{C}_{\mathfrak{a}}$ and $\mathfrak{D}_{\mathfrak{a}}$, and $\mathcal{P}$-$\depth_{R}(\mathfrak{a},M) < \infty$, then
\begin{center}
$\mathcal{P}$-$\depth_{R}(\mathfrak{a},M)$ $+$ $\mathcal{P}$-$\width_{R}(\mathfrak{a},M) \leq \ara(\mathfrak{a}).$
\end{center}
\end{enumerate}
\end{corollary}

\begin{prf}
(i) and (ii): Follows from Theorems \ref{2.2.10} and \ref{2.2.12}.

(iii): Clear by (i), (ii), and Corollary \ref{2.2.14} (iii).
\end{prf}

Now we specialize the results of this section to obtain some characterizations of noetherian local homology modules and artinian local cohomology modules. The following result generalizes \cite[Propositions 7.1, 7.2 and 7.4]{WW} when applied to modules.

\begin{proposition} \label{2.2.16}
Let $\mathfrak{a}=(a_{1},\ldots,a_{n})$ be an ideal of $R$, $\underline{a}=a_{1},\ldots,a_{n}$, and $M$ an $R$-module. Then the following assertions are equivalent for any given $s \geq 0$:
\begin{enumerate}
\item[(i)] $H_{i}(\underline{a};M)$ is a finitely generated $R$-module for every $0 \leq i \leq s$.
\item[(ii)] $\Tor^{R}_{i}(N,M)$ is a finitely generated $R$-module for every finitely generated $R$-module $N$ with $\Supp_{R}(N)\subseteq V(\mathfrak{a})$, and for every $0 \leq i \leq s$.
\item[(iii)] $\Tor^{R}_{i}(N,M)$ is a finitely generated $R$-module for some finitely generated $R$-module $N$ with $\Supp_{R}(N)= V(\mathfrak{a})$, and for every $0 \leq i \leq s$.
\item[(iv)] $H^{\mathfrak{a}}_{i}(M)$ is a finitely generated $\widehat{R}^{\mathfrak{a}}$-module for every $0 \leq i \leq s$.
\end{enumerate}
\end{proposition}

\begin{prf}
Obvious in view of Example \ref{2.2.6} (ii) and Theorem \ref{2.2.10}.
\end{prf}

The following corollary provides a characterization of noetherian local homology modules in its full generality.

\begin{corollary} \label{2.2.17}
Let $\mathfrak{a}=(a_{1},\ldots,a_{n})$ be an ideal of $R$, $\underline{a}=a_{1},\ldots,a_{n}$, and $M$ an $R$-module. Then the following assertions are equivalent:
\begin{enumerate}
\item[(i)] $H^{\mathfrak{a}}_{i}(M)$ is a finitely generated $\widehat{R}^{\mathfrak{a}}$-module for every $i \geq 0$.
\item[(ii)] $H_{i}(\underline{a};M)$ is a finitely generated $R$-module for every $0 \leq i \leq n$.
\item[(iii)] $H_{i}(\underline{a};M)$ is a finitely generated $R$-module for every $0 \leq i \leq \hd(\mathfrak{a},M)$.
\item[(iv)] $\Tor^{R}_{i}(N,M)$ is a finitely generated $R$-module for every finitely generated $R$-module $N$ with $\Supp_{R}(N)\subseteq V(\mathfrak{a})$, and for every $i \geq 0$.
\item[(v)] $\Tor^{R}_{i}(N,M)$ is a finitely generated $R$-module for some finitely generated $R$-module $N$ with $\Supp_{R}(N)= V(\mathfrak{a})$, and for every $0 \leq i \leq \hd(\mathfrak{a},M)$.
\item[(vi)] $\Ext^{i}_{R}(N,M)$ is a finitely generated $R$-module for every finitely generated $R$-module $N$ with $\Supp_{R}(N)\subseteq V(\mathfrak{a})$, and for every $i \geq 0$.
\item[(vii)] $\Ext^{i}_{R}(N,M)$ is a finitely generated $R$-module for some finitely generated $R$-module $N$ with $\Supp_{R}(N)= V(\mathfrak{a})$, and for every $0 \leq i \leq n$.
\end{enumerate}
\end{corollary}

\begin{prf}
(ii) $\Leftrightarrow$ (iv) $\Leftrightarrow$ (vi) $\Leftrightarrow$ (vii): Follows from Corollary \ref{2.2.13}.

(iii) $\Leftrightarrow$ (v): Follows from Proposition \ref{2.2.16} upon setting $s= \hd(\mathfrak{a},M)$.

(i) $\Leftrightarrow$ (iii): Since $H^{\mathfrak{a}}_{i}(M)=0$ for every $i > \hd(\mathfrak{a},M)$, the result follows from Proposition \ref{2.2.16}.

(i) $\Leftrightarrow$ (iv): Follows from Proposition \ref{2.2.16}.
\end{prf}

One should note that a slightly weaker version of the following result has been proved in \cite[Theorem 5.5]{Me2} by using a different method.

\begin{proposition} \label{2.2.18}
Let $\mathfrak{a}=(a_{1},\ldots,a_{n})$ be an ideal of $R$, $\underline{a}=a_{1},\ldots,a_{n}$, and $M$ an $R$-module. Then the following assertions are equivalent for any given $s \geq 0$:
\begin{enumerate}
\item[(i)] $H_{n-i}(\underline{a};M)$ is an artinian $R$-module for every $0 \leq i \leq s$.
\item[(ii)] $\Ext^{i}_{R}(N,M)$ is an artinian $R$-module for every finitely generated $R$-module $N$ with $\Supp_{R}(N)\subseteq V(\mathfrak{a})$, and for every $0 \leq i \leq s$.
\item[(iii)] $\Ext^{i}_{R}(N,M)$ is an artinian $R$-module for some finitely generated $R$-module $N$ with $\Supp_{R}(N)= V(\mathfrak{a})$, and for every $0 \leq i \leq s$.
\item[(iv)] $H^{i}_{\mathfrak{a}}(M)$ is an artinian $R$-module for every $0 \leq i \leq s$.
\end{enumerate}
\end{proposition}

\begin{prf}
Obvious in view of Example \ref{2.2.8} (ii) and Theorem \ref{2.2.12}.
\end{prf}

The following corollary provides a characterization of artinian local cohomology modules in its full generality.

\begin{corollary} \label{2.2.19}
Let $\mathfrak{a}=(a_{1},\ldots,a_{n})$ be an ideal of $R$, $\underline{a}=a_{1},\ldots,a_{n}$, and $M$ an $R$-module. Then the following assertions are equivalent:
\begin{enumerate}
\item[(i)] $H^{i}_{\mathfrak{a}}(M)$ is an artinian $R$-module for every $i \geq 0$.
\item[(ii)] $H_{i}(\underline{a};M)$ is an artinian $R$-module for every $0 \leq i \leq n$.
\item[(iii)] $H_{n-i}(\underline{a};M)$ is an artinian $R$-module for every $0 \leq i \leq \cd(\mathfrak{a},M)$.
\item[(iv)] $\Ext^{i}_{R}(N,M)$ is an artinian $R$-module for every finitely generated $R$-module $N$ with $\Supp_{R}(N)\subseteq V(\mathfrak{a})$, and for every $i \geq 0$.
\item[(v)] $\Ext^{i}_{R}(N,M)$ is an artinian $R$-module for some finitely generated $R$-module $N$ with $\Supp_{R}(N)= V(\mathfrak{a})$, and for every $0 \leq i \leq \cd(\mathfrak{a},M)$.
\item[(vi)] $\Tor^{R}_{i}(N,M)$ is an artinian $R$-module for every finitely generated $R$-module $N$ with $\Supp_{R}(N)\subseteq V(\mathfrak{a})$, and for every $i \geq 0$.
\item[(vii)] $\Tor^{R}_{i}(N,M)$ is an artinian $R$-module for some finitely generated $R$-module $N$ with $\Supp_{R}(N)= V(\mathfrak{a})$, and for every $0 \leq i \leq n$.
\end{enumerate}
\end{corollary}

\begin{prf}
(ii) $\Leftrightarrow$ (iv) $\Leftrightarrow$ (vi) $\Leftrightarrow$ (vii): Follows from Corollary \ref{2.2.13}.

(iii) $\Leftrightarrow$ (v): Follows from Proposition \ref{2.2.18} upon setting $s= \cd(\mathfrak{a},M)$.

(i) $\Leftrightarrow$ (iii): Since $H^{i}_{\mathfrak{a}}(M)=0$ for every $i > \cd(\mathfrak{a},M)$, the result follows from Proposition \ref{2.2.18}.

(i) $\Leftrightarrow$ (iv): Follows from Proposition \ref{2.2.18}.
\end{prf}

Since the property of being zero is a Serre property that satisfies the condition $\mathfrak{D}_{\mathfrak{a}}$, we specialize the results of this section to obtain the following result.

\begin{proposition} \label{2.2.20}
Let $\mathfrak{a}=(a_{1},\ldots,a_{n})$ be an ideal of $R$, $\underline{a}=a_{1},\ldots,a_{n}$, and $M$ an $R$-module. Then the following assertions are equivalent for any given $s \geq 0$:
\begin{enumerate}
\item[(i)] $H_{i}(\underline{a};M)=0$ for every $0 \leq i \leq s$.
\item[(ii)] $\Tor^{R}_{i}(N,M)=0$ for every $R$-module $N$ with $\Supp_{R}(N)\subseteq V(\mathfrak{a})$, and for every $0 \leq i \leq s$.
\item[(iii)] $\Tor^{R}_{i}(N,M)=0$ for some finitely generated $R$-module $N$ with $\Supp_{R}(N)= V(\mathfrak{a})$, and for every $0 \leq i \leq s$.
\item[(iv)] $H^{\mathfrak{a}}_{i}(M)=0$ for every $0 \leq i \leq s$.
\end{enumerate}
\end{proposition}

\begin{prf}
Immediate from Theorem \ref{2.2.10}. For part (ii), note that every module is a direct limit of its finitely generated submodules and Tor functor commutes with direct limits.
\end{prf}

The following result is proved in \cite[Theorem 4.4]{St} using a different method, but here it is an immediate by-product of Proposition \ref{2.2.20}.

\begin{corollary} \label{2.2.21}
Let $M$ be an $R$-module, and $N$ a finitely generated $R$-module. Then the following conditions are equivalent for any given $s \geq 0$:
\begin{enumerate}
\item[(i)] $\Tor^{R}_{i}(N,M)=0$ for every $0 \leq i \leq s$.
\item[(ii)] $\Tor^{R}_{i}\left(R/ \ann_{R}(N),M\right)=0$ for every $0 \leq i \leq s$.
\end{enumerate}
\end{corollary}

\begin{prf}
Immediate from Proposition \ref{2.2.20}.
\end{prf}

We observe that Corollary \ref{2.2.22} below generalizes \cite[Corollary 4.3]{Va}, which states that if $N$ is a faithful finitely generated $R$-module, then $N\otimes_{R}M=0$ if and only if $M=0$.

\begin{corollary} \label{2.2.22}
Let $M$ be an $R$-module, and $N$ a finitely generated $R$-module. Then the following conditions are equivalent:
\begin{enumerate}
\item[(i)] $M\otimes_{R}N=0$.
\item[(ii)] $M= \ann_{R}(N)M$.
\end{enumerate}
In particular, we have
$$\Supp_{R}(M\otimes_{R}N)=\Supp_{R}\left(M/ \ann_{R}(N)M\right).$$
\end{corollary}

\begin{prf}
For the equivalence of (i) and (ii), let $s=0$ in Corollary \ref{2.2.21}. For the second part, note that given any $\mathfrak{p}\in \Spec(R)$,
we have $\ann_{R}(N)_{\frak p}=\ann_{R_{\frak p}}(N_{\frak p})$, and so $(M\otimes_{R}N)_{\mathfrak{p}}=0$ if and only if $\left(M/\ann_{R}(N)M
\right)_{\mathfrak{p}}=0.$
\end{prf}

The support formula in Corollary \ref{2.2.22} generalizes the well-known formula
$$\Supp_{R}(M\otimes_{R}N)=\Supp_{R}(M)\cap \Supp_{R}(N),$$
which holds whenever $M$ and $N$ are both assumed to be finitely generated.

Using Corollary \ref{2.2.15}, we intend to obtain two somewhat known descriptions of the numerical invariant $\width_{R}(\mathfrak{a},M)$ in terms of the Koszul homology and local homology. However, we need some generalizations of the five-term exact sequences. To the best of our knowledge, the only place where one may find such generalizations is \cite[Corollaries 10.32 and 10.34]{Ro}. But, the statements there are not correct and no proof is presented. Hence due to lack of a suitable reference, we deem it appropriate to include the correct statements with proofs for the convenience of the reader.

\begin{lemma} \label{2.2.23}
Let $E^{2}_{p,q} \underset {p} \Rightarrow H_{p+q}$ be a spectral sequence. Then the following assertions hold:
\begin{enumerate}
\item[(i)] If $E^{2}_{p,q} \underset {p} \Rightarrow H_{p+q}$ is first quadrant and there is an integer $n \geq 1$ such that $E^{2}_{p,q}=0$ for every $q \leq n-2$, then there is a five-term exact sequence
$$H_{n+1} \rightarrow E^{2}_{2,n-1} \rightarrow E^{2}_{0,n} \rightarrow H_{n} \rightarrow E^{2}_{1,n-1} \rightarrow 0.$$
\item[(ii)] If $E^{2}_{p,q} \underset {p} \Rightarrow H_{p+q}$ is third quadrant and there is an integer $n \geq 1$ such that $E^{2}_{p,q}=0$ for every $q \geq 2-n$, then there is a five-term exact sequence
$$0 \rightarrow E^{2}_{-1,1-n} \rightarrow H_{-n} \rightarrow E^{2}_{0,-n} \rightarrow E^{2}_{-2,1-n} \rightarrow H_{-n-1}.$$
\end{enumerate}
\end{lemma}

\begin{prf}
(i): Consider the following homomorphisms
$$0=E^{2}_{4,n-2} \xrightarrow{d^{2}_{4,n-2}} E^{2}_{2,n-1} \xrightarrow{d^{2}_{2,n-1}} E^{2}_{0,n} \xrightarrow{d^{2}_{0,n}} E^{2}_{-2,n+1}=0.$$
We thus have
$$E^{3}_{2,n-1} \cong \ker d^{2}_{2,n-1}/ \im d^{2}_{4,n-2} \cong \ker d^{2}_{2,n-1},$$
and
$$E^{3}_{0,n}\cong \ker d^{2}_{0,n}/ \im d^{2}_{2,n-1}= \coker d^{2}_{2,n-1}.$$
Let $r \geq 3$. Consider the following homomorphisms
$$0=E^{r}_{r+2,n-r} \xrightarrow{d^{r}_{r+2,n-r}} E^{r}_{2,n-1} \xrightarrow{d^{r}_{2,n-1}} E^{r}_{2-r,n+r-2}=0.$$
We thus have
$$E^{r+1}_{2,n-1} \cong \ker d^{r}_{2,n-1}/ \im d^{r}_{r+2,n-r} \cong E^{r}_{2,n-1}.$$
Therefore,
$$\ker d^{2}_{2,n-1}\cong E^{3}_{2,n-1} \cong E^{4}_{2,n-1} \cong \cdots \cong E^{\infty}_{2,n-1}.$$
Further, consider the following homomorphisms
$$0=E^{r}_{r,n-r+1} \xrightarrow{d^{r}_{r,n-r+1}} E^{r}_{0,n} \xrightarrow{d^{r}_{0,n}} E^{r}_{-r,n+r-1}=0.$$
We thus have
$$E^{r+1}_{0,n}\cong \ker d^{r}_{0,n}/ \im d^{r}_{r,n-r+1} \cong E^{r}_{0,n}.$$
Therefore,
$$\coker d^{2}_{2,n-1} \cong E^{3}_{0,n} \cong E^{4}_{0,n} \cong \cdots \cong E^{\infty}_{0,n}.$$
Hence, we get the following exact sequence
\begin{equation} \label{eq:2.2.23.1}
0 \rightarrow E^{\infty}_{2,n-1} \rightarrow E^{2}_{2,n-1} \xrightarrow{d^{2}_{2,n-1}} E^{2}_{0,n} \rightarrow E^{\infty}_{0,n} \rightarrow 0.
\end{equation}
There is a finite filtration
$$0=U^{-1} \subseteq U^{0} \subseteq \cdots \subseteq U^{n+1}=H_{n+1},$$
such that $E^{\infty}_{p,n+1-p}\cong U^{p}/U^{p-1}$ for every $0 \leq p \leq n+1$. If $3 \leq p \leq n+1$, then $n+1-p \leq n-2$, so $E^{\infty}_{p,n+1-p}=0$. It follows that
$$U^{2}=U^{3}=\cdots =U^{n+1}.$$
Now since $E^{\infty}_{2,n-1}\cong U^{2}/U^{1}$, we get an exact sequence
$$H_{n+1}\rightarrow E^{\infty}_{2,n-1}\rightarrow 0.$$
Splicing this exact sequence to the exact sequence \eqref{eq:2.2.23.1}, we get the following exact sequence
\begin{equation} \label{eq:2.2.23.2}
H_{n+1} \rightarrow E^{2}_{2,n-1} \xrightarrow{d^{2}_{2,n-1}} E^{2}_{0,n} \rightarrow E^{\infty}_{0,n} \rightarrow 0.
\end{equation}
On the other hand, there is a finite filtration
$$0=V^{-1} \subseteq V^{0} \subseteq \cdots \subseteq V^{n}=H_{n},$$
such that $E^{\infty}_{p,n-p}\cong V^{p}/V^{p-1}$ for every $0 \leq p \leq n$. If $2 \leq p \leq n$, then $n-p \leq n-2$, so $E^{\infty}_{p,n-p}=0$. It follows that
$$V^{1}=V^{2}=\cdots =V^{n}.$$
As $E^{\infty}_{0,n}\cong V^{0}/V^{-1} = V^{0}$ and $E^{\infty}_{1,n-1} \cong V^{1}/V^{0}$, we get the short exact sequence
$$0 \rightarrow E^{\infty}_{0,n} \rightarrow H_{n} \rightarrow E^{\infty}_{1,n-1} \rightarrow 0.$$
Splicing this short exact sequence to the exact sequence \eqref{eq:2.2.23.2}, yields the exact sequence
\begin{equation} \label{eq:2.2.23.3}
H_{n+1} \rightarrow E^{2}_{2,n-1} \xrightarrow{d^{2}_{2,n-1}} E^{2}_{0,n} \rightarrow H_{n} \rightarrow E^{\infty}_{1,n-1} \rightarrow 0.
\end{equation}
Let $r \geq 2$, and consider the following homomorphisms
$$0=E^{r}_{r+1,n-r} \xrightarrow{d^{r}_{r+1,n-r}} E^{r}_{1,n-1} \xrightarrow{d^{r}_{1,n-1}} E^{r}_{1-r,n+r-2}=0.$$
We thus have
$$E^{r+1}_{1,n-1}\cong \ker d^{r}_{1,n-1}/ \im d^{r}_{r+1,n-r} \cong E^{r}_{1,n-1}.$$
Therefore,
$$E^{2}_{1,n-1}\cong E^{3}_{1,n-1}\cong \cdots \cong E^{\infty}_{1,n-1}.$$
Thus from the exact sequence \eqref{eq:2.2.23.3}, we get the desired exact sequence
$$H_{n+1} \rightarrow E^{2}_{2,n-1} \xrightarrow{d^{2}_{2,n-1}} E^{2}_{0,n} \rightarrow H_{n} \rightarrow E^{2}_{1,n-1} \rightarrow 0.$$

(ii): Similar to (i).
\end{prf}

For the next result, we need to recall the notion of coassociated prime ideals. Given an $R$-module $M$, a prime ideal $\mathfrak{p}\in \Spec (R)$ is said to be a coassociated prime ideal of $M$ if $\mathfrak{p}= \ann_{R}(M/N)$ for some submodule $N$ of $M$ such that $M/N$ is an artinian $R$-module. The set of coassociated prime ideals of $M$ is denoted by $\Coass_{R}(M)$.

\begin{corollary} \label{2.2.24}
Let $\mathfrak{a}=(a_{1},\ldots,a_{n})$ be an ideal of $R$, $\underline{a}=a_{1},\ldots,a_{n}$, and $M$ an $R$-module. Then the following assertions hold:
\begin{enumerate}
\item[(i)] $\width_{R}(\mathfrak{a},M) = \inf \left\{i \geq 0 \suchthat H_{i}(\underline{a};M) \neq 0 \right\} = \inf \left\{i \geq 0 \suchthat H_{i}^{\mathfrak{a}}(M) \neq 0 \right\}$.
\item[(ii)] $\Tor^{R}_{\width_{R}(\mathfrak{a},M)}(R/\mathfrak{a},M) \cong (R/\mathfrak{a})\otimes_{R}H_{\width_{R}(\mathfrak{a},M)}^{\mathfrak{a}}(M)$.
\item[(iii)] $\Lambda^{\mathfrak{a}}\left(H_{\width_{R}(\mathfrak{a},M)}^{\mathfrak{a}}(M)\right) \cong \underset{n}\varprojlim \Tor^{R}_{\width_{R}(\mathfrak{a},M)}(R/\mathfrak{a}^{n},M)$.
\item[(iv)] $\Coass_{R}\left(\Tor^{R}_{\width_{R}(\mathfrak{a},M)}(R/\mathfrak{a},M)\right)=\Coass_{R}\left(H_{\width_{R}(\mathfrak{a},M)}^{\mathfrak{a}}(M)\right) \cap V(\mathfrak{a})$.
\end{enumerate}
\end{corollary}

\begin{prf}
(i): Follows immediately from Corollaries \ref{2.2.14} and \ref{2.2.15}.

(ii): By Lemma \ref{2.2.9}, there is a first quadrant spectral sequence
$$E^{2}_{p,q}= \Tor^{R}_{p}\left(R/\mathfrak{a},H^{\mathfrak{a}}_{q}(M)\right) \underset {p} \Rightarrow \Tor^{R}_{p+q}(R/\mathfrak{a},M).$$
Let $n=\width_{R}(\mathfrak{a},M)$. Then by (i), $E^{2}_{p,q}=0$ for every $q\leq n-1$. Now, Lemma \ref{2.2.23} (i) gives the exact sequence
$$\Tor^{R}_{n+1}(R/\mathfrak{a},M) \rightarrow \Tor^{R}_{2}\left(R/\mathfrak{a},H^{\mathfrak{a}}_{n-1}(M)\right) \rightarrow (R/\mathfrak{a})
\otimes_{R} H^{\mathfrak{a}}_{n}(M)\rightarrow $$$$ \Tor^{R}_{n}(R/\mathfrak{a},M) \rightarrow \Tor_{1}^{R}\left(R/\mathfrak{a},
H^{\mathfrak{a}}_{n-1}(M)\right) \rightarrow 0.$$
But
$$\Tor^{R}_{2}\left(R/\mathfrak{a},H^{\mathfrak{a}}_{n-1}(M)\right)=0=\Tor_{1}^{R}\left(R/\mathfrak{a},H^{\mathfrak{a}}_{n-1}(M)\right),$$
so
$$(R/\mathfrak{a})\otimes_{R} H^{\mathfrak{a}}_{n}(M) \cong \Tor^{R}_{n}(R/\mathfrak{a},M).$$

(iii): Using (ii) and the facts that $H_{i}^{\mathfrak{a}^{n}}(M)\cong H_{i}^{\mathfrak{a}}(M)$ and $\width_{R}(\mathfrak{a}^{n},M)=\width_{R}(\mathfrak{a},M)$
for every $n,i\geq 0$, we get:
\begin{equation*}
\begin{split}
\Lambda^{\mathfrak{a}}\left(H_{\width_{R}(\mathfrak{a},M)}^{\mathfrak{a}}(M)\right) & = \underset{n}\varprojlim \left(H_{\width_{R}(\mathfrak{a},M)}^{\mathfrak{a}}(M)/\mathfrak{a}^{n}H_{\width_{R}(\mathfrak{a},M)}^{\mathfrak{a}}(M)\right) \\
 & \cong \underset{n}\varprojlim \Tor^{R}_{\width_{R}(\mathfrak{a},M)}(R/\mathfrak{a}^{n},M).
\end{split}
\end{equation*}

(iv): Using (ii) and \cite[Theorem 1.21]{Ya2}, we have
\begin{equation*}
\begin{split}
\Coass_{R}\left(\Tor^{R}_{\width_{R}(\mathfrak{a},M)}(R/\mathfrak{a},M)\right) & = \Coass_{R}\left((R/\mathfrak{a})\otimes_{R}H_{\width_{R}(\mathfrak{a},M)}^{\mathfrak{a}}(M)\right) \\
 & = \Supp_{R}(R/\mathfrak{a}) \cap \Coass_{R}\left(H_{\width_{R}(\mathfrak{a},M)}^{\mathfrak{a}}(M)\right) \\
 & = V(\mathfrak{a}) \cap \Coass_{R}\left(H_{\width_{R}(\mathfrak{a},M)}^{\mathfrak{a}}(M)\right).
\end{split}
\end{equation*}
\end{prf}

Note that part (iii) of Corollary \ref{2.2.24} is proved in \cite[Proposition 2.5]{Si2} by deploying a different method. On the other hand, in parallel with Corollary \ref{2.2.29} (iii) below, one may wonder if intersecting with $V(\mathfrak{a})$ in part (iv) of Corollary \ref{2.2.24} is redundant. In other words, $\Coass_{R}\left(H_{\width_{R}(\mathfrak{a},M)}^{\mathfrak{a}}(M)\right)$ may be contained in $V(\mathfrak{a})$. However, the following example shows that this is not the case in general.

\begin{example} \label{2.2.25}
Let $R:=\mathbb{Q}[X,Y]_{(X,Y)}$ and $\mathfrak{m}:=(X,Y)R$. Then $\width_{R}(\mathfrak{m},R)=0$ and
$$H^{\frak m}_0(R)\cong \widehat{R}^{\mathfrak{m}}\cong \mathbb{Q}[[X,Y]].$$
For each $n \in \mathbb{Z}$, let $\fp_n:=(X-nY)R$. Then it is easy to see that $R/\fp_n\cong \mathbb{Q}[Y]_{(Y)}$, and so it is not a complete local ring. By \cite [Beispiel 2.4]{Z},
\begin{equation*}
\begin{split}
\Coass_R\left(H^{\frak m}_0(R)\right) & = \Coass_R\left(\widehat{R}^{\mathfrak{m}}\right) \\
 & = \{\fm\} \cup \left\{\mathfrak{p}\in \Spec R \suchthat R/\mathfrak{p} \text{ is not a complete local ring} \right\}.
\end{split}
\end{equation*}
Hence $\Coass_{R}\left(H_0^{\frak m}(R)\right)$ is not a finite set, while
$$\Coass_{R}\left(H_0^{\frak m}(R)\right)\cap V(\mathfrak{m})= \{\mathfrak{m}\}.$$
In particular, $\Coass_{R}\left(H_0^{\frak m}(R)\right) \nsubseteq V(\mathfrak{m})$.
\end{example}

Since the property of being zero is a Serre property that satisfies the condition $\mathfrak{C}_{\mathfrak{a}}$, we obtain the following result.

\begin{proposition} \label{2.2.26}
Let $\mathfrak{a}=(a_{1},\ldots,a_{n})$ be an ideal of $R$, $\underline{a}=a_{1},\ldots,a_{n}$, and $M$ an $R$-module. Then the following assertions are equivalent for any given $s \geq 0$:
\begin{enumerate}
\item[(i)] $H_{n-i}(\underline{a};M)=0$ for every $0 \leq i \leq s$.
\item[(ii)] $\Ext^{i}_{R}(N,M)=0$ for every finitely generated $R$-module $N$ with $\Supp_{R}(N)\subseteq V(\mathfrak{a})$, and for every $0 \leq i \leq s$.
\item[(iii)] $\Ext^{i}_{R}(N,M)=0$ for some finitely generated $R$-module $N$ with $\Supp_{R}(N)= V(\mathfrak{a})$, and for every $0 \leq i \leq s$.
\item[(iv)] $H^{i}_{\mathfrak{a}}(M)=0$ for every $0 \leq i \leq s$.
\end{enumerate}
\end{proposition}

\begin{prf}
Immediate from Theorem \ref{2.2.12}.
\end{prf}

The following result is proved in \cite[Theorem 3.2]{St} using a different method, but here it is an immediate by-product of Proposition \ref{2.2.26}.

\begin{corollary} \label{2.2.27}
Let $M$ be an $R$-module, and $N$ a finitely generated $R$-module. Then the following conditions are equivalent for any given $s \geq 0$:
\begin{enumerate}
\item[(i)] $\Ext^{i}_{R}(N,M)=0$ for every $0 \leq i \leq s$.
\item[(ii)] $\Ext^{i}_{R}\left(R/ \ann_{R}(N),M\right)=0$ for every $0 \leq i \leq s$.
\end{enumerate}
\end{corollary}

\begin{prf}
Immediate from Proposition \ref{2.2.26}.
\end{prf}

The following special case may be of independent interest.

\begin{corollary} \label{2.2.28}
Let $M$ be an $R$-module, and $N$ a finitely generated $R$-module. Then the following conditions are equivalent:
\begin{enumerate}
\item[(i)] $\Hom_{R}(N,M)=0$.
\item[(ii)] $\left(0:_{M} \ann_{R}(N)\right)=0$.
\end{enumerate}
\end{corollary}

\begin{prf}
Let $s=0$ in Corollary \ref{2.2.27}.
\end{prf}

It is easy to deduce from Corollary \ref{2.2.28} that given a finitely generated $R$-module $N$, we have $\Hom_{R}(N,M)\neq 0$ if and only if there are elements $x\in N$ and $0\neq y \in M$ with $\ann_{R}(x) \subseteq \ann_{R}(y)$, which is known as the Hom Vanishing Lemma in \cite[Page 11]{C}.

We state the dual result to Corollary \ref{2.2.24} for the sake of integrity and completeness. Parts (ii) and (iii) of Corollary \ref{2.2.29} below are stated in \cite[Proposition 1.1]{Mar}, and it is only mentioned that part (ii) can be deduced from a spectral sequence. In addition, a proof of this result is offered in \cite[Corollary 2.3]{Me1} by using different techniques.

\begin{corollary} \label{2.2.29}
Let $\mathfrak{a}=(a_{1},\ldots,a_{n})$ be an ideal of $R$, $\underline{a}=a_{1},\ldots,a_{n}$, and $M$ an $R$-module. Then the following assertions hold:
\begin{enumerate}
\item[(i)] $\depth_{R}(\mathfrak{a},M) = \inf \left\{i \geq 0 \suchthat H_{n-i}(\underline{a};M) \neq 0 \right\} = \inf \left\{i \geq 0 \suchthat H_{\mathfrak{a}}^{i}(M) \neq 0 \right\}$.
\item[(ii)] $\Ext^{\depth_{R}(\mathfrak{a},M)}_{R}(R/\mathfrak{a},M) \cong \Hom_{R}\left(R/\mathfrak{a},H^{\depth_{R}(\mathfrak{a},M)}_{\mathfrak{a}}(M)\right)$.
\item[(iii)] $\Ass_{R}\left(\Ext^{\depth_{R}(\mathfrak{a},M)}_{R}(R/\mathfrak{a},M)\right)=\Ass_{R}\left(H^{\depth_{R}(\mathfrak{a},M)}_{\mathfrak{a}}(M)\right)$.
\end{enumerate}
\end{corollary}

\begin{prf}
(i): Follows immediately from Corollaries \ref{2.2.14} and \ref{2.2.15}.

(ii): By Lemma \ref{2.2.9}, there is a third quadrant spectral sequence
$$E^{2}_{p,q}= \Ext^{-p}_{R}\left(R/\mathfrak{a},H^{-q}_{\mathfrak{a}}(M)\right) \underset {p} \Rightarrow \Ext^{-p-q}_{R}(R/\mathfrak{a},M).$$
Now use Lemma \ref{2.2.23} (ii).

(iii): Follows from (ii).
\end{prf}

Finally, we present the following comprehensive vanishing result.

\begin{corollary} \label{2.2.30}
Let $\mathfrak{a}=(a_{1},\ldots,a_{n})$ be an ideal of $R$, $\underline{a}=a_{1},\ldots,a_{n}$, and $M$ an $R$-module. Then the following assertions are equivalent:
\begin{enumerate}
\item[(i)] $H_{i}(\underline{a};M)=0$ for every $0 \leq i \leq n$.
\item[(ii)] $H_{i}(\underline{a};M)=0$ for every $0 \leq i \leq \hd(\mathfrak{a},M)$.
\item[(iii)] $H_{i}(\underline{a};M)=0$ for every $n-\cd(\mathfrak{a},M) \leq i \leq n$.
\item[(iv)] $\Tor^{R}_{i}(N,M)=0$ for every $R$-module $N$ with $\Supp_{R}(N)\subseteq V(\mathfrak{a})$, and for every $i \geq 0$.
\item[(v)] $\Tor^{R}_{i}(N,M)=0$ for some finitely generated $R$-module $N$ with $\Supp_{R}(N)= V(\mathfrak{a})$, and for every $0 \leq i \leq \hd(\mathfrak{a},M)$.
\item[(vi)] $\Ext^{i}_{R}(N,M)=0$ for every finitely generated $R$-module $N$ with $\Supp_{R}(N)\subseteq V(\mathfrak{a})$, and for every $i \geq 0$.
\item[(vii)] $\Ext^{i}_{R}(N,M)=0$ for some finitely generated $R$-module $N$ with $\Supp_{R}(N)= V(\mathfrak{a})$, and for every $0 \leq i \leq \cd(\mathfrak{a},M)$.
\item[(viii)] $H^{\mathfrak{a}}_{i}(M)=0$ for every $i \geq 0$.
\item[(ix)] $H^{i}_{\mathfrak{a}}(M)=0$ for every $i \geq 0$.
\end{enumerate}
\end{corollary}

\begin{prf}
Follows from Propositions \ref{2.2.20} and \ref{2.2.26}.
\end{prf}

The following corollary is proved in \cite[Corollary 1.7]{Si2}. However, it is an immediate consequence of the results thus far obtained.

\begin{corollary} \label{2.2.31}
Let $\mathfrak{a}$ be an ideal of $R$, and $M$ an $R$-module. Then $\depth_{R}(\mathfrak{a},M) < \infty$ if and only if $\width_{R}(\mathfrak{a},M) < \infty$.
Moreover in this case, we have $$\depth_{R}(\mathfrak{a},M)+ \width_{R}(\mathfrak{a},M) \leq \ara(\mathfrak{a}).$$
\end{corollary}

\begin{prf}
Clear by Corollary \ref{2.2.15} (iii).
\end{prf}

\section{Cofiniteness of Modules}

In this section, we investigate the notion of $\mathfrak{a}$-cofiniteness for $R$-modules.

\begin{definition} \label{2.3.1}
Let $\mathfrak{a}$ be an ideal of $R$. An $R$-module $M$ is said to be $\mathfrak{a}$-\textit{cofinite} \index{cofinite module} if $\Supp_{R}(M)
\subseteq \V(\mathfrak{a})$, and $\Ext^{i}_{R}(R/\mathfrak{a},M)$ is a finitely generated $R$-module for every $i \geq 0$.
\end{definition}

The next theorem, whose techniques are mostly developed by Melkersson and Marley, embraces the scattered results in the literature concerning the cofiniteness of the local cohomology modules all at once.

First we need a lemma.

\begin{lemma} \label{2.3.2}
Let $\mathfrak{a}=(a_{1},...,a_{n})$ be an ideal of $R$, $\underline{a}=a_{1},...,a_{n}$, and $M$ an artinian $R$-module. Then $\ell_{R}(H_{n}\left(\underline{a};M)\right) < \infty$ if and only if $\ell_{R}(H_{i}\left(\underline{a};M)\right) < \infty$ for every $0 \leq i \leq n$.
\end{lemma}

\begin{prf}
Suppose that $\ell_{R}(H_{n}\left(\underline{a};M)\right) < \infty$. For any given $0 \leq i \leq n$, let $B_{i}:=B_{i}\left(K^{R}(\underline{a})\otimes_{R}M\right)$ and $Z_{i}:=Z_{i}\left(K^{R}(\underline{a})\otimes_{R}M\right)$. Since $M$ is artinian, $\left(K^{R}(\underline{a})\otimes_{R}M\right)_{i}\cong M^{t_{i}}$ is artinian, and thus $Z_{i}$ is artinian for every $0 \leq i \leq n$. By \cite[Proposition 3]{Ki}, there is an integer $t \geq 0$ such that
$$B_{i}+(0:_{Z_{i}}\mathfrak{a}^{m}) = \left(B_{i}+(0:_{Z_{i}}\mathfrak{a}^{t}):_{Z_{i}}\mathfrak{a}^{m-t}\right)$$
for every $m \geq t$. Let $m=t+1$, so that
$$B_{i}+(0:_{Z_{i}}\mathfrak{a}^{t+1}) = \left(B_{i}+(0:_{Z_{i}}\mathfrak{a}^{t}):_{Z_{i}}\mathfrak{a}\right) \supseteq (B_{i}:_{Z_{i}}\mathfrak{a}).$$
As $H_{i}(\underline{a};M)=Z_{i}/B_{i}$ and $\mathfrak{a}H_{i}(\underline{a};M)=0$, we have $(B_{i}:_{Z_{i}}\mathfrak{a})=Z_{i}$. Hence
$$Z_{i}=(B_{i}:_{Z_{i}}\mathfrak{a}) \subseteq B_{i}+(0:_{Z_{i}}\mathfrak{a}^{t+1}) \subseteq Z_{i},$$
so $Z_{i}= B_{i}+(0:_{Z_{i}}\mathfrak{a}^{t+1})$. Therefore,
\begin{equation*}
\begin{split}
H_{i}(\underline{a};M) & = Z_{i}/B_{i} \\
 & = \frac{B_{i}+(0:_{Z_{i}}\mathfrak{a}^{t+1})}{B_{i}} \\
 & \cong \frac{(0:_{Z_{i}}\mathfrak{a}^{t+1})}{(0:_{Z_{i}}\mathfrak{a}^{t+1}) \cap B_{i}}.
\end{split}
\end{equation*}
Thus $\Att_{R}\left(H_{i}(\underline{a};M)\right) \subseteq \Att_{R}\left((0:_{Z_{i}}\mathfrak{a}^{t+1})\right)$.

On the other hand,
$$\ell_{R}\left((0:_{M}\mathfrak{a})\right)=\ell_{R}\left(H_{n}(\underline{a};M)\right) < \infty$$
by the hypothesis. Hence it is easy to see that $\ell_{R}\left((0:_{M}\mathfrak{a}^{t+1})\right) < \infty$, and thus
$$\ell_{R}\left((0:_{M^{t_{i}}}\mathfrak{a}^{t+1})\right)=\ell_{R}\left((0:_{M}\mathfrak{a}^{t+1})^{t_{i}}\right) < \infty.$$
It follows that $\Att_{R}\left((0:_{M^{t_{i}}}\mathfrak{a}^{t+1})\right) \subseteq \Max(R)$.

Now if $\mathfrak{p} \in \Att_{R}\left(H_{i}(\underline{a};M)\right)$, then $\mathfrak{p} \in \Att_{R}\left((0:_{Z_{i}}\mathfrak{a}^{t+1})\right)$. Since $\left(K^{R}(\underline{a})\otimes_{R}M\right)_{i}\cong M^{t_{i}}$, we have an exact sequence
$$0 \rightarrow Z_{i} \rightarrow M^{t_{i}},$$
which induces the exact sequence
$$0 \rightarrow (0:_{Z_{i}}\mathfrak{a}^{t+1}) \rightarrow (0:_{M^{t_{i}}}\mathfrak{a}^{t+1}).$$
Since $M$ is artinian, $(0:_{M^{t_{i}}}\mathfrak{a}^{t+1})$ is artinian. We thus have
\begin{equation*}
\begin{split}
\bigcap_{\mathfrak{q}\in\Att_{R}\left((0:_{M^{t_{i}}}\mathfrak{a}^{t+1})\right)}\mathfrak{q} & = \sqrt{\ann_{R}\left((0:_{M^{t_{i}}}\mathfrak{a}^{t+1})\right)} \\
 & \subseteq \sqrt{\ann_{R}\left((0:_{Z_{i}}\mathfrak{a}^{t+1})\right)} \\
 & = \bigcap_{\mathfrak{p} \in \Att_{R}\left((0:_{Z_{i}}\mathfrak{a}^{t+1})\right)}\mathfrak{p}.
\end{split}
\end{equation*}
Thus $\mathfrak{p} \supseteq \bigcap_{\mathfrak{q} \in \Att_{R}\left((0:_{M^{t_{i}}}\mathfrak{a}^{t+1})\right)}\mathfrak{q}$. As $\Att_{R}\left((0:_{M^{t_{i}}}\mathfrak{a}^{t+1})\right)$ is finite, we infer that $\mathfrak{p} \supseteq \mathfrak{q}$ for some $\mathfrak{q} \in \Att_{R}\left((0:_{M^{t_{i}}}\mathfrak{a}^{t+1})\right)$. But $\mathfrak{q}$ is a maximal ideal, so $\mathfrak{p}=\mathfrak{q}$, and thus $\mathfrak{p} \in \Max(R)$. It follows that $\Att_{R}\left(H_{i}(\underline{a};M)\right) \subseteq \Max(R)$. As $M$ is artinian, $H_{i}(\underline{a};M)$ is artinian, so $\ell_{R}(H_{i}\left(\underline{a};M)\right) < \infty$.

The converse in trivial.
\end{prf}

\begin{theorem} \label{2.3.3}
Let $\mathfrak{a}$ be an ideal of $R$, and $M$ a finitely generated $R$-module. Suppose that either of the following conditions holds:
\begin{enumerate}
\item[(a)] $j=\dim_{R}(M)$.
\item[(b)] $\cd(\mathfrak{a},M) \leq 1$.
\item[(c)] $\dim_{R}(M) \leq 2$.
\item[(d)] $\dim(R/\mathfrak{a}) \leq 1$.
\end{enumerate}
Then the following assertions hold for any given $j \geq 0$:
\begin{enumerate}
\item[(i)] $H^{j}_{\mathfrak{a}}(M)$ is $\mathfrak{a}$-cofinite.
\item[(ii)] $\Ass_{R}\left(H^{j}_{\mathfrak{a}}(M)\right)$ is a finite set.
\item[(iii)] $\beta^{R}_{i}\left(\mathfrak{p},H^{j}_{\mathfrak{a}}(M)\right)$ is finite for every $\mathfrak{p} \in \Spec(R)$ and every $i \geq 0$.
\item[(iv)] $\mu^{i}_{R}\left(\mathfrak{p},H^{j}_{\mathfrak{a}}(M)\right)$ is finite for every $\mathfrak{p} \in \Spec(R)$ and every $i \geq 0$.
\end{enumerate}
\end{theorem}

\begin{prf}
We first show that (i) implies (ii), (iii), and (iv). Suppose that (i) holds. In particular, $\Hom_{R}\left(R/ \mathfrak{a},H_{\mathfrak{a}}^{j}(M)\right)$ is finitely generated, so $\Ass_{R} \left(\Hom_{R}\left(R/ \mathfrak{a},H_{\mathfrak{a}}^{j}(M)\right)\right)$ is finite.
But
$$\Ass_{R}\left(\Hom_{R}\left(R/\mathfrak{a},H_{\mathfrak{a}}^{j}(M)\right)\right) = V(\mathfrak{a}) \cap \Ass_{R}\left(H_{\mathfrak{a}}^{j}(M)\right) = \Ass_{R}\left(H_{\mathfrak{a}}^{j}(M)\right)$$
since
$$\Ass_{R}\left(H_{\mathfrak{a}}^{j}(M)\right) \subseteq \Supp_{R}\left(H_{\mathfrak{a}}^{j}(M)\right) \subseteq V(\mathfrak{a}).$$
Hence $\Ass_{R}\left(H_{\mathfrak{a}}^{j}(M)\right)$ is finite.
Moreover, we have by definition
$$\beta_{i}^{R}\left(\mathfrak{p},H_{\mathfrak{a}}^{j}(M)\right) :=  \rank_{R_{\mathfrak{p}}/ \mathfrak{p}R_{\mathfrak{p}}}\left(\Tor_{i}^{R_{\mathfrak{p}}}\left(R_{\mathfrak{p}}/ \mathfrak{p}R_\mathfrak{p},H_{\mathfrak{a}}^{j}(M)_{\mathfrak{p}}\right)\right)$$
and
$$\mu_{R}^{i}\left(\mathfrak{p},H_{\mathfrak{a}}^{j}(M)\right) := \rank_{R_{\mathfrak{p}}/ \mathfrak{p}R_{\mathfrak{p}}}\left(\Ext_{R_{\mathfrak{p}}}^i\left(R_{\mathfrak{p}}/ \mathfrak{p}R_{\mathfrak{p}},H_{\mathfrak{a}}^{j}(M)_{\mathfrak{p}}\right)\right).$$
If $\mathfrak{p} \notin V(\mathfrak{a})$, then $\mathfrak{p} \notin \Supp_{R}\left(H^{j}_{\mathfrak{a}}(M)\right)$, so $$\beta_{i}^{R}\left(\mathfrak{p},H_{\mathfrak{a}}^{j}(M)\right) = \mu_{R}^{i}\left(\mathfrak{p},H_{\mathfrak{a}}^{j}(M)\right)=0$$
for every $i \geq 0$.
On the other hand, if $\mathfrak{p} \in V(\mathfrak{a})$, then $\Supp_{R}(R/\mathfrak{p}) = V(\mathfrak{p}) \subseteq V(\mathfrak{a})$. Hence Corollary \ref{2.2.17} implies that
$$\Tor_{i}^{R_{\mathfrak{p}}}\left(R_{\mathfrak{p}}/ \mathfrak{p}R_\mathfrak{p},H_{\mathfrak{a}}^{j}(M)_{\mathfrak{p}}\right) \cong \Tor_{i}^{R}\left(R/ \mathfrak{p},H_{\mathfrak{a}}^{j}(M)\right)_{\mathfrak{p}}$$
and
$$\Ext_{R_{\mathfrak{p}}}^i\left(R_{\mathfrak{p}}/ \mathfrak{p}R_{\mathfrak{p}},H_{\mathfrak{a}}^{j}(M)_{\mathfrak{p}}\right) \cong \Ext_{R}^{i}\left(R/ \mathfrak{p},H_{\mathfrak{a}}^{j}(M)\right)_{\mathfrak{p}}$$
are finitely generated $R_{\mathfrak{p}}$-modules for every $i \geq 0$. It follows that $\beta_{i}^{R}\left(\mathfrak{p},H_{\mathfrak{a}}^{j}(M)\right)$ and $\mu_{R}^{i}\left(\mathfrak{p},H_{\mathfrak{a}}^{j}(M)\right)$ are finite for every $i \geq 0$. Therefore, it suffices to prove (i) in each case. We accordingly proceed as follows.

(a): Let $\mathfrak{a}= (a_{1},...,a_{n})$, and $\underline{a}= a_{1},...,a_{n}$. We argue by induction on $d=\dim_{R}(M)$. If $d=0$, then $H^{d}_{\mathfrak{a}}(M)\cong\Gamma_{\mathfrak{a}}(M)$ is finitely generated and the result obviously holds. Now suppose that $d \geq 1$ and the result holds for $d-1$. Since $d \geq 1$, we may assume that there is an element $a \in \mathfrak{a}$ which is a nonzerodivisor on $M$. The short exact sequence
$$0 \rightarrow M \xrightarrow{a} M \rightarrow M/aM \rightarrow 0$$
yields the exact sequence
\begin{equation} \label{eq:2.3.3.1}
H^{d-1}_{\mathfrak{a}}(M/aM) \rightarrow H^{d}_{\mathfrak{a}}(M) \xrightarrow{a} H^{d}_{\mathfrak{a}}(M) \rightarrow H^{d}_{\mathfrak{a}}(M/aM)=0,
\end{equation}
where the vanishing is due to $\dim_{R}(M/aM) \leq d-1$ in light of Grothendieck's Vanishing Theorem. If $\dim_{R}(M/aM) < d-1$, then $H^{d-1}_{\mathfrak{a}}(M/aM)=0$, and thus $\Ext^{i}_{R}\left(R/\mathfrak{a},H^{d-1}_{\mathfrak{a}}(M/aM)\right)=0$ for every $i \geq 0$. If $\dim_{R}(M/aM) = d-1$, then by the induction hypothesis, $\Ext^{i}_{R}\left(R/\mathfrak{a},H^{d-1}_{\mathfrak{a}}(M/aM)\right)$ is finitely generated for every $i \geq 0$. The short exact sequence
$$0 \rightarrow K \rightarrow H^{d-1}_{\mathfrak{a}}(M/aM) \rightarrow (0:_{H^{d}_{\mathfrak{a}}(M)}a) \rightarrow 0$$
extracted from the exact sequence \eqref{eq:2.3.3.1}, yields the exact sequence
$$0 \rightarrow \Hom_{R}(R/\mathfrak{a},K) \rightarrow \Hom_{R}\left(R/\mathfrak{a},H^{d-1}_{\mathfrak{a}}(M/aM)\right) \rightarrow \Hom_{R}\left(R/\mathfrak{a},(0:_{H^{d}_{\mathfrak{a}}(M)}a)\right)$$ $$\rightarrow \Ext^{1}_{R}(R/\mathfrak{a},K).$$
But as $a \in \mathfrak{a}$, we see that
$$\Hom_{R}\left(R/\mathfrak{a},(0:_{H^{d}_{\mathfrak{a}}(M)}a)\right) \cong \left(0:_{(0:_{H^{d}_{\mathfrak{a}}(M)}a)}\mathfrak{a}\right) = (0:_{H^{d}_{\mathfrak{a}}(M)}\mathfrak{a}) \cong \Hom_{R}\left(R/\mathfrak{a},H^{d}_{\mathfrak{a}}(M)\right),$$
so we get the exact sequence
\begin{equation} \label{eq:2.3.3.2}
0 \rightarrow \Hom_{R}(R/\mathfrak{a},K) \rightarrow \Hom_{R}\left(R/\mathfrak{a},H^{d-1}_{\mathfrak{a}}(M/aM)\right) \rightarrow \Hom_{R}\left(R/\mathfrak{a},H^{d}_{\mathfrak{a}}(M)\right)$$ $$\rightarrow \Ext^{1}_{R}(R/\mathfrak{a},K).
\end{equation}
Since $H^{d-1}_{\mathfrak{a}}(M/aM)$ is artinian, we conclude that $\ell_{R}\left(\Hom_{R}\left(R/\mathfrak{a},H^{d-1}_{\mathfrak{a}}(M/aM)\right)\right) < \infty$. Therefore, the exact sequence \eqref{eq:2.3.3.2} shows that
$$\ell_{R}\left(H_{n}(\underline{a};K)\right) = \ell_{R}\left(\Hom_{R}(R/\mathfrak{a},K)\right) < \infty.$$
Since $K$ is artinian, Lemma \ref{2.3.2} implies that $\ell_{R}\left(H_{i}(\underline{a};K)\right) < \infty$ for every $0 \leq i \leq n$, and thus $\ell_{R}\left(\Ext^{i}_{R}(R/\mathfrak{a},K)\right) < \infty$ for every $i \geq 0$ by Corollaries \ref{2.2.17} and \ref{2.2.19}.
Now from the exact sequence \eqref{eq:2.3.3.2}, we conclude that
$$\ell_{R}\left(H_{n}\left(\underline{a};H^{d}_{\mathfrak{a}}(M)\right)\right) = \ell_{R}\left(\Hom_{R}\left(R/\mathfrak{a},H^{d}_{\mathfrak{a}}(M)\right)\right) < \infty.$$
In a similar fashion, since $H^{d}_{\mathfrak{a}}(M)$ is artinian, Lemma \ref{2.3.2} implies that $\ell_{R}\left(H_{i}\left(\underline{a};H^{d}_{\mathfrak{a}}(M)\right)\right) < \infty$ for every $0 \leq i \leq n$, and thus $\ell_{R}\left(\Ext^{i}_{R}\left(R/\mathfrak{a},H^{d}_{\mathfrak{a}}(M)\right)\right) < \infty$ for every $i \geq 0$ by Corollaries \ref{2.2.17} and \ref{2.2.19}. This shows that $H^{d}_{\mathfrak{a}}(M)$ is $\mathfrak{a}$-cofinite.

(b): We note that $H^{0}_{\mathfrak{a}}(M) \cong \Gamma_{\mathfrak{a}}(M)$ is finitely generated, and $H^{j}_{\mathfrak{a}}(M)=0$ for every $j \geq 2$ by the hypothesis, so we deduce that $\Ext^{i}_{R}\left(R/\mathfrak{a},H^{j}_{\mathfrak{a}}(M)\right)$ is finitely generated for every $j\neq 1$ and every $i \geq 0$. Therefore, it remains to show that $\Ext^{i}_{R}\left(R/\mathfrak{a},H^{1}_{\mathfrak{a}}(M)\right)$ is finitely generated for every $i \geq 0$. Since $H^{1}_{\mathfrak{a}}(M) \cong H^{1}_{\mathfrak{a}}\left(M/\Gamma_{\mathfrak{a}}(M)\right)$, by replacing $M$ with $M/\Gamma_{\mathfrak{a}}(M)$, we may assume that $\Gamma_{\mathfrak{a}}(M)=0$. Consider the spectral sequence
$$E^{2}_{p,q}= \Ext^{-p}_{R}\left(R/\mathfrak{a},H^{-q}_{\mathfrak{a}}(M)\right) \underset {p} \Rightarrow \Ext^{-p-q}_{R}(R/\mathfrak{a},M)$$
from Lemma \ref{2.2.9}.
It follows that $E^{2}_{p,q}=0$ for every $p \in \mathbb{Z}$ and $q \neq -1$. Therefore, the spectral sequence collapses and yields the isomorphism
$$\Ext^{-p}_{R}\left(R/\mathfrak{a},H^{1}_{\mathfrak{a}}(M)\right) \cong \Ext^{-p+1}_{R}(R/\mathfrak{a},M)$$
for every $p \in \mathbb{Z}$. As $\Ext^{-p+1}_{R}(R/\mathfrak{a},M)$ is finitely generated, we are done.

(c): If $\dim_{R}(M)=0$, then $H^{j}_{\mathfrak{a}}(M)=0$ for every $j \geq 1$. Further, $H^{0}_{\mathfrak{a}}(M) \cong \Gamma_{\mathfrak{a}}(M)$ is finitely generated, so we are through. If $\dim_{R}(M)=1$, then $H^{j}_{\mathfrak{a}}(M)=0$ for every $j \geq 2$, and the result follows similar to the proof of (b). Hence suppose that $\dim_{R}(M)=2$. We note that $H^{0}_{\mathfrak{a}}(M) \cong \Gamma_{\mathfrak{a}}(M)$ is finitely generated, and $H^{j}_{\mathfrak{a}}(M)=0$ for every $j \geq 3$, so we deduce that $\Ext^{i}_{R}\left(R/\mathfrak{a},H^{j}_{\mathfrak{a}}(M)\right)$ is finitely generated for every $j \neq 1,2$ and every $i \geq 0$. On the other hand, by (a), $\Ext^{i}_{R}\left(R/\mathfrak{a},H^{2}_{\mathfrak{a}}(M)\right)$ is finitely generated for every $i \geq 0$. Therefore, it remains to show that $\Ext^{i}_{R}\left(R/\mathfrak{a},H^{1}_{\mathfrak{a}}(M)\right)$ is finitely generated for every $i \geq 0$. Consider the spectral sequence
$$E^{2}_{p,q}= \Ext^{-p}_{R}\left(R/\mathfrak{a},H^{-q}_{\mathfrak{a}}(M)\right) \underset {p} \Rightarrow \Ext^{-p-q}_{R}(R/\mathfrak{a},M)$$
from Lemma \ref{2.2.9}.
It follows that $E^{2}_{p,q}$ is finitely generated for every $p \in \mathbb{Z}$ and $q \neq -1$. As $E^{r}_{p,q}$ is a subquotient of $E^{2}_{p,q}$ for every $r \geq 2$, we conclude that $E^{r}_{p,q}$ is finitely generated for every $p \in \mathbb{Z}$, $q \neq -1$, and $r \geq 2$. Consider the following homomorphisms for every $r \geq 2$ and $p,q \in \mathbb{Z}$:
$$E^{r}_{p+r,q-r+1} \xrightarrow{d^{r}_{p+r,q-r+1}} E^{r}_{p,q} \xrightarrow{d^{r}_{p,q}} E^{r}_{p-r,q+r-1}$$
If $q \neq -r$, then $q+r-1 \neq -1$, and thus $E^{r}_{p-r,q+r-1}$ is finitely generated. It follows that $\im d^{r}_{p,q}$ is finitely generated for every $q \neq -r$. If $q=-r$, then $q \neq -1$ as $r \geq 2$, and thus $E^{r}_{p,q}=E^{r}_{p,-r}$ is finitely generated, so $\im d^{r}_{p,q}= \im d^{r}_{p,-r}$ is finitely generated. Therefore, $\im d^{r}_{p,q}$ is finitely generated for every $r \geq 2$ and $p,q \in \mathbb{Z}$. Since $\Ext^{-p-q}_{R}(R/\mathfrak{a},M)$ is finitely generated and $E^{\infty}_{p,q}$ is a subquotient of $\Ext^{-p-q}_{R}(R/\mathfrak{a},M)$ for every $p,q \in \mathbb{Z}$, we infer that $E^{\infty}_{p,q}$ is finitely generated for every $p,q \in \mathbb{Z}$. On the other hand, there is an integer $r_{0} \geq 2$ such that $E^{\infty}_{p,q} = E^{r}_{p,q}$ for every $r \geq r_{0}+1$. It follows that $E^{r_{0}+1}_{p,-1}$ is finitely generated for every $p \in \mathbb{Z}$. But $$E^{r_{0}+1}_{p,-1} \cong \ker d^{r_{0}}_{p,-1}/\im d^{r_{0}}_{p+r_{0},-r_{0}}$$
for every $p \in \mathbb{Z}$, so there is a short exact sequence
$$0 \rightarrow \im d^{r_{0}}_{p+r_{0},-r_{0}} \rightarrow \ker d^{r_{0}}_{p,-1} \rightarrow E^{r_{0}+1}_{p,-1} \rightarrow 0,$$
which implies that $\ker d^{r_{0}}_{p,-1}$ is finitely generated for every $p \in \mathbb{Z}$. Furthermore, the short exact sequence
$$0 \rightarrow \ker d^{r_{0}}_{p,-1} \rightarrow E^{r_{0}}_{p,-1} \rightarrow \im d^{r_{0}}_{p,-1} \rightarrow 0$$
implies that $E^{r_{0}}_{p,-1}$ is finitely generated. Continuing in this fashion, we deduce that $E^{2}_{p,-1} = \Ext^{-p}_{R}\left(R/\mathfrak{a},H^{1}_{\mathfrak{a}}(M)\right)$ is finitely generated for every $p \in \mathbb{Z}$.

(d): We prove by induction on $j$ that if $N$ is an $R$-module with $\Ext^{i}_{R}(R/\mathfrak{a},N)$ finitely generated for every $i \geq 0$, then $\Ext^{i}_{R}\left(R/\mathfrak{a},H^{j}_{\mathfrak{a}}(N)\right)$ is finitely generated for every $i \geq 0$. Note that since $M$ is finitely generated, we get the desired result by letting $N=M$.

First we deal with the base step. Let $j=0$. Let $N$ be an $R$-module such that $\Ext^{i}_{R}(R/\mathfrak{a},N)$ is finitely generated for every $i \geq 0$. The short exact sequence
$$0 \rightarrow \Gamma_{\mathfrak{a}}(N) \rightarrow N \rightarrow N/\Gamma_{\mathfrak{a}}(N) \rightarrow 0$$
yields the exact sequence
\begin{equation} \label{eq:2.3.3.3}
0 \rightarrow \Hom_{R}\left(R/\mathfrak{a},\Gamma_{\mathfrak{a}}(N)\right) \rightarrow \Hom_{R}(R/\mathfrak{a},N) \rightarrow \Hom_{R}\left(R/\mathfrak{a},N/\Gamma_{\mathfrak{a}}(N)\right) $$$$ \rightarrow \Ext_{R}^{1}\left(R/\mathfrak{a},\Gamma_{\mathfrak{a}}(N)\right) \rightarrow \Ext_{R}^{1}(R/\mathfrak{a},N).
\end{equation}
Since $R/\mathfrak{a}$ is $\mathfrak{a}$-torsion and $N/\Gamma_{\mathfrak{a}}(N)$ is $\mathfrak{a}$-torsionfree, we have $\Hom_{R}\left(R/\mathfrak{a},N/\Gamma_{\mathfrak{a}}(N)\right)=0$. It follows from the exact sequence \eqref{eq:2.3.3.3} that $\Hom_{R}\left(R/\mathfrak{a},\Gamma_{\mathfrak{a}}(N)\right)$ and $\Ext_{R}^{1}\left(R/\mathfrak{a},\Gamma_{\mathfrak{a}}(N)\right)$ are finitely generated. We claim that $\Ext_{R}^{i}\left(R/\mathfrak{a},\Gamma_{\mathfrak{a}}(N)\right)$ is finitely generated for every $i \geq 0$.

First suppose that $\dim(R/\mathfrak{a})=0$. Then
$$\Supp_R\left(\Hom_{R}\left(R/\mathfrak{a},\Gamma_{\mathfrak{a}}(N)\right)\right) \subseteq \Supp_{R}\left(\Gamma_{\mathfrak{a}}(N)\right) \subseteq V(\mathfrak{a}) \subseteq \Max(R).$$
But $\Hom_{R}\left(R/\mathfrak{a},\Gamma_{\mathfrak{a}}(N)\right)$ is finitely generated, so
$$\ell_{R}\left(H_{n}\left(\underline{a};\Gamma_{\mathfrak{a}}(N)\right)\right) = \ell_{R}\left(\Hom_{R}\left(R/\mathfrak{a},\Gamma_{\mathfrak{a}}(N)\right)\right) < \infty.$$
In particular, $(0:_{\Gamma_{\mathfrak{a}}(N)}\mathfrak{a}) \cong \Hom_{R}\left(R/\mathfrak{a},\Gamma_{\mathfrak{a}}(N)\right)$ is artinian, so $\Gamma_{\mathfrak{a}}(N)$ is artinian by \cite[Theorem 1.3]{Me3}. It follows that $\ell_{R}\left(H_{i}\left(\underline{a};\Gamma_{\mathfrak{a}}(N)\right)\right) < \infty$ for every $0 \leq i \leq n$ by Lemma \ref{2.3.2}, and consequently, $\ell_{R}\left(\Ext_{R}^{i}\left(R/\mathfrak{a},\Gamma_{\mathfrak{a}}(N)\right)\right) < \infty$ for every $i \geq 0$ by Corollaries \ref{2.2.17} and \ref{2.2.19}.

Next suppose that $\dim(R/\mathfrak{a})=1$. Assume to the contrary that there is an integer $i \geq 2$ for which $\Ext^{i}_{R}\left(R/\mathfrak{a},\Gamma_{\mathfrak{a}}(N)\right)$ is not finitely generated. Let $\mathcal{S}$ be the class of all $R$-modules $L$ with $\Supp_{R}(L) \subseteq V(\mathfrak{a})$, such that $\Ext^{i}_{R}(R/\mathfrak{a},L)$ is finitely generated for $i=0,1$, but $\Ext^{i}_{R}(R/\mathfrak{a},L)$ fails to be finitely generated for some $i \geq 2$. Let $\mathfrak{H}:= \{\ann_{R}(L) \suchthat L \in \mathcal{S} \}$. As $\Gamma_{\mathfrak{a}}(N) \in \mathcal{S}$, we have $\mathfrak{H} \neq \emptyset$. Since $R$ is noetherian, there is an $R$-module $L$ with $\ann_{R}(L)$ maximal in $\mathfrak{H}$.
Let
$$U:=R \setminus \bigcup_{\mathfrak{p}\in \Assh_{R}(R/\mathfrak{a})} \mathfrak{p},$$
and form the localization $U^{-1}R$. Then we have $V(U^{-1}\mathfrak{a}) \subseteq \Max(U^{-1}R)$. Indeed if $U^{-1}\mathfrak{q} \in V(U^{-1}\mathfrak{a})$, then $U^{-1}\mathfrak{a} \subseteq U^{-1}\mathfrak{q}$, so $\mathfrak{a} \subseteq \mathfrak{q} \subseteq \mathfrak{p}$ for some $\mathfrak{p} \in \Assh_{R}(R/\mathfrak{a})$. Since
$$\Assh_{R}(R/\mathfrak{a}) \subseteq \Min_{R}(R/\mathfrak{a}) = \Min(\mathfrak{a}),$$
we see that $\mathfrak{q} = \mathfrak{p}$. Therefore, $U^{-1}\mathfrak{q}=U^{-1}\mathfrak{p} \in \Max(U^{-1}R)$.
It now follows that
\begin{equation*}
\begin{split}
\Supp_{U^{-1}R}\left(\Hom_{U^{-1}R}(U^{-1}R/U^{-1}\mathfrak{a},U^{-1}L)\right) & \subseteq \Supp_{U^{-1}R}(U^{-1}L) \\
 & = \left\{U^{-1}\mathfrak{p} \suchthat \mathfrak{p} \in \Supp_{R}(L), \mathfrak{p} \cap U=\emptyset \right\} \\
 & \subseteq V(U^{-1}\mathfrak{a}) \\
 & \subseteq \Max(U^{-1}R).
\end{split}
\end{equation*}
But
$$\Hom_{U^{-1}R}\left(U^{-1}R/U^{-1}\mathfrak{a},U^{-1}L\right) \cong U^{-1} \Hom_{R}(R/\mathfrak{a},L)$$
is a finitely generated $U^{-1}R$-module, so
$$\ell_{U^{-1}R}\left(H_{n}\left(\underline{a}/1;U^{-1}L\right)\right) = \ell_{U^{-1}R}\left(\Hom_{U^{-1}R}\left(U^{-1}R/U^{-1}\mathfrak{a},U^{-1}L\right)\right)<\infty.$$
In particular,
$$(0:_{U^{-1}L}U^{-1}\mathfrak{a}) \cong \Hom_{U^{-1}R}\left(U^{-1}R/U^{-1}\mathfrak{a},U^{-1}L\right)$$
is an artinian $U^{-1}R$-module, so $U^{-1}L=\Gamma_{U^{-1}\mathfrak{a}}(U^{-1}L)$ is an artinian $U^{-1}R$-module by \cite[Theorem 1.3]{Me3}. It follows from Lemma \ref{2.3.2} that $\ell_{U^{-1}R}\left(H_{i}\left(\underline{a}/1;U^{-1}L\right)\right)< \infty$ for every $0 \leq i \leq n$, and consequently,
$$\ell_{U^{-1}R}\left(\Ext_{U^{-1}R}^{i}\left(U^{-1}R/U^{-1}\mathfrak{a},U^{-1}L\right)\right)<\infty$$
for every $i \geq 0$ by Corollaries \ref{2.2.17} and \ref{2.2.19}. It follows from Corollary \ref{2.2.17} that
$$U^{-1}(L/\mathfrak{a}L) \cong U^{-1}L/(U^{-1}\mathfrak{a})(U^{-1}L) \cong (U^{-1}R/U^{-1}\mathfrak{a}) \otimes_{U^{-1}R} U^{-1}L$$
is a finitely generated $U^{-1}R$-module, say
$$U^{-1}(L/\mathfrak{a}L) = \left\langle(x_{1}+\mathfrak{a}L)/1,...,(x_{1}+\mathfrak{a}L)/1)\right\rangle$$
for some elements $x_{1},...,x_{n}\in L$. Let $K=\langle x_{1},...,x_{n} \rangle \subseteq L$. Then
$$U^{-1}(L/\mathfrak{a}L) = U^{-1}\left((K+\mathfrak{a}L)/\mathfrak{a}L\right),$$
and thus $U^{-1}L=U^{-1}(K+\mathfrak{a}L)$. This implies that
$$U^{-1}(L/K) = (U^{-1}\mathfrak{a})\left(U^{-1}(L/K)\right).$$
As $U^{-1}(L/K) \cong U^{-1}L/U^{-1}K$ is an artinian $U^{-1}R$-module, there is an element $a \in \mathfrak{a}$ such that $U^{-1}(L/K)=aU^{-1}(L/K)$, i.e. $U^{-1}\left(L/(aL+K)\right)=0$. Hence
$$\Supp_{R}\left(L/(aL+K)\right) \subseteq V(\mathfrak{a})\setminus \Assh_{R}(R/\mathfrak{a}) \subseteq \Max(R).$$
The short exact sequence
\begin{equation} \label{eq:2.3.3.4}
0 \rightarrow K \rightarrow L \rightarrow L/K \rightarrow 0
\end{equation}
yields the exact sequence
\begin{equation} \label{eq:2.3.3.5}
\Hom_{R}(R/\mathfrak{a},L) \rightarrow \Hom_{R}(R/\mathfrak{a},L/K) \rightarrow \Ext^{1}_{R}(R/\mathfrak{a},K) \rightarrow \Ext^{1}_{R}(R/\mathfrak{a},L) $$$$ \rightarrow \Ext^{1}_{R}(R/\mathfrak{a},L/K) \rightarrow \Ext^{2}_{R}(R/\mathfrak{a},K).
\end{equation}
Since $K$ is finitely generated, $\Ext^{1}_{R}(R/\mathfrak{a},K)$ and $\Ext^{2}_{R}(R/\mathfrak{a},K)$ are finitely generated. On the other hand, $\Hom_{R}(R/\mathfrak{a},L)$ and $\Ext^{1}_{R}(R/\mathfrak{a},L)$ are finitely generated. Therefore, we deduce form the exact sequence \eqref{eq:2.3.3.5} that $\Hom_{R}(R/\mathfrak{a},L/K)$ and $\Ext^{1}_{R}(R/\mathfrak{a},L/K)$ are finitely generated. Moreover, the short exact sequence \eqref{eq:2.3.3.4} induces the exact sequence
\begin{equation} \label{eq:2.3.3.6}
\Ext^{i}_{R}(R/\mathfrak{a},K) \rightarrow \Ext^{i}_{R}(R/\mathfrak{a},L) \rightarrow \Ext^{i}_{R}(R/\mathfrak{a},L/K).
\end{equation}
Since $K$ is finitely generated, $\Ext^{i}_{R}(R/\mathfrak{a},K)$ is finitely generated for every $i \geq 0$. If $\Ext^{i}_{R}(R/\mathfrak{a},L/K)$ is finitely generated for every $i \geq 0$, then the exact sequence \eqref{eq:2.3.3.6} implies that $\Ext^{i}_{R}(R/\mathfrak{a},L)$ is finitely generated for every $i \geq 0$, which is in contrast with our assumption. Therefore, there is an integer $i \geq 2$ for which $\Ext^{i}_{R}(R/\mathfrak{a},L/K)$ is not finitely generated. Since
$$\Supp_{R}(L/K) \subseteq \Supp_{R}(L) \subseteq V(\mathfrak{a}),$$
we see that $L/K \in \mathcal{S}$, so that $\ann_{R}(L/K) \in \mathfrak{H}$. But $\ann_{R}(L) \subseteq \ann_{R}(L/K)$, and $\ann_{R}(L)$ is a maximal element of $\mathfrak{H}$, so we must have $\ann_{R}(L) = \ann_{R}(L/K)$. As a consequence, we may replace $L$ with $L/K$ and accordingly assume that $K=0$. Now the inclusion $\Supp_{R}\left(L/(aL+K)\right) \subseteq \Max(R)$ simplifies to $\Supp_{R}(L/aL) \subseteq \Max(R)$.
If $aL=0$, then $\Supp_{R}(L) \subseteq \Max(R)$. It follows that
$$\Supp_{R}\left(\Hom_{R}(R/\mathfrak{a},L)\right) \subseteq \Supp_{R}(L) \subseteq \Max(R).$$
But $\Hom_{R}(R/\mathfrak{a},L)$ is finitely generated, so
$$\ell_{R}\left(H_{n}(\underline{a};L)\right) = \ell_{R}\left(\Hom_{R}(R/\mathfrak{a},L)\right) < \infty.$$
In particular, $(0:_{L}\mathfrak{a}) \cong \Hom_{R}(R/\mathfrak{a},L)$ is artinian, so $L= \Gamma_{\mathfrak{a}}(L)$ is artinian by \cite[Theorem 1.3]{Me3}. It follows from Lemma \ref{2.3.2} that $\ell_{R}\left(H_{i}(\underline{a};L)\right) < \infty$ for every $0 \leq i \leq n$, and consequently, $\ell_{R}\left(\Ext_{R}^{i}(R/\mathfrak{a},L)\right) < \infty$ for every $i \geq 0$ by Corollaries \ref{2.2.17} and \ref{2.2.19}, which contradicts the assumption. Therefore, $a \notin \ann_{R}(L)$.

The short exact sequence
\begin{equation} \label{eq:2.3.3.7}
0 \rightarrow aL \rightarrow L \rightarrow L/aL \rightarrow 0
\end{equation}
yields the exact sequence
\begin{equation*}
0 \rightarrow \Hom_{R}(R/\mathfrak{a},aL) \rightarrow \Hom_{R}(R/\mathfrak{a},L),
\end{equation*}
which implies that $\Hom_{R}(R/\mathfrak{a},aL)$ is finitely generated.
The short exact sequence
\begin{equation} \label{eq:2.3.3.8}
0 \rightarrow (0:_{L}a) \rightarrow L \rightarrow aL \rightarrow 0
\end{equation}
yields the exact sequence
\begin{equation} \label{eq:2.3.3.9}
0 \rightarrow \Hom_{R}\left(R/\mathfrak{a},(0:_{L}a)\right) \rightarrow \Hom_{R}(R/\mathfrak{a},L) \rightarrow \Hom_{R}(R/\mathfrak{a},aL) $$$$ \rightarrow \Ext^{1}_{R}\left(R/\mathfrak{a},(0:_{L}a)\right) \rightarrow \Ext^{1}_{R}(R/\mathfrak{a},L).
\end{equation}
The exact sequence \eqref{eq:2.3.3.9} implies that $\Hom_{R}\left(R/\mathfrak{a},(0:_{L}a)\right)$ and $\Ext^{1}_{R}\left(R/\mathfrak{a},(0:_{L}a)\right)$ are finitely generated. We note that $\ann_{R}(L) \subset \ann_{R}\left((0:_{L}a)\right)$, since $a \in \ann_{R}\left((0:_{L}a)\right) \setminus \ann_{R}(L)$. If $\Ext^{i}_{R}\left(R/\mathfrak{a},(0:_{L}a)\right)$ is not finitely generated for some $i \geq 2$, then $\ann_{R}\left((0:_{L}a)\right) \in \mathfrak{H}$ which contradicts the maximality of $\ann_{R}(L)$. Therefore, $\Ext^{i}_{R}\left(R/\mathfrak{a},(0:_{L}a)\right)$ is finitely generated for every $i \geq 0$. Now the short exact sequence \eqref{eq:2.3.3.8} induces the exact sequence
$$\Ext^{1}_{R}(R/\mathfrak{a},L) \rightarrow \Ext^{1}_{R}(R/\mathfrak{a},aL) \rightarrow \Ext^{2}_{R}\left(R/\mathfrak{a},(0:_{L}a)\right),$$
which implies that $\Ext^{1}_{R}(R/\mathfrak{a},aL)$ is finitely generated. The short exact sequence \eqref{eq:2.3.3.7} yields the exact sequence
$$\Hom_{R}(R/\mathfrak{a},L) \rightarrow \Hom_{R}(R/\mathfrak{a},L/aL) \rightarrow \Ext^{1}_{R}(R/\mathfrak{a},aL),$$
which implies that $\Hom_{R}(R/\mathfrak{a},L/aL)$ is finitely generated. But
$$\Supp_{R}\left(\Hom_{R}\left(R/\mathfrak{a},L/aL\right)\right) \subseteq \Supp_{R}(L/aL) \subseteq \Max(R).$$
It follows that
$$\ell_{R}\left(H_{n}\left(\underline{a};L/aL\right)\right) = \ell_{R}\left(\Hom_{R}\left(R/\mathfrak{a},L/aL\right)\right) < \infty.$$
In particular, $(0:_{L/aL}\mathfrak{a}) \cong \Hom_{R}(R/\mathfrak{a},L/aL)$ is artinian. But
$$\Supp_{R}(L/aL) \subseteq \Supp_{R}(L) \subseteq V(\mathfrak{a}),$$
so it follows from \cite[Theorem 1.3]{Me3} that $L/aL= \Gamma_{\mathfrak{a}}(L/aL)$ is artinian. Hence $\ell_{R}\left(H_{i}\left(\underline{a};L/aL\right)\right) < \infty$ for every $0 \leq i \leq n$ by Lemma \ref{2.3.3}, and consequently, $\ell_{R}\left(\Ext_{R}^{i}\left(R/\mathfrak{a},L/aL\right)\right) < \infty$ for every $i \geq 0$ by Corollaries \ref{2.2.17} and \ref{2.2.19}.

Now consider the exact sequence
$$0 \rightarrow (0:_{L}a) \rightarrow L \xrightarrow{a} L \rightarrow L/aL \rightarrow 0,$$
and break it into two short exact sequences
\begin{equation} \label{eq:2.3.3.10}
0 \rightarrow (0:_{L}a) \rightarrow L \xrightarrow{f} aL \rightarrow 0,
\end{equation}
and
\begin{equation} \label{eq:2.3.3.11}
0 \rightarrow aL \xrightarrow{g} L \rightarrow L/aL \rightarrow 0,
\end{equation}
where $gf=a1^{L}$. Fix $i \geq 2$. From the short exact sequence \eqref{eq:2.3.3.10}, we get the exact sequence
$$\Ext^{i}_{R}\left(R/\mathfrak{a},(0:_{L}a)\right) \rightarrow \Ext^{i}_{R}(R/\mathfrak{a},L) \xrightarrow{\Ext^{i}_{R}(R/\mathfrak{a},f)} \Ext^{i}_{R}(R/\mathfrak{a},aL),$$
which in turn implies that $\ker\left(\Ext^{i}_{R}(R/\mathfrak{a},f)\right)$ is finitely generated as it is a quotient of $\Ext^{i}_{R}\left(R/\mathfrak{a},(0:_{L}a)\right)$. From the short exact sequence \eqref{eq:2.3.3.11}, we get the exact sequence
$$\Ext^{i-1}_{R}\left(R/\mathfrak{a},L/aL\right) \rightarrow \Ext^{i}_{R}(R/\mathfrak{a},aL) \xrightarrow{\Ext^{i}_{R}(R/\mathfrak{a},g)} \Ext^{i}_{R}(R/\mathfrak{a},L),$$
which in turn implies that $\ker\left(\Ext^{i}_{R}(R/\mathfrak{a},g)\right)$ is finitely generated as it is a quotient of $\Ext^{i-1}_{R}\left(R/\mathfrak{a},L/aL\right)$. On the other hand, we have
$$\Ext^{i}_{R}(R/\mathfrak{a},g)\Ext^{i}_{R}(R/\mathfrak{a},f)=\Ext^{i}_{R}\left(R/\mathfrak{a},a1^{L}\right).$$
Hence we get the exact sequence
$$0 \rightarrow \ker\left(\Ext^{i}_{R}(R/\mathfrak{a},f)\right) \rightarrow \ker\left(\Ext^{i}_{R}\left(R/\mathfrak{a},a1^{L}\right)\right) \rightarrow \ker\left(\Ext^{i}_{R}(R/\mathfrak{a},g)\right),$$
which implies that $\ker\left(\Ext^{i}_{R}\left(R/\mathfrak{a},a1^{L}\right)\right)$ is finitely generated. However, $\Ext^{i}_{R}\left(R/\mathfrak{a},a1^{L}\right)=0$ since $a \in \mathfrak{a}$. Therefore, $\Ext^{i}_{R}(R/\mathfrak{a},L) = \ker\left(\Ext^{i}_{R}\left(R/\mathfrak{a},a1^{L}\right)\right)$ is finitely generated. This contradiction establishes the base case of the induction.

Next we deal with the inductive step. Let $j \geq 1$, and suppose that the result holds for $j-1$. Let $N$ be an $R$-module such that $\Ext^{i}_{R}(R/\mathfrak{a},N)$ is finitely generated for every $i \geq 0$. By the base step, $\Ext^{i}_{R}\left(R/\mathfrak{a},\Gamma_{\mathfrak{a}}(N)\right)$ is finitely generated for every $i \geq 0$. The short exact sequence
$$0 \rightarrow \Gamma_{\mathfrak{a}}(N) \rightarrow N \rightarrow N/\Gamma_{\mathfrak{a}}(N) \rightarrow 0$$
yields the exact sequence
$$\Ext^{i}_{R}(R/\mathfrak{a},N) \rightarrow \Ext^{i}_{R}\left(R/\mathfrak{a},N/\Gamma_{\mathfrak{a}}(N)\right) \rightarrow \Ext^{i+1}_{R}\left(R/\mathfrak{a},\Gamma_{\mathfrak{a}}(N)\right),$$
which in turn implies that $\Ext^{i}_{R}\left(R/\mathfrak{a},N/\Gamma_{\mathfrak{a}}(N)\right)$ is finitely generated for every $i \geq 0$. Moreover, since $H^{j}_{\mathfrak{a}}\left(N/\Gamma_{\mathfrak{a}}(N)\right) \cong H^{j}_{\mathfrak{a}}(N)$, we may replace $N$ with $N/\Gamma_{\mathfrak{a}}(N)$, and consequently assume that $\Gamma_{\mathfrak{a}}(N)=0$. We then have $\Gamma_{\mathfrak{a}}\left(E_{R}(N)\right) \cong E_{R}\left(\Gamma_{\mathfrak{a}}(N)\right)=0$. In particular, $\Hom_{R}\left(R/\mathfrak{a},E_{R}(N)\right) \cong (0:_{E_{R}(N)}\mathfrak{a})=0$. Thus the short exact sequence
\begin{equation} \label{eq:2.3.3.12}
0 \rightarrow N \rightarrow E_{R}(N) \rightarrow C \rightarrow 0
\end{equation}
yields the exact sequence
\begin{equation} \label{eq:2.3.3.13}
0=\Hom_{R}\left(R/\mathfrak{a},E_{R}(N)\right) \rightarrow \Hom_{R}(R/\mathfrak{a},C) \rightarrow \Ext^{1}_{R}(R/\mathfrak{a},N) \rightarrow \Ext^{1}_{R}\left(R/\mathfrak{a},E_{R}(N)\right) $$$$ \rightarrow \Ext^{1}_{R}(R/\mathfrak{a},C) \rightarrow \Ext^{2}(R/\mathfrak{a},N) \rightarrow \Ext^{2}_{R}\left(R/\mathfrak{a},E_{R}(N)\right) \rightarrow \cdots.
\end{equation}
Since $\Ext^{i}_{R}\left(R/\mathfrak{a},E_{R}(N)\right)=0$ for every $i \geq 1$, we infer from the exact sequence \eqref{eq:2.3.3.13} that $\Ext^{i-1}_{R}(R/\mathfrak{a},C) \cong \Ext^{i}_{R}(R/\mathfrak{a},N)$ for every $i \geq 1$. Therefore, $\Ext^{i-1}_{R}(R/\mathfrak{a},C)$ is finitely generated for every $i \geq 1$. Applying the induction hypothesis to $C$, we conclude that $\Ext^{i}_{R}\left(R/\mathfrak{a},H^{j-1}_{\mathfrak{a}}(C)\right)$ is finitely generated for every $i \geq 0$. From the short exact sequence \eqref{eq:2.3.3.12}, we get the exact sequence
$$0=H^{j-1}_{\mathfrak{a}}\left(E_{R}(N)\right) \rightarrow H^{j-1}_{\mathfrak{a}}(C) \rightarrow H^{j}_{\mathfrak{a}}(N) \rightarrow H^{j}_{\mathfrak{a}}\left(E_{R}(N)\right)=0,$$
which in turn shows that $H^{j-1}_{\mathfrak{a}}(C) \cong H^{j}_{\mathfrak{a}}(N)$ for every $j \geq 1$. Therefore, $\Ext^{i}_{R}\left(R/\mathfrak{a},H^{j}_{\mathfrak{a}}(N)\right)$ is finitely generated for every $i \geq 0$.
\end{prf}

Our next goal in this section is to show that the Serre property of being noetherian satisfies the condition $\mathfrak{D}_{\mathfrak{a}}$ described in the previous section. We discuss it through a chain of interrelated lemmas in the sequel.

\begin{lemma} \label{2.3.4}
Let $\mathfrak{a}$ be an ideal of $R$ and $M$ an $R$-module. If $M/\mathfrak{a}M$ is a finitely generated $R$-module, then $\widehat{M}^{\mathfrak{a}}$
is a finitely generated $\widehat{R}^{\mathfrak{a}}$-module.
\end{lemma}

\begin{prf}
As $M/\mathfrak{a}M$ is finitely generated, there is a finitely generated submodule $N$ of $M$ such that $M=N+\mathfrak{a}M$. Let $\iota:N\rightarrow M$
be the inclusion map. Then by \cite[Lemma 1.2]{Si1}, the $\widehat{R}^{\mathfrak{a}}$-homomorphism $\widehat{\iota}:\widehat{N}^{\mathfrak{a}}\rightarrow
\widehat{M}^{\mathfrak{a}}$ is surjective.
On the other hand, since $\widehat{N}^{\mathfrak{a}} \cong \widehat{R}^{\mathfrak{a}} \otimes_{R}N$, it is clear that $\widehat{N}^{\mathfrak{a}}$ is
a finitely generated $\widehat{R}^{\mathfrak{a}}$-module. It then follows that $\widehat{M}^{\mathfrak{a}}$ is a finitely generated
$\widehat{R}^{\mathfrak{a}}$-module.
\end{prf}

We recall a definition.

\begin{definition} \label{2.3.5}
Let $\mathfrak{a}$ be an ideal of $R$, and $M$ an $R$-module. Then:
\begin{enumerate}
\item[(i)] $M$ is said to be \textit{$\mathfrak{a}$-adically separated} \index{$\mathfrak{a}$-adically separated} if the completion map $\theta_{M}:M\rightarrow \widehat{M}^{\mathfrak{a}}$ is a monomorphism.
\item[(ii)] $M$ is said to be \textit{$\mathfrak{a}$-adically quasi-complete} \index{$\mathfrak{a}$-adically quasi-complete} if the completion map $\theta_{M}:M\rightarrow \widehat{M}^{\mathfrak{a}}$ is an epimorphism.
\item[(iii)] $M$ is said to be \textit{$\mathfrak{a}$-adically complete} \index{$\mathfrak{a}$-adically complete} if the completion map $\theta_{M}:M\rightarrow \widehat{M}^{\mathfrak{a}}$ is an isomorphism.
\end{enumerate}
\end{definition}

The next lemma may be of independent interest.

\begin{lemma} \label{2.3.6}
Let $\mathfrak{a}$ be an ideal of $R$. Then the following assertions hold:
\begin{enumerate}
\item[(i)] Any submodule of an $\mathfrak{a}$-adically separated $R$-module is $\mathfrak{a}$-adically separated.
\item[(ii)] Any homomorphic image of an $\mathfrak{a}$-adically quasi-complete $R$-module is $\mathfrak{a}$-adically quasi-complete.
\item[(iii)] If $f:M\rightarrow N$ is a homomorphism of $\mathfrak{a}$-adically complete $R$-modules, then both $\ker f$ and $\im f$ are
$\mathfrak{a}$-adically complete.
\end{enumerate}
\end{lemma}

\begin{prf}
(i): As $\ker \theta_{M}=\bigcap_{i=1}^{\infty} \mathfrak{a}^{i}M$ for any $R$-module $M$, the assertion is clear.

(ii): If $M$ is $\mathfrak{a}$-adically quasi-complete and $f:M\rightarrow N$ is an epimorphism, then it follows from \cite[Lemma 1.2]{Si1} that $\widehat{f}:\widehat{M}^{\mathfrak{a}}\rightarrow \widehat{N}^{\mathfrak{a}}$ is surjective. Therefore, the commutative diagram
\[\begin{tikzpicture}[every node/.style={midway},]
  \matrix[column sep={3em}, row sep={2.5em}]
  {\node(1) {$M$}; & \node(2) {$N$}; \\
  \node(3) {$\widehat{M}^{\mathfrak{a}}$}; & \node(4) {$\widehat{N}^{\mathfrak{a}}$}; \\};
  \draw[decoration={markings,mark=at position 1 with {\arrow[scale=1.5]{>}}},postaction={decorate},shorten >=0.5pt] (1) -- (2) node[anchor=south] {$f$};
  \draw[decoration={markings,mark=at position 1 with {\arrow[scale=1.5]{>}}},postaction={decorate},shorten >=0.5pt] (3) -- (4) node[anchor=south] {$\widehat{f}$};
  \draw[decoration={markings,mark=at position 1 with {\arrow[scale=1.5]{>}}},postaction={decorate},shorten >=0.5pt] (1) -- (3) node[anchor=west] {$\theta_{M}$};
  \draw[decoration={markings,mark=at position 1 with {\arrow[scale=1.5]{>}}},postaction={decorate},shorten >=0.5pt] (2) -- (4) node[anchor=west] {$\theta_{N}$};
\end{tikzpicture}\]
shows that $\theta_{N}$ is surjective, i.e. $N$ is $\mathfrak{a}$-adically quasi-complete.

(iii): As $N$ is $\mathfrak{a}$-adically complete, it can be seen by inspection that
$$\ker f= \bigcap_{i=1}^{\infty} \left(\ker f + \mathfrak{a}^{i}M\right),$$
i.e. $\ker f$ is a closed submodule of $M$ in the $\mathfrak{a}$-adic topology. It now follows from \cite[Proposition 1.3 (ii)]{Si1} that $\ker f$
is $\mathfrak{a}$-adically complete. On the other hand, $\im f$ is both a submodule of $N$ and a homomorphic image of $M$, so using (i) and (ii),
we infer that $\im f$ is $\mathfrak{a}$-adically complete.
\end{prf}

\begin{lemma} \label{2.3.7}
Let $\mathfrak{a}$ be an ideal of $R$ and $M$ an $R$-module. Suppose that $R$ is $\mathfrak{a}$-adically complete, $M$ is $\mathfrak{a}$-adically
quasi-complete, and $M/\mathfrak{a}M$ is a finitely generated $R$-module. Let $\theta_{M}:M \rightarrow \widehat{M}^{\mathfrak{a}}$ be the
completion map with $K= \ker \theta_{M}$. Then $\widehat{K}^{\mathfrak{a}}=0=H^{\mathfrak{a}}_{0}(K)$.
\end{lemma}

\begin{prf}
Since $M/\mathfrak{a}M$ is a finitely generated $R$-module, Lemma \ref{2.3.4} implies that $\widehat{M}^{\mathfrak{a}}$ is a finitely generated
$R$-module. From the short exact sequence
$$0\rightarrow K \rightarrow M \xrightarrow {\theta_{M}} \widehat{M}^{\mathfrak{a}} \rightarrow 0,$$
we get the exact sequence
$$\Tor^{R}_{1}\left(R/\mathfrak{a},\widehat{M}^{\mathfrak{a}}\right)\rightarrow K/\mathfrak{a}K \rightarrow M/\mathfrak{a}M,$$
which implies that $K/\mathfrak{a}K$ is a finitely generated $R$-module. Subsequently, from the short exact sequence
$$0\rightarrow K/\mathfrak{a}K \rightarrow M/\mathfrak{a}K \rightarrow \widehat{M}^{\mathfrak{a}} \rightarrow 0,$$
we deduce that $M/\mathfrak{a}K$ is a finitely generated $R$-module. It then follows that the zero submodule of $M/\mathfrak{a}K$ has a minimal
primary decomposition
$$\mathfrak{a}K/\mathfrak{a}K= \bigcap_{i=1}^{n}Q_{i}/\mathfrak{a}K,$$
which in turn gives a minimal primary decomposition
$$\mathfrak{a}K= \bigcap_{i=1}^{n}Q_{i}$$
of $\mathfrak{a}K$.
We prove that $K=\mathfrak{a}K$, by showing that $K \subseteq Q_{i}$ for every $1 \leq i \leq n$. Assume to the contrary that there is an integer
$1 \leq j \leq n$ and an element $x\in K \backslash Q_{j}$. Then we have $\mathfrak{a}x \subseteq \mathfrak{a}K \subseteq Q_{j}$. This means that
the homomorphism $M/Q_{j} \xrightarrow {a} M/Q_{j}$ is not injective for every $a\in \mathfrak{a}$. As $Q_{j}$ is a primary submodule of $M$, we
conclude that $\mathfrak{a} \subseteq \sqrt{(Q_{j}:_{R}M)}$. Therefore, there is an integer $t \geq 1$ such that $\mathfrak{a}^{t}M \subseteq Q_{j}$.
But $$K=\bigcap_{i=1}^{\infty}\mathfrak{a}^{i}M \subseteq \mathfrak{a}^{t}M \subseteq Q_{j},$$
which is a contradiction. Therefore, $K=\mathfrak{a}K$, and so by \cite[Lemma 5.1 (ii)]{Si1}, $\widehat{K}^{\mathfrak{a}}=0=H^{\mathfrak{a}}_{0}(K)$.
\end{prf}

The next lemma may be of independent interest.

\begin{lemma} \label{2.3.8}
Let $\mathfrak{a}$ be an ideal of $R$ and $M$ an $R$-module. Then for any $j\geq 0$ we have:
\[
    H_{i}^{\mathfrak{a}}\left(H_{j}^{\mathfrak{a}}(M)\right)\cong
\begin{dcases}
    H_{j}^{\mathfrak{a}}(M) & \text{if } i= 0\\
    0              & \text{if } i \geq 1.
\end{dcases}
\]
\end{lemma}

\begin{prf}
Let
$$F= \cdots \rightarrow F_{2} \xrightarrow {\partial^{F}_{2}} F_{1} \xrightarrow {\partial^{F}_{1}} F_{0} \rightarrow 0$$
be a free resolution of $M$. Then by definition
$$H^{\mathfrak{a}}_{j}(M)= \ker \widehat{\partial^{F}_{j}}/ \im  \widehat{\partial^{F}_{j+1}},$$
for any $j\geq 0$. By Lemma \ref{2.3.6} (iii), both $\im \widehat{\partial^{F}_{j+1}}$ and $\ker\widehat{\partial^{F}_{j}}$ are $\mathfrak{a}$-adically complete. Therefore, invoking \cite[Lemmas 5.1 (i), and 5.2 (i)]{Si1}, there are natural isomorphisms $$\alpha: H_{0}^{\mathfrak{a}}\left(\im \widehat{\partial^{F}_{j+1}}\right) \rightarrow \im \widehat{\partial^{F}_{j+1}},$$and$$\beta: H_{0}^{\mathfrak{a}}\left(\ker \widehat{\partial^{F}_{j}}\right) \rightarrow \ker \widehat{\partial^{F}_{j}}.$$ Now the short exact sequence
\begin{equation} \label{eq:2.3.8.1}
0 \rightarrow  \im  \widehat{\partial^{F}_{j+1}} \xrightarrow {\iota} \ker \widehat{\partial^{F}_{j}} \xrightarrow {\pi}
H^{\mathfrak{a}}_{j}(M) \rightarrow 0,
\end{equation}
yields the following commutative diagram with exact rows
\[\begin{tikzpicture}[every node/.style={midway},]
  \matrix[column sep={3em}, row sep={2.5em}]
  {\node(1) {}; & \node(2) {$H_{0}^{\mathfrak{a}}\left(\im  \widehat{\partial^{F}_{j+1}}\right)$}; & \node(3) {$H_{0}^{\mathfrak{a}}
  \left(\ker  \widehat{\partial^{F}_{j}}\right)$}; & \node(4) {$H_{0}^{\mathfrak{a}}\left(H^{\mathfrak{a}}_{j}(M)\right)$}; &
  \node(5) {$0$}; \\
  \node(6) {$0$}; & \node(7) {$\im  \widehat{\partial^{F}_{j+1}}$}; & \node(8) {$\ker \widehat{\partial^{F}_{j}}$}; & \node(9)
  {$H^{\mathfrak{a}}_{j}(M)$}; & \node(10) {$0$}; \\};
  \draw[decoration={markings,mark=at position 1 with {\arrow[scale=1.5]{>}}},postaction={decorate},shorten >=0.5pt] (2) -- (3)
  node[anchor=south] {$H_{0}^{\mathfrak{a}}(\iota)$};
  \draw[decoration={markings,mark=at position 1 with {\arrow[scale=1.5]{>}}},postaction={decorate},shorten >=0.5pt] (3) -- (4)
  node[anchor=south] {$H_{0}^{\mathfrak{a}}(\pi)$};
  \draw[decoration={markings,mark=at position 1 with {\arrow[scale=1.5]{>}}},postaction={decorate},shorten >=0.5pt] (4) -- (5)
  node[anchor=south] {};
  \draw[decoration={markings,mark=at position 1 with {\arrow[scale=1.5]{>}}},postaction={decorate},shorten >=0.5pt] (6) -- (7)
  node[anchor=south] {};
  \draw[decoration={markings,mark=at position 1 with {\arrow[scale=1.5]{>}}},postaction={decorate},shorten >=0.5pt] (7) -- (8)
  node[anchor=south] {$\iota$};
  \draw[decoration={markings,mark=at position 1 with {\arrow[scale=1.5]{>}}},postaction={decorate},shorten >=0.5pt] (8) -- (9)
  node[anchor=south] {$\pi$};
  \draw[decoration={markings,mark=at position 1 with {\arrow[scale=1.5]{>}}},postaction={decorate},shorten >=0.5pt] (9) -- (10)
  node[anchor=south] {};
  \draw[decoration={markings,mark=at position 1 with {\arrow[scale=1.5]{>}}},postaction={decorate},shorten >=0.5pt] (2) -- (7)
  node[anchor=west] {$\cong$} node[anchor=east] {$\alpha$};
  \draw[decoration={markings,mark=at position 1 with {\arrow[scale=1.5]{>}}},postaction={decorate},shorten >=0.5pt] (3) -- (8)
  node[anchor=west] {$\cong$} node[anchor=east] {$\beta$};
  \draw[dashed,decoration={markings,mark=at position 1 with {\arrow[scale=1.5]{>}}},postaction={decorate},shorten >=0.5pt] (4) --
  (9) node[anchor=west] {$\gamma$} node[anchor=east] {$\exists$};
\end{tikzpicture}\]
from which we deduce that $H_{0}^{\mathfrak{a}}(\iota)$ is injective and $\gamma$ is an isomorphism.
On the other hand, the short exact sequence \eqref{eq:2.3.8.1} yields an exact sequence
$$0= H_{i+1}^{\mathfrak{a}}\left(\ker \widehat{\partial^{F}_{j}}\right) \rightarrow H_{i+1}^{\mathfrak{a}}
\left(H^{\mathfrak{a}}_{j}(M)\right) \rightarrow H_{i}^{\mathfrak{a}}\left(\im  \widehat{\partial^{F}_{j+1}}\right) =0,$$
for every $i \geq 1$, where the vanishing follows from \cite[Lemma 5.2 (i)]{Si1}.
It follows that $H_{i+1}^{\mathfrak{a}}\left(H^{\mathfrak{a}}_{j}(M)\right)=0$ for every $i \geq 1$. It further yields the exact
sequence
$$0= H_{1}^{\mathfrak{a}}\left(\ker \widehat{\partial^{F}_{j}}\right) \rightarrow H_{1}^{\mathfrak{a}}\left(H^{\mathfrak{a}}_{j}(M)
\right) \rightarrow H_{0}^{\mathfrak{a}}\left(\im  \widehat{\partial^{F}_{j+1}}\right) \xrightarrow {H_{0}^{\mathfrak{a}}(\iota)}
H_{0}^{\mathfrak{a}}\left(\ker \widehat{\partial^{F}_{j}}\right).$$
As $H_{0}^{\mathfrak{a}}(\iota)$ is injective, we conclude that $H_{1}^{\mathfrak{a}}\left(H^{\mathfrak{a}}_{j}(M)\right)=0$.
\end{prf}

\begin{lemma} \label{2.3.9}
Let $\mathfrak{a}$ be an ideal of $R$, $M$ an $R$-module, and $s\geq 0$ an integer. Suppose that $R$ is $\mathfrak{a}$-adically complete, and
$(R/\mathfrak{a})\otimes_{R}H^{\mathfrak{a}}_{s}(M)$ is a finitely generated $R$-module. Then $H^{\mathfrak{a}}_{s}(M)$ is a finitely
generated $R$-module.
\end{lemma}

\begin{prf}
Let $L=H^{\mathfrak{a}}_{s}(M)$. As $L/\mathfrak{a}L$ is finitely generated, Lemma \ref{2.3.4} implies that $\widehat{L}^{\mathfrak{a}}$
is a finitely generated $R$-module. Lemma \ref{2.3.6} yields that $L$ is $\mathfrak{a}$-adically quasi-complete, so Lemma \ref{2.3.7} implies
that $H_{0}^{\mathfrak{a}}(\ker \theta_{L})=0$, where $\theta_{L}:L \rightarrow \widehat{L}^{\mathfrak{a}}$
is the completion map. Hence from the short exact sequence
$$0 \rightarrow \ker \theta_{L} \rightarrow L \xrightarrow {\theta_{L}} \widehat{L}^{\mathfrak{a}} \rightarrow 0,$$
we get the exact sequence
$$0=H_{0}^{\mathfrak{a}}(\ker \theta_{L}) \rightarrow H_{0}^{\mathfrak{a}}(L) \rightarrow H_{0}^{\mathfrak{a}}\left(\widehat{L}^
{\mathfrak{a}}\right) \rightarrow 0,$$
implying that $H_{0}^{\mathfrak{a}}(L)\cong H_{0}^{\mathfrak{a}}\left(\widehat{L}^{\mathfrak{a}}\right)$. By \cite[Lemma 5.2 (i)]{Si1},
we have $H_{0}^{\mathfrak{a}}\left(\widehat{L}^{\mathfrak{a}}\right)\cong \widehat{L}^{\mathfrak{a}}$. On the other hand, Lemma \ref{2.3.8}
implies that $H_{0}^{\mathfrak{a}}(L)\cong L$. Hence $L\cong \widehat{L}^{\mathfrak{a}}$ is a finitely generated $R$-module.
\end{prf}

Our final goal in this section is to show that given an ideal $\mathfrak{a}$ of $R$ with $\cd(\frak{a},R)\leq 1$, the subcategory $\mathcal{M}(R,\mathfrak{a})_{cof}$ of $\mathcal{M}(R)$ is abelian. As mentioned earlier, this fact is proved in \cite[Theorem 2.4]{PAB}, under the extra assumption that $R$ is local. The tool deployed here is the local homology functor.

\begin{theorem} \label{2.3.10}
Let $\mathfrak{a}$ be an ideal of $R$. Then the following assertions hold:
\begin{enumerate}
\item[(i)] An $R$-module $M$ with $\Supp_{R}(M)\subseteq V(\mathfrak{a})$ is $\mathfrak{a}$-cofinite if and only if $H^{\mathfrak{a}}_{i}(M)$ is a finitely generated $\widehat{R}^{\mathfrak{a}}$-module for every $0 \leq i \leq \cd(\mathfrak{a},R)$.
\item[(ii)] If $\cd(\mathfrak{a},R)\leq 1$, then $\mathcal{M}(R,\mathfrak{a})_{cof}$ is an abelian subcategory of $\mathcal{M}(R)$.
\end{enumerate}
\end{theorem}

\begin{prf}
(i): By \cite[Theorem 2.5 and Corollary 3.2]{GM}, $H^{\mathfrak{a}}_{i}(M)=0$ for every $i> \cd(\mathfrak{a},R)$. Therefore, the assertion follows from Corollary \ref{2.2.17}.

(ii): Let $M$ and $N$ be two $\mathfrak{a}$-cofinite $R$-modules and $f:M\rightarrow N$ an $R$-homomorphism. The short exact sequence
\begin{equation} \label{eq:2.3.10.1}
0 \rightarrow \ker f \rightarrow M \rightarrow \im f \rightarrow 0,
\end{equation}
gives the exact sequence
$$H^{\mathfrak{a}}_{0}(M) \rightarrow H^{\mathfrak{a}}_{0}\left(\im f \right) \rightarrow 0,$$
which in turn implies that $H^{\mathfrak{a}}_{0}\left(\im f \right)$ is finitely generated $\widehat{R}^{\mathfrak{a}}$-module
since $H^{\mathfrak{a}}_{0}(M)$ is so. The short exact sequence
\begin{equation} \label{eq:2.3.10.2}
0 \rightarrow \im f \rightarrow N \rightarrow \coker f \rightarrow 0,
\end{equation}
gives the exact sequence
\begin{equation} \label{eq:2.3.10.3}
H^{\mathfrak{a}}_{1}(N) \rightarrow H^{\mathfrak{a}}_{1}\left(\coker f \right) \rightarrow H^{\mathfrak{a}}_{0}\left(\im f \right)
\rightarrow H^{\mathfrak{a}}_{0}(N) \rightarrow H^{\mathfrak{a}}_{0}\left(\coker f \right) \rightarrow 0.
\end{equation}
As $H^{\mathfrak{a}}_{0}(N)$, $H^{\mathfrak{a}}_{0}\left(\im f \right)$, and $H^{\mathfrak{a}}_{1}(N)$ are finitely generated $\widehat{R}^{\mathfrak{a}}$-modules, the exact sequence \eqref{eq:2.3.10.3} shows that $H^{\mathfrak{a}}_{0}\left(\coker f \right)$ and $H^{\mathfrak{a}}_{1}\left(\coker f \right)$ are finitely generated  $\widehat{R}^{\mathfrak{a}}$-modules, and thus $\coker f$ is $\mathfrak{a}$-cofinite by (i). From the short exact sequence \eqref{eq:2.3.10.2}, we conclude that $\im f$ is $\mathfrak{a}$-cofinite, and from the short exact sequence \eqref{eq:2.3.10.1}, we infer that $\ker f$ is $\mathfrak{a}$-cofinite. It follows that $\mathcal{M}(R,\mathfrak{a})_{cof}$ is an abelian subcategory of $\mathcal{M}(R)$.
\end{prf}

\begin{corollary} \label{2.3.11}
Let $\mathfrak{a}$ be an ideal of $R$ such that either $\cd(\mathfrak{a},R) \leq 1$, or $\dim(R) \leq 2$, or $\dim(R/\mathfrak{a}) \leq 1$. Then $H^{i}_{\mathfrak{a}}(M)$ is $\mathfrak{a}$-cofinite for every finitely generated $R$-module $M$ and every $i \geq 0$, and $\mathcal{M}(R,\mathfrak{a})_{cof}$ is an abelian subcategory of $\mathcal{M}(R)$.
\end{corollary}

\begin{prf}
Follows from Theorems \ref{2.1.5} and \ref{2.3.10}.
\end{prf}

\section{Cofiniteness of Complexes}

In this section, we first present some background material on complexes which will be used efficaciously in the rest of the work. For more information, refer to \cite{AF}, \cite{Ha2}, \cite{Fo}, \cite{Li}, and \cite{Sp}. In what follows, $\mathcal{C}(R)$ denotes the category of $R$-complexes.

The theory of derived category is the ultimate formulation of homological algebra. The derived category $\mathcal{D}(R)$ \index{derived category} is defined as the localization of the homotopy category $\mathcal{K}(R)$ with respect to the multiplicative system of quasi-isomorphisms. Simply put, an object in $\mathcal{D}(R)$ is an $R$-complex $X$ displayed in the standard homological style
$$X= \cdots \rightarrow X_{i+1} \xrightarrow {\partial^{X}_{i+1}} X_{i} \xrightarrow {\partial^{X}_{i}} X_{i-1} \rightarrow \cdots,$$
and a morphism $\varphi:X\rightarrow Y$ in $\mathcal{D}(R)$ is given by the equivalence class of a pair of morphisms
$X \xleftarrow{g} U \xrightarrow{f} Y$ in $\mathcal{C}(R)$ with $g$ a quasi-isomorphism, under the equivalence relation that identifies two such pairs $X \xleftarrow{g} U \xrightarrow{f} Y$ and $X \xleftarrow{g^{\prime}} U^{\prime} \xrightarrow{f^{\prime}} Y$, whenever there is a diagram in $\mathcal{C}(R)$ as follows which commutes up to homotopy:
\[\begin{tikzpicture}[every node/.style={midway}]
  \matrix[column sep={3em}, row sep={3em}]
  {\node(1) {$$}; & \node(2) {$U$}; & \node(3) {$$};\\
  \node(4) {$X$}; & \node(5) {$V$}; & \node(6) {$Y$};\\
  \node(7) {$$}; & \node(8) {$U^{\prime}$}; & \node(9) {$$};\\};
  \draw[decoration={markings,mark=at position 1 with {\arrow[scale=1.5]{>}}},postaction={decorate},shorten >=0.5pt] (2) -- (4) node[anchor=south] {$g$} node[anchor=west] {$\simeq$};
  \draw[decoration={markings,mark=at position 1 with {\arrow[scale=1.5]{>}}},postaction={decorate},shorten >=0.5pt] (2) -- (6) node[anchor=south] {$f$};
  \draw[decoration={markings,mark=at position 1 with {\arrow[scale=1.5]{>}}},postaction={decorate},shorten >=0.5pt] (8) -- (4) node[anchor=north] {$g^{\prime}$} node[anchor=west] {$\simeq$};
  \draw[decoration={markings,mark=at position 1 with {\arrow[scale=1.5]{>}}},postaction={decorate},shorten >=0.5pt] (8) -- (6) node[anchor=north] {$f^{\prime}$};
  \draw[decoration={markings,mark=at position 1 with {\arrow[scale=1.5]{>}}},postaction={decorate},shorten >=0.5pt] (5) -- (2) node[anchor=east] {$$};
  \draw[decoration={markings,mark=at position 1 with {\arrow[scale=1.5]{>}}},postaction={decorate},shorten >=0.5pt] (5) -- (4) node[anchor=south] {$\simeq$};
  \draw[decoration={markings,mark=at position 1 with {\arrow[scale=1.5]{>}}},postaction={decorate},shorten >=0.5pt] (5) -- (6) node[anchor=east] {$$};
  \draw[decoration={markings,mark=at position 1 with {\arrow[scale=1.5]{>}}},postaction={decorate},shorten >=0.5pt] (5) -- (8) node[anchor=east] {$$};
\end{tikzpicture}\]
The isomorphisms in $\mathcal{D}(R)$ are marked by the symbol $\simeq$.

The derived category $\mathcal{D}(R)$ is triangulated. A distinguished triangle \index{distinguished triangle} in $\mathcal{D}(R)$ is a triangle that is isomorphic to a triangle of the form
$$X \xrightarrow {\mathfrak{L}(f)} Y \xrightarrow{\mathfrak{L}(\varepsilon)} \Cone(f) \xrightarrow{\mathfrak{L}(\varpi)} \Sigma X,$$
for some morphism $f:X \rightarrow Y$ in $\mathcal{C}(R)$ with the mapping cone sequence
$$0 \rightarrow Y \xrightarrow{\varepsilon} \Cone(f) \xrightarrow{\varpi} \Sigma X \rightarrow 0,$$
in which $\mathfrak{L}:\mathcal{C}(R) \rightarrow \mathcal{D}(R)$ is the canonical functor that is defined as $\mathfrak{L}(X)=X$ for every $R$-complex $X$, and $\mathfrak{L}(f)=\varphi$ where $\varphi$ is represented by the morphisms $X \xleftarrow{1^{X}} X \xrightarrow{f} Y$ in $\mathcal{C}(R)$. We note that if $f$ is a quasi-isomorphism in $\mathcal{C}(R)$, then $\mathfrak{L}(f)$ is an isomorphism in $\mathcal{D}(R)$. We sometimes use the shorthand notation
$$X \rightarrow Y \rightarrow Z \rightarrow$$
for a distinguished triangle.

We let $\mathcal{D}_{\sqsubset}(R)$ (res. $\mathcal{D}_{\sqsupset}(R)$) denote the full subcategory of $\mathcal{D}(R)$ consisting of $R$-complexes $X$ with $H_{i}(X)=0$ for $i \gg 0$ (res. $i \ll 0$), and  $D_{\square}(R):=\mathcal{D}_{\sqsubset}(R)\cap \mathcal{D}_{\sqsupset}(R)$. We further let $\mathcal{D}^{f}(R)$ denote the full subcategory of $\mathcal{D}(R)$ consisting of $R$-complexes $X$ with finitely generated homology modules. We also feel free to use any combination of the subscripts and the superscript as in $\mathcal{D}^{f}_{\square}(R)$, with the obvious meaning of the intersection of the two subcategories involved. Given an $R$-complex $X$, the standard notions
$$\sup(X) = \sup \left\{i \in \mathbb{Z} \suchthat H_{i}(X) \neq 0 \right\}$$
and
$$\inf(X) = \inf \left\{i \in \mathbb{Z} \suchthat H_{i}(X) \neq 0 \right\}$$
are frequently used, with the convention that $\sup(\emptyset) =- \infty$ and $\inf(\emptyset) = +\infty$.

An $R$-complex $P$ of projective modules is said to be semi-projective if the functor $\Hom_{R}(P,-)$ preserves quasi-isomorphisms. By a semi-projective resolution of an $R$-complex $X$, we mean a quasi-isomorphism $P\xrightarrow {\simeq} X$ in which $P$ is a semi-projective $R$-complex. Dually, an $R$-complex
$I$ of injective modules is said to be semi-injective if the functor $\Hom_{R}(-,I)$ preserves quasi-isomorphisms. By a semi-injective resolution of an $R$-complex $X$, we mean a quasi-isomorphism $X\xrightarrow {\simeq} I$ in which $I$ is a semi-injective $R$-complex. Semi-projective and semi-injective resolutions exist for any $R$-complex. Moreover, any right-bounded $R$-complex of projective modules is semi-projective, and any left-bounded $R$-complex of injective modules is semi-injective.

Let $X$ and $Y$ be two $R$-complexes. Then each of the functors $\Hom_{R}(X,-)$ and $\Hom_{R}(-,Y)$ on $\mathcal{C}(R)$ enjoys a right total derived functor on $\mathcal{D}(R)$, together with a balance property, in the sense that ${\bf R}\Hom_{R}(X,Y)$ can be computed by
$${\bf R}\Hom_{R}(X,Y)\simeq \Hom_{R}(P,Y) \simeq \Hom_{R}(X,I),$$
where $P\xrightarrow {\simeq} X$ is any semi-projective resolution of $X$, and $Y\xrightarrow {\simeq} I$ is any semi-injective resolution of $Y$. In addition, these functors turn out to be triangulated, in the sense that they preserve shifts and distinguished triangles. Moreover, we let $$\Ext^{i}_{R}(X,Y):=H_{-i}\left({\bf R}\Hom_{R}(X,Y)\right)$$
for every $i \in \mathbb{Z}$.

Likewise, each of the functors $X\otimes_{R}-$ and $-\otimes_{R}Y$ on $\mathcal{C}(R)$ enjoys a left total derived functor on $\mathcal{D}(R)$, together with a balance property, in the sense that $X\otimes_{R}^{\bf L}Y$ can be computed by
$$X\otimes_{R}^{\bf L}Y \simeq P\otimes_{R}Y \simeq X\otimes_{R}Q,$$
where $P\xrightarrow {\simeq} X$ is any semi-projective resolution of $X$, and $Q\xrightarrow {\simeq} Y$ is any semi-projective resolution of $Y$. Besides, these functors turn out to be triangulated. Moreover, we let $$\Tor^{R}_{i}(X,Y):=H_{i}\left(X\otimes_{R}^{\bf L}Y\right)$$
for every $i \in \mathbb{Z}$.

We frequently use the following case of the Hom Evaluation Isomorphism. The natural morphism
$$X\otimes_{R}^{\bf L}{\bf R}\Hom_{R}(Y,Z)\rightarrow {\bf R}\Hom_{R}\left({\bf R}\Hom_{R}(X,Y),Z\right)$$
is an isomorphism in $\mathcal{D}(R)$ provided that $X \in \mathcal{D}^{f}_{\sqsupset}(R)$, $Y\in \mathcal{D}_{\sqsubset}(R)$, and $\id_{R}(Z) < \infty$.
We also use the following case of the Tensor Evaluation Isomorphism. The natural morphism
$${\bf R}\Hom_{R}(X,Y)\otimes_{R}^{\bf L}Z \rightarrow {\bf R}\Hom_{R}\left(X,Y\otimes_{R}^{\bf L}Z\right)$$
is an isomorphism in $\mathcal{D}(R)$ provided that $X \in \mathcal{D}^{f}_{\sqsupset}(R)$, $Y\in \mathcal{D}_{\sqsubset}(R)$, and $\fd_{R}(Z) < \infty$.

Now let $\mathfrak{a}$ be an ideal of $R$. We let $\Gamma_{\mathfrak{a}}(M):=\left\{x\in M \suchthat \mathfrak{a}^{t}x=0 \text{ for some } t\geq 1 \right\}$ for an $R$-module $M$, and $\Gamma_{\mathfrak{a}}(f):=f|_{\Gamma_{\mathfrak{a}}(M)}$ for an $R$-homomorphism $f:M\rightarrow N$. This provides us with the so-called torsion functor $\Gamma_{\mathfrak{a}}(-)$ on the category of $R$-modules, which extends naturally to a functor on the category of $R$-complexes. The extended functor preserves homotopy equivalences, and thus enjoys a right derived functor ${\bf R}\Gamma_{\mathfrak{a}}(-):\mathcal{D}(R)\rightarrow \mathcal{D}(R)$, that can be computed by
${\bf R}\Gamma_{\mathfrak{a}}(X)\simeq \Gamma_{\mathfrak{a}}(I)$, where $X \xrightarrow {\simeq} I$ is any semi-injective resolution of $X$. Besides, we define the $i$th local cohomology module of $X$ to be
$$H^{i}_{\mathfrak{a}}(X):= H_{-i}\left({\bf R}\Gamma_{\mathfrak{a}}(X)\right)$$
for every $i\in\mathbb{Z}$. The functor ${\bf R}\Gamma_{\mathfrak{a}}(-)$ turns out to be triangulated.

Similarly, we let $\Lambda^{\mathfrak{a}}(M):=\widehat{M}^{\mathfrak{a}}= \varprojlim (M/\mathfrak{a}^{n}M)$ for an $R$-module $M$, and $\Lambda^{\mathfrak{a}}(f):=\widehat{f}$ for an $R$-homomorphism $f:M\rightarrow N$. This provides us with the so-called completion functor $\Lambda^{\mathfrak{a}}(-)$ on the category of $R$-modules, which extends naturally to a functor on the category of $R$-complexes. The extended functor preserves homotopy equivalences, and thus enjoys a left derived functor ${\bf L}\Lambda^{\mathfrak{a}}(-):\mathcal{D}(R)\rightarrow \mathcal{D}(R)$, which can be equally well considered as ${\bf L}\Lambda^{\mathfrak{a}}(-):\mathcal{D}(R)\rightarrow \mathcal{D}\left(\widehat{R}^{\mathfrak{a}}\right)$. This functor can be computed by ${\bf L}\Lambda^{\mathfrak{a}}(X)\simeq \Lambda^{\mathfrak{a}}(P)$, where $P \xrightarrow {\simeq} X$ is any semi-projective resolution of $X$. Moreover, we define the $i$th local homology module of $X$ to be
$$H^{\mathfrak{a}}_{i}(X):= H_{i}\left({\bf L}\Lambda^{\mathfrak{a}}(X)\right)$$
for every $i\in\mathbb{Z}$. The functor ${\bf L}\Lambda^{\mathfrak{a}}(-)$ turns out to be triangulated.

The derived torsion and completion functors are intimately linked with the \v{C}ech complex. Let $\mathfrak{a}=(a_{1},...,a_{n})$ and $\check{C}(\underline{a})$ denote the \v{C}ech complex on the sequence of elements $\underline{a}=a_{1},...,a_{n}\in R$. Then one has the following natural isomorphisms in $\mathcal{D}(R)$:
$${\bf R}\Gamma_{\mathfrak{a}}(X)\simeq \check{C}(\underline{a})\otimes_{R}^{\bf L} X,$$
and
$${\bf L}\Lambda^{\mathfrak{a}}(X)\simeq {\bf R}\Hom_{R}\left(\check{C}(\underline{a}),X\right).$$

We are now ready to investigate the notion of $\mathfrak{a}$-cofiniteness for $R$-complexes.

\begin{definition} \label{2.4.1}
Let $\mathfrak{a}$ be an ideal of $R$. An $R$-complex $X$ is said to be $\mathfrak{a}$-\textit{cofinite} \index{cofinite complex} if $\Supp_{R}(X)\subseteq \V(\mathfrak{a})$, and ${\bf R}\Hom_{R}(R/\mathfrak{a},X) \in \mathcal{D}^{f}(R)$.
\end{definition}

Some special cases of the following lemma is more or less proved in \cite[Propositions 7.1, 7.2, and 7.4]{WW} and \cite[Theorem 3.10 and Proposition 3.13]{PSY2}. However, we include it here with a different and shorter proof due to its pivotal role in the theory of cofiniteness.

We recall that an $R$-complex $D\in \mathcal{D}^{f}_{\square}(R)$ is said to be a dualizing complex \index{dualizing complex} for $R$ if the homothety morphism $R \rightarrow {\bf R}\Hom_{R}(D,D)$ is an isomorphism in $\mathcal{D}(R)$, and $\id_{R}(D)<\infty$.

\begin{lemma} \label{2.4.2}
Let $\mathfrak{a}$ be an ideal of $R$ and $X\in \mathcal{D}_{\square}(R)$. Then the following assertions are equivalent:
\begin{enumerate}
\item[(i)] ${\bf R}\Hom_{R}(R/\mathfrak{a},X) \in \mathcal{D}^{f}(R)$.
\item[(ii)] ${\bf R}\Hom_{R}(Y,X) \in \mathcal{D}^{f}(R)$ for every $Y\in \mathcal{D}^{f}_{\square}(R)$ with $\Supp_{R}(Y)\subseteq V(\mathfrak{a})$.
\item[(iii)] $(R/\mathfrak{a})\otimes_{R}^{\bf L}X \in \mathcal{D}^{f}(R)$.
\item[(iv)] $Y\otimes_{R}^{\bf L}X \in \mathcal{D}^{f}(R)$ for every $Y\in \mathcal{D}^{f}_{\square}(R)$ with $\Supp_{R}(Y)\subseteq V(\mathfrak{a})$.
\item[(v)] ${\bf L}\Lambda^{\mathfrak{a}}(X) \in \mathcal{D}^{f}_{\square}\left(\widehat{R}^{\mathfrak{a}}\right)$.
\item[(vi)] ${\bf R}\Gamma_{\mathfrak{a}}(X)\simeq {\bf R}\Gamma_{\mathfrak{a}}(Z)$ for some $Z \in \mathcal{D}^{f}_{\square}\left(\widehat{R}^{\mathfrak{a}}\right)$.
\item[(vii)] ${\bf R}\Gamma_{\mathfrak{a}}(X)\simeq {\bf R}\Hom_{\widehat{R}^{\mathfrak{a}}}\left(Y,{\bf R}\Gamma_{\mathfrak{a}}(D)\right)$ for some
$Y \in \mathcal{D}^{f}_{\square}\left(\widehat{R}^{\mathfrak{a}}\right)$, provided that $\widehat{R}^{\mathfrak{a}}$ enjoys a dualizing complex $D$.
\end{enumerate}
\end{lemma}

\begin{prf}
For the equivalences (i) $\Leftrightarrow$ (ii) $\Leftrightarrow$ (iii) $\Leftrightarrow$ (iv) $\Leftrightarrow$ (v), see \cite[Propositions 7.1, 7.2, and Theorem 7.4]{WW}.

(i) $\Rightarrow$ (vi): By \cite[Propositions 7.4]{WW}, $Z:={\bf L}\Lambda^{\mathfrak{a}}(X) \in \mathcal{D}^{f}_{\square}\left(\widehat{R}^{\mathfrak{a}}\right)$. Then by \cite[Corollary after (0.3)$^{\ast}$]{AJL}, we have
$${\bf R}\Gamma_{\mathfrak{a}}(Z)\simeq {\bf R}\Gamma_{\mathfrak{a}}\left({\bf L}\Lambda^{\mathfrak{a}}(X)\right)\simeq {\bf R}\Gamma_{\mathfrak{a}}(X).$$

(vi) $\Rightarrow$ (vii): Set $Y:={\bf R}\Hom_{\widehat{R}^{\mathfrak{a}}}(Z,D)$. If $\id_{\widehat{R}^{\mathfrak{a}}}(D) = n$, then there is a semi-injective resolution $D \xrightarrow {\simeq} I$ of $D$ such that $I_{i}=0$ for every $i > \sup D$ or $i < -n$. In particular, $I$ is bounded. On the other hand, $Z \in \mathcal{D}^{f}_{\square}\left(\widehat{R}^{\mathfrak{a}}\right)$, so there is a bounded $\widehat{R}^{\mathfrak{a}}$-complex $Z^{\prime}$ such that $Z\simeq Z^{\prime}$. Therefore,
$$Y= {\bf R}\Hom_{\widehat{R}^{\mathfrak{a}}}(Z,D) \simeq {\bf R}\Hom_{\widehat{R}^{\mathfrak{a}}}(Z^{\prime},D) \simeq \Hom_{\widehat{R}^{\mathfrak{a}}}(Z^{\prime},I).$$
But it is obvious that $\Hom_{\widehat{R}^{\mathfrak{a}}}(Z^{\prime},I)$ is bounded, so $Y\in \mathcal{D}^{f}_{\square}\left(\widehat{R}^{\mathfrak{a}}\right)$.
Now, let $\check{C}(\underline{a})$ denote the \v{C}ech complex on a sequence of elements $\underline{a}=a_{1},...,a_{n}\in R$ that generate $\mathfrak{a}$.
For any $R$-complex $W$, \cite[Proposition 3.1.2]{Li} yields that ${\bf R}\Gamma_{\mathfrak{a}}(W)\simeq \check{C}(\underline{a})\otimes_{R}^{\bf L} W$. We thus have
\begin{equation*}
\begin{split}
{\bf R}\Gamma_{\mathfrak{a}}(X) & \simeq {\bf R}\Gamma_{\mathfrak{a}}(Z)\\
& \simeq {\bf R}\Gamma_{\mathfrak{a}}\left({\bf R}\Hom_{\widehat{R}^{\mathfrak{a}}}\left({\bf R}\Hom_{\widehat{R}^{\mathfrak{a}}}(Z,D),D\right)\right)\\
& \simeq {\bf R}\Gamma_{\mathfrak{a}}\left({\bf R}\Hom_{\widehat{R}^{\mathfrak{a}}}\left(Y,D\right)\right)\\
& \simeq \check{C}(\underline{a})\otimes_{R}^{\bf L} {\bf R}\Hom_{\widehat{R}^{\mathfrak{a}}}\left(Y,D\right)\\
& \simeq {\bf R}\Hom_{\widehat{R}^{\mathfrak{a}}}\left(Y, \check{C}(\underline{a})\otimes_{R}^{\bf L} D\right)\\
& \simeq {\bf R}\Hom_{\widehat{R}^{\mathfrak{a}}}\left(Y,{\bf R}\Gamma_{\mathfrak{a}}(D)\right).\\
\end{split}
\end{equation*}
The second isomorphism is due to the fact that $D$ is a dualizing complex for $\widehat{R}^{\mathfrak{a}}$, and the fifth isomorphism follows from the application of the Tensor Evaluation Isomorphism. The other isomorphisms are straightforward.

(vii) $\Rightarrow$ (v): Similar to the argument of the implication (vi) $\Rightarrow$ (vii), we conclude that ${\bf R}\Hom_{\widehat{R}^{\mathfrak{a}}}(Y,D)\in \mathcal{D}^{f}_{\square}\left(\widehat{R}^{\mathfrak{a}}\right)$. We further have
\begin{equation*}
\begin{split}
{\bf L}\Lambda^{\mathfrak{a}}(X) & \simeq {\bf L}\Lambda^{\mathfrak{a}}\left({\bf R}\Gamma_{\mathfrak{a}}(X)\right)\\
& \simeq {\bf L}\Lambda^{\mathfrak{a}}\left({\bf R}\Hom_{\widehat{R}^{\mathfrak{a}}}\left(Y,{\bf R}\Gamma_{\mathfrak{a}}(D)\right)\right)\\
& \simeq {\bf L}\Lambda^{\mathfrak{a}}\left({\bf R}\Gamma_{\mathfrak{a}}\left({\bf R}\Hom_{\widehat{R}^{\mathfrak{a}}}(Y,D)\right)\right)\\
& \simeq {\bf L}\Lambda^{\mathfrak{a}}\left({\bf R}\Hom_{\widehat{R}^{\mathfrak{a}}}(Y,D)\right)\\
& \simeq {\bf R}\Hom_{\widehat{R}^{\mathfrak{a}}}(Y,D)\in \mathcal{D}^{f}_{\square}\left(\widehat{R}^{\mathfrak{a}}\right).\\
\end{split}
\end{equation*}
The first isomorphism uses \cite[Corollary after (0.3)$^{\ast}$]{AJL}, the third isomorphism follows from the application of the Tensor Evaluation Isomorphism
just as in the previous paragraph, the fourth isomorphism uses \cite[Corollary after (0.3)$^{\ast}$]{AJL}, and the fifth isomorphism follows from \cite[Theorem 1.21]{PSY2}, noting that as ${\bf R}\Hom_{\widehat{R}^{\mathfrak{a}}}(Y,D)\in \mathcal{D}^{f}_{\square}\left(\widehat{R}^{\mathfrak{a}}\right)$, its homology modules are $\mathfrak{a}$-adically complete $\widehat{R}^{\mathfrak{a}}$-modules. Now the results follows from \cite[Propositions 7.4]{WW}.
\end{prf}

\begin{corollary} \label{2.4.3}
Let $\mathfrak{a}$ be an ideal of $R$ for which $R$ is $\mathfrak{a}$-adically complete and $X\in \mathcal{D}_{\square}(R)$. Then the following assertions are equivalent:
\begin{enumerate}
\item[(i)] $X$ is $\mathfrak{a}$-cofinite.
\item[(ii)] $X\simeq {\bf R}\Gamma_{\mathfrak{a}}(Z)$ for some $Z\in \mathcal{D}^{f}_{\square}(R)$.
\item[(iii)] $X\simeq {\bf R}\Hom_{R}\left(Y,{\bf R}\Gamma_{\mathfrak{a}}(D)\right)$ for some $Y\in \mathcal{D}^{f}_{\square}(R)$, provided that $R$ enjoys a dualizing complex $D$.
\end{enumerate}
\end{corollary}

\begin{prf}
For any two $R$-complexes $V\in \mathcal{D}^{f}_{\square}(R)$ and $W\in \mathcal{D}_{\square}(R)$, one may easily see that $$\Supp_{R}
\left({\bf R}\Hom_{R}\left(V,{\bf R}\Gamma_{\mathfrak{a}}(W)\right)\right) \subseteq V(\mathfrak{a}).$$ Also, for any $U\in \mathcal{D}(R)$,
\cite[Corollary 3.2.1]{Li} yields that $\Supp_{R}(U)\subseteq V(\mathfrak{a})$ if and only if ${\bf R}\Gamma_{\mathfrak{a}}(U)\simeq U$.
Hence the assertions follow from Lemma \ref{2.4.2}.
\end{prf}

We next collect some basic properties of $\mathfrak{a}$-cofinite $R$-complexes in the following result.

\begin{lemma} \label{2.4.4}
Let $\mathfrak{a}$ be an ideal of $R$ and $X\in \mathcal{D}(R)$. Then the following assertions hold:
\begin{enumerate}
\item[(i)] If $X\in \mathcal{D}_{\sqsubset}(R)$ and $H_{i}(X)$ is $\mathfrak{a}$-cofinite for every $i\in \mathbb{Z}$, then $X$ is $\mathfrak{a}$-cofinite.
\item[(ii)] If $X\in \mathcal{D}^{f}_{\sqsubset}(R)$, then ${\bf R}\Gamma_{\mathfrak{a}}(X)$ is $\mathfrak{a}$-cofinite.
\item[(iii)] If $X\in \mathcal{D}_{\sqsubset}(R)$ and $X$ is $\mathfrak{a}$-cofinite, then the Bass number $\mu^{i}_{R}(\mathfrak{p},X)$ is finite for every $\mathfrak{p}\in \Spec (R)$ and every $i \in \mathbb{Z}$.
\item[(iv)] If $X\in \mathcal{D}_{\square}(R)$ and $X$ is $\mathfrak{a}$-cofinite, then the Betti number $\beta_{i}^{R}(\mathfrak{p},X)$ is finite for every $\mathfrak{p}\in \Spec (R)$ and every $i \in \mathbb{Z}$.
\end{enumerate}
\end{lemma}

\begin{prf}
(i): Since $H_{i}(X)$ is $\mathfrak{a}$-cofinite for every $i\in \mathbb{Z}$, we have
$$\Supp_{R}(X)=\bigcup_{i\in \mathbb{Z}}\Supp_{R}\left(H_{i}(X)\right) \subseteq V(\mathfrak{a}).$$
The spectral sequence
$$E_{p,q}^{2}=\Ext_{R}^{p}\left(R/\mathfrak{a},H_{-q}(X)\right) \underset{p}\Rightarrow \Ext_{R}^{p+q}(R/\mathfrak{a},X)$$
from the proof of \cite[Proposition 6.2]{Ha1}, together with the assumption that $E_{p,q}^{2}$ is finitely generated for every $p,q\in \mathbb{Z}$, conspire to imply that $\Ext_{R}^{p+q}(R/\mathfrak{a},X)$ is finitely generated, i.e. $X$ is $\mathfrak{a}$-cofinite.

(ii): It is clear that $\Supp_{R}\left({\bf R}\Gamma_{\mathfrak{a}}(X)\right)\subseteq V(\mathfrak{a})$. Since $X$ is homologically left-bounded, there exists a left-bounded semi-injective resolution $X \xrightarrow {\simeq} I$ of $X$.
As $I$ is an $R$-complex of injective modules, we see that $\Gamma_{\mathfrak{a}}(I)$ is a left-bounded $R$-complex of injective modules, and thus $\Gamma_{\mathfrak{a}}(I)$ is semi-injective. Hence one has
\begin{equation*}
\begin{split}
{\bf R}\Hom_{R}\left(R/\mathfrak{a},{\bf R}\Gamma_{\mathfrak{a}}(X)\right) & \simeq {\bf R}\Hom_{R}\left(R/\mathfrak{a},\Gamma_{\mathfrak{a}}(I)\right) \\
 & \simeq \Hom_{R}\left(R/\mathfrak{a},\Gamma_{\mathfrak{a}}(I)\right) \\
 & \simeq \Hom_{R}(R/\mathfrak{a},I) \\
 & \simeq {\bf R}\Hom_{R}(R/\mathfrak{a},X).
\end{split}
\end{equation*}
But $X\in \mathcal{D}^{f}_{\sqsubset}(R)$, so ${\bf R}\Hom_{R}(R/\mathfrak{a},X) \in \mathcal{D}^{f}_{\sqsubset}(R)$, and thus the assertion follows.

(iii): For every $i\in \mathbb{Z}$, we have
$$\mu^{i}_{R}(\mathfrak{p},X):= \rank_{R_{\mathfrak{p}}/\mathfrak{p}R_{\mathfrak{p}}}\left(\Ext_{R_{\mathfrak{p}}}^{i}\left(R_{\mathfrak{p}}/\mathfrak{p}R_{\mathfrak{p}},
X_{\mathfrak{p}}\right)\right).$$
If $\mathfrak{p}\not\in V(\mathfrak{a})$,
then $\mathfrak{p}\not\in \Supp_{R}(X)$, so $\mu^{i}_{R}(\mathfrak{p},X)=0$. If $\mathfrak{p}\in V(\mathfrak{a})$, then $\V(\mathfrak{p}) \subseteq V(\mathfrak{a})$, so by \cite[Proposition 7.1]{WW}, ${\bf R}\Hom_{R}(R/\mathfrak{p},X)\in \mathcal{D}^{f}(R)$, whence
$\mu^{i}_{R}(\mathfrak{p},X)< \infty$.

(iv): For every $i\in \mathbb{Z}$, we have
$$\beta_{i}^{R}(\mathfrak{p},X):= \rank_{R_{\mathfrak{p}}/\mathfrak{p}R_{\mathfrak{p}}}\left(\Tor_{i}^{R_{\mathfrak{p}}}\left(R_{\mathfrak{p}}/\mathfrak{p}R_{\mathfrak{p}},
X_{\mathfrak{p}}\right)\right).$$
If $\mathfrak{p}\not\in V(\mathfrak{a})$, then as in (iii), $\beta_{i}^{R}(\mathfrak{p},X)=0$. If $\mathfrak{p}\in
V(\mathfrak{a})$, then \cite[Proposition 7.4]{WW} implies that $(R/\mathfrak{p})\otimes_{R}^{\bf L} X\in \mathcal{D}^{f}(R)$, thereby
$\beta_{i}^{R}(\mathfrak{p},X)<\infty$.
\end{prf}

The next proposition deals with the change of rings principle for cofiniteness of complexes.

\begin{proposition} \label{2.4.5}
Let $\mathfrak{a}$ an ideal of $R$. Let $f:R\rightarrow S$ be a faithfully flat ring homomorphism. Then the following assertions hold:
\begin{enumerate}
\item[(i)] For each $R$-module $M$, $\Supp_{R}(M)\subseteq \V(\mathfrak{a})$ if and only if $\Supp_{S}\left(S\otimes_{R}M\right)\subseteq
\V\left(\mathfrak{a}S\right)$.
\item[(ii)] A complex $X\in\mathcal{D}_{\square}(R)$ is $\mathfrak{a}$-cofinite if and only if the complex $S\otimes_{R}X \in\mathcal{D}_{\square}
\left(S\right)$ is $\mathfrak{a}S$-cofinite.
\item[(iii)] If $\mathcal{M}\left(S,\mathfrak{a}S\right)_{cof}$ is an abelian category, then
$\mathcal{M}(R,\mathfrak{a})_{cof}$ is an abelian category as well.
\end{enumerate}
\end{proposition}

\begin{prf}
(i): First suppose that $\Supp_{R}(M)\subseteq \V(\mathfrak{a})$, and take $\mathfrak{q}\in \Supp_{S}\left(S\otimes_{R}M\right)$. Let
$\mathfrak{p}=\mathfrak{q}\cap R$. Since
$$M_{\mathfrak{p}}\otimes_{R_{\mathfrak{p}}} S_{\mathfrak{q}} \cong \left(M\otimes_{R}S\right)_{\mathfrak{q}}\neq 0,$$
we deduce that $M_{\mathfrak{p}} \neq 0$, i.e. $\mathfrak{p} \in \Supp_{R}(M)$, so
$\mathfrak{p}\supseteq \mathfrak{a}$. Now we have
$$\mathfrak{a}S \subseteq \mathfrak{p}S = (\mathfrak{q}\cap R)S\subseteq \mathfrak{q},$$
i.e. $\mathfrak{q} \in \V\left(\mathfrak{a}S\right)$.

Next, suppose that $\Supp_{S}\left(S\otimes_{R}M\right)\subseteq \V\left(\mathfrak{a}S\right)$, and take $\mathfrak{p}\in \Supp_{R}(M)$.
Since $f:R \rightarrow S$ is faithfully flat, it has the lying-over property, so there is a $\mathfrak{q} \in \Spec\left(S\right)$ such
that $\mathfrak{q}\cap R=\mathfrak{p}$. As the induced ring homomorphism $\tilde{f}:R_{\mathfrak{p}}\rightarrow S_{\mathfrak{q}}$ is faithfully
flat, we conclude that $$\left(M\otimes_{R}S\right)_{\mathfrak{q}}\cong M_{\mathfrak{p}}\otimes_{R_{\mathfrak{p}}}S_{\mathfrak{q}}
\neq 0,$$ i.e. $\mathfrak{q}\in \Supp_{S}\left(S\otimes_{R}M\right)$, so $\mathfrak{q}\supseteq \mathfrak{a}S$. It follows that
$$\mathfrak{a}\subseteq \mathfrak{a}S \cap R \subseteq \mathfrak{q}\cap R = \mathfrak{p},$$
i.e. $\mathfrak{p} \in \V(\mathfrak{a})$.

(ii): We note that $\Supp_{R}(X)= \Supp_{R}\left(\bigoplus_{i\in\mathbb{Z}}H_{i}(X)\right)$, and
\begin{equation*}
\begin{split}
\Supp_{R}\left(S\otimes_{R}X\right) & =\Supp_{R}\left(\bigoplus_{i\in\mathbb{Z}}H_{i}\left(S \otimes_{R}X\right)\right) \\
 & =\Supp_{R}\left(S\otimes_{R}\left(\bigoplus_{i\in\mathbb{Z}}H_{i}(X)\right)\right). \\
\end{split}
\end{equation*}
Therefore, by part (i), $\Supp_{R}(X)\subseteq \V(\mathfrak{a})$ if and only if $\Supp_{S}\left(S\otimes_{R}X\right)\subseteq \V\left(\mathfrak{a}S\right)$. On the other hand,
using the Tensor Evaluation and Adjointness Isomorphisms we have
\begin{equation*}
\begin{split}
{\bf R}\Hom_{R}(R/\mathfrak{a},X)\otimes_{R}^{\bf L}S & \simeq {\bf R}\Hom_{R}\left(R/\mathfrak{a},X \otimes_{R}^{\bf L}S\right) \\
 & \simeq {\bf R}\Hom_{R}\left(R/\mathfrak{a},{\bf R}\Hom_{S}\left(S,X \otimes_{R}^{\bf L}S\right)\right) \\
 & \simeq {\bf R}\Hom_{S}\left((R/\mathfrak{a})\otimes_{R}^{\bf L}S,X \otimes_{R}^{\bf L}S\right) \\
 & \simeq {\bf R}\Hom_{S}\left(S/\mathfrak{a}S,X \otimes_{R}^{\bf L}S\right), \\
\end{split}
\end{equation*}
and thus taking the homology, we are left with the isomorphism
$$\Ext^{i}_{R}(R/\mathfrak{a},X)\otimes_{R}S\cong \Ext^{i}_{S}\left(S/\mathfrak{a}S,X\otimes_{R} S\right)$$
for every $i \in \mathbb{Z}$.
As $f:R\rightarrow S$ is faithfully flat, we infer that $\Ext^{i}_{R}(R/\mathfrak{a},X)$ is a finitely generated
$R$-module if and only if $\Ext^{i}_{S}\left(S/\mathfrak{a}S,X\otimes_{R}S\right)$ is a finitely generated $T$-module. The conclusion now follows.

(iii) Let $f:M\rightarrow N$ be an $R$-homomorphism in which $M$ and $N$ are $\mathfrak{a}$-cofinite $R$-modules. Then by (ii), the $S$-modules
$S\otimes_RM$ and $S\otimes_R N$ are $\fa S$-cofinite. Since the category $\mathcal{M}\left(S,\mathfrak{a} S\right)_{cof}$ is abelian, we conclude
that $$\ker \left(S\otimes_{R}f\right)\cong S\otimes_{R} \ker f$$ and $$\coker \left(S\otimes_{R}f\right)\cong S\otimes_{R} \coker f$$ are
$\mathfrak{a} S$-cofinite, so by part (ii), $\ker f$ and $\coker f$ are $\mathfrak{a}$-cofinite. Therefore, $\mathcal{M}(R,\mathfrak{a})_{cof}$
is an abelian category.
\end{prf}

Next, we exploit the technique of way-out functors as the main tool to depart from modules to complexes. We first need the definitions.

\begin{definition} \label{2.4.6}
Let $R$ and $S$ be two rings, and $\mathcal{F}: \mathcal{D}(R) \rightarrow \mathcal{D}(S)$ a covariant functor. We say that
\begin{enumerate}
\item[(i)] $\mathcal{F}$ is \textit{way-out left} \index{way-out left functor on derived category} if for every $n \in \mathbb{Z}$, there is an $m \in \mathbb{Z}$, such that for any $R$-complex $X$ with $\sup X \leq m$, we have $\sup \mathcal{F}(X) \leq n$.
\item[(ii)] $\mathcal{F}$ is \textit{way-out right} \index{way-out right functor on derived category} if for every $n \in \mathbb{Z}$, there is an $m \in \mathbb{Z}$, such that for any $R$-complex $X$ with $\inf X \geq m$, we have $\inf \mathcal{F}(X) \geq n$.
\item[(iii)] $\mathcal{F}$ is \textit{way-out} \index{way-out functor on derived category} if it is both way-out left and way-out right.
\end{enumerate}
\end{definition}

The Way-out Lemma \index{way-out lemma} appears in \cite[Ch. I, Proposition 7.3]{Ha2}. However, we need a refined version which is tailored to our needs. Since the proof of the original result in \cite[Ch. I, Proposition 7.3]{Ha2} is left to the reader, we deem it appropriate to include a proof of our refined version for the convenience of the reader as well as bookkeeping.

\begin{lemma} \label{2.4.7}
Let $R$ and $S$ be two rings, and $\mathcal{F}: \mathcal{D}(R) \rightarrow \mathcal{D}(S)$ a triangulated covariant functor. Let $\mathcal{A}$ be an additive subcategory of $\mathcal{M}(R)$, and $\mathcal{B}$ an abelian subcategory of $\mathcal{M}(S)$ which is closed under extensions. Suppose that $H_{i}\left(\mathcal{F}(M)\right) \in \mathcal{B}$ for every $M \in \mathcal{A}$ and every $i \in \mathbb{Z}$. Then the following assertions hold:
\begin{enumerate}
\item[(i)] If $X \in \mathcal{D}_{\square}(R)$ with $H_{i}(X) \in \mathcal{A}$ for every $i \in \mathbb{Z}$, then $H_{i}\left(\mathcal{F}(X)\right) \in \mathcal{B}$ for every $i \in \mathbb{Z}$.
\item[(ii)] If $\mathcal{F}$ is way-out left and $X \in \mathcal{D}_{\sqsubset}(R)$ with $H_{i}(X) \in \mathcal{A}$ for every $i \in \mathbb{Z}$, then $H_{i}\left(\mathcal{F}(X)\right) \in \mathcal{B}$ for every $i \in \mathbb{Z}$.
\item[(iii)] If $\mathcal{F}$ is way-out right and $X \in \mathcal{D}_{\sqsupset}(R)$ with $H_{i}(X) \in \mathcal{A}$ for every $i \in \mathbb{Z}$, then $H_{i}\left(\mathcal{F}(X)\right) \in \mathcal{B}$ for every $i \in \mathbb{Z}$.
\item[(iv)] If $\mathcal{F}$ is way-out and $X \in \mathcal{D}(R)$ with $H_{i}(X) \in \mathcal{A}$ for every $i \in \mathbb{Z}$, then $H_{i}\left(\mathcal{F}(X)\right) \in \mathcal{B}$ for every $i \in \mathbb{Z}$.
\end{enumerate}
\end{lemma}

\begin{prf}
(i): Let $s= \sup(X)$. Since $\amp(X) < \infty$, we argue by induction on $n= \amp(X)$. If $n=0$, then $X \simeq \Sigma^{s} H_{s}(X)$. Therefore,
$$H_{i}\left(\mathcal{F}(X)\right) \cong H_{i}\left(\mathcal{F}\left(\Sigma^{s} H_{s}(X)\right)\right) \cong H_{i-s}\left(\mathcal{F}\left(H_{s}(X)\right)\right) \in \mathcal{B},$$
as $H_{s}(X) \in \mathcal{A}$. Now, let $n \geq 1$ and assume that the result holds for amplitude less than $n$. Since $X \simeq X_{s \subset}$, there is a distinguished triangle
\begin{equation} \label{eq:2.4.7.1}
\Sigma^{s}H_{s}(X) \rightarrow X \rightarrow X_{s-1\subset} \rightarrow.
\end{equation}
It is clear that the two $R$-complexes $\Sigma^{s}H_{s}(X)$ and $X_{s-1\subset}$ have all their homology modules in $\mathcal{A}$
and their amplitudes are less than $n$. Therefore, the induction hypothesis implies that $H_{i}\left(\mathcal{F}\left(\Sigma^{s}H_{s}(X)\right)\right) \in \mathcal{B}$ and $H_{i}\left(\mathcal{F}\left(X_{s-1\subset})\right)\right) \in \mathcal{B}$ for every $i \in \mathbb{Z}$. Applying the functor $\mathcal{F}$ to the distinguished triangle \eqref{eq:2.4.7.1}, we get the distinguished triangle
$$\mathcal{F}\left(\Sigma^{s}H_{s}(X)\right) \rightarrow \mathcal{F}(X) \rightarrow \mathcal{F}(X_{s-1\subset}) \rightarrow,$$
which in turn yields the long exact homology sequence
$$\cdots \rightarrow H_{i+1}\left(\mathcal{F}(X_{s-1\subset})\right) \rightarrow H_{i}\left(\mathcal{F}\left(\Sigma^{s}H_{s}(X)\right)\right) \rightarrow H_{i}\left(\mathcal{F}(X)\right) \rightarrow $$$$ H_{i}\left(\mathcal{F}(X_{s-1\subset})\right) \rightarrow H_{i-1}\left(\mathcal{F}\left(\Sigma^{s}H_{s}(X)\right)\right) \rightarrow \cdots.$$
We break the displayed part of the above exact sequence into the following exact sequences
$$H_{i+1}\left(\mathcal{F}(X_{s-1\subset})\right) \rightarrow H_{i}\left(\mathcal{F}\left(\Sigma^{s}H_{s}(X)\right)\right)\rightarrow K\rightarrow 0,$$
$$0\rightarrow K\rightarrow H_{i}\left(\mathcal{F}(X)\right)\rightarrow L\rightarrow 0,$$
$$0\rightarrow L \rightarrow H_{i}\left(\mathcal{F}(X_{s-1\subset})\right)\rightarrow H_{i-1}\left(\mathcal{F}\left(\Sigma^{s}H_{s}(X)\right)\right).$$
Since the subcategory $\mathcal{B}$ is abelian, we conclude from the first and the third exact sequences above that $K,L \in \mathcal{B}$. Since $\mathcal{B}$ is closed under extensions, the second exact sequence above implies that $H_{i}\left(\mathcal{F}(X)\right) \in \mathcal{B}$ for every $i \in \mathbb{Z}$.

(ii): Let $i \in \mathbb{Z}$. Since $\mathcal{F}$ is way-out left, we can choose an integer $j \in \mathbb{Z}$ corresponding to $i-1$. Apply the functor $\mathcal{F}$ to the distinguished triangle
$$X_{\supset j+1} \rightarrow X \rightarrow X_{j \subset} \rightarrow,$$
to get the distinguished triangle
$$\mathcal{F}(X_{\supset j+1}) \rightarrow \mathcal{F}(X) \rightarrow \mathcal{F}(X_{j \subset}) \rightarrow.$$
From the associated long exact homology sequence, we get
$$0= H_{i+1}\left(\mathcal{F}(X_{j \subset})\right) \rightarrow H_{i}\left(\mathcal{F}(X_{\supset j+1})\right) \rightarrow H_{i}\left(\mathcal{F}(X)\right) \rightarrow H_{i}\left(\mathcal{F}(X_{j \subset})\right) =0,$$
where the vanishing is due to the choice of $j$. Since $X_{\supset j+1} \in \mathcal{D}_{\square}(R)$ with $H_{i}(X_{\supset j+1}) \in \mathcal{A}$ for every $i \in \mathbb{Z}$, it follows from (i) that $H_{i}\left(\mathcal{F}(X_{\supset j+1})\right) \in \mathcal{B}$ for every $i \in \mathbb{Z}$, and as a consequence, $H_{i}\left(\mathcal{F}(X)\right) \in \mathcal{B}$ for every $i \in \mathbb{Z}$.

(iii): Let $i \in \mathbb{Z}$. Since $\mathcal{F}$ is way-out right, we can choose an integer $j \in \mathbb{Z}$ corresponding to $i+1$. Apply the functor $\mathcal{F}$ to the distinguished triangle
$$X_{\supset j} \rightarrow X \rightarrow X_{j-1 \subset} \rightarrow,$$
to get the distinguished triangle
$$\mathcal{F}(X_{\supset j}) \rightarrow \mathcal{F}(X) \rightarrow \mathcal{F}(X_{j-1 \subset}) \rightarrow.$$
From the associated long exact homology sequence, we get
$$0= H_{i}\left(\mathcal{F}(X_{\supset j})\right) \rightarrow H_{i}\left(\mathcal{F}(X)\right) \rightarrow H_{i}\left(\mathcal{F}(X_{j-1 \subset})\right) \rightarrow H_{i-1}\left(\mathcal{F}(X_{\supset j})\right) =0,$$
where the vanishing is due to the choice of $j$. Since $X_{j-1 \subset} \in \mathcal{D}_{\square}(R)$ with $H_{i}(X_{j-1 \subset}) \in \mathcal{A}$ for every $i \in \mathbb{Z}$, it follows from (i) that $H_{i}\left(\mathcal{F}(X_{j-1 \subset})\right) \in \mathcal{B}$ for every $i \in \mathbb{Z}$, and as a consequence, $H_{i}\left(\mathcal{F}(X)\right) \in \mathcal{B}$ for every $i \in \mathbb{Z}$.

(iv): Apply the functor $\mathcal{F}$ to the distinguished triangle
$$X_{\supset 1} \rightarrow X \rightarrow X_{0 \subset} \rightarrow,$$
to get the distinguished triangle
$$\mathcal{F}(X_{\supset 1}) \rightarrow \mathcal{F}(X) \rightarrow \mathcal{F}(X_{0 \subset}) \rightarrow.$$
Since $X_{0 \subset} \in \mathcal{D}_{\sqsubset}(R)$ and $X_{\supset 1} \in \mathcal{D}_{\sqsupset}(R)$ with $H_{i}(X_{0 \subset}), H_{i}(X_{\supset 1}) \in \mathcal{A}$ for every $i \in \mathbb{Z}$, we deduce from (ii) and (iii) that $H_{i}\left(\mathcal{F}(X_{0 \subset})\right), H_{i}\left(\mathcal{F}(X_{\supset 1})\right) \in \mathcal{B}$ for every $i \in \mathbb{Z}$. Using the associated long exact homology sequence, an argument similar to (i) yields that $H_{i}\left(\mathcal{F}(X)\right) \in \mathcal{B}$ for every $i \in \mathbb{Z}$.
\end{prf}

The next result provides us with a suitable transition device from modules to complexes when dealing with cofiniteness.

\begin{theorem} \label{2.4.8}
If $\mathfrak{a}$ is an ideal of $R$, then the functor ${\bf R}\Gamma_{\mathfrak{a}}(-): \mathcal{D}(R) \rightarrow \mathcal{D}(R)$ is triangulated and way-out.
As a consequence, if $H^{i}_{\mathfrak{a}}(M)$ is $\mathfrak{a}$-cofinite for every finitely generated $R$-module $M$ and every $i \geq 0$, and $\mathcal{M}(R,\mathfrak{a})_{cof}$ is an abelian category, then $H^{i}_{\mathfrak{a}}(X)$ is $\mathfrak{a}$-cofinite for every $X \in \mathcal{D}^{f}(R)$ and every $i \in \mathbb{Z}$.
\end{theorem}

\begin{prf} By \cite[Corollary 3.1.4]{Li}, the functor ${\bf R}\Gamma_{\mathfrak{a}}(-): \mathcal{D}(R) \rightarrow \mathcal{D}(R)$ is triangulated and way-out. Now, let $\mathcal{A}$ be the subcategory of finitely generated $R$-modules, and let $\mathcal{B}:= \mathcal{M}(R,\mathfrak{a})_{cof}$. It can be easily seen that $\mathcal{B}$ is closed under extensions. It now follows from Lemma \ref{2.4.7} that $H^{i}_{\mathfrak{a}}(X) = H_{-i}\left({\bf R}\Gamma_{\mathfrak{a}}(X)\right) \in \mathcal{B}$ for every $X \in \mathcal{D}^{f}(R)$ and every $i \in \mathbb{Z}$.
\end{prf}

\begin{corollary} \label{2.4.9}
Let $\mathfrak{a}$ be an ideal of $R$ such that either $\cd(\fa,R)\leq 1$, or $\dim(R)\leq 2$, or $\dim R/\fa\leq 1$. Then the following assertions hold:
\begin{enumerate}
\item[(i)] $H^{i}_{\mathcal{\fa}}(X)$ is $\mathfrak{a}$-cofinite for every $X \in \mathcal{D}^{f}(R)$ and every $i \in \mathbb{Z}$.
\item[(ii)] $\Ass_{R}\left(H^{i}_{\mathfrak{a}}(X)\right)$ is finite for every $i \in \mathbb{Z}$.
\item[(iii)] The Bass number $\mu^{j}_{R}\left(\mathfrak{p},H^{i}_{\mathfrak{a}}(X)\right)$ is finite for every $\mathfrak{p}\in \Spec(R)$ and every $i,j \in \mathbb{Z}$.
\item[(iv)] The Betti number $\beta_{j}^{R}\left(\mathfrak{p},H^{i}_{\mathfrak{a}}(X)\right)$ is finite for every $\mathfrak{p}\in \Spec(R)$ and every $i,j \in \mathbb{Z}$.
\end{enumerate}
\end{corollary}

\begin{prf}
(i): Obvious in light of Corollaries \ref{2.3.11} and \ref{2.4.8}.

(ii): Follows from (i).

(iii): Follows from (i).

(iv): Follows from (i).
\end{prf}

Now we probe the connection between Hartshorne's questions as highlighted in the Introduction.

\begin{theorem} \label{2.4.10}
Let $\mathfrak{a}$ be an ideal of $R$. Consider the following assertions:
\begin{enumerate}
\item[(i)] $H^{i}_{\mathfrak{a}}(M)$ is $\mathfrak{a}$-cofinite for every finitely generated $R$-module $M$ and every $i \geq 0$, and $\mathcal{M}(R,\mathfrak{a})_{cof}$ is an abelian subcategory of $\mathcal{M}(R)$.
\item[(ii)] $H^{i}_{\mathfrak{a}}(X)$ is $\mathfrak{a}$-cofinite for every $X\in \mathcal{D}^{f}(R)$ and every $i \in \mathbb{Z}$.
\item[(iii)] An $R$-complex $X\in \mathcal{D}_{\square}(R)$ is $\mathfrak{a}$-cofinite if and only if $H_{i}(X)$ is $\mathfrak{a}$-cofinite for every $i \in \mathbb{Z}$.
\end{enumerate}
Then the implications $(i)\Rightarrow (ii)$ and $(iii)\Rightarrow (i)$ hold. Furthermore, if $R$ is $\mathfrak{a}$-adically complete, then all three assertions are equivalent.
\end{theorem}

\begin{prf}
(i) $\Rightarrow$ (ii): Follows from Theorem \ref{2.4.8}.

(iii) $\Rightarrow$ (i): Let $M$ be a finitely generated $R$-module. Since $H^{i}_{\mathfrak{a}}(M)=0$ for every $i<0$ or $i> \ara(\mathfrak{a})$, we have
${\bf R}\Gamma_{\mathfrak{a}}(M)\in \mathcal{D}_{\square}(R)$. However, \cite[Proposition 3.2.2]{Li} implies that
$${\bf R}\Hom_{R}\left(R/\mathfrak{a},{\bf R}\Gamma_{\mathfrak{a}}(M)\right) \simeq {\bf R}\Hom_{R}(R/\mathfrak{a},M),$$
showing that ${\bf R}\Gamma_{\mathfrak{a}}(M)$ is $\mathfrak{a}$-cofinite. The hypothesis now implies that $H^{i}_{\mathfrak{a}}(M)=H_{-i}({\bf R}\Gamma_{\mathfrak{a}}\left(M)\right)$ is $\mathfrak{a}$-cofinite for every $i\geq 0$.

Now, let $M$ and $N$ be two $\mathfrak{a}$-cofinite $R$-modules and $f:M\rightarrow N$ an $R$-homomorphism. Let $\varphi:M \rightarrow N$ be the morphism in $\mathcal{D}(R)$ represented by the diagram $M \xleftarrow{1^{M}} M \xrightarrow{f} N$. From the long exact homology sequence associated to the distinguished triangle
\begin{equation} \label{eq:2.4.10.1}
M \xrightarrow {\varphi} N \rightarrow \Cone(f) \rightarrow,
\end{equation}
we deduce that $\Supp_{R}\left(\Cone(f)\right) \subseteq \V(\mathfrak{a})$.
In addition, applying the functor ${\bf R}\Hom_{R}(R/\mathfrak{a},-)$ to \eqref{eq:2.4.10.1}, gives the distinguished triangle
$${\bf R}\Hom_{R}(R/\mathfrak{a},M) \rightarrow {\bf R}\Hom_{R}(R/\mathfrak{a},N) \rightarrow {\bf R}\Hom_{R}\left(R/\mathfrak{a},\Cone(f)\right)
\rightarrow,$$
whose associated long exact homology sequence shows that
$${\bf R}\Hom_{R}\left(R/\mathfrak{a},\Cone(f)\right)\in \mathcal{D}^{f}(R).$$
Hence $\Cone(f)$ is $\mathfrak{a}$-cofinite. However, we have
$$\Cone(f)= \cdots \rightarrow 0 \rightarrow M \xrightarrow {f} N \rightarrow 0 \rightarrow \cdots,$$
so $\Cone(f)\in \mathcal{D}_{\square}(R)$. Thus the hypothesis implies that $H_{i}\left(\Cone(f)\right)$ is $\mathfrak{a}$-cofinite for every
$i\in \mathbb{Z}$. It follows that $\ker f$ and $\coker f$ are $\mathfrak{a}$-cofinite, and as a consequence $\mathcal{M}(R,\mathfrak{a})_{cof}$ is
an abelian subcategory of $\mathcal{M}(R)$.

Now, suppose that $R$ is $\mathfrak{a}$-adically complete.

(ii) $\Rightarrow$ (iii): Let $X\in \mathcal{D}_{\square}(R)$. Suppose that $H_{i}(X)$ is $\mathfrak{a}$-cofinite for every $i \in \mathbb{Z}$. The spectral sequence
$$E_{p,q}^{2}=\Ext_{R}^{p}\left(R/\mathfrak{a},H_{-q}(X)\right) \underset{p}\Rightarrow \Ext_{R}^{p+q}(R/\mathfrak{a},X)$$
from the proof of \cite[Proposition 6.2]{Ha1}, together with the assumption that $E_{p,q}^{2}$ is finitely generated for every $p,q\in \mathbb{Z}$, conspire to imply that $\Ext_{R}^{p+q}(R/\mathfrak{a},X)$ is finitely generated, i.e. $X$ is $\mathfrak{a}$-cofinite. Conversely, assume that $X$ is $\mathfrak{a}$-cofinite. Then by Corollary \ref{2.4.3}, $X\simeq {\bf R}\Gamma_{\mathfrak{a}}(Z)$ for some $Z\in \mathcal{D}^{f}_{\square}(R)$. Thus the hypothesis implies that $$H_{i}(X)
\cong H_{i}\left({\bf R}\Gamma_{\mathfrak{a}}(Z)\right)=H^{-i}_{\mathfrak{a}}(Z)$$ is $\mathfrak{a}$-cofinite for every $i \in \mathbb{Z}$.
\end{prf}

The next result answers Hartshorne's third question.

\begin{corollary} \label{2.4.11}
Let $\mathfrak{a}$ be an ideal of $R$ for which $R$ is $\mathfrak{a}$-adically complete. Suppose that either $\cd(\mathfrak{a},R) \leq 1$, or $\dim(R) \leq 2$, or $\dim(R/\mathfrak{a}) \leq 1$. Then an $R$-complex $X\in \mathcal{D}_{\square}(R)$ is $\mathfrak{a}$-cofinite if and only
if $H_{i}(X)$ is $\mathfrak{a}$-cofinite for every $i \in \mathbb{Z}$.
\end{corollary}

\begin{prf}
Obvious in light of Corollary \ref{2.4.9} and Theorem \ref{2.4.10}.
\end{prf}

The next result provide a relatively complete answer to Hartshorn's third question from another perspective.

\begin{corollary} \label{2.4.12}
Suppose that $R$ admits a dualizing complex $D$, and $\mathfrak{a}$ is an ideal of $R$. Further, suppose that $H^{i}_{\mathfrak{a}}(M)$ is $\mathfrak{a}$-cofinite for every finitely generated $R$-module $M$ and every $i \geq 0$, and $\mathcal{M}(R,\mathfrak{a})_{cof}$ is an abelian subcategory of $\mathcal{M}(R)$. Let $Y\in \mathcal{D}_{\sqsupset}^{f}(R)$, and $X := {\bf R}\Hom_{R}\left(Y,{\bf R}\Gamma_{\mathfrak{a}}(D)\right)$. Then $H_{i}(X)$ is $\mathfrak{a}$-cofinite for every $i \in \mathbb{Z}$.
\end{corollary}

\begin{prf}
Set $Z:= {\bf R}\Hom_{R}(Y,D)$. Then clearly, $Z \in \mathcal{D}_{\sqsubset}^{f}(R)$. Let $\check{C}(\underline{a})$ denote the \v{C}ech complex on a sequence of elements $\underline{a}=a_{1},...,a_{n}\in R$ that generates $\mathfrak{a}$. Now, one has
\begin{equation*}
\begin{split}
X & = {\bf R}\Hom_{R}\left(Y,{\bf R}\Gamma_{\mathfrak{a}}(D)\right)\\
& \simeq \check{C}(\underline{a})\otimes_{R}^{\bf L} {\bf R}\Hom_{\widehat{R}^{\mathfrak{a}}}(Y,D)\\
& \simeq \check{C}(\underline{a})\otimes_{R}^{\bf L} Z\\
& \simeq {\bf R}\Gamma_{\mathfrak{a}}(Z).\\
\end{split}
\end{equation*}
Hence $H_{i}(X) \cong H^{-i}_{\mathfrak{a}}(Z)$ for every $i \in \mathbb{Z}$. Now the result follows.
\end{prf}

\section{Cofiniteness with Respect to Stable Under Specialization Sets}

In this section, we study the notion of cofiniteness for generalized local cohomology modules with respect to stable under specialization sets.

A subset $\mathcal{Z} \subseteq \Spec(R)$ is said to be stable under specialization \index{stable under specialization set} if $V(\fp)\subseteq \mathcal{Z}$ for every $\fp\in \mathcal{Z}$. For such a subset $\mathcal{Z}$, we set
$$F(\mathcal{Z}):=\left\{\mathfrak{a}\lhd R \suchthat V(\mathfrak{a})\subseteq \mathcal{Z}\right\}.$$

If $M$ is an $R$-module, then clearly $\Supp_{R}(M)$ is a stable under specialization subset of $\Spec(R)$. Conversely, given any stable under specialization subset $\mathcal{Z}$ of $\Spec(R)$, one readily checks out that $\mathcal{Z} = \Supp_{R}\left(\bigoplus_{\mathfrak{a} \in F(\mathcal{Z})} R/\mathfrak{a}\right)$. In other words, subsets of $\Spec(R)$ which are stable under specialization are precisely supports of $R$-modules. In particular, any subset of $\Max(R)$ and $V(\mathfrak{a})$ for any ideal $\mathfrak{a}$ of $R$, are stable under specialization subsets of $\Spec(R)$.

Given a stable under specialization subset $\mathcal{Z}$ of $\Spec(R)$, we let
$$\Gamma_{\mathcal{Z}}(M) := \left\{x\in M \suchthat \Supp_{R}(Rx)\subseteq \mathcal{Z}\right\}$$
for an $R$-module $M$, and $\Gamma_{\mathcal{Z}}(f) := f|_{\Gamma_{\mathcal{Z}}(M)}$ for an $R$-homomorphism $f:M\rightarrow N$. This provides us with the so-called $\mathcal{Z}$-torsion functor $\Gamma_{\mathcal{Z}}(-)$ on $\mathcal{M}(R)$, which extends naturally to a functor on $\mathcal{C}(R)$. The extended functor enjoys a right total derived functor ${\bf R}\Gamma_{\mathcal{Z}}(-)$ on $\mathcal{D}(R)$, that can be computed by ${\bf R}\Gamma_{\mathcal{Z}}(X)\simeq \Gamma_{\mathcal{Z}}(I)$, where $X \xrightarrow {\simeq} I$ is any semi-injective resolution of $X$. The functor ${\bf R}\Gamma_{\mathcal{Z}}(-)$ turns out to be triangulated. Besides, we define the $i$th local cohomology module of $X$ with support in $\mathcal{Z}$ as
$$H^{i}_{\mathcal{Z}}(X):= H_{-i}\left({\bf R}\Gamma_{\mathcal{Z}}(X)\right)$$
for every $i\in\mathbb{Z}$. It is clear that upon setting $\mathcal{Z}=V(\mathfrak{a})$ for some ideal $\mathfrak{a}$ of $R$, we recover the usual local cohomology module with respect to $\mathfrak{a}$.

It is straightforward to see that the set $F(\mathcal{Z})$ is a directed poset under reverse inclusion. Let $X$ be an $R$-complex and $X \xrightarrow {\simeq} I$ a semi-injective resolution of $X$. Then one can see by inspection that
$$\Gamma_{\mathcal{Z}}(I_{i})=\bigcup_{\mathfrak{a}\in F(\mathcal{Z})}\Gamma_{\mathfrak{a}}(I_{i}) \cong \underset{\fa\in F(\mathcal{Z})}{\varinjlim}\Gamma_{\mathfrak{a}}(I_{i})$$
for every $i\in \mathbb{Z}$, which in turn implies that $\Gamma_{\mathcal{Z}}(I) \cong \underset{\fa\in F(\mathcal{Z})}{\varinjlim}\Gamma_{\mathfrak{a}}(I)$. Therefore, we have

\begin{equation*}
\begin{split}
H_{\mathcal{Z}}^i(X) & \cong H_{-i}\left({\bf R}\Gamma_{\mathcal{Z}}(X)\right) \\
 & \cong H_{-i}\left(\Gamma_{\mathcal{Z}}(I)\right) \\
 & \cong H_{-i}\left(\underset{\fa\in F(\mathcal{Z})}{\varinjlim}\Gamma_{\mathfrak{a}}(I)\right) \\
 & \cong \underset{\fa\in F(\mathcal{Z})}{\varinjlim} H_{-i}\left(\Gamma_{\mathfrak{a}}(I)\right) \\
 & \cong \underset{\fa\in F(\mathcal{Z})}{\varinjlim} H_{-i}\left({\bf R}\Gamma_{\mathfrak{a}}(X)\right) \\
 & \cong \underset{\fa\in F(\mathcal{Z})}{\varinjlim}H_{\mathfrak{a}}^i(X)
\end{split}
\end{equation*}
for every $i \in \mathbb{Z}$.

Now, we are ready to define the general notion of $\mathcal{Z}$-cofiniteness.

\begin{definition} \label{2.5.1}
Let $\mathcal{Z}$ be a stable under specialization subset of $\Spec(R)$. An $R$-complex $X\in \mathcal{D}(R)$ is said to be $\mathcal{Z}$\textit{-cofinite} \index{cofiniteness of complexes with respect to a stable under specialization sets} if $\Supp_{R}(X)\subseteq \mathcal{Z}$, and ${\bf R}\Hom_{R}(R/\mathfrak{a},X)\in \mathcal{D}^{f}(R)$ for every $\mathfrak{a}\in F(\mathcal{Z})$.
\end{definition}

We denote the full subcategory of $\mathcal{M}(R)$ consisting of $\mathcal{Z}$-cofinite $R$-modules by $\mathcal{M}(R,\mathcal{Z})_{cof}$. The next result lays on some characterizations of $\mathcal{Z}$-cofinite complexes.

\begin{lemma} \label{2.5.2}
Let $\mathcal{Z}$ be a stable under specialization subset of $\Spec(R)$, and $X \in \mathcal{D}(R)$ with $\Supp_{R}(X)\subseteq \mathcal{Z}$. Consider the following conditions:
\begin{enumerate}
\item[(a)] $X$ is $\mathcal{Z}$-cofinite.
\item[(b)] ${\bf R}\Hom_{R}(Y,X)\in \mathcal{D}^{f}(R)$ for every $Y \in \mathcal{D}^{f}_{\square}(R)$ with $\Supp_{R}(Y)\subseteq \mathcal{Z}$.
\item[(c)] $(R/\mathfrak{a})\otimes_{R}^{\bf L} X\in \mathcal{D}^{f}(R)$ for every $\mathfrak{a}\in F(\mathcal{Z})$.
\item[(d)] $Y\otimes_{R}^{\bf L} X\in \mathcal{D}^{f}(R)$ for every $Y \in \mathcal{D}^{f}_{\square}(R)$ with $\Supp_{R}(Y)\subseteq \mathcal{Z}$.
\end{enumerate}
Then the following assertions hold:
\begin{enumerate}
\item[(i)] If $X \in \mathcal{D}_{\sqsubset}(R)$, then (a) and (b) are equivalent.
\item[(ii)] If $X \in \mathcal{D}_{\sqsupset}(R)$, then (c) and (d) are equivalent.
\item[(iii)] If $X \in \mathcal{D}_{\square}(R)$, then all the assertions are equivalent.
\end{enumerate}
\end{lemma}

\begin{prf}
(i): Suppose that (a) holds and $Y \in \mathcal{D}^{f}_{\square}(R)$ with $\Supp_{R}(Y)\subseteq \mathcal{Z}$. Let $\mathfrak{a}=\ann_{R}\left(\oplus_{i\in \mathbb{Z}}H_{i}(Y)\right)$.
Then $V(\mathfrak{a})= \Supp_{R}(Y) \subseteq \mathcal{Z}$, so $\mathfrak{a}\in F(\mathcal{Z})$. Therefore, ${\bf R}\Hom_{R}(R/\mathfrak{a},X)\in\mathcal{D}^{f}(R)$.
Now, \cite[Proposition 7.2]{WW} implies that ${\bf R}\Hom_{R}(Y,X)\in \mathcal{D}^{f}(R)$. The converse is clear.

(ii): Similar to (i) using \cite[Proposition 7.1]{WW}.

(iii): Fix $\mathfrak{a}\in F(\mathcal{Z})$. Then \cite[Proposition 7.4]{WW} yields that ${\bf R}\Hom_{R}(R/\mathfrak{a},X)\in \mathcal{D}^{f}(R)$ if and only if $(R/\mathfrak{a})\otimes_{R}^{\bf L} X\in \mathcal{D}^{f}(R)$.
\end{prf}

We collect some basic properties of $\mathcal{Z}$-cofinite $R$-complexes in the following result. Its first part indicates that in the case where $\mathcal{Z}= V(\mathfrak{a})$ for some ideal $\mathfrak{a}$ of $R$, our definition of $\mathcal{Z}$-cofiniteness coincides with the notion of $\mathfrak{a}$-cofiniteness defined in the previous sections. In other words, $\mathcal{Z}$-cofiniteness is a generalization of $\mathfrak{a}$-cofiniteness.

\begin{lemma} \label{2.5.3}
Let $\mathcal{Z}$ be a stable under specialization subset of $\Spec(R)$ and $X\in \mathcal{D}(R)$. Then the following assertions hold:
\begin{enumerate}
\item[(i)] If $\mathfrak{a}$ is an ideal of $R$ and $X \in \mathcal{D}_{\sqsubset}(R)$, then $X$ is $V(\mathfrak{a})$-cofinite if and only if $\Supp_{R}(X)\subseteq V(\mathfrak{a})$ and ${\bf R}\Hom_{R}(R/\mathfrak{a},X)\in \mathcal{D}^{f}(R)$.
\item[(ii)] If $X\in \mathcal{D}_{\sqsubset}(R)$ and $H_{i}(X)$ is $\mathcal{Z}$-cofinite for every $i\in \mathbb{Z}$, then $X$ is $\mathcal{Z}$-cofinite.
\item[(iii)] If $X\in \mathcal{D}^{f}_{\sqsubset}(R)$, then ${\bf R}\Gamma_{\mathcal{Z}}(X)$ is $\mathcal{Z}$-cofinite.
\item[(iv)] If $X\in \mathcal{D}_{\sqsubset}(R)$ and $X$ is $\mathcal{Z}$-cofinite, then the Bass number $\mu^{i}_{R}(\mathfrak{p},X)$ is finite for every $\mathfrak{p}\in \Spec (R)$ and every $i \in \mathbb{Z}$.
\item[(v)] If $X\in \mathcal{D}_{\square}(R)$ and $X$ is $\mathcal{Z}$-cofinite, then the Betti number $\beta_{i}^{R}(\mathfrak{p},X)$ is finite for every $\mathfrak{p}\in \Spec (R)$ and every $i \in \mathbb{Z}$.
\end{enumerate}
\end{lemma}

\begin{prf}
(i): Suppose that $\Supp_{R}(X)\subseteq V(\mathfrak{a})$ and ${\bf R}\Hom_{R}(R/\mathfrak{a},X)\in \mathcal{D}^{f}(R)$. Let $\mathfrak{b}\in F\left(V(\mathfrak{a})\right)$. It follows that $\Supp_{R}(R/\mathfrak{b}) \subseteq V(\mathfrak{a})$. Now, by \cite[Proposition 7.2]{WW} we are through. The converse is clear.

(ii): Since $H_{i}(X)$ is $\mathcal{Z}$-cofinite for every $i\in \mathbb{Z}$, we have
$$\Supp_{R}(X)=\bigcup_{i\in \mathbb{Z}}\Supp_{R}\left(H_{i}(X)\right) \subseteq \mathcal{Z}.$$
Let $\mathfrak{a}\in F(\mathcal{Z})$.
The spectral sequence
$$E_{p,q}^{2}=\Ext_{R}^{p}\left(R/\mathfrak{a},H_{-q}(X)\right) \underset{p}\Rightarrow \Ext_{R}^{p+q}(R/\mathfrak{a},X)$$
from the proof of \cite[Proposition 6.2]{Ha1}, together with the assumption that $E_{p,q}^{2}$ is finitely generated for every $p,q\in \mathbb{Z}$, conspire to imply that $\Ext_{R}^{p+q}(R/\mathfrak{a},X)$ is finitely generated, i.e. $X$ is $\mathcal{Z}$-cofinite.

(iii): It is clear that $\Supp_{R}\left({\bf R}\Gamma_{\mathcal{Z}}(X)\right)\subseteq \mathcal{Z}$. Since $X$ is homologically left-bounded, there exists a left-bounded semi-injective resolution $X \xrightarrow {\simeq} I$ of $X$.
As $I$ is an $R$-complex of injective modules, and
$$\Gamma_{\mathcal{Z}}(I_{i}) \cong \underset{\fa\in F(\mathcal{Z})}{\varinjlim}\Gamma_{\mathfrak{a}}(I_{i})$$
for every $i \in \mathbb{Z}$, we see that $\Gamma_{\mathcal{Z}}(I)$ is a left-bounded $R$-complex of injective modules, and thus $\Gamma_{\mathcal{Z}}(I)$ is semi-injective.
Fix $\mathfrak{a}\in F(\mathcal{Z})$. For every $R$-module $\Hom_{R}\left(R/\mathfrak{a},\Gamma_{\mathcal{Z}}(I_{i})\right)$ in the $R$-complex $\Hom_{R}\left(R/\mathfrak{a},\Gamma_{\mathcal{Z}}(I)\right)$, we have
$$\Hom_{R}\left(R/\mathfrak{a},\Gamma_{\mathcal{Z}}(I_{i})\right) \cong \Hom_{R}\left(R/\mathfrak{a},\Gamma_{\mathfrak{a}}(I_{i})\right) \cong \Hom_{R}\left(R/\mathfrak{a},I_{i}\right).$$
Indeed, if $N$ is a submodule of $M$ containing $\Gamma_{\mathfrak{a}}(M)$, then every $R$-homomorphism $f: R/\mathfrak{a} \rightarrow N$ has its image in $\Gamma_{\mathfrak{a}}(M)$.
Hence one has
\begin{equation*}
\begin{split}
{\bf R}\Hom_{R}\left(R/\mathfrak{a},{\bf R}\Gamma_{\mathcal{Z}}(X)\right) & \simeq {\bf R}\Hom_{R}\left(R/\mathfrak{a},\Gamma_{\mathcal{Z}}(I)\right) \\
 & \simeq \Hom_{R}\left(R/\mathfrak{a},\Gamma_{\mathcal{Z}}(I)\right) \\
 & \simeq \Hom_{R}(R/\mathfrak{a},I) \\
 & \simeq {\bf R}\Hom_{R}(R/\mathfrak{a},X).
\end{split}
\end{equation*}
But $X\in \mathcal{D}^{f}_{\sqsubset}(R)$, so ${\bf R}\Hom_{R}(R/\mathfrak{a},X) \in \mathcal{D}^{f}_{\sqsubset}(R)$, and thus the assertion follows.

(iv): For every $i\in \mathbb{Z}$, we have
$$\mu^{i}_{R}(\mathfrak{p},X):= \rank_{R_{\mathfrak{p}}/\mathfrak{p}R_{\mathfrak{p}}}\left(\Ext_{R_{\mathfrak{p}}}^{i}\left(R_{\mathfrak{p}}/\mathfrak{p}R_{\mathfrak{p}},
X_{\mathfrak{p}}\right)\right).$$
If $\mathfrak{p}\not\in \mathcal{Z}$, then $\mathfrak{p}\not\in \Supp_{R}(X)$, so $\mu^{i}_{R}(\mathfrak{p},X)=0$. If $\mathfrak{p}\in \mathcal{Z}$, then $\V(\mathfrak{p}) \subseteq \mathcal{Z}$, so by definition, ${\bf R}\Hom_{R}(R/\mathfrak{p},X)\in \mathcal{D}^{f}(R)$, whence
$\mu^{i}_{R}(\mathfrak{p},X)< \infty$.

(v): For every $i\in \mathbb{Z}$, we have
$$\beta_{i}^{R}(\mathfrak{p},X):= \rank_{R_{\mathfrak{p}}/\mathfrak{p}R_{\mathfrak{p}}}\left(\Tor_{i}^{R_{\mathfrak{p}}}\left(R_{\mathfrak{p}}/\mathfrak{p}R_{\mathfrak{p}},
X_{\mathfrak{p}}\right)\right).$$
If $\mathfrak{p}\not\in \mathcal{Z}$, then as in (iv), $\beta_{i}^{R}(\mathfrak{p},X)=0$. If $\mathfrak{p}\in \mathcal{Z}$, then Lemma \ref{2.5.2} implies that $(R/\mathfrak{p})\otimes_{R}^{\bf L} X\in \mathcal{D}^{f}(R)$, thereby
$\beta_{i}^{R}(\mathfrak{p},X)< \infty$.
\end{prf}

The following result may be of independent interest.

\begin{lemma} \label{2.5.4}
Let $\mathcal{Z}$ be a stable under specialization subset of $\Spec(R)$, and $M$ an $R$-module. Then for any $\mathfrak{a} \in F(\mathcal{Z})$, there is a spectral sequence
$$E^{2}_{p,q}=\Ext^{-p}_{R}\left(R/\mathfrak{a},H^{-q}_{\mathcal{Z}}(M)\right) \underset {p} \Rightarrow \Ext^{-p-q}_{R}(R/\mathfrak{a},M).$$
\end{lemma}

\begin{prf}
Let $\mathcal{F}=\Hom_{R}(R/\mathfrak{a},-)$, and $\mathcal{G}=\Gamma_{\mathcal{Z}}(-)$. Then $\mathcal{F}$ is left exact, and
$\mathcal{G}(I)$ is right $\mathcal{F}$-acyclic for every injective $R$-module $I$, since
$$\mathcal{G}(I)=\Gamma_{\mathcal{Z}}(I)\cong \underset{\fa\in F(\mathcal{Z})}{\varinjlim}\Gamma_{\mathfrak{a}}(I)$$
is injective. Therefore, by \cite[Theorem 10.47]{Ro}, there is a Grothendieck spectral sequence
$$E^{2}_{p,q}= R^{-p}\mathcal{F}\left(R^{-q}\mathcal{G}(M)\right) \underset {p} \Rightarrow R^{-p-q}(\mathcal{FG})(M).$$
If $N$ is an $R$-module, then as observed before, $\Hom_{R}\left(R/\mathfrak{a},\Gamma_{\mathcal{Z}}(N)\right)\cong \Hom_{R}(R/\mathfrak{a},N)$,
whence $R^{-p-q}(\mathcal{FG})(M) \cong \Ext^{-p-q}_{R}(R/\mathfrak{a},M)$.
\end{prf}

Now we proceed to extend the results of the previous sections on $\mathfrak{a}$-cofiniteness to $\mathcal{Z}$-cofiniteness.

We first deal with the case $\cd(\mathcal{Z},X) \leq 1$. Given a stable under specialization subset $\mathcal{Z}$ of $\Spec(R)$, we define the cohomological dimension \index{cohomological dimension with respect to a stable under specialization set} of an $R$-complex $X$ with respect to $\mathcal{Z}$ as
$$\cd(\mathcal{Z},X):=\sup \left\{i \in \mathbb{Z} \suchthat H^{i}_{\mathcal{Z}}(X)\neq 0 \right\}.$$

\begin{theorem} \label{2.5.5}
Let $\mathcal{Z}$ be a stable under specialization subset of $\Spec(R)$, and $M$ a finitely generated $R$-module. If $\cd(\mathcal{Z},M)\leq 1$, then $H^{i}_{\mathcal{Z}}(M)$
is $\mathcal{Z}$-cofinite for every $i \geq 0$.
\end{theorem}

\begin{prf}
Clearly, one has $\Supp_{R}\left(H^{i}_{\mathcal{Z}}(M)\right) \subseteq \mathcal{Z}$. Let $\mathfrak{a} \in F(\mathcal{Z})$. Then the $R$-module $\Ext^{j}_{R}\left(R/\mathfrak{a},\Gamma_{\mathcal{Z}}(M)\right)$ is finitely generated for every $j \geq 0$, and $H^{i}_{\mathcal{Z}}(M)=0$ for every $i \geq 2$ by the assumption. Therefore, it remains to show that $\Ext^{j}_{R}\left(R/\mathfrak{a},H^{1}_{\mathcal{Z}}(M)\right)$ is finitely generated for every $j \geq 0$. Since $\Gamma_{\mathcal{Z}}(M)$ is $\mathcal{Z}$-torsion, $\Gamma_{\mathcal{Z}}(M)$ has an injective resolution whose modules are $\mathcal{Z}$-torsion.
It follows that $H^{i}_{\mathcal{Z}}\left(\Gamma_{\mathcal{Z}}(M)\right) = 0$ for every $i \geq 1$, and thus the long exact local cohomology sequence associated with the short exact sequence
$$0 \rightarrow \Gamma_{\mathcal{Z}}(M) \rightarrow M \rightarrow M/\Gamma_{\mathcal{Z}}(M) \rightarrow 0,$$
shows that $H^{i}_{\mathcal{Z}}(M) \cong H^{i}_{\mathcal{Z}}\left(M/\Gamma_{\mathcal{Z}}(M)\right)$ for every $i \geq 1$. Therefore, replacing $M$ by $M/\Gamma_{\mathcal{Z}}(M)$, we may assume that $\Gamma_{\mathcal{Z}}(M)=0$. As a consequence, the spectral sequence
$$E^{2}_{-p,-q}=\Ext^{-p}_{R}\left(R/\mathfrak{a},H^{-q}_{\mathcal{Z}}(M)\right) \underset {p} \Rightarrow \Ext^{-p-q}_{R}(R/\mathfrak{a},M)$$
of Lemma \ref{2.5.4} collapses on the row $q=-1$, and thus yields an isomorphism
$$\Ext^{j}_{R}\left(R/\mathfrak{a},H^{1}_{\mathcal{Z}}(M)\right) \cong \Ext^{j+1}_{R}(R/\mathfrak{a},M)$$
for every $j \geq 0$, which in turn establishes the desired conclusion.
\end{prf}

\begin{lemma} \label{2.5.6}
Let $\mathcal{Z}$ be a stable under specialization subset of $\Spec(R)$. Let $M$ be an $R$-module such that $\Ass_{R}(M) \cap \mathcal{Z} \cap \Max(R)$ is a finite set. Then the following assertions are equivalent:
\begin{enumerate}
\item[(i)] $(0:_{M}\mathfrak{a})$ is an artinian $R$-module for every $\mathfrak{a} \in F(\mathcal{Z})$.
\item[(ii)] $\Gamma_{\mathfrak{a}}(M)$ is an artinian $R$-module for every $\mathfrak{a} \in F(\mathcal{Z})$.
\item[(iii)] $\Gamma_{\mathcal{Z}}(M)$ is an artinian $R$-module.
\end{enumerate}
\end{lemma}

\begin{prf}
(i) $\Rightarrow$ (ii): Let $\mathfrak{a}\in F(\mathcal{Z})$. Then $\Gamma_{\mathfrak{a}}(M)$ is $\mathfrak{a}$-torsion, and
$(0:_{\Gamma_{\mathfrak{a}}(M)}\mathfrak{a}) = (0:_{M}\mathfrak{a})$
is artinian, so using \cite[Theorem 1.3]{Me3}, we conclude that $\Gamma_{\mathfrak{a}}(M)$ is artinian.

(ii) $\Rightarrow$ (iii): Set
$$\mathfrak{J}:= \bigcap_{\mathfrak{m}\in \Ass_{R}(M)\cap \mathcal{Z} \cap \Max(R)} \mathfrak{m}.$$
Since $\Ass_{R}(M)\cap \mathcal{Z} \cap \Max(R)$ is a finite set, we see that
$$V(\mathfrak{J})= \Ass_{R}(M)\cap \mathcal{Z} \cap \Max(R),$$
so that $\mathfrak{J} \in F(\mathcal{Z})$. Therefore, $\Gamma_{\mathfrak{J}}(M)$ is an artinian $R$-module. Now, let $\mathfrak{a} \in F(\mathcal{Z})$. Then $\Gamma_{\mathfrak{a}}(M)$ is an artinian $R$-module by the assumption. Let $x \in \Gamma_{\mathfrak{a}}(M)$. As $Rx$ is artinian, it turns out that
\begin{equation*}
\begin{split}
\Ass_{R}(Rx) & \subseteq \Ass_{R}(M) \cap V(\mathfrak{a}) \cap \Max(R) \\
 & \subseteq \Ass_{R}(M) \cap \mathcal{Z} \cap \Max(R) \\
 & = V(\mathfrak{J}).
\end{split}
\end{equation*}
Hence $x \in \Gamma_{\mathfrak{J}}(M)$. This yields that
$$\Gamma_{\mathcal{Z}}(M) = \bigcup_{\mathfrak{a}\in F(\mathcal{Z})} \Gamma_{\mathfrak{a}}(M) \subseteq \Gamma_{\mathfrak{J}}(M) \subseteq \Gamma_{\mathcal{Z}}(M),$$
thereby $\Gamma_{\mathcal{Z}}(M)=\Gamma_{\mathfrak{J}}(M)$ is artinian.

(iii) $\Rightarrow$ (i): Clear, since $(0:_{M}\mathfrak{a}) \subseteq \Gamma_{\mathcal{Z}}(M)$ for every $\mathfrak{a} \in F(\mathcal{Z})$.
\end{prf}

\begin{lemma} \label{2.5.7}
Let $\mathcal{Z}$ be a stable under specialization subset of $\Spec(R)$. Let $M$ be an $R$-module, and $r \geq 0$ an integer. Consider the following conditions:
\begin{enumerate}
\item[(a)] $H^{i}_{\mathcal{Z}}(M)$ is an artinian $R$-module for every $0 \leq i \leq r$.
\item[(b)] $\Ext^{i}_{R}(R/\mathfrak{a},M)$ is an artinian $R$-module for every $\mathfrak{a} \in F(\mathcal{Z})$ and for every $0 \leq i \leq r$.
\end{enumerate}
Then the following assertions hold:
\begin{enumerate}
\item[(i)] (a) implies (b).
\item[(ii)] If $\Supp_{R}(M)\cap \mathcal{Z} \cap \Max(R)$ is a finite set, then (a) and (b) are equivalent.
\end{enumerate}
\end{lemma}

\begin{prf}
Let
$$I: 0\rightarrow I_{0} \xrightarrow {\partial^{I}_{0}} I_{-1} \xrightarrow {\partial^{I}_{-1}} I_{-2} \rightarrow \cdots$$
be a minimal injective resolution of $M$. Given any $\mathfrak{a}\in F(\mathcal{Z})$, consider the two $R$-complexes
$$\Hom_{R}(R/\mathfrak{a},I): 0\rightarrow \Hom_{R}(R/\mathfrak{a},I_{0}) \xrightarrow {\Hom_{R}\left(R/\mathfrak{a},\partial^{I}_{0}\right)} \Hom_{R}(R/\mathfrak{a},I_{-1}) $$$$ \xrightarrow {\Hom_{R}\left(R/\mathfrak{a},\partial^{I}_{-1}\right)} \Hom_{R}(R/\mathfrak{a},I_{-2}) \rightarrow \cdots,$$
and
$$\Gamma_{\mathcal{Z}}(I): 0\rightarrow \Gamma_{\mathcal{Z}}(I_{0}) \xrightarrow {\Gamma_{\mathcal{Z}}\left(\partial^{I}_{0}\right)} \Gamma_{\mathcal{Z}}(I_{-1}) \xrightarrow {\Gamma_{\mathcal{Z}}\left(\partial^{I}_{-1}\right)} \Gamma_{\mathcal{Z}}(I_{-2}) \rightarrow \cdots.$$
One can easily check that $\ker \left(\Hom_{R}\left(R/\mathfrak{a},\partial^{I}_{-i}\right)\right)$ is an essential submodule of $\Hom(R/\mathfrak{a},I_{-i})$, and $\ker \left(\Gamma_{\mathcal{Z}}\left(\partial^{I}_{-i}\right)\right)$ is an essential submodule of $\Gamma_{\mathcal{Z}}(I_{-i})$ for every $i \geq 0$.

Let $$X: 0 \rightarrow X_{0} \xrightarrow {\partial^{X}_{0}} X_{-1} \xrightarrow {\partial^{X}_{-1}} X_{-2} \rightarrow \cdots$$
be an $R$-complex such that $\ker \partial^{X}_{-i}$ is an essential submodule of $X_{-i}$ for every $i \geq 0$. For any given
$r\geq 0$, \cite[Lemma 5.4]{Me2} yields that $X_{-i}$ is an artinian $R$-module for every $0 \leq i \leq r$ if and only if $H_{-i}(X)$ is an artinian $R$-module for every $0 \leq i \leq r$. In the remainder of the proof, we apply this twice.

Now, we prove the following:

(i): Let $\mathfrak{a} \in F(\mathcal{Z})$. Applying the discussion above to the $R$-complex $\Gamma_{\mathcal{Z}}(I)$, we
see that $\Gamma_{\mathcal{Z}}(I_{-i})$ is artinian for every $0 \leq i \leq r$. Since
$$\Hom_{R}(R/\mathfrak{a},I_{-i}) \cong (0:_{I_{-i}}\mathfrak{a}) \subseteq \Gamma_{\mathcal{Z}}(I_{-i}),$$
it is obvious that $\Hom_{R}(R/\mathfrak{a},I_{-i})$ is artinian for every $0 \leq i \leq r$.
This shows that $\Ext^{i}_{R}(R/\mathfrak{a},M)$ is artinian for every $0 \leq i \leq r$.

(ii): Another application of the discussion above to the $R$-complex $\Hom_{R}(R/\mathfrak{a},I)$, yields that $\Hom_{R}(R/\mathfrak{a},I_{-i})$ is artinian for every $\mathfrak{a} \in F(\mathcal{Z})$ and $0 \leq i \leq r$. Since $I$ is a minimal injective resolution of $M$, one can see that $\Supp_{R}(I_{-i}) \subseteq \Supp_{R}(M)$ for every $i \geq 0$. Therefore, we can use Lemma \ref{2.5.6} to deduce that $\Gamma_{\mathcal{Z}}(I_{-i})$ is artinian for every $0 \leq i \leq r$. A fortiori, $H^{i}_{\mathcal{Z}}(M)$ is artinian for
every $0 \leq i \leq r$.
\end{prf}

The next example shows that the finiteness condition of the sets $\Ass_{R}(M)\cap \mathcal{Z} \cap \Max(R)$ and $\Supp_{R}(M)\cap \mathcal{Z} \cap \Max(R)$ cannot be removed from Lemmas \ref{2.5.6} and \ref{2.5.7}, respectively.

\begin{example} \label{2.5.8}
Let $R$ be a ring with infinitely many maximal ideals. Let $\mathcal{Z}:= \Max(R)$, and $M:= \bigoplus_{\mathfrak{m} \in \mathcal{Z}} R/\mathfrak{m}$. Then $F(\mathcal{Z})= \{\mathfrak{a} \lhd R \suchthat \dim(R/\mathfrak{a})\leq 0 \}$. Thus for any given proper ideal $\mathfrak{a} \in F(\mathcal{Z})$, there are finitely many maximal ideals in $V(\mathfrak{a})$, say $\mathfrak{m}_{1},..., \mathfrak{m}_{n}$. It follows that $\Gamma_{\mathfrak{a}}(M)=\bigoplus_{i=1}^{n} R/\mathfrak{m}_{i}$. Now $\Gamma_{\mathfrak{a}}(M)$ is artinian, while $\Gamma_{\mathcal{Z}}(M)=M$ fails to be artinian as it contains infinitely
many direct summands.
\end{example}

\begin{corollary} \label{2.5.9}
Let $R$ be a semilocal ring with Jacobson radical $\mathfrak{J}$, $\mathcal{Z}$ a stable under specialization subset of $\Spec(R)$, and $M$ an $R$-module. Let $(-)^{\vee}:= \Hom_{R}\left(-,E_{R}\left(R/\mathfrak{J}\right)\right)$ be the Matlis duality functor. Then the following assertions are equivalent for any given $r \geq 0$:
\begin{enumerate}
\item[(i)] $H^{i}_{\mathcal{Z}}\left(M^{\vee}\right)$ is an artinian $R$-module for every $0 \leq i \leq r$.
\item[(ii)] $\Ext^{i}_{R}\left(R/\mathfrak{a},M^{\vee}\right)$ is an artinian $R$-module for every $\mathfrak{a} \in F(\mathcal{Z})$ and every $0 \leq i \leq r$.
\item[(iii)] $\Tor^{R}_{i}(R/\mathfrak{a},M)$ is a finitely generated $R$-module for every $\mathfrak{a} \in F(\mathcal{Z})$ and every $0 \leq i \leq r$.
\end{enumerate}
\end{corollary}

\begin{prf}
(i) $\Leftrightarrow$ (ii): Follows from Lemma \ref{2.5.7}.

One can easily see that
$$E_{R}\left(R/\mathfrak{J}\right) \cong \bigoplus_{\mathfrak{m}\in \Max(R)}E_{R}(R/\mathfrak{m})$$
is an artinian injective cogenerator for $R$. Fix $\mathfrak{a} \in F(\mathcal{Z})$ and $0 \leq i \leq r$, and let $N:=
\Tor^{R}_{i}(R/\mathfrak{a},M)$ for the rest of the proof. It is easy to see that $N^{\vee}\cong  \Ext^{i}_{R}
\left(R/\mathfrak{a},M^{\vee}\right)$.

(ii) $\Rightarrow$ (iii): Let $\Max(R)=\{\mathfrak{m}_1,\dots ,\mathfrak{m}_n\}$ and set $T:=\widehat{R}^{\mathfrak{J}}$
and $T_j:=\widehat{R_{\mathfrak{m}_{j}}}^{\mathfrak{m}_{j}}$ for every $j=1,\dots, n$. We know that $T\cong \prod_{j=1}^nT_j$ and $T$ is a
$\mathfrak{J}T$-adically complete semilocal ring with
$$\Max(T)=\{\mathfrak{m}_{j}T \suchthat j=1,\dots, n\}.$$
Any $\mathfrak{J}$-torsion $R$-module possesses a $T$-module structure in such a way a subset is an $R$-submodule if and only
if it is a $T$-submodule. In particular, $N^{\vee}$ and $E_{R}\left(R/\mathfrak{J}\right)$ are artinian $T$-modules.
Moreover, one may easily check that the two $T$-modules $E_{R}\left(R/\mathfrak{J}\right)$ and $E_{T}\left(T/\mathfrak{J}T\right)$
are isomorphic and $\mathfrak{J}T$ is the Jacobson radical of $T$. Putting everything together, we obtain
\begin{equation*}
\begin{split}
N^{\vee} & \cong \Hom_{R}\left(N,E_{R}\left(R/\mathfrak{J}\right)\right) \\
 & \cong \Hom_{R}\left(N,E_{T}\left(T/\mathfrak{J}T\right)\right) \\
 & \cong \Hom_{R}\left(N,\Hom_{T}\left(T,E_{T}\left(T/\mathfrak{J}T\right)\right)\right)  \\
 & \cong \Hom_{T}\left(N\otimes_{R} T,E_{T}\left(T/\mathfrak{J}T\right)\right). \\
 \end{split}
\end{equation*}
Applying the Matlis Duality Theorem over the ring $T$ \cite[Proposition 4 (c)]{CW}, we deduce that $N\otimes_{R}T$ is a
finitely generated $T$-module and by the faithfully flatness of the completion map $\theta_{R}: R \rightarrow T$,
we infer that $N$ is a finitely generated  $R$-module.

(iii) $\Rightarrow$ (ii): There is an exact sequence $$R^{n}\rightarrow N \rightarrow 0,$$ which yields the exact sequence
$$0 \rightarrow N^{\vee} \rightarrow E_{R}\left(R/\mathfrak{J}\right)^{n}.$$ It follows that $\Ext^{i}_{R}\left(R/\mathfrak{a},
M^{\vee}\right) \cong N^{\vee}$ is artinian.
\end{prf}

\begin{corollary} \label{2.5.10}
Let $R$ be a semilocal ring, $\mathcal{Z}$ a stable under specialization subset of $\Spec(R)$, and $M$ an $R$-module with
$\Supp_{R}(M) \subseteq \mathcal{Z}$. Then the following conditions are equivalent:
\begin{enumerate}
\item[(i)] $M$ is $\mathcal{Z}$-cofinite.
\item[(ii)] $\Tor^{R}_{i}(R/\mathfrak{a},M)$ is finitely generated for every $\mathfrak{a} \in F(\mathcal{Z})$ and every $i \geq 0$.
\item[(iii)] $H^{i}_{\mathcal{Z}}\left(M^{\vee}\right)$ is artinian for every $i \geq 0$.
\item[(iv)] $H^{i}_{\mathcal{Z}}\left(M^{\vee}\right)$ is artinian for every $0 \leq i \leq \cd(\mathcal{Z},M^{\vee})$.
\item[(v)] $\Ext^{i}_{R}\left(R/\mathfrak{a},M^{\vee}\right)$ is artinian for every $\mathfrak{a} \in F(\mathcal{Z})$ and every $0 \leq i \leq \cd(\mathcal{Z},M^{\vee})$.
\item[(vi)] $\Tor^{R}_{i}(R/\mathfrak{a},M)$ is finitely generated for every $\mathfrak{a} \in F(\mathcal{Z})$ and every $0 \leq i \leq \cd(\mathcal{Z},M^{\vee})$.
\end{enumerate}
\end{corollary}

\begin{prf}
(i) $\Leftrightarrow$ (ii): Follows from Lemma \ref{2.5.2}.

(ii) $\Leftrightarrow$ (iii): Follows from Corollary \ref{2.5.9}.

(iii) $\Leftrightarrow$ (iv): Obvious.

(iv) $\Leftrightarrow$ (v) $\Leftrightarrow$ (vi): Follows from Corollary \ref{2.5.9}.
\end{prf}

The next lemma generalizes \cite[Theorem 2.2]{DNT}, which may sound appealing in its sake.

\begin{lemma} \label{2.5.11}
Let $\mathcal{Z}$ be a stable under specialization subset of $\Spec(R)$, $M$ an $R$-module, and $N$ a finitely generated $R$-module. If $\Supp_{R}(M)\subseteq \Supp_{R}(N)$, then $\cd(\mathcal{Z},M) \leq \cd(\mathcal{Z},N)$.
\end{lemma}

\begin{prf}
Since $H^{i}_{\mathcal{Z}}(-)$ commutes with direct limits and $M$ can be written as a direct limit of its finitely generated submodules, we may assume that $M$ is finitely generated. Now, the proof is a straightforward adaptation of the argument given in \cite[Theorem 2.2]{DNT}.
\end{prf}

In the sequel, we use the straightforward observation that if any two modules in a short exact sequence are $\mathcal{Z}$-cofinite, then so is the third.

\begin{theorem} \label{2.5.12}
Let $R$ be a semilocal ring, $\mathcal{Z}$ a stable under specialization subset of $\Spec(R)$ such that $\cd(\mathcal{Z},R) \leq 1$. Then $\mathcal{M}(R,\mathcal{Z})_{cof}$ is an abelian subcategory of $\mathcal{M}(R)$.
\end{theorem}

\begin{prf}
First of all note that Lemma \ref{2.5.11} implies that $\cd(\mathcal{Z},L) \leq 1$ for every $R$-module $L$.
Let $M$ and $N$ be two $\mathcal{Z}$-cofinite $R$-modules and $f:M\rightarrow N$ an $R$-homomorphism. Let $\mathfrak{a}\in F(\mathcal{Z})$. The short exact sequence
\begin{equation} \label{eq:2.5.12.1}
0 \rightarrow \ker f \rightarrow M \rightarrow \im f \rightarrow 0,
\end{equation}
gives the exact sequence
$$(R/\mathfrak{a})\otimes_{R}M \rightarrow (R/\mathfrak{a})\otimes_{R} \im f \rightarrow 0,$$
which in turn implies that $(R/\mathfrak{a})\otimes_{R} \im f$ is finitely generated, since $(R/\mathfrak{a})\otimes_{R}M$ is finitely generated by Corollary \ref{2.5.10}. The short exact sequence
\begin{equation} \label{eq:2.5.12.2}
0 \rightarrow \im f \rightarrow N \rightarrow \coker f \rightarrow 0,
\end{equation}
gives the exact sequence
\begin{equation} \label{eq:2.5.12.3}
\Tor_{1}^{R}(R/\mathfrak{a},N) \rightarrow \Tor_{1}^{R}(R/\mathfrak{a},\coker f) \rightarrow (R/\mathfrak{a})\otimes_{R} \im f \rightarrow (R/\mathfrak{a})\otimes_{R} N \rightarrow $$$$(R/\mathfrak{a})\otimes_{R} \coker f \rightarrow 0.
\end{equation}
As $(R/\mathfrak{a})\otimes_{R} N$ and $\Tor_{1}^{R}(R/\mathfrak{a},N)$ are finitely generated by Corollary \ref{2.5.10}, the exact sequence \eqref{eq:2.5.12.3} shows that $(R/\mathfrak{a})\otimes_{R} \coker f$ and $\Tor_{1}^{R}(R/\mathfrak{a},\coker f)$ are finitely generated, and thus Corollary \ref{2.5.10} implies that $\coker f$ is $\mathcal{Z}$-cofinite. From the short exact sequence \eqref{eq:2.5.12.2}, we conclude that $\im f$ is $\mathcal{Z}$-cofinite, and from the short exact sequence \eqref{eq:2.5.12.1}, we infer that $\ker f$ is $\mathcal{Z}$-cofinite. It follows that $\mathcal{M}(R,\mathcal{Z})_{cof}$ is an abelian subcategory of $\mathcal{M}(R)$.
\end{prf}

In Theorem \ref{2.5.12}, the assumption that $R$ is semilocal is somehow not desirable. Accordingly, we pose the following question.

\begin{question} \label{2.5.13}
Let $\mathcal{Z}$ be a stable under specialization subset of $\Spec(R)$ such that $\cd(\mathcal{Z},R) \leq 1$. Is $\mathcal{M}(R,\mathcal{Z})_{cof}$ an abelian subcategory of $\mathcal{M}(R)$?
\end{question}

Next we take care of the case $\dim(\mathcal{Z}) \leq 1$. Given a stable under specialization subset $\mathcal{Z}$ of $\Spec(R)$, we define the dimension of $\mathcal{Z}$ as
$$\dim(\mathcal{Z}):=\sup \left\{\dim(R/\mathfrak{a}) \suchthat \mathfrak{a}\in F(\mathcal{Z})\right\}.$$

\begin{lemma} \label{2.5.14}
Let $\mathcal{Z}$ be a subset of $\Max(R)$. Let $M$ be an $R$-module such that $\Supp_{R}(M)$ is a finite subset of $\mathcal{Z}$. If $\Hom_{R}(R/\mathfrak{a},M)$ is
finitely generated for every $\mathfrak{a} \in F(\mathcal{Z})$, then $M$ is artinian and $\mathcal{Z}$-cofinite.
\end{lemma}

\begin{prf}
For every $\mathfrak{a} \in F(\mathcal{Z})$, one has
$$\Supp_{R}\left(\Hom_{R}(R/\mathfrak{a},M)\right) \subseteq \Supp_{R}(M)\subseteq \Max(R).$$
Therefore, $\Hom_{R}(R/\mathfrak{a},M)$ has finite length for every $\mathfrak{a} \in F(\mathcal{Z})$. Hence, Lemma \ref{2.5.6} implies that $M=\Gamma_{\mathcal{Z}}(M)$ is artinian. Let $\mathfrak{a} \in F(\mathcal{Z})$. There are finitely many maximal ideals $\mathfrak{m}_{1},\ldots,\mathfrak{m}_{n}$ of $R$ such that $\mathfrak{m}_{1} \cap \cdots \cap \mathfrak{m}_{n} = \sqrt{\mathfrak{a}}$. Given any $i \geq 0$, $\Ext^{i}_{R}(R/\mathfrak{a},M)$ is an artinian $R$-module and there exists an integer $t \geq 0$ such that
$$(\mathfrak{m}_{1} \cap \cdots \cap \mathfrak{m}_{n})^{t} \Ext^{i}_{R}(R/\mathfrak{a},M)=0,$$
so $\Ext^{i}_{R}(R/\mathfrak{a},M)$ is finitely generated. Therefore, $M$ is $\mathcal{Z}$-cofinite.
\end{prf}

The following example shows that the finiteness assumption on $\Supp_{R}(M)$ in Lemma \ref{2.5.14} cannot be removed. This example also demonstrates that unlike $\mathfrak{a}$-cofinite modules, a $\mathcal{Z}$-cofinite module can have infinitely many associated prime ideals.

\begin{example} \label{2.5.15}
Let $R$ be a Gorenstein ring of dimension $d$ such that
$$\mathcal{Z}:= \{\mathfrak{m}\in \Max(R) \suchthat \Ht(\mathfrak{m})=d \}$$
is an infinite set. Then by \cite[Remark 2.12]{HD}, it turns out that
\[
    H_{\mathcal{Z}}^{i}(R)\cong
\begin{dcases}
    \bigoplus_{\mathfrak{m}\in \mathcal{Z}}E_{R}(R/\mathfrak{m}) & \text{if } i= d\\
    0              & \text{if } i \neq d.
\end{dcases}
\]
Thus $\Ext^{i}_{R}\left(R/\mathfrak{a},H^{j}_{\mathcal{Z}}(R)\right)$ is finitely generated for every $\mathfrak{a} \in F(\mathcal{Z})$ and every $i,j \geq 0$. It follows that $H^{j}_{\mathcal{Z}}(R)$ is $\mathcal{Z}$-cofinite for every $j \geq 0$, whereas $H^{d}_{\mathcal{Z}}(R)$ is not artinian and $\Ass_R\left(H^{d}_{\mathcal{Z}}(R)\right)$ is not finite.
\end{example}

\begin{lemma} \label{2.5.16}
Let $\mathcal{Z}$ be a stable under specialization subset of $\Spec(R)$ with $\dim(\mathcal{Z})\leq 1$, and $M$ an $R$-module such that $\Supp_{R}(M) \subseteq \mathcal{Z}$. Then the following assertions are equivalent:
\begin{enumerate}
\item[(i)] $M$ is $\mathcal{Z}$-cofnite.
\item[(ii)] $\Hom_{R}(R/\mathfrak{a},M)$ and $\Ext^{1}_{R}(R/\mathfrak{a},M)$ are finitely generated for every $\mathfrak{a} \in F(\mathcal{Z})$.
\end{enumerate}
\end{lemma}

\begin{prf}
(i) $\Rightarrow$ (ii): Clear.

(ii) $\Rightarrow$ (i): Since $\Supp_{R}(M) \subseteq \mathcal{Z}$, the hypothesis implies that $\dim(M) \leq 1$. On the other hand, for every $\mathfrak{a} \in F(\mathcal{Z})$, we have $\dim(R/\mathfrak{a}) \leq 1$. Now, \cite[Theorem 2.5]{BNS2} establishes the result.
\end{prf}

\begin{lemma} \label{2.5.17}
Let $\mathcal{Z}$ be a stable under specialization subset of $\Spec(R)$ with $\dim(\mathcal{Z}) \leq 1$. Let $M$ be an $R$-module such that $\Ext^{i}_{R}(R/\mathfrak{a},M)$ is a finitely generated $R$-module for every $\mathfrak{a} \in F(\mathcal{Z})$ and every $i \geq 0$. Then $H^{i}_{\mathcal{Z}}(M)$ is $\mathcal{Z}$-cofinite for every $i \geq 0$.
\end{lemma}

\begin{prf}
Clearly, $\Supp_{R}\left(H^{i}_{\mathcal{Z}}(M)\right) \subseteq \mathcal{Z}$ for every $i \geq 0$. Thus, it remains to show that $\Ext^{j}_{R}\left(R/\mathfrak{a},H^{i}_{\mathcal{Z}}(M)\right)$ is finitely generated for every $\mathfrak{a} \in F(\mathcal{Z})$ and every $i,j \geq 0$.

By induction on $i$, we show that $\Ext^{j}_{R}\left(R/\mathfrak{b},H^{i}_{\mathcal{Z}}(M)\right)$ is finitely generated for every $\mathfrak{b}\in F(\mathcal{Z})$ and every $j \geq 0$. The short exact sequence
$$0 \rightarrow \Gamma_{\mathcal{Z}}(M) \rightarrow M \rightarrow M/\Gamma_{\mathcal{Z}}(M) \rightarrow 0,$$
yields the exact sequence
$$0 \rightarrow \Hom_{R}\left(R/\mathfrak{b},\Gamma_{\mathcal{Z}}(M)\right) \rightarrow \Hom_{R}(R/\mathfrak{b},M) \rightarrow \Hom_{R}\left(R/\mathfrak{b},M/\Gamma_{\mathcal{Z}}(M)\right) \rightarrow $$$$ \Ext_{R}^{1}\left(R/\mathfrak{b},\Gamma_{\mathcal{Z}}(M)\right) \rightarrow \Ext_{R}^{1}(R/\mathfrak{b},M).$$
It can be seen by inspection that $\Hom_{R}\left(R/\mathfrak{b},M/\Gamma_{\mathcal{Z}}(M)\right)=0$ for every $\mathfrak{b}\in F(\mathcal{Z})$. Hence the above exact sequence shows that the $R$-modules $\Hom_{R}\left(R/\mathfrak{b},\Gamma_{\mathcal{Z}}(M)\right)$ and $\Ext_{R}^{1}\left(R/\mathfrak{b},\Gamma_{\mathcal{Z}}(M)\right)$ are finitely generated for every $\mathfrak{b}\in F(\mathcal{Z})$. Therefore, by Lemma \ref{2.5.16} the case $i=0$ holds true.

Now, suppose that $i \geq 1$ and make the obvious induction hypothesis. From the exact sequence
$$\Ext_{R}^{j}(R/\mathfrak{b},M) \rightarrow \Ext_{R}^{j}\left(R/\mathfrak{b},M/\Gamma_{\mathcal{Z}}(M)\right) \rightarrow \Ext_{R}^{j+1}\left(R/\mathfrak{b},\Gamma_{\mathcal{Z}}(M)\right),$$
using the base case $i=0$, we deduce that $\Ext_{R}^{j}\left(R/\mathfrak{b},M/\Gamma_{\mathcal{Z}}(M)\right)$ is finitely generated for every for every $\mathfrak{b}\in F(\mathcal{Z})$ and every $j \geq 0$. Since $H^{i}_{\mathcal{Z}}(M) \cong H^{i}_{\mathcal{Z}}\left(M/\Gamma_{\mathcal{Z}}(M)\right)$ for every $i \geq 1$, we may assume that $\Gamma_{\mathcal{Z}}(M)=0$. Let $E:= E_{R}(M)$ and $N:=E/M$. We have $\Gamma_{\mathcal{Z}}(E) \cong E_{R}\left(\Gamma_{\mathcal{Z}}(M)\right)=0$, and $\Hom_{R}(R/\mathfrak{b},E)=0$ for every $\mathfrak{b}\in F(\mathcal{Z})$. Then from the short exact sequence
$$0 \rightarrow M \rightarrow E \rightarrow N \rightarrow 0,$$
we conclude that $H^{k}_{\mathcal{Z}}(M) \cong H^{k-1}_{\mathcal{Z}}(N)$ and $\Ext^{k}_{R}(R/\mathfrak{b},M) \cong \Ext^{k-1}_{R}(R/\mathfrak{b},N)$ for every $\mathfrak{b}\in F(\mathcal{Z})$ and every $k \geq 1$. Hence the assumption is satisfied by $N$, and thus $H^{i-1}_{\mathcal{Z}}(N)$ is $\mathcal{Z}$-cofinite by the induction hypothesis, and we are through.
\end{prf}

\begin{theorem} \label{2.5.18}
Let $\mathcal{Z}$ be a stable under specialization subset of $\Spec (R)$ with $\dim(\mathcal{Z}) \leq 1$, and $M$ a finitely generated $R$-module. Then $H^{i}_{\mathcal{Z}}(M)$ is $\mathcal{Z}$-cofinite for every $i \geq 0$.
\end{theorem}

\begin{prf}
Immediate from Lemma \ref{2.5.17}.
\end{prf}

\begin{theorem} \label{2.5.19}
Let $\mathcal{Z}$ be a stable under specialization subset of $\Spec(R)$ with $\dim(\mathcal{Z}) \leq 1$. Then $\mathcal{M}(R,\mathcal{Z})_{cof}$ is an abelian subcategory of $\mathcal{M}(R)$.
\end{theorem}

\begin{prf}
Let $M$ and $N$ be two $\mathcal{Z}$-cofinite $R$-modules and let $f:M\rightarrow N$ be an $R$-homomorphism. Consider the short exact sequences
\begin{equation} \label{eq:2.5.19.1}
0 \rightarrow \ker f \rightarrow M \rightarrow \im f \rightarrow 0,
\end{equation}
and
\begin{equation} \label{eq:2.5.19.2}
0 \rightarrow \im f \rightarrow N \rightarrow \coker f \rightarrow 0.
\end{equation}
From the short exact sequence \eqref{eq:2.5.19.1}, we infer that $\Hom_{R}(R/\mathfrak{a},\ker f)$ is finitely generated for every $\mathfrak{a} \in F(\mathcal{Z})$. From the short exact sequence \eqref{eq:2.5.19.2}, we deduce that $\Hom_{R}(R/\mathfrak{a},\im f)$ is finitely generated for every $\mathfrak{a} \in F(\mathcal{Z})$. Thus from the short exact sequence \eqref{eq:2.5.19.1}, we conclude that $\Ext^{1}_{R}(R/\mathfrak{a},\ker f)$ is finitely generated. Now, Lemma \ref{2.5.16} implies that $\ker f$ is $\mathcal{Z}$-cofinite. It follows that $\coker f$ is $\mathcal{Z}$-cofinite as well, which means that $\mathcal{M}(R,\mathcal{Z})_{cof}$ is an abelian subcategory of $\mathcal{M}(R)$.
\end{prf}

We further investigate the case $\dim(R) \leq 2$.

Given a stable under specialization subset $\mathcal{Z}$ of $\Spec(R)$ and a finitely generated $R$-module $M$, we remind that
$$\depth_{R}(\mathcal{Z},M):= \inf \left\{\depth_{R}(\mathfrak{a},M) \suchthat \mathfrak{a} \in F(\mathcal{Z}) \right\}.$$

\begin{lemma} \label{2.5.20}
Let $\mathcal{Z}$ be a stable under specialization subset of $\Spec(R)$, and $M$ a finitely generated $R$-module. Then
$$\Supp_{R}(M)\cap \mathcal{Z}=\bigcup_{i=0}^{\infty} \Supp_{R}\left(H^{i}_{\mathcal{Z}}(M)\right).$$
\end{lemma}

\begin{prf}
Let $\mathfrak{p} \in \Supp_{R}(M)\cap \mathcal{Z}$. It is straightforward to see that the set
$$\mathcal{Z}_{\mathfrak{p}}:= \left\{\mathfrak{q}R_{\mathfrak{p}} \suchthat \mathfrak{q}\in \mathcal{Z} \text{ and } \mathfrak{q}\subseteq \mathfrak{p}\right\}$$
is a stable under specialization subset of $\Spec (R)_{\mathfrak{p}}$. It is clear that $\mathfrak{p}R_{\mathfrak{p}} \in \mathcal{Z}_{\mathfrak{p}} \cap \Supp_{R_{\mathfrak{p}}}\left(M_{\mathfrak{p}}\right)$, so $\depth_{R_{\mathfrak{p}}}\left(\mathfrak{p}R_{\mathfrak{p}},M_{\mathfrak{p}}\right) < \infty$, and thus $s:= \depth_{R_{\mathfrak{p}}}\left(\mathcal{Z}_{\mathfrak{p}},M_{\mathfrak{p}}\right) < \infty$. However by \cite[Proposition 5.5]{B},
$$\depth_{R_{\mathfrak{p}}}\left(\mathcal{Z}_{\mathfrak{p}},M_{\mathfrak{p}}\right)= \inf \left\{i \in \mathbb{Z} \suchthat H^{i}_{\mathcal{Z_{\mathfrak{p}}}}\left(M_{\mathfrak{p}}\right) \neq 0 \right\},$$
so $H^{s}_{\mathcal{Z_{\mathfrak{p}}}}\left(M_{\mathfrak{p}}\right) \neq 0$. One may check that $H^{s}_{\mathcal{Z}}(M)_{\mathfrak{p}} \cong H^{s}_{\mathcal{Z_{\mathfrak{p}}}}\left(M_{\mathfrak{p}}\right)$, and so $\mathfrak{p}\in \bigcup_{i=0}^{\infty} \Supp_{R}\left(H^{i}_{\mathcal{Z}}(M)\right)$. The reverse inclusion is immediate.
\end{prf}

\begin{lemma} \label{2.5.21}
Let $S$ be a module-finite $R$-algebra. Let $\mathfrak{a}$ be an ideal of $R$, and $M$ an $S$-module. Then the $R$-module $\Ext^{i}_{R}(R/\mathfrak{a},M)$ is
finitely generated for every $i\geq 0$ if and only if the $S$-module $\Ext^{i}_{S}(S/\mathfrak{a}S,M)$ is finitely generated for every $i\geq 0$.
\end{lemma}

\begin{prf}
The proof of \cite[Proposition 2]{DM} establishes the claim. Note that the assumption on the supports is not used in that proof.
\end{prf}

\begin{theorem} \label{2.5.22}
Let $\mathcal{Z}$ be a stable under specialization subset of $\Spec(R)$, and $M$ a finitely generated $R$-module. Suppose that $\dim\left(\Supp_{R}\left(H^{i}_{\mathcal{Z}}(M)\right)\right) \leq 1$ for every $i\geq 0$. Then $H^{i}_{\mathcal{Z}}(M)$ is $\mathcal{Z}$-cofinite for every $i\geq 0$.
\end{theorem}

\begin{prf}
By the assumption and Lemma \ref{2.5.20}, we have $\dim\left(\Supp_{R}(M)\cap \mathcal{Z}\right) \leq 1$. Set
$$\widetilde{\mathcal{Z}}:= \left\{\frac{\mathfrak{p}}{\ann_{R}(M)} \suchthat \mathfrak{p} \in \Supp_{R}(M)\cap \mathcal{Z} \right\},$$
and $S:= R/ \ann_{R}(M)$. Then, it is straightforward to see that $\widetilde{\mathcal{Z}}$ is a stable under specialization subset of
$\Spec(S)$ with $\dim\left(\widetilde{\mathcal{Z}}\right) \leq 1$, and
$F\left(\widetilde{\mathcal{Z}}\right)=\left\{\mathfrak{a}S \suchthat \mathfrak{a} \in F(\mathcal{Z}) \right\}$.
In addition, we have
\begin{equation*}
\begin{split}
H^{i}_{\mathcal{Z}}(M) & \cong \underset{\fa\in F(\mathcal{Z})}{\varinjlim}H_{\mathfrak{a}}^i(M) \\
 & \cong \underset{\fa\in F(\mathcal{Z})}{\varinjlim}H_{\mathfrak{a}S}^i(M) \\
 & \cong \underset{\fb\in F\left(\widetilde{\mathcal{Z}}\right)}{\varinjlim}H_{\mathfrak{b}}^i(M) \\
 & \cong H^{i}_{\widetilde{\mathcal{Z}}}(M)
\end{split}
\end{equation*}
for every $i \geq 0$.
Hence by Theorem \ref{2.5.18}, $N:= H^{i}_{\mathcal{Z}}(M) \cong H^{i}_{\widetilde{\mathcal{Z}}}(M)$ is a $\widetilde{\mathcal{Z}}$-cofinite $S$-module for every $i \geq 0$. Hence, the $S$-module $\Ext^{n}_{S}(S/\mathfrak{a}S,N)$ is finitely generated for every $\mathfrak{a}\in F(\mathcal{Z})$ and for every $n\geq 0$. Now, Lemma \ref{2.5.21} implies that the $R$-module $\Ext^{n}_{R}(R/\mathfrak{a},N)$ is finitely generated for every $\mathfrak{a}\in F(\mathcal{Z})$ and every $n\geq 0$, thereby, $H^{i}_{\mathcal{Z}}(M)$ is $\mathcal{Z}$-cofinite for every $i\geq 0$.
\end{prf}

\begin{lemma} \label{2.5.23}
Let $\mathcal{Z}$ be a stable under specialization subset of $\Spec(R)$, and $M$ a finitely generated $R$-module. Then $$\dim\left(\Supp_{R}\left(H^{i}_{\mathcal{Z}}(M)\right)\right) \leq \dim_{R}(M)-1$$
for every $i \geq 1$.
\end{lemma}

\begin{prf}
Clearly, we may assume that $\dim_{R}(M) < \infty$. Since $H_{\mathcal{Z}}^i(M) \cong \underset{\fa\in F(\mathcal{Z})}{\varinjlim}H_{\mathfrak{a}}^i(M)$ for every $i \geq 0$, it suffices to show that
$$\dim \left(\Supp_{R}\left(H_{\mathfrak{a}}^i(M)\right)\right) \leq \dim_{R}(M)-1$$
for every $\fa\in F(\mathcal{Z})$ and every $i \geq 1$. As $H^{i}_{\mathfrak{a}}(M) \cong H^{i}_{\mathfrak{a}}\left(M/\Gamma_{\mathfrak{a}}(M)\right)$ for every $i \geq 1$, we may assume that $\Gamma_{\mathfrak{a}}(M)=0$. Consequently, we conclude that $\mathfrak{a}$ contains a nonzerodivisor $r$ on $M$. As $\Supp_{R}\left(H^{i}_{\mathfrak{a}}(M)\right) \subseteq \Supp_{R}(M/\mathfrak{a}M)$, we have
\begin{equation*}
\begin{split}
\dim\left(\Supp_{R}\left(H^{i}_{\mathfrak{a}}(M)\right)\right) & \leq \dim_{R}(M/\mathfrak{a}M) \\
 & \leq \dim_{R}(M/rM) \\
 & \leq \dim_{R}(M)-1.
\end{split}
\end{equation*}
\end{prf}

\begin{corollary} \label{2.5.24}
Let $\mathcal{Z}$ be a stable under specialization subset of $\Spec(R)$, and $M$ a finitely generated $R$-module with $\dim_{R}(M) \leq 2$. Then $H^{i}_{\mathcal{Z}}(M)$ is $\mathcal{Z}$-cofinite for every $i \geq 0$.
\end{corollary}

\begin{prf}
Clearly, $\Gamma_{\mathcal{Z}}(M)$ is $\mathcal{Z}$-cofinite.  So by replacing $M$ with $M/\Gamma_{\mathcal{Z}}(M)$, we can assume that $\Gamma_{\mathcal{Z}}(M)=0$.
Now, Lemma \ref{2.5.23} implies that
$$\dim\left(\Supp_{R}\left(H^{i}_{\mathcal{Z}}(M)\right)\right) \leq 1$$
for every $i\geq 0$, thereby Theorem \ref{2.5.22} completes the argument.
\end{prf}

Similar to the previous section, we apply the technique of way-out functors to depart from  modules to complexes. The next result provides us with a suitable transition device from modules to complexes when dealing with cofiniteness.

\begin{theorem} \label{2.5.26}
Let $\mathcal{Z}$ a stable under specialization subset of $\Spec (R)$. Then the functor ${\bf R}\Gamma_{\mathcal{Z}}(-): \mathcal{D}(R) \rightarrow \mathcal{D}(R)$ is triangulated and way-out left. As a consequence, if $H^{i}_{\mathcal{Z}}(M)$ is $\mathcal{Z}$-cofinite for every finitely generated $R$-module $M$ and every $i \geq 0$, and $\mathcal{M}(R,\mathcal{Z})_{cof}$ is an abelian category, then $H^{i}_{\mathcal{Z}}(X)$ is $\mathcal{Z}$-cofinite for every $X \in \mathcal{D}^{f}_{\sqsubset}(R)$ and every $i \in \mathbb{Z}$.
\end{theorem}

\begin{prf}
It is folklore that if a functor on $\mathcal{D}(R)$ extends from a functor on $\mathcal{M}(R)$, then it commutes with mapping cones. Hence, it can be easily verified that the functor ${\bf R}\Gamma_{\mathcal{Z}}(-): \mathcal{D}(R) \rightarrow \mathcal{D}(R)$ is triangulated. Moreover, if $\sup X \leq n$, then there is a semi-injective resolution $X \xrightarrow{\simeq} I$ of $X$ such that $I_{i}=0$ for every $i \geq n$. Applying the functor $\Gamma_{\mathcal{Z}}(-)$ to $I$ and taking homology, we see that $\sup {\bf R}\Gamma_{\mathcal{Z}}(X) \leq n$. It follows that the functor ${\bf R}\Gamma_{\mathcal{Z}}(-)$ is way-out left. Now, let $\mathcal{A}$ be the subcategory of finitely generated $R$-modules, and let $\mathcal{B}:= \mathcal{M}(R,\mathcal{Z})_{cof}$. It can be easily seen that $\mathcal{B}$ is closed under extensions. It now follows from Lemma \ref{2.4.7} that $H^{i}_{\mathfrak{a}}(X) = H_{-i}\left({\bf R}\Gamma_{\mathfrak{a}}(X)\right) \in \mathcal{B}$ for every $X \in \mathcal{D}^{f}(R)$ and every $i \in \mathbb{Z}$.
\end{prf}

\begin{corollary} \label{2.5.27}
Let $\mathcal{Z}$ be a stable under specialization subset of $\Spec(R)$ and $X \in \mathcal{D}^{f}_{\sqsubset}(R)$. Suppose that either $R$ is semilocal and $\cd(\mathcal{Z},R) \leq 1$, or $\dim(\mathcal{Z}) \leq 1$. Then $H^{i}_{\mathcal{Z}}(X)$ is $\mathcal{Z}$-cofinite for every $i \in \mathcal{Z}$.
\end{corollary}

\begin{prf}
Clear in view of Theorems \ref{2.5.5}, \ref{2.5.12}, \ref{2.5.18}, \ref{2.5.19}, and \ref{2.5.26}.
\end{prf}

\begin{theorem} \label{2.5.28}
Let $\mathcal{Z}$ a stable under specialization subset of $\Spec(R)$, and $X \in \mathcal{D}(R)$ with
$\dim\left(\Supp_{R}(X)\right) \leq 2$. If $X \in \mathcal{D}_{\sqsubset}^{f}(R)$, then $H^{i}_{\mathcal{Z}}(X)$ is $\mathcal{Z}$-cofinite for every $i \in \mathbb{Z}$.
\end{theorem}

\begin{prf}
Let
$$\mathcal{A}:= \left\{M \in \mathcal{M}(R) \suchthat M \text{ is finitely generated and } \dim_{R}(M) \leq 2 \right\},$$
and
$$\mathcal{B}:= \left\{M \in \mathcal{M}(R) \suchthat M \text{ is } \mathcal{Z} \text{-cofinite and } \dim_{R}\left(\Supp_{R}(M)\right) \leq 1 \right\}.$$
By the argument of the proof of \cite[Theorem 2.7]{BNS1}, $\mathcal{B}$ is an abelian subcategory of $\mathcal{M}(R)$. In addition, it is closed under extensions. Now, by Lemma \ref{2.5.23} and Corollary \ref{2.5.24}, we have $H^{i}_{\mathcal{Z}}(M) \in \mathcal{B}$ for any $M \in \mathcal{A}$. Considering the triangulated way-out left functor ${\bf R}\Gamma_{\mathcal{Z}}(-)$, Lemma \ref{2.4.7} implies that $H^{i}_{\mathcal{Z}}(X)=H_{-i}\left({\bf R}\Gamma_{\mathcal{Z}}(X)\right) \in \mathcal{B}$ for every $X \in \mathcal{D}_{\sqsubset}(R)$ with $H_{i}(X) \in \mathcal{A}$ and for every $i \in \mathbb{Z}$.
\end{prf}

\begin{remark} \label{2.5.29}
Given a stable under specialization subset $\mathcal{Z}$ of $\Spec(R)$, the local cohomology module $H^{i}_{\mathcal{Z}}(X)$ of an $R$-complex $X$ with support in $\mathcal{Z}$, is an all-in-one generalization of the previously known generalized local cohomology modules as outlined in the following discussion.
\begin{enumerate}
\item[(i)] Let $\mathfrak{a}$ be an ideal of $R$, and $M$ and $N$ two $R$-modules. The generalized local cohomology module $H^{i}_{\mathfrak{a}}(M,N)$ is defined in \cite{He} as
    $$H^{i}_{\mathfrak{a}}(M,N):= \underset{n}{\varinjlim} \Ext^{i}_{R}(M/\mathfrak{a}^{n}M,N)$$
    for every $i\geq0$. It is shown in \cite{Ya1} that if $M$ is finitely generated, then we have $H^{i}_{\mathfrak{a}}(M,N)=H^{i}_{\mathcal{Z}}(X)$ for every $i \geq 0$, where $\mathcal{Z}= V(\mathfrak{a})$ and $X={\bf R}\Hom_{R}(M,N)$.
\item[(ii)] Let $\mathfrak{a}$ and $\mathfrak{b}$ be two ideals of $R$, and $M$ an $R$-module. Let
    $$W(\mathfrak{a},\mathfrak{b}):=\left\{\mathfrak{p}\in \Spec (R) \suchthat \mathfrak{a}^{n} \subseteq \mathfrak{p}+\mathfrak{b} \text{ for some integer } n\geq 1\right\}.$$
    Define a functor $\Gamma_{\mathfrak{a},\mathfrak{b}}(-)$ on $\mathcal{M}(R)$ by setting
    $$\Gamma_{\mathfrak{a},\mathfrak{b}}(M):= \left\{x\in M \suchthat \Supp_{R}(Rx) \subseteq W(\mathfrak{a},\mathfrak{b}) \right\},$$
    for an $R$-module $M$, and $\Gamma_{\mathfrak{a},\mathfrak{b}}(f):= f|_{\Gamma_{\mathfrak{a},\mathfrak{b}}(M)}$ for an $R$-homomorphism $f:M\rightarrow N$. The generalized local cohomology module $H^{i}_{\mathfrak{a},\mathfrak{b}}(M)$ is defined in \cite{TYY} to be
    $H^{i}_{\mathfrak{a},\mathfrak{b}}(M):= R^{i}\Gamma_{\mathfrak{a},\mathfrak{b}}(M)$ for every $i\geq0$.
    It is clear that $H^{i}_{\mathfrak{a},\mathfrak{b}}(M)=H^{i}_{\mathcal{Z}}(X)$ for every $i \geq 0$, where $\mathcal{Z}= W(\mathfrak{a},\mathfrak{b})$
    and $X=M$.
\item[(iii)] Let $\Phi$ be a directed poset. By a system of ideals $\varphi$, we mean a family $\varphi= \{\mathfrak{a}_{\alpha}\}_{\alpha\in\Phi}$ of ideals of $R$, such that $\mathfrak{a}_{\alpha}\subseteq \mathfrak{a}_{\beta}$ whenever $\alpha\geq \beta$, and for any $\alpha,\beta \in \Phi$, there is a $\gamma \in \Phi$ with $\mathfrak{a}_{\gamma} \subseteq \mathfrak{a}_{\alpha}\mathfrak{a}_{\beta}$. Given a system of ideals $\varphi$, define a functor $\Gamma_{\varphi}(-)$ on $\mathcal{M}(R)$ by setting
    $$\Gamma_{\varphi}(M):= \left\{x\in M \suchthat \mathfrak{a}x=0  \text{ for some } \mathfrak{a} \in \varphi \right\},$$
    for an $R$-module $M$, and $\Gamma_{\varphi}(f):= f|_{\Gamma_{\varphi}(M)}$ for an $R$-homomorphism $f:M\rightarrow N$.
    Then the generalized local cohomology module $H^{i}_{\varphi}(M)$ is defined in \cite[Notation 2.2.2]{BS} to be
    $H^{i}_{\varphi}(M):=R^{i}\Gamma_{\varphi}(M)$ for every $i\geq0$.
    It is easy to see that $H^{i}_{\varphi}(M)=H^{i}_{\mathcal{Z}}(X)$ for every $i \geq 0$, where $\mathcal{Z}=\bigcup_{\mathfrak{a}\in \varphi}V(\mathfrak{a})$ and $X=M$.
\item[(iv)]   Yoshino and Yoshizawa \cite[Theorem 2.10]{YY} have shown that for any abstract local cohomology functor $\delta:\mathcal{D}_{\sqsubset}(R)
\rightarrow \mathcal{D}_{\sqsubset}(R)$, there is a stable under specialization subset $\mathcal{Z}$ of $\Spec(R)$ such that $\delta\cong {\bf R}\Gamma_{\mathcal{Z}}(-)$.
\end{enumerate}
Accordingly, our results in this section generalize the following results:
\begin{enumerate}
\item[(a)] \cite[Proposition 3.6, Corollaries 3.9, 3.10, 3.11, and 3.12]{HV}.
\item[(b)] \cite[Theorem 2.5]{DH}.
\item[(c)] \cite[Theorem 1.3]{DS}.
\item[(d)] \cite[Theorems 1.1 and 1.3]{TGV}.
\item[(e)] \cite[Proposition 6.1, Corollary 6.3, Proposition 7.6, Corollary 7.7]{Ha1}.
\item[(f)] \cite[Theorem 2.1]{Ka1}.
\item[(g)] \cite[Corollary 2.8]{BNS1}.
\end{enumerate}
\end{remark}

\chapter{Greenlees-May Duality in a Nutshell}

\section{Introduction}

Throughout this chapter, all rings are assumed to be commutative and noetherian with identity.

In his algebraic geometry seminars of 1961-2, Grothendieck founded the theory of local cohomology as an indispensable tool in both algebraic geometry and commutative algebra. Given an ideal $\mathfrak{a}$ of $R$, the local cohomology functor $H^{i}_{\mathfrak{a}}(-)$ is defined as the $i$th right derived functor of the $\mathfrak{a}$-torsion functor $\Gamma_{\mathfrak{a}}(-) \cong \varinjlim \Hom_{R}(R/\mathfrak{a}^{n},-)$. Among a myriad of outstanding results, he proved the Local Duality Theorem.

\begin{theorem} \label{3.1.1}
Let $(R,\mathfrak{m})$ be a local ring with a dualizing module $\omega_{R}$, and $M$ a finitely generated $R$-module. Then
$$H^{i}_{\mathfrak{m}}(M) \cong \Ext^{\dim(R)-i}_{R}(M,\omega_{R})^{\vee}$$
for every $i \geq 0$.
\end{theorem}

The dual theory to local cohomology, i.e. local homology, was initiated by Matlis \cite{Mat2} in 1974, and its study was continued by Simon in \cite{Si1} and \cite{Si2}. Given an ideal $\mathfrak{a}$ of $R$, the local homology functor $H^{\mathfrak{a}}_{i}(-)$ is defined as the $i$th left derived functor of the $\mathfrak{a}$-adic completion functor $\Lambda^{\mathfrak{a}}(-) \cong \varprojlim (R/\mathfrak{a}^{n}\otimes_{R}-)$.

The existence of a dualizing module in Theorem \ref{3.1.1} is rather restrictive as it forces $R$ to be Cohen-Macaulay. To proceed further and generalize Theorem \ref{3.1.1}, Greenlees and May \cite[Propositions 3.1 and 3.8]{GM}, established a spectral sequence
\begin{equation} \label{eq:3.1.1.1}
E^{2}_{p,q}= \Ext^{-p}_{R}\left(H^{q}_{\mathfrak{a}}(R),M\right) \underset {p} \Rightarrow H^{\mathfrak{a}}_{p+q}(M)
\end{equation}
for any $R$-module $M$.
One can also settle the dual spectral sequence
\begin{equation} \label{eq:3.1.1.2}
E^{2}_{p,q}= \Tor^{R}_{p}\left(H^{q}_{\mathfrak{a}}(R),M\right) \underset {p} \Rightarrow H^{p+q}_{\mathfrak{a}}(M)
\end{equation}
for any $R$-module $M$.
It is by and large more palatable to have isomorphisms rather than spectral sequences. But the problem is that the category of $R$-modules $\mathcal{M}(R)$ is not rich enough to allow for the coveted isomorphisms. We need to enlarge this category to the category of $R$-complexes $\mathcal{C}(R)$, and even enrich it further, to the derived category $\mathcal{D}(R)$. The derived category $\mathcal{D}(R)$ is privileged with extreme maturity to accommodate the sought isomorphisms. As a matter of fact, the spectral sequence \eqref{eq:3.1.1.1} turns into the isomorphism
\begin{equation} \label{eq:3.1.1.3}
{\bf R}\Hom_{R}\left({\bf R}\Gamma_{\mathfrak{a}}(R),X\right)\simeq {\bf L}\Lambda^{\mathfrak{a}}(X),
\end{equation}
and the spectral sequence \eqref{eq:3.1.1.2} turns into the isomorphism
\begin{equation} \label{eq:3.1.1.4}
{\bf R}\Gamma_{\mathfrak{a}}(R)\otimes_{R}^{\bf L}X \simeq {\bf R}\Gamma_{\mathfrak{a}}(X)
\end{equation}
in $\mathcal{D}(R)$ for any $R$-complex $X$. Patching the two isomorphisms \eqref{eq:3.1.1.3} and \eqref{eq:3.1.1.4} together, we are blessed with the celebrated Greenlees-May Duality.
\begin{theorem} \label{3.1.2}
Let $\mathfrak{a}$ be an ideal of $R$, and $X,Y \in \mathcal{D}(R)$. Then there is a natural isomorphism
$${\bf R}\Hom_{R}\left({\bf R}\Gamma_{\mathfrak{a}}(X),Y\right) \simeq {\bf R}\Hom_{R}\left(X,{\bf L}\Lambda^{\mathfrak{a}}(Y)\right)$$
in $\mathcal{D}(R)$.
\end{theorem}
This was first proved by Alonso Tarr\'{i}o, Jerem\'{i}as L\'{o}pez and Lipman in \cite{AJL}. Theorem \ref{3.1.2} is a far-reaching generalization of Theorem \ref{3.1.1} and indeed extends it to its full generality. This theorem also demonstrates perfectly some sort of adjointness between derived local cohomology and homology.

Despite its incontrovertible impact on the theory of derived local homology and cohomology, we regretfully notice that there is no comprehensive and accessible treatment of the Greenlees-May Duality in the literature. There are some papers that touch on the subject, each from a different perspective, but none of them present a clear-cut and thorough proof that is fairly readable for non-experts; see for example \cite{GM}, \cite{AJL}, \cite{PSY1}, and \cite{Sc}. To remedy this defect, we commence on probing this theorem by providing the prerequisites from scratch and build upon a well-documented rigorous proof which is basically presented in layman's terms. In the course of our proof, some arguments are familiar while some others are novel. However, all the details are fully worked out so as to set forth a satisfactory exposition of the subject. We finally depict the highly non-trivial fact that the Greenlees-May Duality generalizes the Local Duality in simple and traceable steps.

\section{Module Prerequisites}

In this section, we embark on providing the requisite tools on modules which are to be recruited in Section 3.4.

First we recall the notion of a $\delta$-functor which will be used as a powerful tool to establish natural isomorphisms.

\begin{definition} \label{3.2.1}
Let $R$ and $S$ be two rings. Then:
\begin{enumerate}
\item[(i)] A \textit{homological covariant $\delta$-functor} \index{homological covariant $\delta$-functor} is a sequence $\left(\mathcal{F}_{i}:\mathcal{M}(R)\rightarrow \mathcal{M}(S)\right)_{i \geq 0}$ of additive covariant functors with the property that every short exact sequence
$$0 \rightarrow M' \rightarrow M \rightarrow M'' \rightarrow 0$$
of $R$-modules gives rise to a long exact sequence
$$\cdots \rightarrow \mathcal{F}_{2}(M'') \xrightarrow{\delta_{2}} \mathcal{F}_{1}(M') \rightarrow \mathcal{F}_{1}(M) \rightarrow \mathcal{F}_{1}(M'') \xrightarrow{\delta_{1}} \mathcal{F}_{0}(M') \rightarrow \mathcal{F}_{0}(M) \rightarrow \mathcal{F}_{0}(M'') \rightarrow 0$$
of $S$-modules, such that the connecting morphisms $\delta_{i}$'s are natural in the sense that any commutative diagram
\[\begin{tikzpicture}[every node/.style={midway},]
  \matrix[column sep={2.5em}, row sep={2.5em}]
  {\node(1) {$0$}; & \node(2) {$M'$}; & \node(3) {$M$}; & \node(4) {$M''$}; & \node(5) {$0$};\\
  \node(6) {$0$}; & \node(7) {$N'$}; & \node(8) {$N$}; & \node(9) {$N''$}; & \node(10) {$0$};\\};
  \draw[decoration={markings,mark=at position 1 with {\arrow[scale=1.5]{>}}},postaction={decorate},shorten >=0.5pt] (4) -- (9) node[anchor=east] {};
  \draw[decoration={markings,mark=at position 1 with {\arrow[scale=1.5]{>}}},postaction={decorate},shorten >=0.5pt] (3) -- (8) node[anchor=east] {};
  \draw[decoration={markings,mark=at position 1 with {\arrow[scale=1.5]{>}}},postaction={decorate},shorten >=0.5pt] (2) -- (7) node[anchor=east] {};
  \draw[decoration={markings,mark=at position 1 with {\arrow[scale=1.5]{>}}},postaction={decorate},shorten >=0.5pt] (1) -- (2) node[anchor=south] {};
  \draw[decoration={markings,mark=at position 1 with {\arrow[scale=1.5]{>}}},postaction={decorate},shorten >=0.5pt] (2) -- (3) node[anchor=south] {};
  \draw[decoration={markings,mark=at position 1 with {\arrow[scale=1.5]{>}}},postaction={decorate},shorten >=0.5pt] (3) -- (4) node[anchor=south] {};
  \draw[decoration={markings,mark=at position 1 with {\arrow[scale=1.5]{>}}},postaction={decorate},shorten >=0.5pt] (4) -- (5) node[anchor=south] {};
  \draw[decoration={markings,mark=at position 1 with {\arrow[scale=1.5]{>}}},postaction={decorate},shorten >=0.5pt] (6) -- (7) node[anchor=south] {};
  \draw[decoration={markings,mark=at position 1 with {\arrow[scale=1.5]{>}}},postaction={decorate},shorten >=0.5pt] (7) -- (8) node[anchor=south] {};
  \draw[decoration={markings,mark=at position 1 with {\arrow[scale=1.5]{>}}},postaction={decorate},shorten >=0.5pt] (8) -- (9) node[anchor=south] {};
  \draw[decoration={markings,mark=at position 1 with {\arrow[scale=1.5]{>}}},postaction={decorate},shorten >=0.5pt] (9) -- (10) node[anchor=south] {};
\end{tikzpicture}\]
of $R$-modules with exact rows induces a commutative diagram
\[\begin{tikzpicture}[every node/.style={midway},]
  \matrix[column sep={1em}, row sep={2.5em}]
  {\node(1) {$\cdots$}; & \node(2) {$\mathcal{F}_{2}(M'')$}; & \node(3) {$\mathcal{F}_{1}(M')$}; & \node(4) {$\mathcal{F}_{1}(M)$}; & \node(5) {$\mathcal{F}_{1}(M'')$}; & \node(6) {$\mathcal{F}_{0}(M')$}; & \node(7) {$\mathcal{F}_{0}(M)$}; & \node(8) {$\mathcal{F}_{0}(M'')$}; & \node(9) {$0$};\\
  \node(10) {$\cdots$}; & \node(11) {$\mathcal{F}_{2}(N'')$}; & \node(12) {$\mathcal{F}_{1}(N')$}; & \node(13) {$\mathcal{F}_{1}(N)$}; & \node(14) {$\mathcal{F}_{1}(N'')$}; & \node(15) {$\mathcal{F}_{0}(N')$}; & \node(16) {$\mathcal{F}_{0}(N)$}; & \node(17) {$\mathcal{F}_{0}(N'')$}; & \node(18) {$0$};\\};
  \draw[decoration={markings,mark=at position 1 with {\arrow[scale=1.5]{>}}},postaction={decorate},shorten >=0.5pt] (1) -- (2) node[anchor=east] {};
  \draw[decoration={markings,mark=at position 1 with {\arrow[scale=1.5]{>}}},postaction={decorate},shorten >=0.5pt] (2) -- (3) node[anchor=south] {$\delta_{2}$};
  \draw[decoration={markings,mark=at position 1 with {\arrow[scale=1.5]{>}}},postaction={decorate},shorten >=0.5pt] (3) -- (4) node[anchor=east] {};
  \draw[decoration={markings,mark=at position 1 with {\arrow[scale=1.5]{>}}},postaction={decorate},shorten >=0.5pt] (4) -- (5) node[anchor=east] {};
  \draw[decoration={markings,mark=at position 1 with {\arrow[scale=1.5]{>}}},postaction={decorate},shorten >=0.5pt] (5) -- (6) node[anchor=south] {$\delta_{1}$};
  \draw[decoration={markings,mark=at position 1 with {\arrow[scale=1.5]{>}}},postaction={decorate},shorten >=0.5pt] (6) -- (7) node[anchor=east] {};
  \draw[decoration={markings,mark=at position 1 with {\arrow[scale=1.5]{>}}},postaction={decorate},shorten >=0.5pt] (7) -- (8) node[anchor=east] {};
  \draw[decoration={markings,mark=at position 1 with {\arrow[scale=1.5]{>}}},postaction={decorate},shorten >=0.5pt] (8) -- (9) node[anchor=east] {};
  \draw[decoration={markings,mark=at position 1 with {\arrow[scale=1.5]{>}}},postaction={decorate},shorten >=0.5pt] (10) -- (11) node[anchor=east] {};
  \draw[decoration={markings,mark=at position 1 with {\arrow[scale=1.5]{>}}},postaction={decorate},shorten >=0.5pt] (11) -- (12) node[anchor=south] {$\Delta_{2}$};
  \draw[decoration={markings,mark=at position 1 with {\arrow[scale=1.5]{>}}},postaction={decorate},shorten >=0.5pt] (12) -- (13) node[anchor=east] {};
  \draw[decoration={markings,mark=at position 1 with {\arrow[scale=1.5]{>}}},postaction={decorate},shorten >=0.5pt] (13) -- (14) node[anchor=east] {};
  \draw[decoration={markings,mark=at position 1 with {\arrow[scale=1.5]{>}}},postaction={decorate},shorten >=0.5pt] (14) -- (15) node[anchor=south] {$\Delta_{1}$};
  \draw[decoration={markings,mark=at position 1 with {\arrow[scale=1.5]{>}}},postaction={decorate},shorten >=0.5pt] (15) -- (16) node[anchor=east] {};
  \draw[decoration={markings,mark=at position 1 with {\arrow[scale=1.5]{>}}},postaction={decorate},shorten >=0.5pt] (16) -- (17) node[anchor=east] {};
  \draw[decoration={markings,mark=at position 1 with {\arrow[scale=1.5]{>}}},postaction={decorate},shorten >=0.5pt] (17) -- (18) node[anchor=east] {};
  \draw[decoration={markings,mark=at position 1 with {\arrow[scale=1.5]{>}}},postaction={decorate},shorten >=0.5pt] (2) -- (11) node[anchor=east] {};
  \draw[decoration={markings,mark=at position 1 with {\arrow[scale=1.5]{>}}},postaction={decorate},shorten >=0.5pt] (3) -- (12) node[anchor=east] {};
  \draw[decoration={markings,mark=at position 1 with {\arrow[scale=1.5]{>}}},postaction={decorate},shorten >=0.5pt] (4) -- (13) node[anchor=east] {};
  \draw[decoration={markings,mark=at position 1 with {\arrow[scale=1.5]{>}}},postaction={decorate},shorten >=0.5pt] (5) -- (14) node[anchor=east] {};
  \draw[decoration={markings,mark=at position 1 with {\arrow[scale=1.5]{>}}},postaction={decorate},shorten >=0.5pt] (6) -- (15) node[anchor=east] {};
  \draw[decoration={markings,mark=at position 1 with {\arrow[scale=1.5]{>}}},postaction={decorate},shorten >=0.5pt] (7) -- (16) node[anchor=east] {};
  \draw[decoration={markings,mark=at position 1 with {\arrow[scale=1.5]{>}}},postaction={decorate},shorten >=0.5pt] (8) -- (17) node[anchor=east] {};
\end{tikzpicture}\]
of $S$-modules with exact rows.
\item[(ii)] A \textit{cohomological covariant $\delta$-functor} \index{cohomological covariant $\delta$-functor} is a sequence $\left(\mathcal{F}^{i}:\mathcal{M}(R)\rightarrow \mathcal{M}(S)\right)_{i \geq 0}$ of additive covariant functors with the property that every short exact sequence
$$0 \rightarrow M' \rightarrow M \rightarrow M'' \rightarrow 0$$
of $R$-modules gives rise to a long exact sequence
$$0 \rightarrow \mathcal{F}^{0}(M') \rightarrow \mathcal{F}^{0}(M) \rightarrow \mathcal{F}^{0}(M'') \xrightarrow{\delta^{0}} \mathcal{F}^{1}(M') \rightarrow \mathcal{F}^{1}(M) \rightarrow \mathcal{F}^{1}(M'') \xrightarrow{\delta^{1}} \mathcal{F}^{2}(M') \rightarrow \cdots$$
of $S$-modules, such that the connecting morphisms $\delta^{i}$'s are natural in the sense that any commutative diagram
\[\begin{tikzpicture}[every node/.style={midway},]
  \matrix[column sep={2.5em}, row sep={2.5em}]
  {\node(1) {$0$}; & \node(2) {$M'$}; & \node(3) {$M$}; & \node(4) {$M''$}; & \node(5) {$0$};\\
  \node(6) {$0$}; & \node(7) {$N'$}; & \node(8) {$N$}; & \node(9) {$N''$}; & \node(10) {$0$};\\};
  \draw[decoration={markings,mark=at position 1 with {\arrow[scale=1.5]{>}}},postaction={decorate},shorten >=0.5pt] (4) -- (9) node[anchor=east] {};
  \draw[decoration={markings,mark=at position 1 with {\arrow[scale=1.5]{>}}},postaction={decorate},shorten >=0.5pt] (3) -- (8) node[anchor=east] {};
  \draw[decoration={markings,mark=at position 1 with {\arrow[scale=1.5]{>}}},postaction={decorate},shorten >=0.5pt] (2) -- (7) node[anchor=east] {};
  \draw[decoration={markings,mark=at position 1 with {\arrow[scale=1.5]{>}}},postaction={decorate},shorten >=0.5pt] (1) -- (2) node[anchor=south] {};
  \draw[decoration={markings,mark=at position 1 with {\arrow[scale=1.5]{>}}},postaction={decorate},shorten >=0.5pt] (2) -- (3) node[anchor=south] {};
  \draw[decoration={markings,mark=at position 1 with {\arrow[scale=1.5]{>}}},postaction={decorate},shorten >=0.5pt] (3) -- (4) node[anchor=south] {};
  \draw[decoration={markings,mark=at position 1 with {\arrow[scale=1.5]{>}}},postaction={decorate},shorten >=0.5pt] (4) -- (5) node[anchor=south] {};
  \draw[decoration={markings,mark=at position 1 with {\arrow[scale=1.5]{>}}},postaction={decorate},shorten >=0.5pt] (6) -- (7) node[anchor=south] {};
  \draw[decoration={markings,mark=at position 1 with {\arrow[scale=1.5]{>}}},postaction={decorate},shorten >=0.5pt] (7) -- (8) node[anchor=south] {};
  \draw[decoration={markings,mark=at position 1 with {\arrow[scale=1.5]{>}}},postaction={decorate},shorten >=0.5pt] (8) -- (9) node[anchor=south] {};
  \draw[decoration={markings,mark=at position 1 with {\arrow[scale=1.5]{>}}},postaction={decorate},shorten >=0.5pt] (9) -- (10) node[anchor=south] {};
\end{tikzpicture}\]
of $R$-modules with exact rows induces a commutative diagram
\[\begin{tikzpicture}[every node/.style={midway},]
  \matrix[column sep={1em}, row sep={2.5em}]
  {\node(1) {$0$}; & \node(2) {$\mathcal{F}^{0}(M')$}; & \node(3) {$\mathcal{F}^{0}(M)$}; & \node(4) {$\mathcal{F}^{0}(M'')$}; & \node(5) {$\mathcal{F}^{1}(M')$}; & \node(6) {$\mathcal{F}^{1}(M)$}; & \node(7) {$\mathcal{F}^{1}(M'')$}; & \node(8) {$\mathcal{F}^{2}(M')$}; & \node(9) {$\cdots$};\\
  \node(10) {$0$}; & \node(11) {$\mathcal{F}^{0}(N')$}; & \node(12) {$\mathcal{F}^{0}(N)$}; & \node(13) {$\mathcal{F}^{0}(N'')$}; & \node(14) {$\mathcal{F}^{1}(N')$}; & \node(15) {$\mathcal{F}^{1}(N)$}; & \node(16) {$\mathcal{F}^{1}(N'')$}; & \node(17) {$\mathcal{F}^{2}(N')$}; & \node(18) {$\cdots$};\\};
  \draw[decoration={markings,mark=at position 1 with {\arrow[scale=1.5]{>}}},postaction={decorate},shorten >=0.5pt] (1) -- (2) node[anchor=east] {};
  \draw[decoration={markings,mark=at position 1 with {\arrow[scale=1.5]{>}}},postaction={decorate},shorten >=0.5pt] (2) -- (3) node[anchor=south] {};
  \draw[decoration={markings,mark=at position 1 with {\arrow[scale=1.5]{>}}},postaction={decorate},shorten >=0.5pt] (3) -- (4) node[anchor=east] {};
  \draw[decoration={markings,mark=at position 1 with {\arrow[scale=1.5]{>}}},postaction={decorate},shorten >=0.5pt] (4) -- (5) node[anchor=south] {$\delta^{0}$};
  \draw[decoration={markings,mark=at position 1 with {\arrow[scale=1.5]{>}}},postaction={decorate},shorten >=0.5pt] (5) -- (6) node[anchor=south] {};
  \draw[decoration={markings,mark=at position 1 with {\arrow[scale=1.5]{>}}},postaction={decorate},shorten >=0.5pt] (6) -- (7) node[anchor=east] {};
  \draw[decoration={markings,mark=at position 1 with {\arrow[scale=1.5]{>}}},postaction={decorate},shorten >=0.5pt] (7) -- (8) node[anchor=south] {$\delta^{1}$};
  \draw[decoration={markings,mark=at position 1 with {\arrow[scale=1.5]{>}}},postaction={decorate},shorten >=0.5pt] (8) -- (9) node[anchor=east] {};
  \draw[decoration={markings,mark=at position 1 with {\arrow[scale=1.5]{>}}},postaction={decorate},shorten >=0.5pt] (10) -- (11) node[anchor=east] {};
  \draw[decoration={markings,mark=at position 1 with {\arrow[scale=1.5]{>}}},postaction={decorate},shorten >=0.5pt] (11) -- (12) node[anchor=south] {};
  \draw[decoration={markings,mark=at position 1 with {\arrow[scale=1.5]{>}}},postaction={decorate},shorten >=0.5pt] (12) -- (13) node[anchor=east] {};
  \draw[decoration={markings,mark=at position 1 with {\arrow[scale=1.5]{>}}},postaction={decorate},shorten >=0.5pt] (13) -- (14) node[anchor=south] {$\Delta^{0}$};
  \draw[decoration={markings,mark=at position 1 with {\arrow[scale=1.5]{>}}},postaction={decorate},shorten >=0.5pt] (14) -- (15) node[anchor=south] {};
  \draw[decoration={markings,mark=at position 1 with {\arrow[scale=1.5]{>}}},postaction={decorate},shorten >=0.5pt] (15) -- (16) node[anchor=east] {};
  \draw[decoration={markings,mark=at position 1 with {\arrow[scale=1.5]{>}}},postaction={decorate},shorten >=0.5pt] (16) -- (17) node[anchor=south] {$\Delta^{1}$};
  \draw[decoration={markings,mark=at position 1 with {\arrow[scale=1.5]{>}}},postaction={decorate},shorten >=0.5pt] (17) -- (18) node[anchor=east] {};
  \draw[decoration={markings,mark=at position 1 with {\arrow[scale=1.5]{>}}},postaction={decorate},shorten >=0.5pt] (2) -- (11) node[anchor=east] {};
  \draw[decoration={markings,mark=at position 1 with {\arrow[scale=1.5]{>}}},postaction={decorate},shorten >=0.5pt] (3) -- (12) node[anchor=east] {};
  \draw[decoration={markings,mark=at position 1 with {\arrow[scale=1.5]{>}}},postaction={decorate},shorten >=0.5pt] (4) -- (13) node[anchor=east] {};
  \draw[decoration={markings,mark=at position 1 with {\arrow[scale=1.5]{>}}},postaction={decorate},shorten >=0.5pt] (5) -- (14) node[anchor=east] {};
  \draw[decoration={markings,mark=at position 1 with {\arrow[scale=1.5]{>}}},postaction={decorate},shorten >=0.5pt] (6) -- (15) node[anchor=east] {};
  \draw[decoration={markings,mark=at position 1 with {\arrow[scale=1.5]{>}}},postaction={decorate},shorten >=0.5pt] (7) -- (16) node[anchor=east] {};
  \draw[decoration={markings,mark=at position 1 with {\arrow[scale=1.5]{>}}},postaction={decorate},shorten >=0.5pt] (8) -- (17) node[anchor=east] {};
\end{tikzpicture}\]
of $S$-modules with exact rows.
\end{enumerate}
\end{definition}

\begin{example} \label{3.2.2}
Let $R$ and $S$ be two rings, and $\mathcal{F}:\mathcal{M}(R) \rightarrow \mathcal{M}(S)$ an additive covariant functor. Then $\left(L_{i}\mathcal{F}:\mathcal{M}(R) \rightarrow \mathcal{M}(S)\right)_{i \geq 0}$ is a homological covariant $\delta$-functor and $\left(R^{i}\mathcal{F}:\mathcal{M}(R) \rightarrow \mathcal{M}(S)\right)_{i \geq 0}$ is a cohomological covariant $\delta$-functor.
\end{example}

\begin{definition} \label{3.2.3}
Let $R$ and $S$ be two rings. Then:
\begin{enumerate}
\item[(i)] A \textit{morphism} \index{morphism of homological covariant $\delta$-functors}
$$\tau: \left(\mathcal{F}_{i}:\mathcal{M}(R) \rightarrow \mathcal{M}(S)\right)_{i \geq 0} \rightarrow \left(\mathcal{G}_{i}:\mathcal{M}(R) \rightarrow \mathcal{M}(S)\right)_{i \geq 0}$$
of homological covariant $\delta$-functors is a sequence $\tau= \left(\tau_{i}:\mathcal{F}_{i} \rightarrow \mathcal{G}_{i}\right)_{i \geq 0}$ of natural transformations of functors, such that any short exact sequence
$$0 \rightarrow M' \rightarrow M \rightarrow M'' \rightarrow 0$$
of $R$-modules induces a commutative diagram
\[\begin{tikzpicture}[every node/.style={midway},]
  \matrix[column sep={1em}, row sep={2.5em}]
  {\node(1) {$\cdots$}; & \node(2) {$\mathcal{F}_{2}(M'')$}; & \node(3) {$\mathcal{F}_{1}(M')$}; & \node(4) {$\mathcal{F}_{1}(M)$}; & \node(5) {$\mathcal{F}_{1}(M'')$}; & \node(6) {$\mathcal{F}_{0}(M')$}; & \node(7) {$\mathcal{F}_{0}(M)$}; & \node(8) {$\mathcal{F}_{0}(M'')$}; & \node(9) {$0$};\\
  \node(10) {$\cdots$}; & \node(11) {$\mathcal{G}_{2}(M'')$}; & \node(12) {$\mathcal{G}_{1}(M')$}; & \node(13) {$\mathcal{G}_{1}(M)$}; & \node(14) {$\mathcal{G}_{1}(M'')$}; & \node(15) {$\mathcal{G}_{0}(M')$}; & \node(16) {$\mathcal{G}_{0}(M)$}; & \node(17) {$\mathcal{G}_{0}(M'')$}; & \node(18) {$0$};\\};
  \draw[decoration={markings,mark=at position 1 with {\arrow[scale=1.5]{>}}},postaction={decorate},shorten >=0.5pt] (1) -- (2) node[anchor=east] {};
  \draw[decoration={markings,mark=at position 1 with {\arrow[scale=1.5]{>}}},postaction={decorate},shorten >=0.5pt] (2) -- (3) node[anchor=south] {$\delta_{2}$};
  \draw[decoration={markings,mark=at position 1 with {\arrow[scale=1.5]{>}}},postaction={decorate},shorten >=0.5pt] (3) -- (4) node[anchor=east] {};
  \draw[decoration={markings,mark=at position 1 with {\arrow[scale=1.5]{>}}},postaction={decorate},shorten >=0.5pt] (4) -- (5) node[anchor=east] {};
  \draw[decoration={markings,mark=at position 1 with {\arrow[scale=1.5]{>}}},postaction={decorate},shorten >=0.5pt] (5) -- (6) node[anchor=south] {$\delta_{1}$};
  \draw[decoration={markings,mark=at position 1 with {\arrow[scale=1.5]{>}}},postaction={decorate},shorten >=0.5pt] (6) -- (7) node[anchor=east] {};
  \draw[decoration={markings,mark=at position 1 with {\arrow[scale=1.5]{>}}},postaction={decorate},shorten >=0.5pt] (7) -- (8) node[anchor=east] {};
  \draw[decoration={markings,mark=at position 1 with {\arrow[scale=1.5]{>}}},postaction={decorate},shorten >=0.5pt] (8) -- (9) node[anchor=east] {};
  \draw[decoration={markings,mark=at position 1 with {\arrow[scale=1.5]{>}}},postaction={decorate},shorten >=0.5pt] (10) -- (11) node[anchor=east] {};
  \draw[decoration={markings,mark=at position 1 with {\arrow[scale=1.5]{>}}},postaction={decorate},shorten >=0.5pt] (11) -- (12) node[anchor=south] {$\Delta_{2}$};
  \draw[decoration={markings,mark=at position 1 with {\arrow[scale=1.5]{>}}},postaction={decorate},shorten >=0.5pt] (12) -- (13) node[anchor=east] {};
  \draw[decoration={markings,mark=at position 1 with {\arrow[scale=1.5]{>}}},postaction={decorate},shorten >=0.5pt] (13) -- (14) node[anchor=east] {};
  \draw[decoration={markings,mark=at position 1 with {\arrow[scale=1.5]{>}}},postaction={decorate},shorten >=0.5pt] (14) -- (15) node[anchor=south] {$\Delta_{1}$};
  \draw[decoration={markings,mark=at position 1 with {\arrow[scale=1.5]{>}}},postaction={decorate},shorten >=0.5pt] (15) -- (16) node[anchor=east] {};
  \draw[decoration={markings,mark=at position 1 with {\arrow[scale=1.5]{>}}},postaction={decorate},shorten >=0.5pt] (16) -- (17) node[anchor=east] {};
  \draw[decoration={markings,mark=at position 1 with {\arrow[scale=1.5]{>}}},postaction={decorate},shorten >=0.5pt] (17) -- (18) node[anchor=east] {};
  \draw[decoration={markings,mark=at position 1 with {\arrow[scale=1.5]{>}}},postaction={decorate},shorten >=0.5pt] (2) -- (11) node[anchor=west] {$\tau_{2}(M'')$};
  \draw[decoration={markings,mark=at position 1 with {\arrow[scale=1.5]{>}}},postaction={decorate},shorten >=0.5pt] (3) -- (12) node[anchor=west] {$\tau_{1}(M')$};
  \draw[decoration={markings,mark=at position 1 with {\arrow[scale=1.5]{>}}},postaction={decorate},shorten >=0.5pt] (4) -- (13) node[anchor=west] {$\tau_{1}(M)$};
  \draw[decoration={markings,mark=at position 1 with {\arrow[scale=1.5]{>}}},postaction={decorate},shorten >=0.5pt] (5) -- (14) node[anchor=west] {$\tau_{1}(M'')$};
  \draw[decoration={markings,mark=at position 1 with {\arrow[scale=1.5]{>}}},postaction={decorate},shorten >=0.5pt] (6) -- (15) node[anchor=west] {$\tau_{0}(M')$};
  \draw[decoration={markings,mark=at position 1 with {\arrow[scale=1.5]{>}}},postaction={decorate},shorten >=0.5pt] (7) -- (16) node[anchor=west] {$\tau_{0}(M)$};
  \draw[decoration={markings,mark=at position 1 with {\arrow[scale=1.5]{>}}},postaction={decorate},shorten >=0.5pt] (8) -- (17) node[anchor=west] {$\tau_{0}(M'')$};
\end{tikzpicture}\]
of $S$-modules with exact rows. If in particular, $\tau_{i}$ is an isomorphism for every $i \geq 0$, then $\tau$ is called an \textit{isomorphism} \index{isomorphism of homological covariant $\delta$-functors} of $\delta$-functors.
\item[(ii)] A \textit{morphism} \index{morphism of cohomological covariant $\delta$-functors}
$$\tau: \left(\mathcal{F}^{i}:\mathcal{M}(R) \rightarrow \mathcal{M}(S)\right)_{i \geq 0} \rightarrow \left(\mathcal{G}^{i}:\mathcal{M}(R) \rightarrow \mathcal{M}(S)\right)_{i \geq 0}$$
of cohomological covariant $\delta$-functors is a sequence $\tau= \left(\tau^{i}:\mathcal{F}^{i} \rightarrow \mathcal{G}^{i}\right)_{i \geq 0}$ of natural transformations of functors, such that any short exact sequence
$$0 \rightarrow M' \rightarrow M \rightarrow M'' \rightarrow 0$$
of $R$-modules induces a commutative diagram
\[\begin{tikzpicture}[every node/.style={midway},]
  \matrix[column sep={1em}, row sep={2.5em}]
  {\node(1) {$0$}; & \node(2) {$\mathcal{F}^{0}(M')$}; & \node(3) {$\mathcal{F}^{0}(M)$}; & \node(4) {$\mathcal{F}^{0}(M'')$}; & \node(5) {$\mathcal{F}^{1}(M')$}; & \node(6) {$\mathcal{F}^{1}(M)$}; & \node(7) {$\mathcal{F}^{1}(M'')$}; & \node(8) {$\mathcal{F}^{2}(M')$}; & \node(9) {$\cdots$};\\
  \node(10) {$0$}; & \node(11) {$\mathcal{G}^{0}(M')$}; & \node(12) {$\mathcal{G}^{0}(M)$}; & \node(13) {$\mathcal{G}^{0}(M'')$}; & \node(14) {$\mathcal{G}^{1}(M')$}; & \node(15) {$\mathcal{G}^{1}(M)$}; & \node(16) {$\mathcal{G}^{1}(M'')$}; & \node(17) {$\mathcal{G}^{2}(M')$}; & \node(18) {$\cdots$};\\};
  \draw[decoration={markings,mark=at position 1 with {\arrow[scale=1.5]{>}}},postaction={decorate},shorten >=0.5pt] (1) -- (2) node[anchor=east] {};
  \draw[decoration={markings,mark=at position 1 with {\arrow[scale=1.5]{>}}},postaction={decorate},shorten >=0.5pt] (2) -- (3) node[anchor=south] {};
  \draw[decoration={markings,mark=at position 1 with {\arrow[scale=1.5]{>}}},postaction={decorate},shorten >=0.5pt] (3) -- (4) node[anchor=east] {};
  \draw[decoration={markings,mark=at position 1 with {\arrow[scale=1.5]{>}}},postaction={decorate},shorten >=0.5pt] (4) -- (5) node[anchor=south] {$\delta^{0}$};
  \draw[decoration={markings,mark=at position 1 with {\arrow[scale=1.5]{>}}},postaction={decorate},shorten >=0.5pt] (5) -- (6) node[anchor=south] {};
  \draw[decoration={markings,mark=at position 1 with {\arrow[scale=1.5]{>}}},postaction={decorate},shorten >=0.5pt] (6) -- (7) node[anchor=east] {};
  \draw[decoration={markings,mark=at position 1 with {\arrow[scale=1.5]{>}}},postaction={decorate},shorten >=0.5pt] (7) -- (8) node[anchor=south] {$\delta^{1}$};
  \draw[decoration={markings,mark=at position 1 with {\arrow[scale=1.5]{>}}},postaction={decorate},shorten >=0.5pt] (8) -- (9) node[anchor=east] {};
  \draw[decoration={markings,mark=at position 1 with {\arrow[scale=1.5]{>}}},postaction={decorate},shorten >=0.5pt] (10) -- (11) node[anchor=east] {};
  \draw[decoration={markings,mark=at position 1 with {\arrow[scale=1.5]{>}}},postaction={decorate},shorten >=0.5pt] (11) -- (12) node[anchor=south] {};
  \draw[decoration={markings,mark=at position 1 with {\arrow[scale=1.5]{>}}},postaction={decorate},shorten >=0.5pt] (12) -- (13) node[anchor=east] {};
  \draw[decoration={markings,mark=at position 1 with {\arrow[scale=1.5]{>}}},postaction={decorate},shorten >=0.5pt] (13) -- (14) node[anchor=south] {$\Delta^{0}$};
  \draw[decoration={markings,mark=at position 1 with {\arrow[scale=1.5]{>}}},postaction={decorate},shorten >=0.5pt] (14) -- (15) node[anchor=south] {};
  \draw[decoration={markings,mark=at position 1 with {\arrow[scale=1.5]{>}}},postaction={decorate},shorten >=0.5pt] (15) -- (16) node[anchor=east] {};
  \draw[decoration={markings,mark=at position 1 with {\arrow[scale=1.5]{>}}},postaction={decorate},shorten >=0.5pt] (16) -- (17) node[anchor=south] {$\Delta^{1}$};
  \draw[decoration={markings,mark=at position 1 with {\arrow[scale=1.5]{>}}},postaction={decorate},shorten >=0.5pt] (17) -- (18) node[anchor=east] {};
  \draw[decoration={markings,mark=at position 1 with {\arrow[scale=1.5]{>}}},postaction={decorate},shorten >=0.5pt] (2) -- (11) node[anchor=west] {$\tau^{0}(M')$};
  \draw[decoration={markings,mark=at position 1 with {\arrow[scale=1.5]{>}}},postaction={decorate},shorten >=0.5pt] (3) -- (12) node[anchor=west] {$\tau^{0}(M)$};
  \draw[decoration={markings,mark=at position 1 with {\arrow[scale=1.5]{>}}},postaction={decorate},shorten >=0.5pt] (4) -- (13) node[anchor=west] {$\tau^{0}(M'')$};
  \draw[decoration={markings,mark=at position 1 with {\arrow[scale=1.5]{>}}},postaction={decorate},shorten >=0.5pt] (5) -- (14) node[anchor=west] {$\tau^{1}(M')$};
  \draw[decoration={markings,mark=at position 1 with {\arrow[scale=1.5]{>}}},postaction={decorate},shorten >=0.5pt] (6) -- (15) node[anchor=west] {$\tau^{1}(M)$};
  \draw[decoration={markings,mark=at position 1 with {\arrow[scale=1.5]{>}}},postaction={decorate},shorten >=0.5pt] (7) -- (16) node[anchor=west] {$\tau^{1}(M'')$};
  \draw[decoration={markings,mark=at position 1 with {\arrow[scale=1.5]{>}}},postaction={decorate},shorten >=0.5pt] (8) -- (17) node[anchor=west] {$\tau^{2}(M')$};
\end{tikzpicture}\]
of $S$-modules with exact rows. If in particular, $\tau^{i}$ is an isomorphism for every $i \geq 0$, then $\tau$ is called an \textit{isomorphism} \index{isomorphism of homological covariant $\delta$-functors} of $\delta$-functors.
\end{enumerate}
\end{definition}

The following remarkable theorem due to Grothendieck provides hands-on conditions that ascertain the existence of isomorphisms between $\delta$-functors.

\begin{theorem} \label{3.2.4}
Let $R$ and $S$ be two rings. Then the following assertions hold:
\begin{enumerate}
\item[(i)] Assume that $\left(\mathcal{F}_{i}:\mathcal{M}(R) \rightarrow \mathcal{M}(S)\right)_{i \geq 0}$ and $\left(\mathcal{G}_{i}:\mathcal{M}(R) \rightarrow \mathcal{M}(S)\right)_{i \geq 0}$ are two homological covariant $\delta$-functors such that $\mathcal{F}_{i}(F)=0=\mathcal{G}_{i}(F)$ for every free $R$-module $F$ and every $i \geq 1$. If there is a natural transformation $\eta: \mathcal{F}_{0} \rightarrow \mathcal{G}_{0}$ of functors which is an isomorphism on free $R$-modules, then there is a unique isomorphism $\tau:(\mathcal{F}_{i})_{i \geq 0} \rightarrow (\mathcal{G}_{i})_{i \geq 0}$ of $\delta$-functors such that $\tau_{0}=\eta$.
\item[(ii)] Assume that $\left(\mathcal{F}^{i}:\mathcal{M}(R) \rightarrow \mathcal{M}(S)\right)_{i \geq 0}$ and $\left(\mathcal{G}^{i}:\mathcal{M}(R) \rightarrow \mathcal{M}(S)\right)_{i \geq 0}$ are two cohomological covariant $\delta$-functors such that $\mathcal{F}^{i}(I)=0=\mathcal{G}^{i}(I)$ for every injective $R$-module $I$ and every $i \geq 1$. If there is a natural transformation $\eta: \mathcal{F}^{0} \rightarrow \mathcal{G}^{0}$ of functors which is an isomorphism on injective $R$-modules, then there is a unique isomorphism $\tau:(\mathcal{F}^{i})_{i \geq 0} \rightarrow (\mathcal{G}^{i})_{i \geq 0}$ of $\delta$-functors such that $\tau^{0}=\eta$.
\end{enumerate}
\end{theorem}

\begin{prf}
The proof is standard and can be found in almost every book on homological algebra. For example, see \cite[Corollaries 6.34 and 6.49]{Ro}. One should note that the above version is somewhat stronger than what is normally recorded in the books. However, the same proof can be modified in a suitable way to imply the above version.
\end{prf}

The following corollary sets forth a special case of Theorem \ref{3.2.4} which frequently occurs in practice.

\begin{corollary} \label{3.2.5}
Let $R$ and $S$ be two rings. Then the following assertions hold:
\begin{enumerate}
\item[(i)] Assume that $\mathcal{F}:\mathcal{M}(R) \rightarrow \mathcal{M}(S)$ is an additive covariant functor, and $\left(\mathcal{F}_{i}:\mathcal{M}(R) \rightarrow \mathcal{M}(S)\right)_{i \geq 0}$ is a homological covariant $\delta$-functor such that $\mathcal{F}_{i}(F)=0$ for every free $R$-module $F$ and every $i \geq 1$. If there is a natural transformation $\eta: L_{0}\mathcal{F} \rightarrow \mathcal{F}_{0}$ of functors which is an isomorphism on free $R$-modules, then there is a unique isomorphism $\tau:(L_{i}\mathcal{F})_{i \geq 0} \rightarrow (\mathcal{F}_{i})_{i \geq 0}$ of $\delta$-functors such that $\tau_{0}=\eta$.
\item[(ii)] Assume that $\mathcal{F}:\mathcal{M}(R) \rightarrow \mathcal{M}(S)$ is an additive covariant functor, and $\left(\mathcal{F}^{i}:\mathcal{M}(R) \rightarrow \mathcal{M}(S)\right)_{i \geq 0}$ is a cohomological covariant $\delta$-functor such that $\mathcal{F}^{i}(I)=0$ for every injective $R$-module $I$ and every $i \geq 1$. If there is a natural transformation $\eta: \mathcal{F}^{0} \rightarrow R^{0}\mathcal{F}$ of functors which is an isomorphism on injective $R$-modules, then there is a unique isomorphism $\tau: (\mathcal{F}^{i})_{i \geq 0} \rightarrow (R^{i}\mathcal{F})_{i \geq 0}$ of $\delta$-functors such that $\tau^{0}=\eta$.
\end{enumerate}
\end{corollary}

\begin{prf}
(i): We note that $(L_{i}\mathcal{F})(F)=0$ for every $i \geq 1$ and every free $R$-module $F$. Now the result follows from Theorem \ref{3.2.4} (i).

(ii): We note that $(R^{i}\mathcal{F})(I)=0$ for every $i \geq 1$ and every injective $R$-module $I$. Now the result follows from Theorem \ref{3.2.4} (ii).
\end{prf}

We next need direct and inverse systems of Koszul complexes and Koszul homologies.

\begin{remark} \label{3.2.6}
We have:
\begin{enumerate}
\item[(i)] Given an element $a \in R$, we define a morphism $\varphi^{k,l}_{a}:K^{R}(a^{k}) \rightarrow K^{R}(a^{l})$ of $R$-complexes for every $k \leq l$ as follows:
\[\begin{tikzpicture}[every node/.style={midway},]
  \matrix[column sep={3em}, row sep={3em}]
  {\node(1) {$0$}; & \node(2) {$R$}; & \node(3) {$R$}; & \node(4) {$0$};\\
  \node(5) {$0$}; & \node(6) {$R$}; & \node(7) {$R$}; & \node(8) {$0$};\\};
  \draw[decoration={markings,mark=at position 1 with {\arrow[scale=1.5]{>}}},postaction={decorate},shorten >=0.5pt] (1) -- (2) node[anchor=east] {};
  \draw[decoration={markings,mark=at position 1 with {\arrow[scale=1.5]{>}}},postaction={decorate},shorten >=0.5pt] (2) -- (3) node[anchor=south] {$a^{k}$};
  \draw[decoration={markings,mark=at position 1 with {\arrow[scale=1.5]{>}}},postaction={decorate},shorten >=0.5pt] (3) -- (4) node[anchor=east] {};
  \draw[decoration={markings,mark=at position 1 with {\arrow[scale=1.5]{>}}},postaction={decorate},shorten >=0.5pt] (5) -- (6) node[anchor=south] {};
  \draw[decoration={markings,mark=at position 1 with {\arrow[scale=1.5]{>}}},postaction={decorate},shorten >=0.5pt] (6) -- (7) node[anchor=south] {$a^{l}$};
  \draw[decoration={markings,mark=at position 1 with {\arrow[scale=1.5]{>}}},postaction={decorate},shorten >=0.5pt] (7) -- (8) node[anchor=south] {};
  \draw[double distance=1.5pt] (2) -- (6) node[anchor=west] {};
  \draw[decoration={markings,mark=at position 1 with {\arrow[scale=1.5]{>}}},postaction={decorate},shorten >=0.5pt] (3) -- (7) node[anchor=west] {$a^{l-k}$};
\end{tikzpicture}\]
It is easily seen that $\left\{K^{R}(a^{k}),\varphi^{k,l}_{a}\right\}_{k \in \mathbb{N}}$ is a direct system of $R$-complexes. Given elements $\underline{a}=a_{1},...,a_{n} \in R$, we let $\underline{a}^{k}=a_{1}^{k},...,a_{n}^{k}$ for every $k \geq 1$. Now
$$K^{R}(\underline{a}^{k})=K^{R}(a^{k}_{1}) \otimes_{R} \cdots \otimes_{R} K^{R}(a^{k}_{n}),$$
and we let
$$\varphi^{k,l}:=\varphi^{k,l}_{a_{1}}\otimes_{R}\cdots\otimes_{R}\varphi^{k,l}_{a_{n}}.$$
It follows that $\left\{K^{R}(\underline{a}^{k}),\varphi^{k,l}\right\}_{k \in \mathbb{N}}$ is a direct system of $R$-complexes. \index{direct system of Koszul complexes} It is also clear that $\left\{H_{i}\left(\underline{a}^{k};M\right),H_{i}\left(\varphi^{k,l}\otimes_{R}M\right)\right\}_{k \in \mathbb{N}}$ is a direct system of $R$-modules for every $i \in \mathbb{Z}$. \index{direct system of Koszul homologies}

\item[(ii)] Given an element $a \in R$, we define a morphism $\psi^{k,l}_{a}:K^{R}(a^{k}) \rightarrow K^{R}(a^{l})$ of $R$-complexes for every $k \geq l$ as follows:
\[\begin{tikzpicture}[every node/.style={midway},]
  \matrix[column sep={3em}, row sep={3em}]
  {\node(1) {$0$}; & \node(2) {$R$}; & \node(3) {$R$}; & \node(4) {$0$};\\
  \node(5) {$0$}; & \node(6) {$R$}; & \node(7) {$R$}; & \node(8) {$0$};\\};
  \draw[decoration={markings,mark=at position 1 with {\arrow[scale=1.5]{>}}},postaction={decorate},shorten >=0.5pt] (1) -- (2) node[anchor=east] {};
  \draw[decoration={markings,mark=at position 1 with {\arrow[scale=1.5]{>}}},postaction={decorate},shorten >=0.5pt] (2) -- (3) node[anchor=south] {$a^{k}$};
  \draw[decoration={markings,mark=at position 1 with {\arrow[scale=1.5]{>}}},postaction={decorate},shorten >=0.5pt] (3) -- (4) node[anchor=east] {};
  \draw[decoration={markings,mark=at position 1 with {\arrow[scale=1.5]{>}}},postaction={decorate},shorten >=0.5pt] (5) -- (6) node[anchor=south] {};
  \draw[decoration={markings,mark=at position 1 with {\arrow[scale=1.5]{>}}},postaction={decorate},shorten >=0.5pt] (6) -- (7) node[anchor=south] {$a^{l}$};
  \draw[decoration={markings,mark=at position 1 with {\arrow[scale=1.5]{>}}},postaction={decorate},shorten >=0.5pt] (7) -- (8) node[anchor=south] {};
  \draw[double distance=1.5pt] (3) -- (7) node[anchor=west] {};
  \draw[decoration={markings,mark=at position 1 with {\arrow[scale=1.5]{>}}},postaction={decorate},shorten >=0.5pt] (2) -- (6) node[anchor=west] {$a^{k-l}$};
\end{tikzpicture}\]
It is easily seen that $\left\{K^{R}(a^{k}),\varphi^{k,l}_{a}\right\}_{k \in \mathbb{N}}$ is an inverse system of $R$-complexes. Given elements $\underline{a}=a_{1},...,a_{n} \in R$, we let $\underline{a}^{k}=a_{1}^{k},...,a_{n}^{k}$ for every $k \geq 1$. Now
$$K^{R}(\underline{a}^{k})=K^{R}(a^{k}_{1}) \otimes_{R} \cdots \otimes_{R} K^{R}(a^{k}_{n}),$$
and we let
$$\psi^{k,l}:=\psi^{k,l}_{a_{1}}\otimes_{R}\cdots\otimes_{R}\psi^{k,l}_{a_{n}}.$$
It follows that $\left\{K^{R}(\underline{a}^{k}),\psi^{k,l}\right\}_{k \in \mathbb{N}}$ is an inverse system of $R$-complexes. \index{inverse system of Koszul complexes} It is also clear that $\left\{H_{i}\left(\underline{a}^{k};M\right),H_{i}\left(\psi^{k,l}\otimes_{R}M\right)\right\}_{k \in \mathbb{N}}$ is an inverse system of $R$-modules for every $i \in \mathbb{Z}$. \index{inverse system of Koszul homologies}
\end{enumerate}
\end{remark}

Recall that an inverse system $\left\{M_{\alpha},\varphi_{\alpha,\beta}\right\}_{\alpha \in \mathbb{N}}$ of $R$-modules is said to satisfy the trivial Mittag-Leffler condition \index{trivial Mittag-Leffler condition} if for every $\beta \in \mathbb{N}$, there is an $\alpha \geq \beta$ such that $\varphi_{\alpha \beta}=0$. Besides, the inverse system $\left\{M_{\alpha},\varphi_{\alpha,\beta}\right\}_{\alpha \in \mathbb{N}}$ of $R$-modules is said to satisfy the Mittag-Leffler condition \index{Mittag-Leffler condition} if for every $\beta \in \mathbb{N}$, there is an $\alpha_{0} \geq \beta$ such that $\im \varphi_{\alpha \beta}=\im \varphi_{\alpha_{0} \beta}$ for every $\alpha \geq \alpha_{0} \geq \beta$. It is straightforward to verify that the trivial Mittag-Leffler condition implies the Mittag-Leffler condition.

The following lemma reveals a significant feature of Koszul homology and lies at the heart of the proof of Greenlees-May Duality. The idea of the proof is taken from \cite{Sc}.

\begin{lemma} \label{3.2.7}
Let $\underline{a}=a_{1},...,a_{n} \in R$, and $\underline{a}^{k}=a_{1}^{k},...,a_{n}^{k}$ for every $k \geq 1$. Then the inverse system $\left\{H_{i}\left(\underline{a}^{k};R\right)\right\}_{k \in \mathbb{N}}$ satisfies the trivial Mittag-Leffler condition for every $i \geq 1$.
\end{lemma}

\begin{prf}
Let $a\in R$ and $M$ a finitely generated $R$-module. The transition maps of the inverse system $\left\{K^{R}(a^{k})\otimes_{R}M\right\}_{k \in \mathbb{N}}$ can be identified with the following morphisms of $R$-complexes for every $k \geq l$:
\[\begin{tikzpicture}[every node/.style={midway},]
  \matrix[column sep={3em}, row sep={3em}]
  {\node(1) {$0$}; & \node(2) {$M$}; & \node(3) {$M$}; & \node(4) {$0$};\\
  \node(5) {$0$}; & \node(6) {$M$}; & \node(7) {$M$}; & \node(8) {$0$};\\};
  \draw[decoration={markings,mark=at position 1 with {\arrow[scale=1.5]{>}}},postaction={decorate},shorten >=0.5pt] (1) -- (2) node[anchor=east] {};
  \draw[decoration={markings,mark=at position 1 with {\arrow[scale=1.5]{>}}},postaction={decorate},shorten >=0.5pt] (2) -- (3) node[anchor=south] {$a^{k}$};
  \draw[decoration={markings,mark=at position 1 with {\arrow[scale=1.5]{>}}},postaction={decorate},shorten >=0.5pt] (3) -- (4) node[anchor=east] {};
  \draw[decoration={markings,mark=at position 1 with {\arrow[scale=1.5]{>}}},postaction={decorate},shorten >=0.5pt] (5) -- (6) node[anchor=south] {};
  \draw[decoration={markings,mark=at position 1 with {\arrow[scale=1.5]{>}}},postaction={decorate},shorten >=0.5pt] (6) -- (7) node[anchor=south] {$a^{l}$};
  \draw[decoration={markings,mark=at position 1 with {\arrow[scale=1.5]{>}}},postaction={decorate},shorten >=0.5pt] (7) -- (8) node[anchor=south] {};
  \draw[double distance=1.5pt] (3) -- (7) node[anchor=west] {};
  \draw[decoration={markings,mark=at position 1 with {\arrow[scale=1.5]{>}}},postaction={decorate},shorten >=0.5pt] (2) -- (6) node[anchor=west] {$a^{k-l}$};
\end{tikzpicture}\]
Since $H_{1}\left(a^{k};M\right)=\left(0:_{M}a^{k}\right)$, the transition maps of the inverse system $\left\{H_{1}\left(a^{k};M\right)\right\}_{k \in \mathbb{N}}$ can be identified with the $R$-homomorphisms
$$\left(0:_{M}a^{k}\right) \xrightarrow{a^{k-l}} \left(0:_{M}a^{l}\right)$$
for every $k \geq l$. Fix $l \in \mathbb{N}$. Since $R$ is noetherian and $M$ is finitely generated, the ascending chain
$$\left(0:_{M}a\right) \subseteq \left(0:_{M}a^{2}\right) \subseteq \cdots$$
of submodules of $M$ stabilizes, i.e. there is an integer $t \geq 1$ such that
$$\left(0:_{M}a^{t}\right) = \left(0:_{M}a^{t+1}\right) = \cdots.$$
Set $k:=t+l$. Then the transition map $\left(0:_{M}a^{k}\right) \xrightarrow{a^{k-l}} \left(0:_{M}a^{l}\right)$ is zero. Indeed, if $x \in \left(0:_{M}a^{k}\right)$, then since
$$\left(0:_{M}a^{k}\right) = \left(0:_{M}a^{t+l}\right) = \left(0:_{M}a^{t}\right),$$
we have $x \in \left(0:_{M}a^{t}\right)$, so $a^{k-l}x=a^{t}x=0$. This shows that the inverse system $\left\{H_{1}\left(a^{k};M\right)\right\}_{k \in \mathbb{N}}$ satisfies the trivial Mittag-Leffler condition. But $H_{i}\left(a^{k};M\right)=0$ for every $i \geq 2$, so the inverse system $\left\{H_{i}\left(a^{k};M\right)\right\}_{k \in \mathbb{N}}$ satisfies the trivial Mittag-Leffler condition for every $i \geq 1$.

Now we argue by induction on $n$. If $n=1$, then the inverse system $\left\{H_{i}\left(a_{1}^{k};R\right)\right\}_{k \in \mathbb{N}}$ satisfies the trivial Mittag-Leffler condition for every $i \geq 1$ by the discussion above. Now assume that $n \geq 2$, and make the obvious induction hypothesis. There is an exact sequence of inverse systems
\begin{equation} \label{eq:3.2.7.1}
\left\{H_{i}\left(a_{1}^{k},...,a_{n-1}^{k};R\right)\right\}_{k \in \mathbb{N}} \rightarrow \left\{H_{0}\left(a_{n}^{k};H_{i}\left(a_{1}^{k},...,a_{n-1}^{k};R\right)\right)\right\}_{k \in \mathbb{N}} \rightarrow 0
\end{equation}
of $R$-modules for every $i \geq 0$. By the induction hypothesis, the inverse system $\left\{H_{i}\left(a_{1}^{k},...,a_{n-1}^{k};R\right)\right\}_{k \in \mathbb{N}}$ satisfies the trivial Mittag-Leffler condition for every $i \geq 1$, so the exact sequence \eqref{eq:3.2.7.1} shows that the inverse system $\left\{H_{0}\left(a_{n}^{k};H_{i}\left(a_{1}^{k},...,a_{n-1}^{k};R\right)\right)\right\}_{k \in \mathbb{N}}$ satisfies the Mittag-Leffler condition for every $i \geq 1$. On the other hand, there is a short exact sequence of inverse systems
\begin{equation} \label{eq:3.2.7.2}
0 \rightarrow \left\{H_{0}\left(a_{n}^{k};H_{i}\left(a_{1}^{k},...,a_{n-1}^{k};R\right)\right)\right\}_{k \in \mathbb{N}} \rightarrow \left\{H_{i}\left(a_{1}^{k},...,a_{n}^{k};R\right)\right\}_{k \in \mathbb{N}} \rightarrow $$$$ \left\{H_{1}\left(a_{n}^{k};H_{i-1}\left(a_{1}^{k},...,a_{n-1}^{k};R\right)\right)\right\}_{k \in \mathbb{N}} \rightarrow 0
\end{equation}
of $R$-modules for every $i \geq 0$.
Since $H_{i-1}\left(a_{1}^{k},...,a_{n-1}^{k};R\right)$ is a finitely generated $R$-module for every $i \geq 1$, the discussion above shows that
$\left\{H_{1}\left(a_{n}^{k};H_{i-1}\left(a_{1}^{k},...,a_{n-1}^{k};R\right)\right)\right\}_{k \in \mathbb{N}}$ satisfies the Mittag-Leffler condition for every $i \geq 1$. Therefore, the short exact sequence \eqref{eq:3.2.7.2} shows that the inverse system $\left\{H_{i}\left(a_{1}^{k},...,a_{n}^{k};R\right)\right\}_{k \in \mathbb{N}}$ satisfies the trivial Mittag-Leffler condition for every $i \geq 1$.
\end{prf}

The category $\mathcal{C}(R)$ of $R$-complexes enjoys direct limits and inverse limits. However, the derived category $\mathcal{D}(R)$ does not support the notions of direct limits and inverse limits. But this situation is remedied by the existence of homotopy direct limits and homotopy inverse limits as defined in triangulated categories with countable products and coproducts.

\begin{remark} \label{3.2.8}
Let $\left\{X^{\alpha},\varphi^{\alpha \beta}\right\}_{\alpha \in \mathbb{N}}$ be a direct system of $R$-complexes, and $\left\{Y^{\alpha},\psi^{\alpha \beta}\right\}_{\alpha \in \mathbb{N}}$ an inverse system of $R$-complexes. Then we have:
\begin{enumerate}
\item[(i)] The direct limit \index{direct limit of complexes} of the direct system $\left\{X^{\alpha},\varphi^{\alpha \beta}\right\}_{\alpha \in \mathbb{N}}$ is an $R$-complex $\varinjlim X^{\alpha}$ given by $\left(\varinjlim X^{\alpha}\right)_{i}=\varinjlim X^{\alpha}_{i}$ and $\partial^{\varinjlim X^{\alpha}}_{i}=\varinjlim \partial^{X^{\alpha}}_{i}$ for every $i \in \mathbb{Z}$. Indeed, it is easy to see that $\varinjlim X^{\alpha}$ satisfies the universal property of direct limits in a category.
\item[(ii)] The homotopy direct limit \index{homotopy direct limit of complexes} of the direct system $\left\{X^{\alpha},\varphi^{\alpha \beta}\right\}_{\alpha \in \mathbb{N}}$ is given by $\hocolim X^{\alpha} = \Cone(\vartheta)$, where the morphism $\vartheta: \bigoplus_{\alpha=1}^{\infty}X^{\alpha} \rightarrow \bigoplus_{\alpha=1}^{\infty}X^{\alpha}$ is given by $\vartheta_{i}\left((x_{i}^{\alpha})\right)= \iota_{i}^{\alpha}(x_{i}^{\alpha})-\iota_{i}^{\alpha+1}\left(\varphi_{i}^{\alpha,\alpha+1}(x_{i}^{\alpha})\right)$ for every $i \in \mathbb{Z}$. Indeed, it is easy to see that the morphism $\vartheta$ fits into a distinguished triangle
    $$\bigoplus_{\alpha=1}^{\infty}X^{\alpha} \rightarrow \bigoplus_{\alpha=1}^{\infty}X^{\alpha} \rightarrow \hocolim X^{\alpha} \rightarrow.$$
\item[(iii)] The inverse limit \index{inverse limit of complexes} of the inverse system $\left\{Y^{\alpha},\psi^{\alpha \beta}\right\}_{\alpha \in \mathbb{N}}$ is an $R$-complex $\varprojlim Y^{\alpha}$ given by $\left(\varprojlim Y^{\alpha}\right)_{i}=\varprojlim Y^{\alpha}_{i}$ and $\partial^{\varprojlim Y^{\alpha}}_{i}=\varprojlim \partial^{Y^{\alpha}}_{i}$ for every $i \in \mathbb{Z}$. Indeed, it is easy to see that $\varprojlim X^{\alpha}$ satisfies the universal property of inverse limits in a category.
\item[(iv)] The homotopy inverse limit \index{homotopy inverse limit of complexes} of the inverse system $\left\{Y^{\alpha},\psi^{\alpha \beta}\right\}_{\alpha \in \mathbb{N}}$ is given by $\holim Y^{\alpha} = \Sigma^{-1}\Cone(\varpi)$, where the morphism $\varpi: \prod_{\alpha=1}^{\infty}Y^{\alpha} \rightarrow \prod_{\alpha=1}^{\infty}Y^{\alpha}$ is given by $\varpi_{i}\left((y_{i}^{\alpha})\right)= \left(y_{i}^{\alpha}-\psi_{i}^{\alpha+1,\alpha}(y_{i}^{\alpha+1})\right)$ for every $i \in \mathbb{Z}$. Indeed, it is easy to see that the morphism $\varpi$ fits into a distinguished triangle
    $$\holim Y^{\alpha} \rightarrow \prod_{\alpha=1}^{\infty}Y^{\alpha} \rightarrow \prod_{\alpha=1}^{\infty}Y^{\alpha} \rightarrow.$$
\end{enumerate}
\end{remark}

The Mittag-Leffler condition forces many limits to be zero.

\begin{lemma} \label{3.2.9}
Let $\left\{M_{\alpha},\varphi_{\alpha \beta}\right\}_{\alpha \in \mathbb{N}}$ be an inverse system of $R$-modules that satisfies the trivial Mittag-Leffler condition, and $\mathcal{F}:\mathcal{M}(R) \rightarrow \mathcal{M}(R)$ an additive contravariant functor. Then the following assertions hold:
\begin{enumerate}
\item[(i)] $\varprojlim M_{\alpha} = 0 = \varprojlim^{1} M_{\alpha}$.
\item[(ii)] $\varinjlim \mathcal{F}(M_{\alpha})=0$.
\end{enumerate}
\end{lemma}

\begin{prf}
(i): Let $\varpi:\prod_{\alpha \in \mathbb{N}}M_{\alpha} \rightarrow \prod_{\alpha \in \mathbb{N}}M_{\alpha}$ be an $R$-homomorphism given by $\varpi\left((x_{\alpha})\right)=\left(x_{\alpha}-\varphi_{\alpha+1,\alpha}(x_{\alpha+1})\right)$. We show that $\varpi$ is an isomorphism. Let $(x_{\alpha}) \in \prod_{\alpha \in \mathbb{N}}M_{\alpha}$ be such that $x_{\alpha}=\varphi_{\alpha+1,\alpha}(x_{\alpha+1})$ for every $\alpha \in \mathbb{N}$. Fix $\alpha \in \mathbb{N}$, and by the trivial Mittag-Leffler condition choose $\gamma \geq \alpha$ such that $\varphi_{\gamma \alpha}=0$. Then we have
\begin{equation*}
\begin{split}
 x_{\alpha} & = \varphi_{\alpha+1,\alpha}(x_{\alpha+1}) \\
 & = \varphi_{\alpha+1,\alpha}\left(\varphi_{\alpha+2,\alpha+1}\left(...\left(\varphi_{\gamma,\gamma-1} (x_{\gamma})\right)\right)\right)\\
 & = \varphi_{\gamma \alpha}(x_{\gamma})\\
 & = 0.\\
\end{split}
\end{equation*}
Hence $(x_{\alpha})=0$, and thus $\varpi$ is injective. Now let $(y_{\alpha}) \in \prod_{\alpha \in \mathbb{N}}M_{\alpha}$. For any $\beta \in \mathbb{N}$, we set $x_{\beta}:= \sum_{\alpha=\beta}^{\infty}\varphi_{\alpha \beta}(y_{\alpha})$ which is a finite sum by the trivial Mittag-Leffler condition. Then we have
\begin{equation*}
\begin{split}
 \varphi_{\beta+1,\beta}(x_{\beta+1}) & = \varphi_{\beta+1,\beta}\left(\sum_{\alpha=\beta+1}^{\infty}\varphi_{\alpha, \beta+1}(y_{\alpha})\right) \\
 & = \sum_{\alpha=\beta+1}^{\infty} \varphi_{\alpha \beta}(y_{\alpha})\\
 & = \sum_{\alpha=\beta}^{\infty} \varphi_{\alpha \beta}(y_{\alpha})-\varphi_{\beta \beta}(y_{\beta})\\
 & = x_{\beta}-y_{\beta}.\\
\end{split}
\end{equation*}
Therefore, we have
$$\varpi\left((x_{\alpha})\right)=\left(x_{\alpha}-\varphi_{\alpha+1,\alpha}(x_{\alpha+1})\right)=(y_{\alpha}),$$
so $\varpi$ is surjective. It follows that $\varpi$ is an isomorphism. Therefore, $\varprojlim M_{\alpha} \cong \ker \varpi=0$ and $\varprojlim^{1} M_{\alpha} \cong \coker \varpi=0$.

(ii): First we note that $\left\{\mathcal{F}(M_{\alpha}),\psi_{\beta\alpha}:=\mathcal{F}(\varphi_{\alpha\beta})\right\}_{\alpha \in \mathbb{N}}$ is a direct system of $R$-modules. Let $\psi_{\alpha}:\mathcal{F}(M_{\alpha}) \rightarrow \varinjlim \mathcal{F}(M_{\alpha})$ be the natural injection of direct limit for every $\alpha \in \mathbb{N}$. We know that an arbitrary element of $\varinjlim \mathcal{F}(M_{\alpha})$ is of the form $\psi_{t}(y)$ for some $t \in \mathbb{N}$ and some $y \in \mathcal{F}(M_{t})$. By the trivial Mittag-Leffler condition, there is an integer $s \geq t$ such that $\varphi_{st}=0$, so that $\psi_{ts}=\mathcal{F}(\varphi_{st})=0$. Then $\psi_{t}(y)=\psi_{s}\left(\psi_{ts}(y)\right)=0$. Hence $\varinjlim \mathcal{F}(M_{\alpha})=0$.
\end{prf}

The next proposition collects some information on the homology of limits.

\begin{proposition} \label{3.2.10}
Let $\left\{X^{\alpha},\varphi^{\alpha \beta}\right\}_{\alpha \in \mathbb{N}}$ be a direct system of $R$-complexes, and $\left\{Y^{\alpha},\psi^{\alpha \beta}\right\}_{\alpha \in \mathbb{N}}$ an inverse system of $R$-complexes. Then the following assertions hold for every $i \in \mathbb{Z}$:
\begin{enumerate}
\item[(i)] There is a natural isomorphism $H_{i}\left(\varinjlim X^{\alpha}\right) \cong \varinjlim H_{i}(X^{\alpha})$.
\item[(ii)] There is a natural isomorphism $H_{i}\left(\hocolim X^{\alpha}\right) \cong \varinjlim H_{i}(X^{\alpha})$.
\item[(iii)] If the inverse system $\left\{Y^{\alpha}_{i},\psi^{\alpha \beta}_{i}\right\}_{\alpha \in \mathbb{N}}$ of $R$-modules satisfies the Mittag-Leffler condition for every $i \in \mathbb{Z}$, then there is a short exact sequence
\begin{center}
$0 \rightarrow \varprojlim^{1} H_{i+1}(Y^{\alpha}) \rightarrow H_{i}\left(\varprojlim Y^{\alpha}\right) \rightarrow \varprojlim H_{i}(Y^{\alpha}) \rightarrow 0$
\end{center}
of $R$-modules.
\item[(iv)] There is a short exact sequence
\begin{center}
$0 \rightarrow \varprojlim^{1} H_{i+1}(Y^{\alpha}) \rightarrow H_{i}\left(\holim Y^{\alpha}\right) \rightarrow \varprojlim H_{i}(Y^{\alpha}) \rightarrow 0$
\end{center}
of $R$-modules.
\end{enumerate}
\end{proposition}

\begin{prf}
(i): See \cite[Theorem 4.2.4]{Se}.

(ii): See the paragraph before \cite[Lemma 0.1]{GM}.

(iii): See \cite[Theorem 3.5.8]{We}.

(iv): See the paragraph after \cite[Lemma 0.1]{GM}.
\end{prf}

Now we are ready to present the following definitions.

\begin{definition} \label{3.2.11}
Let $\underline{a}=a_{1},...,a_{n} \in R$. Then:
\begin{enumerate}
\item[(i)] Define the \textit{\v{C}ech complex} \index{\v{C}ech complex} on the elements $\underline{a}$ to be $\check{C}(\underline{a}):= \varinjlim \Sigma^{-n}K^{R}(\underline{a}^{k})$.
\item[(ii)] Define the \textit{stable \v{C}ech complex} \index{stable \v{C}ech complex} on the elements $\underline{a}$ to be $\check{C}_{\infty}(\underline{a}):= \hocolim \Sigma^{-n}K^{R}(\underline{a}^{k})$.
\end{enumerate}
\end{definition}

We note that $\check{C}(\underline{a})$ is a bounded $R$-complex of flat modules concentrated in degrees $0,-1,...,-n$, and $\check{C}_{\infty}(\underline{a})$ is a bounded $R$-complex of free modules concentrated in degrees $1,0,...,-n$. Moreover, it can be shown that there is a quasi-isomorphism $\check{C}_{\infty}(\underline{a}) \xrightarrow{\simeq} \check{C}(\underline{a})$, which in turn implies that $\check{C}_{\infty}(\underline{a}) \simeq \check{C}(\underline{a})$ in $\mathcal{D}(R)$. Therefore, $\check{C}_{\infty}(\underline{a})$ is a semi-projective approximation of the semi-flat $R$-complex $\check{C}(\underline{a})$.

The next proposition investigates the relation between local cohomology and local homology with \v{C}ech complex and stable \v{C}ech complex, and provides the first essential step towards the Greenlees-May Duality.

\begin{proposition} \label{3.2.12}
Let $\mathfrak{a}= (a_{1},...,a_{n})$ be an ideal of $R$, $\underline{a}=a_{1},...,a_{n}$, and $M$ an $R$-module. Then there are natural isomorphisms for every $i \geq 0$:
\begin{enumerate}
\item[(i)] $H^{i}_{\mathfrak{a}}(M) \cong H_{-i}\left(\check{C}(\underline{a})\otimes_{R}M\right) \cong H_{-i}\left(\check{C}_{\infty}(\underline{a})\otimes_{R}M\right)$.
\item[(ii)] $H^{\mathfrak{a}}_{i}(M) \cong H_{i}\left(\Hom_{R}\left(\check{C}_{\infty}(\underline{a}),M\right)\right)$.
\end{enumerate}
\end{proposition}

\begin{prf}
(i): Let $\mathcal{F}^{i}=H_{-i}\left(\check{C}(\underline{a})\otimes_{R}-\right):\mathcal{M}(R) \rightarrow \mathcal{M}(R)$ for every $i \geq 0$. Given a short exact sequence
$$0 \rightarrow M' \rightarrow M \rightarrow M'' \rightarrow 0$$
of $R$-modules, since $\check{C}(\underline{a})$ is an $R$-complex of flat modules, the functor $\check{C}(\underline{a})\otimes_{R}-:\mathcal{C}(R) \rightarrow \mathcal{C}(R)$ is exact, whence we get a short exact sequence
$$0 \rightarrow \check{C}(\underline{a})\otimes_{R}M' \rightarrow \check{C}(\underline{a})\otimes_{R}M \rightarrow \check{C}(\underline{a})\otimes_{R}M'' \rightarrow 0$$
of $R$-complexes, which in turn yields a long exact homology sequence in a functorial way. This shows that $\left(\mathcal{F}^{i}:\mathcal{M}(R) \rightarrow \mathcal{M}(R)\right)_{i \geq 0}$ is a cohomological covariant $\delta$-functor. Moreover, using Proposition \ref{3.2.10} (i), we have
\begin{equation} \label{eq:3.2.12.1}
\begin{split}
 \mathcal{F}^{i} & = H_{-i}\left(\check{C}(\underline{a})\otimes_{R}-\right) \\
 & = H_{-i}\left(\left(\varinjlim \Sigma^{-n}K^{R}(\underline{a}^{k})\right)\otimes_{R}-\right) \\
 & \cong \varinjlim H_{n-i}\left(K^{R}(\underline{a}^{k})\otimes_{R}-\right) \\
 & \cong \varinjlim H_{n-i}\left(\underline{a}^{k};-\right)\\
\end{split}
\end{equation}
for every $i \geq 0$.

Let $I$ be an injective $R$-module. Then by the display \eqref{eq:3.2.12.1}, we have
\begin{equation} \label{eq:3.2.12.2}
\begin{split}
 \mathcal{F}^{i}(I) & = \varinjlim H_{n-i}\left(\underline{a}^{k};I\right)\\
 & \cong \varinjlim H_{-i}\left(\Hom_{R}\left(K^{R}(\underline{a}^{k}),I\right)\right) \\
 & \cong \varinjlim \Hom_{R}\left(H_{i}\left(K^{R}(\underline{a}^{k})\right),I\right) \\
 & \cong \varinjlim \Hom_{R}\left(H_{i}\left(\underline{a}^{k};R\right),I\right).\\
\end{split}
\end{equation}
By Lemma \ref{3.2.7}, the inverse system $\left\{H_{i}\left(\underline{a}^{k};R\right)\right\}_{k \in \mathbb{N}}$ satisfies the trivial Mittag-Leffler condition for every $i \geq 1$. Now Lemma \ref{3.2.9} (ii) implies that $\varinjlim \Hom_{R}\left(H_{i}\left(\underline{a}^{k};R\right),I\right)=0$, thereby the display \eqref{eq:3.2.12.2} shows that $\mathcal{F}^{i}(I)=0$ for every $i \geq 1$.

Let $M$ be an $R$-module. Then by the display \eqref{eq:3.2.12.1}, we have the natural isomorphisms
\begin{equation*}
\begin{split}
 \mathcal{F}^{0}(M) & \cong \varinjlim H_{n}\left(\underline{a}^{k};M\right) \\
 & \cong \varinjlim \left(0:_{M}(\underline{a}^{k})\right) \\
 & \cong \varinjlim \Hom_{R}\left(R/(\underline{a}^{k}),M\right) \\
 & \cong \varinjlim \Hom_{R}\left(R/\mathfrak{a}^{k},M\right)\\
 & \cong \Gamma_{\mathfrak{a}}(M)\\
 & \cong H^{0}_{\mathfrak{a}}(M).\\
\end{split}
\end{equation*}
It follows from Corollary \ref{3.2.5} (ii) that $H^{i}_{\mathfrak{a}}(-) \cong \mathcal{F}^{i}$ for every $i \geq 0$.

For the second isomorphism, using the display \eqref{eq:3.2.12.1} and Proposition \ref{3.2.10} (ii), we have the natural isomorphisms
\begin{equation*}
\begin{split}
 H^{i}_{\mathfrak{a}}(M) & \cong \mathcal{F}^{i}(M) \\
 & \cong \varinjlim H_{n-i}\left(\underline{a}^{k};M\right) \\
 & \cong \varinjlim H_{n-i}\left(K^{R}(\underline{a}^{k})\otimes_{R}M\right) \\
 & \cong H_{n-i}\left(\hocolim\left(K^{R}(\underline{a}^{k})\otimes_{R}M\right)\right) \\
 & \cong H_{-i}\left(\left(\hocolim \Sigma^{-n}K^{R}(\underline{a})\right)\otimes_{R}M\right) \\
 & \cong H_{-i}\left(\check{C}_{\infty}(\underline{a})\otimes_{R}M\right) \\
\end{split}
\end{equation*}
for every $i \geq 0$.

(ii): Let $\mathcal{F}_{i}=H_{i}\left(\Hom_{R}\left(\check{C}_{\infty}(\underline{a}),-\right)\right):\mathcal{M}(R) \rightarrow \mathcal{M}(R)$ for every $i \geq 0$. Given a short exact sequence
\[0 \rightarrow M' \rightarrow M \rightarrow M'' \rightarrow 0\]
$R$-modules, since $\check{C}_{\infty}(\underline{a})$ is an $R$-complex of free modules, the functor $\Hom_{R}\left(\check{C}_{\infty}(\underline{a}),-\right):\mathcal{C}(R) \rightarrow \mathcal{C}(R)$ is exact, whence we get a short exact sequence
$$0 \rightarrow \Hom_{R}\left(\check{C}_{\infty}(\underline{a}),M'\right) \rightarrow \Hom_{R}\left(\check{C}_{\infty}(\underline{a}),M\right) \rightarrow \Hom_{R}\left(\check{C}_{\infty}(\underline{a}),M''\right) \rightarrow 0$$
of $R$-complexes, which in turn yields a long exact homology sequence in a functorial way. It follows that $\left(\mathcal{F}_{i}:\mathcal{M}(R) \rightarrow \mathcal{M}(R)\right)_{i \geq 0}$ is a homological covariant $\delta$-functor. Moreover, using the self-duality property of Koszul complex, we have
\begin{equation} \label{eq:3.2.12.3}
\begin{split}
 \mathcal{F}_{i} & = H_{i}\left(\Hom_{R}\left(\check{C}_{\infty}(\underline{a}),-\right)\right) \\
 & = H_{i}\left(\Hom_{R}\left(\hocolim \Sigma^{-n}K^{R}(\underline{a}^{k}),-\right)\right) \\
 & \cong H_{i}\left(\holim \Sigma^{n}\Hom_{R}\left(K^{R}(\underline{a}^{k}),-\right)\right)\\
 & \cong H_{i}\left(\holim \left(K^{R}(\underline{a}^{k})\otimes_{R}-\right)\right)\\
\end{split}
\end{equation}
for every $i \geq 0$.

Let $M$ be an $R$-module. By Proposition \ref{3.2.10} (iv), we get a short exact sequence
\begin{center}
$0 \rightarrow \varprojlim ^{1} H_{i+1}\left(K^{R}(\underline{a}^{k})\otimes_{R}M\right) \rightarrow H_{i}\left(\holim\left(K^{R}(\underline{a}^{k})\otimes_{R}M\right)\right) \rightarrow \varprojlim H_{i}\left(K^{R}(\underline{a}^{k})\otimes_{R}M\right) \rightarrow 0,$
\end{center}
which implies the short exact sequence
\begin{center}
$0 \rightarrow \varprojlim^{1} H_{i+1}\left(\underline{a}^{k};M\right) \rightarrow \mathcal{F}_{i}(M) \rightarrow \varprojlim H_{i}\left(\underline{a}^{k};M\right) \rightarrow 0$
\end{center}
of $R$-modules for every $i \geq 0$.

Let $F$ be a free $R$-module. If $i \geq 1$, then the inverse system $\left\{H_{i}\left(\underline{a}^{k};R\right)\right\}_{k \in \mathbb{N}}$ satisfies the trivial Mittag-Leffler condition by Lemma \ref{3.2.7}. But $H_{i}\left(\underline{a}^{k};F\right) \cong H_{i}\left(\underline{a}^{k};R\right)\otimes_{R}F$, so it straightforward to see that the inverse system $\left\{H_{i}\left(\underline{a}^{k};F\right)\right\}_{k \in \mathbb{N}}$ satisfies the trivial Mittag-Leffler condition for every $i \geq 1$. Therefore, Lemma \ref{3.2.9} (i) implies that
\begin{center}
$\varprojlim^{1} H_{i}\left(\underline{a}^{k};F\right)=0=\varprojlim H_{i}\left(\underline{a}^{k};F\right)$
\end{center}
for every $i \geq 1$. It follows from the above short exact sequence that $\mathcal{F}_{i}(F)=0$ for every $i \geq 1$.

Upon setting $i=0$, the above short exact sequence yields
\begin{center}
$0 = \varprojlim^{1} H_{1}\left(\underline{a}^{k};F\right) \rightarrow \mathcal{F}_{0}(F) \rightarrow \varprojlim H_{0}\left(\underline{a}^{k};F\right) \rightarrow 0.$
\end{center}
Thus we get the natural isomorphisms
\begin{equation*}
\begin{split}
 \mathcal{F}_{0}(F) & \cong \varprojlim H_{0}\left(\underline{a}^{k};F\right) \\
 & \cong \varprojlim F/(\underline{a}^{k})F \\
 & \cong \varprojlim F/\mathfrak{a}^{k}F \\
 & = \widehat{F}^{\mathfrak{a}} \\
 & \cong H^{\mathfrak{a}}_{0}(F).\\
\end{split}
\end{equation*}
It now follows from Corollary \ref{3.2.5} (i) that $H^{\mathfrak{a}}_{i}(-) \cong \mathcal{F}_{i}$ for every $i \geq 0$.
\end{prf}

\begin{remark} \label{3.2.13}
One should note that $H^{\mathfrak{a}}_{i}(M) \ncong H_{i}\left(\Hom_{R}\left(\check{C}(\underline{a}),M\right)\right)$.
\end{remark}

\section{Complex Prerequisites}

In this section, we commence on developing the requisite tools on complexes which are to be deployed in Section 3.4. For more information on the material in this section, refer to \cite{AF}, \cite{Ha2}, \cite{Fo}, \cite{Li}, and \cite{Sp}.

The derived category $\mathcal{D}(R)$ \index{derived category} is defined as the localization of the homotopy category $\mathcal{K}(R)$ with respect to the multiplicative system of quasi-isomorphisms. Simply put, an object in $\mathcal{D}(R)$ is an $R$-complex $X$ displayed in the standard homological style
$$X= \cdots \rightarrow X_{i+1} \xrightarrow {\partial^{X}_{i+1}} X_{i} \xrightarrow {\partial^{X}_{i}} X_{i-1} \rightarrow \cdots,$$
and a morphism $\varphi:X\rightarrow Y$ in $\mathcal{D}(R)$ is given by the equivalence class of a pair $(f,g)$ of morphisms
$X \xleftarrow{g} U \xrightarrow{f} Y$ in $\mathcal{C}(R)$ with $g$ a quasi-isomorphism, under the equivalence relation that identifies two such pairs $(f,g)$ and $(f^{\prime},g^{\prime})$, whenever there is a diagram in $\mathcal{C}(R)$ as follows which commutes up to homotopy:
\[\begin{tikzpicture}[every node/.style={midway}]
  \matrix[column sep={3em}, row sep={3em}]
  {\node(1) {$$}; & \node(2) {$U$}; & \node(3) {$$};\\
  \node(4) {$X$}; & \node(5) {$V$}; & \node(6) {$Y$};\\
  \node(7) {$$}; & \node(8) {$U^{\prime}$}; & \node(9) {$$};\\};
  \draw[decoration={markings,mark=at position 1 with {\arrow[scale=1.5]{>}}},postaction={decorate},shorten >=0.5pt] (2) -- (4) node[anchor=south] {$g$} node[anchor=west] {$\simeq$};
  \draw[decoration={markings,mark=at position 1 with {\arrow[scale=1.5]{>}}},postaction={decorate},shorten >=0.5pt] (2) -- (6) node[anchor=south] {$f$};
  \draw[decoration={markings,mark=at position 1 with {\arrow[scale=1.5]{>}}},postaction={decorate},shorten >=0.5pt] (8) -- (4) node[anchor=north] {$g^{\prime}$} node[anchor=west] {$\simeq$};
  \draw[decoration={markings,mark=at position 1 with {\arrow[scale=1.5]{>}}},postaction={decorate},shorten >=0.5pt] (8) -- (6) node[anchor=north] {$f^{\prime}$};
  \draw[decoration={markings,mark=at position 1 with {\arrow[scale=1.5]{>}}},postaction={decorate},shorten >=0.5pt] (5) -- (2) node[anchor=east] {$$};
  \draw[decoration={markings,mark=at position 1 with {\arrow[scale=1.5]{>}}},postaction={decorate},shorten >=0.5pt] (5) -- (4) node[anchor=south] {$\simeq$};
  \draw[decoration={markings,mark=at position 1 with {\arrow[scale=1.5]{>}}},postaction={decorate},shorten >=0.5pt] (5) -- (6) node[anchor=east] {$$};
  \draw[decoration={markings,mark=at position 1 with {\arrow[scale=1.5]{>}}},postaction={decorate},shorten >=0.5pt] (5) -- (8) node[anchor=east] {$$};
\end{tikzpicture}\]
The isomorphisms in $\mathcal{D}(R)$ are marked by the symbol $\simeq$.

The derived category $\mathcal{D}(R)$ is triangulated. A distinguished triangle \index{distinguished triangle} in $\mathcal{D}(R)$ is a triangle that is isomorphic to a triangle of the form
$$X \xrightarrow {\mathfrak{L}(f)} Y \xrightarrow{\mathfrak{L}(\varepsilon)} \Cone(f) \xrightarrow{\mathfrak{L}(\varpi)} \Sigma X,$$
for some morphism $f:X \rightarrow Y$ in $\mathcal{C}(R)$ with the mapping cone sequence
$$0 \rightarrow Y \xrightarrow{\varepsilon} \Cone(f) \xrightarrow{\varpi} \Sigma X \rightarrow 0,$$
in which $\mathfrak{L}:\mathcal{C}(R) \rightarrow \mathcal{D}(R)$ is the canonical functor that is defined as $\mathfrak{L}(X)=X$ for every $R$-complex $X$, and $\mathfrak{L}(f)=\varphi$ where $\varphi$ is represented by the morphisms $X \xleftarrow{1^{X}} X \xrightarrow{f} Y$ in $\mathcal{C}(R)$. We note that if $f$ is a quasi-isomorphism in $\mathcal{C}(R)$, then $\mathfrak{L}(f)$ is an isomorphism in $\mathcal{D}(R)$. We sometimes use the shorthand notation
$$X \rightarrow Y \rightarrow Z \rightarrow$$
for a distinguished triangle.

We let $\mathcal{D}_{\sqsubset}(R)$ (res. $\mathcal{D}_{\sqsupset}(R)$) denote the full subcategory of $\mathcal{D}(R)$ consisting of $R$-complexes $X$ with $H_{i}(X)=0$ for $i \gg 0$ (res. $i \ll 0$), and  $D_{\square}(R):=\mathcal{D}_{\sqsubset}(R)\cap \mathcal{D}_{\sqsupset}(R)$. We further let $\mathcal{D}^{f}(R)$ denote the full subcategory of $\mathcal{D}(R)$ consisting of $R$-complexes $X$ with finitely generated homology modules. We also feel free to use any combination of the subscripts and the superscript as in $\mathcal{D}^{f}_{\square}(R)$, with the obvious meaning of the intersection of the two subcategories involved.

We recall the resolutions of complexes.

\begin{definition} \label{3.3.1}
We have:
\begin{enumerate}
\item[(i)] An $R$-complex $P$ of projective modules is said to be \textit{semi-projective} \index{semi-projective complex} if the functor $\Hom_{R}(P,-)$ preserves quasi-isomorphisms. By a \textit{semi-projective resolution} \index{semi-projective resolution} of an $R$-complex $X$, we mean a quasi-isomorphism $P\xrightarrow {\simeq} X$ in which $P$ is a semi-projective $R$-complex.
\item[(ii)] An $R$-complex
$I$ of injective modules is said to be \textit{semi-injective} \index{semi-injective complex} if the functor $\Hom_{R}(-,I)$ preserves quasi-isomorphisms. By a \textit{semi-injective resolution} \index{semi-injective resolution} of an $R$-complex $X$, we mean a quasi-isomorphism $X\xrightarrow {\simeq} I$ in which $I$ is a semi-injective $R$-complex.
\item[(iii)] An $R$-complex
$F$ of flat modules is said to be \textit{semi-flat} \index{semi-flat complex} if the functor $F\otimes_{R}-$ preserves quasi-isomorphisms. By a \textit{semi-flat resolution} \index{semi-flat resolution} of an $R$-complex $X$, we mean a quasi-isomorphism $F\xrightarrow {\simeq} X$ in which $F$ is a semi-flat $R$-complex.
\end{enumerate}
\end{definition}

Semi-projective, semi-injective, and semi-flat resolutions exist for any $R$-complex. Moreover, any right-bounded $R$-complex of projective (flat) modules is semi-projective (semi-flat), and any left-bounded $R$-complex of injective modules is semi-injective.

We now remind the total derived functors that we need.

\begin{remark} \label{3.3.2}
Let $\mathfrak{a}$ be an ideal of $R$, and $X$ and $Y$ two $R$-complexes. Then we have:
\begin{enumerate}
\item[(i)] Each of the functors $\Hom_{R}(X,-)$ and $\Hom_{R}(-,Y)$ on $\mathcal{C}(R)$ enjoys a right total derived functor on $\mathcal{D}(R)$, together with a balance property, in the sense that ${\bf R}\Hom_{R}(X,Y)$ can be computed by
$${\bf R}\Hom_{R}(X,Y)\simeq \Hom_{R}(P,Y) \simeq \Hom_{R}(X,I),$$
where $P\xrightarrow {\simeq} X$ is any semi-projective resolution of $X$, and $Y\xrightarrow {\simeq} I$ is any semi-injective resolution of $Y$. In addition, these functors turn out to be triangulated, in the sense that they preserve shifts and distinguished triangles. Moreover, we let $$\Ext^{i}_{R}(X,Y):=H_{-i}\left({\bf R}\Hom_{R}(X,Y)\right)$$
for every $i \in \mathbb{Z}$.
\item[(ii)] Each of the functors $X\otimes_{R}-$ and $-\otimes_{R}Y$ on $\mathcal{C}(R)$ enjoys a left total derived functor on $\mathcal{D}(R)$, together with a balance property, in the sense that $X\otimes_{R}^{\bf L}Y$ can be computed by
$$X\otimes_{R}^{\bf L}Y \simeq P\otimes_{R}Y \simeq X\otimes_{R}Q,$$
where $P\xrightarrow {\simeq} X$ is any semi-projective resolution of $X$, and $Q\xrightarrow {\simeq} Y$ is any semi-projective resolution of $Y$. Besides, these functors turn out to be triangulated. Moreover, we let $$\Tor^{R}_{i}(X,Y):=H_{i}\left(X\otimes_{R}^{\bf L}Y\right)$$
for every $i \in \mathbb{Z}$.
\item[(iii)] The functor $\Gamma_{\mathfrak{a}}(-)$ on $\mathcal{M}(R)$ extends naturally to a functor on $\mathcal{C}(R)$. The extended functor enjoys a right total derived functor ${\bf R}\Gamma_{\mathfrak{a}}(-):\mathcal{D}(R)\rightarrow \mathcal{D}(R)$, that can be computed by
${\bf R}\Gamma_{\mathfrak{a}}(X)\simeq \Gamma_{\mathfrak{a}}(I)$, where $X \xrightarrow {\simeq} I$ is any semi-injective resolution of $X$. Besides, we define the $i$th local cohomology module of $X$ to be
$$H^{i}_{\mathfrak{a}}(X):= H_{-i}\left({\bf R}\Gamma_{\mathfrak{a}}(X)\right)$$
for every $i\in\mathbb{Z}$. The functor ${\bf R}\Gamma_{\mathfrak{a}}(-)$ turns out to be triangulated.
\item[(iv)] The functor $\Lambda^{\mathfrak{a}}(-)$ on $\mathcal{M}(R)$ extends naturally to a functor on $\mathcal{C}(R)$. The extended functor enjoys a left total derived functor ${\bf L}\Lambda^{\mathfrak{a}}(-):\mathcal{D}(R)\rightarrow \mathcal{D}(R)$, that can be computed by ${\bf L}\Lambda^{\mathfrak{a}}(X)\simeq \Lambda^{\mathfrak{a}}(P)$, where $P \xrightarrow {\simeq} X$ is any semi-projective resolution of $X$. Moreover, we define the $i$th local homology module of $X$ to be
$$H^{\mathfrak{a}}_{i}(X):= H_{i}\left({\bf L}\Lambda^{\mathfrak{a}}(X)\right)$$
for every $i\in\mathbb{Z}$. The functor ${\bf L}\Lambda^{\mathfrak{a}}(-)$ turns out to be triangulated.
\end{enumerate}
\end{remark}

We further need the notion of way-out functors for functors between the category of complexes.

\begin{definition} \label{3.3.3}
Let $R$ and $S$ be two rings, and $\mathcal{F}: \mathcal{C}(R) \rightarrow \mathcal{C}(S)$ a covariant functor. Then:
\begin{enumerate}
\item[(i)] $\mathcal{F}$ is said to be \textit{way-out left} \index{way-out left functor on category of complexes} if for every $n \in \mathbb{Z}$, there is an $m \in \mathbb{Z}$, such that for any $R$-complex $X$ with $X_{i}=0$ for every $i>m$, we have $\mathcal{F}(X)_{i}=0$ for every $i>n$.
\item[(ii)] $\mathcal{F}$ is said to be \textit{way-out right} \index{way-out right functor on category of complexes} if for every $n \in \mathbb{Z}$, there is an $m \in \mathbb{Z}$, such that for any $R$-complex $X$ with $X_{i}=0$ for every $i<m$, we have $\mathcal{F}(X)_{i}=0$ for every $i<n$.
\item[(iii)] $\mathcal{F}$ is said to be \textit{way-out} \index{way-out functor on category of complexes} if it is both way-out left and way-out right.
\end{enumerate}
\end{definition}

The following lemma is the Way-out Lemma for functors between the category of complexes. We include a proof since there is no account of this version in the literature.

\begin{lemma} \label{3.3.4}
Let $R$ and $S$ be two rings, and $\mathcal{F},\mathcal{G}: \mathcal{C}(R) \rightarrow \mathcal{C}(S)$ two additive covariant functors that commute with shift and preserve the exactness of degreewise split short exact sequences of $R$-complexes. Let $\sigma : \mathcal{F} \rightarrow \mathcal{G}$ be a natural transformation of functors. Then the following assertions hold:
\begin{enumerate}
\item[(i)] If $X$ is a bounded $R$-complex such that $\sigma^{X_{i}}:\mathcal{F}(X_{i}) \rightarrow \mathcal{G}(X_{i})$ is a quasi-isomorphism for every $i \in \mathbb{Z}$, then $\sigma^{X}: \mathcal{F}(X) \rightarrow \mathcal{G}(X)$ is a quasi-isomorphism.
\item[(ii)] If $\mathcal{F}$ and $\mathcal{G}$ are way-out left, and $X$ is a left-bounded $R$-complex such that $\sigma^{X_{i}}:\mathcal{F}(X_{i}) \rightarrow \mathcal{G}(X_{i})$ is a quasi-isomorphism for every $i \in \mathbb{Z}$, then $\sigma^{X}: \mathcal{F}(X) \rightarrow \mathcal{G}(X)$ is a quasi-isomorphism.
\item[(iii)] If $\mathcal{F}$ and $\mathcal{G}$ are way-out right, and $X$ is a right-bounded $R$-complex such that $\sigma^{X_{i}}:\mathcal{F}(X_{i}) \rightarrow \mathcal{G}(X_{i})$ is a quasi-isomorphism for every $i \in \mathbb{Z}$, then $\sigma^{X}: \mathcal{F}(X) \rightarrow \mathcal{G}(X)$ is a quasi-isomorphism.
\item[(iv)] If $\mathcal{F}$ and $\mathcal{G}$ are way-out, and $X$ is an $R$-complex such that $\sigma^{X_{i}}:\mathcal{F}(X_{i}) \rightarrow \mathcal{G}(X_{i})$ is a quasi-isomorphism for every $i \in \mathbb{Z}$, then $\sigma^{X}: \mathcal{F}(X) \rightarrow \mathcal{G}(X)$ is a quasi-isomorphism.
\end{enumerate}
\end{lemma}

\begin{prf}
(i): Without loss of generality we may assume that
\[X: 0\rightarrow X_{n} \xrightarrow{\partial_{n}^{X}} X_{n-1} \rightarrow \cdots \rightarrow X_{1} \xrightarrow{\partial_{1}^{X}} X_{0} \rightarrow 0.\]
Let
\[Y: 0\rightarrow X_{n-1} \xrightarrow{\partial_{n-1}^{X}} X_{n-2} \rightarrow \cdots \rightarrow X_{1} \xrightarrow{\partial_{1}^{X}} X_{0} \rightarrow 0.\]
Consider the degreewise split short exact sequence
\[0\rightarrow Y \rightarrow X \rightarrow \Sigma^{n}X_{n} \rightarrow 0\]
of $R$-complexes, and apply $\mathcal{F}$ and $\mathcal{G}$ to get the following commutative diagram of $S$-complexes with exact rows:
\[\begin{tikzpicture}[every node/.style={midway},]
  \matrix[column sep={2.5em}, row sep={2.5em}]
  {\node(1) {$0$}; & \node(2) {$\mathcal{F}(Y)$}; & \node(3) {$\mathcal{F}(X)$}; & \node(4) {$\Sigma^{n}\mathcal{F}(X_{n})$}; & \node(5) {$0$};\\
  \node(6) {$0$}; & \node(7) {$\mathcal{G}(Y)$}; & \node(8) {$\mathcal{G}(X)$}; & \node(9) {$\Sigma^{n}\mathcal{G}(X_{n})$}; & \node(10) {$0$};\\};
  \draw[decoration={markings,mark=at position 1 with {\arrow[scale=1.5]{>}}},postaction={decorate},shorten >=0.5pt] (4) -- (9) node[anchor=west] {$\Sigma^{n}\sigma_{X_{n}}$};
  \draw[decoration={markings,mark=at position 1 with {\arrow[scale=1.5]{>}}},postaction={decorate},shorten >=0.5pt] (3) -- (8) node[anchor=west] {$\sigma_{X}$};
  \draw[decoration={markings,mark=at position 1 with {\arrow[scale=1.5]{>}}},postaction={decorate},shorten >=0.5pt] (2) -- (7) node[anchor=west] {$\sigma_{Y}$};
  \draw[decoration={markings,mark=at position 1 with {\arrow[scale=1.5]{>}}},postaction={decorate},shorten >=0.5pt] (1) -- (2) node[anchor=south] {};
  \draw[decoration={markings,mark=at position 1 with {\arrow[scale=1.5]{>}}},postaction={decorate},shorten >=0.5pt] (2) -- (3) node[anchor=south] {};
  \draw[decoration={markings,mark=at position 1 with {\arrow[scale=1.5]{>}}},postaction={decorate},shorten >=0.5pt] (3) -- (4) node[anchor=south] {};
  \draw[decoration={markings,mark=at position 1 with {\arrow[scale=1.5]{>}}},postaction={decorate},shorten >=0.5pt] (4) -- (5) node[anchor=south] {};
  \draw[decoration={markings,mark=at position 1 with {\arrow[scale=1.5]{>}}},postaction={decorate},shorten >=0.5pt] (6) -- (7) node[anchor=south] {};
  \draw[decoration={markings,mark=at position 1 with {\arrow[scale=1.5]{>}}},postaction={decorate},shorten >=0.5pt] (7) -- (8) node[anchor=south] {};
  \draw[decoration={markings,mark=at position 1 with {\arrow[scale=1.5]{>}}},postaction={decorate},shorten >=0.5pt] (8) -- (9) node[anchor=south] {};
  \draw[decoration={markings,mark=at position 1 with {\arrow[scale=1.5]{>}}},postaction={decorate},shorten >=0.5pt] (9) -- (10) node[anchor=south] {};
\end{tikzpicture}\]
Note that $\Sigma^{n}\sigma_{X_{n}}$ is a quasi-isomorphism by the assumption. Hence to prove that $\sigma_{X}$ is a quasi-isomorphism, it suffices to show that $\sigma_{Y}$ is a quasi-isomorphism. Since $Y$ is bounded, by continuing this process with $Y$, we reach at a level that we need $\sigma_{X_{0}}$ to be a quasi-isomorphism, which holds by the assumption. Therefore, we are done.

(ii): Without loss of generality we may assume that
\[X: 0\rightarrow X_{n} \xrightarrow{\partial_{n}^{X}} X_{n-1} \rightarrow \cdots.\]
Let $i \in \mathbb{Z}$. We show that $H_{i}(\sigma_{X}):H_{i}\left(\mathcal{F}(X)\right) \rightarrow H_{i}\left(\mathcal{G}(X)\right)$ is an isomorphism. Since $\mathcal{F}$ and $\mathcal{G}$ are way-out left, we can choose an integer $j \in \mathbb{Z}$ corresponding to $i-2$. Let
\[Z: 0\rightarrow X_{n} \xrightarrow{\partial^{X}_{n}} X_{n-1} \rightarrow \cdots \rightarrow X_{j+1} \xrightarrow{\partial^{X}_{j+1}} X_{j} \rightarrow 0\]
and
\[Y: 0 \rightarrow X_{j-1} \xrightarrow{\partial^{X}_{j-1}} X_{j-2} \rightarrow \cdots.\]
Then there is a degreewise split short exact sequence
\[0 \rightarrow Y \rightarrow X \rightarrow Z \rightarrow 0\]
of $R$-complexes. Apply $\mathcal{F}$ and $\mathcal{G}$ to get the following commutative diagram with exact rows:
\[\begin{tikzpicture}[every node/.style={midway},]
  \matrix[column sep={2.5em}, row sep={2.5em}]
  {\node(1) {$0$}; & \node(2) {$\mathcal{F}(Y)$}; & \node(3) {$\mathcal{F}(X)$}; & \node(4) {$\mathcal{F}(Z)$}; & \node(5) {$0$};\\
  \node(6) {$0$}; & \node(7) {$\mathcal{G}(Y)$}; & \node(8) {$\mathcal{G}(X)$}; & \node(9) {$\mathcal{G}(Z)$}; & \node(10) {$0$};\\};
  \draw[decoration={markings,mark=at position 1 with {\arrow[scale=1.5]{>}}},postaction={decorate},shorten >=0.5pt] (4) -- (9) node[anchor=west] {$\sigma_{Z}$};
  \draw[decoration={markings,mark=at position 1 with {\arrow[scale=1.5]{>}}},postaction={decorate},shorten >=0.5pt] (3) -- (8) node[anchor=west] {$\sigma_{X}$};
  \draw[decoration={markings,mark=at position 1 with {\arrow[scale=1.5]{>}}},postaction={decorate},shorten >=0.5pt] (2) -- (7) node[anchor=west] {$\sigma_{Y}$};
  \draw[decoration={markings,mark=at position 1 with {\arrow[scale=1.5]{>}}},postaction={decorate},shorten >=0.5pt] (1) -- (2) node[anchor=south] {};
  \draw[decoration={markings,mark=at position 1 with {\arrow[scale=1.5]{>}}},postaction={decorate},shorten >=0.5pt] (2) -- (3) node[anchor=south] {};
  \draw[decoration={markings,mark=at position 1 with {\arrow[scale=1.5]{>}}},postaction={decorate},shorten >=0.5pt] (3) -- (4) node[anchor=south] {};
  \draw[decoration={markings,mark=at position 1 with {\arrow[scale=1.5]{>}}},postaction={decorate},shorten >=0.5pt] (4) -- (5) node[anchor=south] {};
  \draw[decoration={markings,mark=at position 1 with {\arrow[scale=1.5]{>}}},postaction={decorate},shorten >=0.5pt] (6) -- (7) node[anchor=south] {};
  \draw[decoration={markings,mark=at position 1 with {\arrow[scale=1.5]{>}}},postaction={decorate},shorten >=0.5pt] (7) -- (8) node[anchor=south] {};
  \draw[decoration={markings,mark=at position 1 with {\arrow[scale=1.5]{>}}},postaction={decorate},shorten >=0.5pt] (8) -- (9) node[anchor=south] {};
  \draw[decoration={markings,mark=at position 1 with {\arrow[scale=1.5]{>}}},postaction={decorate},shorten >=0.5pt] (9) -- (10) node[anchor=south] {};
\end{tikzpicture}\]
From the above diagram, we get the following commutative diagram of $S$-modules with exact rows:
\[\begin{tikzpicture}[every node/.style={midway},]
  \matrix[column sep={2.5em}, row sep={2.5em}]
  {\node(1) {$0=H_{i}\left(\mathcal{F}(Y)\right)$}; & \node(2) {$H_{i}\left(\mathcal{F}(X)\right)$}; & \node(3) {$H_{i}\left(\mathcal{F}(Z)\right)$}; & \node(4) {$H_{i-1}\left(\mathcal{F}(Y)\right)=0$}; \\
  \node(5) {$0=H_{i}\left(\mathcal{G}(Y)\right)$}; & \node(6) {$H_{i}\left(\mathcal{G}(X)\right)$}; & \node(7) {$H_{i}\left(\mathcal{G}(Z)\right)$}; & \node(8) {$H_{i-1}\left(\mathcal{G}(Y)\right)=0$};\\};
  \draw[decoration={markings,mark=at position 1 with {\arrow[scale=1.5]{>}}},postaction={decorate},shorten >=0.5pt] (2) -- (6) node[anchor=west] {$H_{i}(\sigma_{X})$};
  \draw[decoration={markings,mark=at position 1 with {\arrow[scale=1.5]{>}}},postaction={decorate},shorten >=0.5pt] (3) -- (7) node[anchor=west] {$H_{i}(\sigma_{Z})$};
  \draw[decoration={markings,mark=at position 1 with {\arrow[scale=1.5]{>}}},postaction={decorate},shorten >=0.5pt] (1) -- (2) node[anchor=south] {};
  \draw[decoration={markings,mark=at position 1 with {\arrow[scale=1.5]{>}}},postaction={decorate},shorten >=0.5pt] (2) -- (3) node[anchor=south] {};
  \draw[decoration={markings,mark=at position 1 with {\arrow[scale=1.5]{>}}},postaction={decorate},shorten >=0.5pt] (3) -- (4) node[anchor=south] {};
  \draw[decoration={markings,mark=at position 1 with {\arrow[scale=1.5]{>}}},postaction={decorate},shorten >=0.5pt] (5) -- (6) node[anchor=south] {};
  \draw[decoration={markings,mark=at position 1 with {\arrow[scale=1.5]{>}}},postaction={decorate},shorten >=0.5pt] (6) -- (7) node[anchor=south] {};
  \draw[decoration={markings,mark=at position 1 with {\arrow[scale=1.5]{>}}},postaction={decorate},shorten >=0.5pt] (7) -- (8) node[anchor=south] {};
\end{tikzpicture}\]
where the vanishing is due to the choice of $j$. Since $Z$ is bounded, it follows from (i) that $H_{i}(\sigma_{Z})$ is an isomorphism, and as a consequence, $H_{i}(\sigma_{X})$ is an isomorphism.

(iii): Without loss of generality we may assume that
\[X: \cdots \rightarrow X_{n+1} \xrightarrow{\partial_{n+1}^{X}} X_{n} \rightarrow 0.\]
Let $i \in \mathbb{Z}$. We show that $H_{i}(\sigma_{X}):H_{i}\left(\mathcal{F}(X)\right) \rightarrow H_{i}\left(\mathcal{G}(X)\right)$ is an isomorphism. Since $\mathcal{F}$ and $\mathcal{G}$ are way-out right, we can choose an integer $j \in \mathbb{Z}$ corresponding to $i+2$. Let
\[Y: 0\rightarrow X_{j-1} \xrightarrow{\partial^{X}_{j-1}} X_{j-2} \rightarrow \cdots \rightarrow X_{n+1} \xrightarrow{\partial^{X}_{n+1}} X_{n} \rightarrow 0\]
and
\[Z: \cdots \rightarrow X_{j+1} \xrightarrow{\partial^{X}_{j+1}} X_{j} \rightarrow 0.\]
Then there is a degreewise split short exact sequence
\[0 \rightarrow Y \rightarrow X \rightarrow Z \rightarrow 0\]
of $R$-complexes. Apply $\mathcal{F}$ and $\mathcal{G}$ to get the following commutative diagram of $S$-complexes with exact rows:
\[\begin{tikzpicture}[every node/.style={midway},]
  \matrix[column sep={2.5em}, row sep={2.5em}]
  {\node(1) {$0$}; & \node(2) {$\mathcal{F}(Y)$}; & \node(3) {$\mathcal{F}(X)$}; & \node(4) {$\mathcal{F}(Z)$}; & \node(5) {$0$};\\
  \node(6) {$0$}; & \node(7) {$\mathcal{G}(Y)$}; & \node(8) {$\mathcal{G}(X)$}; & \node(9) {$\mathcal{G}(Z)$}; & \node(10) {$0$};\\};
  \draw[decoration={markings,mark=at position 1 with {\arrow[scale=1.5]{>}}},postaction={decorate},shorten >=0.5pt] (4) -- (9) node[anchor=west] {$\sigma_{Z}$};
  \draw[decoration={markings,mark=at position 1 with {\arrow[scale=1.5]{>}}},postaction={decorate},shorten >=0.5pt] (3) -- (8) node[anchor=west] {$\sigma_{X}$};
  \draw[decoration={markings,mark=at position 1 with {\arrow[scale=1.5]{>}}},postaction={decorate},shorten >=0.5pt] (2) -- (7) node[anchor=west] {$\sigma_{Y}$};
  \draw[decoration={markings,mark=at position 1 with {\arrow[scale=1.5]{>}}},postaction={decorate},shorten >=0.5pt] (1) -- (2) node[anchor=south] {};
  \draw[decoration={markings,mark=at position 1 with {\arrow[scale=1.5]{>}}},postaction={decorate},shorten >=0.5pt] (2) -- (3) node[anchor=south] {};
  \draw[decoration={markings,mark=at position 1 with {\arrow[scale=1.5]{>}}},postaction={decorate},shorten >=0.5pt] (3) -- (4) node[anchor=south] {};
  \draw[decoration={markings,mark=at position 1 with {\arrow[scale=1.5]{>}}},postaction={decorate},shorten >=0.5pt] (4) -- (5) node[anchor=south] {};
  \draw[decoration={markings,mark=at position 1 with {\arrow[scale=1.5]{>}}},postaction={decorate},shorten >=0.5pt] (6) -- (7) node[anchor=south] {};
  \draw[decoration={markings,mark=at position 1 with {\arrow[scale=1.5]{>}}},postaction={decorate},shorten >=0.5pt] (7) -- (8) node[anchor=south] {};
  \draw[decoration={markings,mark=at position 1 with {\arrow[scale=1.5]{>}}},postaction={decorate},shorten >=0.5pt] (8) -- (9) node[anchor=south] {};
  \draw[decoration={markings,mark=at position 1 with {\arrow[scale=1.5]{>}}},postaction={decorate},shorten >=0.5pt] (9) -- (10) node[anchor=south] {};
\end{tikzpicture}\]
From the above diagram, we get the following commutative diagram of $S$-modules with exact rows:
\[\begin{tikzpicture}[every node/.style={midway},]
  \matrix[column sep={2.5em}, row sep={2.5em}]
  {\node(1) {$0=H_{i+1}\left(\mathcal{F}(Z)\right)$}; & \node(2) {$H_{i}\left(\mathcal{F}(Y)\right)$}; & \node(3) {$H_{i}\left(\mathcal{F}(X)\right)$}; & \node(4) {$H_{i}\left(\mathcal{F}(Z)\right)=0$}; \\
  \node(5) {$0=H_{i+1}\left(\mathcal{G}(Z)\right)$}; & \node(6) {$H_{i}\left(\mathcal{G}(Y)\right)$}; & \node(7) {$H_{i}\left(\mathcal{G}(X)\right)$}; & \node(8) {$H_{i}\left(\mathcal{G}(Z)\right)=0$};\\};
  \draw[decoration={markings,mark=at position 1 with {\arrow[scale=1.5]{>}}},postaction={decorate},shorten >=0.5pt] (2) -- (6) node[anchor=west] {$H_{i}(\sigma_{Y})$};
  \draw[decoration={markings,mark=at position 1 with {\arrow[scale=1.5]{>}}},postaction={decorate},shorten >=0.5pt] (3) -- (7) node[anchor=west] {$H_{i}(\sigma_{X})$};
  \draw[decoration={markings,mark=at position 1 with {\arrow[scale=1.5]{>}}},postaction={decorate},shorten >=0.5pt] (1) -- (2) node[anchor=south] {};
  \draw[decoration={markings,mark=at position 1 with {\arrow[scale=1.5]{>}}},postaction={decorate},shorten >=0.5pt] (2) -- (3) node[anchor=south] {};
  \draw[decoration={markings,mark=at position 1 with {\arrow[scale=1.5]{>}}},postaction={decorate},shorten >=0.5pt] (3) -- (4) node[anchor=south] {};
  \draw[decoration={markings,mark=at position 1 with {\arrow[scale=1.5]{>}}},postaction={decorate},shorten >=0.5pt] (5) -- (6) node[anchor=south] {};
  \draw[decoration={markings,mark=at position 1 with {\arrow[scale=1.5]{>}}},postaction={decorate},shorten >=0.5pt] (6) -- (7) node[anchor=south] {};
  \draw[decoration={markings,mark=at position 1 with {\arrow[scale=1.5]{>}}},postaction={decorate},shorten >=0.5pt] (7) -- (8) node[anchor=south] {};
\end{tikzpicture}\]
where the vanishing is due to the choice of $j$. Since $Y$ is bounded, it follows from (i) that $H_{i}(\sigma_{Y})$ is an isomorphism, and as a consequence, $H_{i}(\sigma_{X})$ is an isomorphism.

(iv): Let
\[Y: 0 \rightarrow X_{0} \xrightarrow{\partial^{X}_{0}} X_{-1} \rightarrow \cdots\]
and
\[Z: \cdots \rightarrow X_{2} \xrightarrow{\partial^{X}_{2}} X_{1} \rightarrow 0.\]
Then there is a degreewise split short exact sequence
\[0 \rightarrow Y \rightarrow X \rightarrow Z \rightarrow 0\]
of $R$-complexes.
Applying $\mathcal{F}$ and $\mathcal{G}$, we get the following commutative diagram of $S$-complexes with exact rows:
\[\begin{tikzpicture}[every node/.style={midway},]
  \matrix[column sep={2.5em}, row sep={2.5em}]
  {\node(1) {$0$}; & \node(2) {$\mathcal{F}(Y)$}; & \node(3) {$\mathcal{F}(X)$}; & \node(4) {$\mathcal{F}(Z)$}; & \node(5) {$0$};\\
  \node(6) {$0$}; & \node(7) {$\mathcal{G}(Y)$}; & \node(8) {$\mathcal{G}(X)$}; & \node(9) {$\mathcal{G}(Z)$}; & \node(10) {$0$};\\};
  \draw[decoration={markings,mark=at position 1 with {\arrow[scale=1.5]{>}}},postaction={decorate},shorten >=0.5pt] (4) -- (9) node[anchor=west] {$\sigma_{Z}$};
  \draw[decoration={markings,mark=at position 1 with {\arrow[scale=1.5]{>}}},postaction={decorate},shorten >=0.5pt] (3) -- (8) node[anchor=west] {$\sigma_{X}$};
  \draw[decoration={markings,mark=at position 1 with {\arrow[scale=1.5]{>}}},postaction={decorate},shorten >=0.5pt] (2) -- (7) node[anchor=west] {$\sigma_{Y}$};
  \draw[decoration={markings,mark=at position 1 with {\arrow[scale=1.5]{>}}},postaction={decorate},shorten >=0.5pt] (1) -- (2) node[anchor=south] {};
  \draw[decoration={markings,mark=at position 1 with {\arrow[scale=1.5]{>}}},postaction={decorate},shorten >=0.5pt] (2) -- (3) node[anchor=south] {};
  \draw[decoration={markings,mark=at position 1 with {\arrow[scale=1.5]{>}}},postaction={decorate},shorten >=0.5pt] (3) -- (4) node[anchor=south] {};
  \draw[decoration={markings,mark=at position 1 with {\arrow[scale=1.5]{>}}},postaction={decorate},shorten >=0.5pt] (4) -- (5) node[anchor=south] {};
  \draw[decoration={markings,mark=at position 1 with {\arrow[scale=1.5]{>}}},postaction={decorate},shorten >=0.5pt] (6) -- (7) node[anchor=south] {};
  \draw[decoration={markings,mark=at position 1 with {\arrow[scale=1.5]{>}}},postaction={decorate},shorten >=0.5pt] (7) -- (8) node[anchor=south] {};
  \draw[decoration={markings,mark=at position 1 with {\arrow[scale=1.5]{>}}},postaction={decorate},shorten >=0.5pt] (8) -- (9) node[anchor=south] {};
  \draw[decoration={markings,mark=at position 1 with {\arrow[scale=1.5]{>}}},postaction={decorate},shorten >=0.5pt] (9) -- (10) node[anchor=south] {};
\end{tikzpicture}\]
Since $Y$ is left-bounded, $\sigma_{Y}$ is a quasi-isomorphism by (ii), and since $Z$ is right-bounded, $\sigma_{Z}$ is a quasi-isomorphism by (iii). Therefore, $\sigma_{X}$ is a quasi-isomorphism.
\end{prf}

Although $\check{C}_{\infty}(\underline{a})$ is suitable in Proposition \ref{3.2.12}, it is not applicable in the next proposition due to the fact that it is concentrated in degrees $1,0,...,-n$. What we really need here is a semi-projective approximation of $\check{C}(\underline{a})$ of the same length, i.e. concentrated in degrees $0,-1,...,-n$. We proceed as follows.

Given an element $a \in R$, consider the following commutative diagram:
\[\begin{tikzpicture}[every node/.style={midway},]
  \matrix[column sep={3em}, row sep={3em}]
  {\node(1) {$0$}; & \node(2) {$R[X]\oplus R$}; & \node(3) {$R[X]$}; & \node(4) {$0$};\\
  \node(5) {$0$}; & \node(6) {$R$}; & \node(7) {$R_{a}$}; & \node(8) {$0$};\\};
  \draw[decoration={markings,mark=at position 1 with {\arrow[scale=1.5]{>}}},postaction={decorate},shorten >=0.5pt] (1) -- (2) node[anchor=east] {};
  \draw[decoration={markings,mark=at position 1 with {\arrow[scale=1.5]{>}}},postaction={decorate},shorten >=0.5pt] (2) -- (3) node[anchor=south] {$f_{a}$};
  \draw[decoration={markings,mark=at position 1 with {\arrow[scale=1.5]{>}}},postaction={decorate},shorten >=0.5pt] (3) -- (4) node[anchor=east] {};
  \draw[decoration={markings,mark=at position 1 with {\arrow[scale=1.5]{>}}},postaction={decorate},shorten >=0.5pt] (5) -- (6) node[anchor=south] {};
  \draw[decoration={markings,mark=at position 1 with {\arrow[scale=1.5]{>}}},postaction={decorate},shorten >=0.5pt] (6) -- (7) node[anchor=south] {$\lambda_{R}^{a}$};
  \draw[decoration={markings,mark=at position 1 with {\arrow[scale=1.5]{>}}},postaction={decorate},shorten >=0.5pt] (7) -- (8) node[anchor=south] {};
  \draw[decoration={markings,mark=at position 1 with {\arrow[scale=1.5]{>}}},postaction={decorate},shorten >=0.5pt] (2) -- (6) node[anchor=west] {$\pi$};
  \draw[decoration={markings,mark=at position 1 with {\arrow[scale=1.5]{>}}},postaction={decorate},shorten >=0.5pt] (3) -- (7) node[anchor=west] {$g_{a}$};
\end{tikzpicture}\]
in which, $f_{a}\left(p(X),b\right)=(aX-1)p(X)+b$, $\pi\left(p(X),b\right)=b$, $\lambda^{a}_{R}$ is the localization map, and $g_{a}\left(p(X)\right)=\frac{b_{k}}{a^{k}}+\cdots+\frac{b_{1}}{a}+\frac{b_{0}}{1}$ where $p(X)=b_{k}X^{k}+\cdots+b_{1}X+b_{0} \in R[X]$. Let $L^{R}(a)$ denote the $R$-complex in the first row of the diagram above concentrated in degrees $0,-1$. Since the second row is isomorphic to $\check{C}(a)$, it can be seen that the diagram above provides a quasi-isomorphism $L^{R}(a) \xrightarrow{\simeq} \check{C}(a)$. Hence $L^{R}(a) \xrightarrow{\simeq} \check{C}(a)$ is a semi-projective resolution of $\check{C}(a)$. Now for the elements $\underline{a}=a_{1},...,a_{n}\in R$, let
$$L^{R}(\underline{a})=L^{R}(a_{1})\otimes_{R}\cdots\otimes_{R}L^{R}(a_{n}).$$
Then $L^{R}(\underline{a})$ is an $R$-complex of free modules concentrated in degrees $0,-1,...,-n$, and $L^{R}(\underline{a}) \xrightarrow{\simeq} \check{C}(\underline{a})$ is a semi-projective resolution of $\check{C}(\underline{a})$.

The next proposition inspects the relation between derived torsion functor and derived completion functor with \v{C}ech complex, and provides the second crucial step towards the Greenlees-May Duality.

\begin{proposition} \label{3.3.5}
Let $\mathfrak{a}= (a_{1},...,a_{n})$ be an ideal of $R$, $\underline{a}=a_{1},...,a_{n}$, and $X \in \mathcal{D}(R)$. Then there are natural isomorphisms in $\mathcal{D}(R)$:
\begin{enumerate}
\item[(i)] ${\bf R}\Gamma_{\mathfrak{a}}(X) \simeq \check{C}(\underline{a})\otimes_{R}^{\bf L}X \simeq \check{C}_{\infty}(\underline{a})\otimes_{R}^{\bf L}X$.
\item[(ii)] ${\bf L}\Lambda^{\mathfrak{a}}(X) \simeq {\bf R}\Hom_{R}\left(\check{C}(\underline{a}),X\right) \simeq {\bf R}\Hom_{R}\left(\check{C}_{\infty}(\underline{a}),X\right)$.
\end{enumerate}
\end{proposition}

\begin{prf}
(i): Let $X \xrightarrow{\simeq} I$ be a semi-injective resolution of $X$. Then ${\bf R}\Gamma_{\mathfrak{a}}(X) \simeq \Gamma_{\mathfrak{a}}(I)$, and
$$\check{C}(\underline{a})\otimes_{R}^{\bf L}X \simeq \check{C}(\underline{a})\otimes_{R}^{\bf L}I \simeq \check{C}(\underline{a})\otimes_{R}I,$$
since $\check{C}(\underline{a})$ is a semi-flat $R$-complex. Hence it suffices to establish a quasi-isomorphism $\Gamma_{\mathfrak{a}}(I) \rightarrow \check{C}(\underline{a})\otimes_{R}I$.

Let $Y$ be an $R$-complex and $i \in \mathbb{Z}$. Let $\sigma_{i}^{Y}: \Gamma_{\mathfrak{a}}(Y)_{i} \rightarrow \left(\check{C}(\underline{a})\otimes_{R}Y\right)_{i}$ be the composition of the following natural $R$-homomorphisms:
$$\Gamma_{\mathfrak{a}}(Y)_{i} = \Gamma_{\mathfrak{a}}(Y_{i}) \xrightarrow{\cong} H^{0}_{\mathfrak{a}}(Y_{i}) \xrightarrow{\cong} H_{0}\left(\check{C}(\underline{a})\otimes_{R}Y_{i}\right)= \ker \left(\partial_{0}^{\check{C}(\underline{a})}\otimes_{R}Y_{i}\right) \rightarrow \check{C}(\underline{a})_{0}\otimes_{R}Y_{i} $$$$ \rightarrow \bigoplus_{s+t=i}\left(\check{C}(\underline{a})_{s}\otimes_{R}Y_{t}\right) = \left(\check{C}(\underline{a})\otimes_{R}Y\right)_{i}$$
We note that the second isomorphism above comes from Proposition \ref{3.2.12} (i). One can easily see that $\sigma^{Y}=(\sigma_{i}^{Y})_{i \in \mathbb{Z}}:\Gamma_{\mathfrak{a}}(Y) \rightarrow \check{C}(\underline{a})\otimes_{R}Y$ is a natural morphism of $R$-complexes.

Since $I_{i}$ is an injective $R$-module for any $i \in \mathbb{Z}$, using Proposition \ref{3.2.12} (i), we get
$$H_{-j}\left(\check{C}(\underline{a})\otimes_{R}I_{i}\right) \cong H^{j}_{\mathfrak{a}}(I_{i})=0$$
for every $j \geq 1$. It follows that $\sigma^{I_{i}}:\Gamma_{\mathfrak{a}}(I_{i}) \rightarrow \check{C}(\underline{a})\otimes_{R}I_{i}$ is a quasi-isomorphism:
\[\begin{tikzpicture}[every node/.style={midway}]
  \matrix[column sep={3em}, row sep={3em}]
  {\node(1) {$0$}; & \node(2) {$\Gamma_{\mathfrak{a}}(I_{i})$}; & \node(3) {$0$}; & \node(4) {$\cdots$}; & \node(6) {$0$}; & \node(7) {$0$};\\
  \node(8) {$0$}; & \node(9) {$\check{C}(\underline{a})_{0}\otimes_{R}I_{i}$}; & \node(10) {$\check{C}(\underline{a})_{-1}\otimes_{R}I_{i}$}; & \node(11) {$\cdots$}; & \node(13) {$\check{C}(\underline{a})_{-n}\otimes_{R}I_{i}$}; & \node(14) {$0$};\\};
  \draw[decoration={markings,mark=at position 1 with {\arrow[scale=1.5]{>}}},postaction={decorate},shorten >=0.5pt] (2) -- (9) node[anchor=west] {$\sigma^{I_{i}}_{0}$};
  \draw[decoration={markings,mark=at position 1 with {\arrow[scale=1.5]{>}}},postaction={decorate},shorten >=0.5pt] (3) -- (10) node[anchor=west] {};
  \draw[decoration={markings,mark=at position 1 with {\arrow[scale=1.5]{>}}},postaction={decorate},shorten >=0.5pt] (6) -- (13) node[anchor=west] {};
  \draw[decoration={markings,mark=at position 1 with {\arrow[scale=1.5]{>}}},postaction={decorate},shorten >=0.5pt] (1) -- (2) node[anchor=south] {};
  \draw[decoration={markings,mark=at position 1 with {\arrow[scale=1.5]{>}}},postaction={decorate},shorten >=0.5pt] (8) -- (9) node[anchor=south] {};
  \draw[decoration={markings,mark=at position 1 with {\arrow[scale=1.5]{>}}},postaction={decorate},shorten >=0.5pt] (2) -- (3) node[anchor=south] {};
  \draw[decoration={markings,mark=at position 1 with {\arrow[scale=1.5]{>}}},postaction={decorate},shorten >=0.5pt] (9) -- (10) node[anchor=south] {};
  \draw[decoration={markings,mark=at position 1 with {\arrow[scale=1.5]{>}}},postaction={decorate},shorten >=0.5pt] (3) -- (4) node[anchor=south] {};
  \draw[decoration={markings,mark=at position 1 with {\arrow[scale=1.5]{>}}},postaction={decorate},shorten >=0.5pt] (10) -- (11) node[anchor=south] {};
  \draw[decoration={markings,mark=at position 1 with {\arrow[scale=1.5]{>}}},postaction={decorate},shorten >=0.5pt] (4) -- (6) node[anchor=south] {};
  \draw[decoration={markings,mark=at position 1 with {\arrow[scale=1.5]{>}}},postaction={decorate},shorten >=0.5pt] (11) -- (13) node[anchor=south] {};
  \draw[decoration={markings,mark=at position 1 with {\arrow[scale=1.5]{>}}},postaction={decorate},shorten >=0.5pt] (6) -- (7) node[anchor=south] {};
  \draw[decoration={markings,mark=at position 1 with {\arrow[scale=1.5]{>}}},postaction={decorate},shorten >=0.5pt] (13) -- (14) node[anchor=south] {};
\end{tikzpicture}\]
In addition, it is easily seen that the functors $\Gamma_{\mathfrak{a}}(-):\mathcal{C}(R) \rightarrow \mathcal{C}(R)$ and $\check{C}(\underline{a})\otimes_{R}-:\mathcal{C}(R) \rightarrow \mathcal{C}(R)$ are additive way-out functors that commute with shift and preserve the exactness of degreewise split short exact sequences of $R$-complexes. Hence by Lemma \ref{3.3.4} (iv), we conclude that $\sigma^{I}:\Gamma_{\mathfrak{a}}(I) \rightarrow \check{C}(\underline{a})\otimes_{R}I$ is a quasi-isomorphism.

The second isomorphism is immediate since $\check{C}(\underline{a}) \simeq \check{C}_{\infty}(\underline{a})$ and $-\otimes_{R}^{\bf L}X$ is a functor on $\mathcal{D}(R)$.

(ii): We know that $L^{R}(\underline{a}) \simeq \check{C}(\underline{a}) \simeq \check{C}_{\infty}(\underline{a})$. Let $P \xrightarrow{\simeq} X$ be a semi-projective resolution of $X$. Then ${\bf L}\Lambda^{\mathfrak{a}}(X) \simeq \Lambda^{\mathfrak{a}}(P)$, and
$${\bf R}\Hom_{R}\left(\check{C}(\underline{a}),X\right) \simeq {\bf R}\Hom_{R}\left(L^{R}(\underline{a}),P\right) \simeq \Hom_{R}\left(L^{R}(\underline{a}),P\right),$$
since $L^{R}(\underline{a})$ is a semi-projective $R$-complex. Moreover, we have
$$\Hom_{R}\left(L^{R}(\underline{a}),P\right) \simeq {\bf R}\Hom_{R}\left(L^{R}(\underline{a}),P\right) \simeq {\bf R}\Hom_{R}\left(\check{C}_{\infty}(\underline{a}),P\right)\simeq \Hom_{R}\left(\check{C}_{\infty}(\underline{a}),P\right),$$
since $\check{C}_{\infty}(\underline{a})$ is a semi-projective $R$-complex. In particular, we get
\begin{equation} \label{eq:3.3.5.1}
H_{i}\left(\Hom_{R}\left(L^{R}(\underline{a}),P\right)\right) \cong H_{i}\left(\Hom_{R}\left(\check{C}_{\infty}(\underline{a}),P\right)\right)
\end{equation}
for every $i \in \mathbb{Z}$. Now it suffices to establish a natural quasi-isomorphism $\Hom_{R}\left(L^{R}(\underline{a}),P\right) \rightarrow \Lambda^{\mathfrak{a}}(P)$.

Let $Y$ be an $R$-complex and $i \in \mathbb{Z}$. Let $\varsigma^{Y}_{i}: \Hom_{R}\left(L^{R}(\underline{a}),Y\right)_{i} \rightarrow \Lambda^{\mathfrak{a}}(Y)_{i}$ be the composition of the following natural $R$-homomorphisms:
$$\Hom_{R}\left(L^{R}(\underline{a}),Y\right)_{i} = \prod_{s \in \mathbb{Z}}\Hom_{R}\left(L^{R}(\underline{a})_{s},Y_{s+i}\right) \rightarrow \Hom_{R}\left(L^{R}(\underline{a})_{0},Y_{i}\right) $$$$ \rightarrow \frac{\Hom_{R}\left(L^{R}(\underline{a})_{0},Y_{i}\right)}{\im \left(\Hom_{R}\left(\partial_{0}^{L^{R}(\underline{a})},Y_{i}\right)\right)} = H_{0}\left(\Hom_{R}\left(L^{R}(\underline{a}),Y_{i}\right)\right) \xrightarrow{\cong} H_{0}\left(\Hom_{R}\left(\check{C}_{\infty}(\underline{a}),Y_{i}\right)\right) $$$$ \xrightarrow{\cong} H^{\mathfrak{a}}_{0}(Y_{i}) \rightarrow \Lambda^{\mathfrak{a}}(Y_{i}) = \Lambda^{\mathfrak{a}}(Y)_{i}$$
We note that the first isomorphism above comes from the isomorphism \eqref{eq:3.3.5.1} and the second comes from Proposition \ref{3.2.12} (ii). One can easily see that
$\varsigma^{Y}=(\varsigma^{Y}_{i})_{i \in \mathbb{Z}} : \Hom_{R}\left(L^{R}(\underline{a}),Y\right)\rightarrow \Lambda^{\mathfrak{a}}(Y)$ is a natural morphism of $R$-complexes.

Since $P_{i}$ is a projective $R$-module for any $i \in \mathbb{Z}$, using the isomorphism \eqref{eq:3.3.5.1} and Proposition \ref{3.2.12} (ii), we get
$$H_{j}\left(\Hom_{R}\left(L^{R}(\underline{a}),P_{i}\right)\right) \cong H_{j}\left(\Hom_{R}\left(\check{C}_{\infty}(\underline{a}),P_{i}\right)\right) \cong H^{\mathfrak{a}}_{j}(P_{i})=0$$
for every $j \geq 1$. It follows that $\varsigma^{P_{i}}: \Hom_{R}\left(L^{R}(\underline{a}),P_{i}\right) \rightarrow \Lambda^{\mathfrak{a}}(P_{i})$ is a quasi-isomorphism:
\[\begin{tikzpicture}[every node/.style={midway}]
  \matrix[column sep={1em}, row sep={3em}]
  {\node(1) {$0$}; & \node(2) {$\Hom_{R}\left(L^{R}(\underline{a})_{-n},P_{i}\right)$}; & \node(3) {$\cdots$}; & \node(4) {$\Hom_{R}\left(L^{R}(\underline{a})_{-1},P_{i}\right)$}; & \node(6) {$\Hom_{R}\left(L^{R}(\underline{a})_{0},P_{i}\right)$}; & \node(7) {$0$};\\
  \node(8) {$0$}; & \node(9) {$0$}; & \node(10) {$\cdots$}; & \node(11) {$0$}; & \node(13) {$\Lambda^{\mathfrak{a}}(P_{i})$}; & \node(14) {$0$};\\};
  \draw[decoration={markings,mark=at position 1 with {\arrow[scale=1.5]{>}}},postaction={decorate},shorten >=0.5pt] (2) -- (9) node[anchor=west] {};
  \draw[decoration={markings,mark=at position 1 with {\arrow[scale=1.5]{>}}},postaction={decorate},shorten >=0.5pt] (4) -- (11) node[anchor=west] {};
  \draw[decoration={markings,mark=at position 1 with {\arrow[scale=1.5]{>}}},postaction={decorate},shorten >=0.5pt] (6) -- (13) node[anchor=west] {$\varsigma^{P_{i}}_{0}$};
  \draw[decoration={markings,mark=at position 1 with {\arrow[scale=1.5]{>}}},postaction={decorate},shorten >=0.5pt] (1) -- (2) node[anchor=south] {};
  \draw[decoration={markings,mark=at position 1 with {\arrow[scale=1.5]{>}}},postaction={decorate},shorten >=0.5pt] (8) -- (9) node[anchor=south] {};
  \draw[decoration={markings,mark=at position 1 with {\arrow[scale=1.5]{>}}},postaction={decorate},shorten >=0.5pt] (2) -- (3) node[anchor=south] {};
  \draw[decoration={markings,mark=at position 1 with {\arrow[scale=1.5]{>}}},postaction={decorate},shorten >=0.5pt] (9) -- (10) node[anchor=south] {};
  \draw[decoration={markings,mark=at position 1 with {\arrow[scale=1.5]{>}}},postaction={decorate},shorten >=0.5pt] (3) -- (4) node[anchor=south] {};
  \draw[decoration={markings,mark=at position 1 with {\arrow[scale=1.5]{>}}},postaction={decorate},shorten >=0.5pt] (10) -- (11) node[anchor=south] {};
  \draw[decoration={markings,mark=at position 1 with {\arrow[scale=1.5]{>}}},postaction={decorate},shorten >=0.5pt] (4) -- (6) node[anchor=south] {};
  \draw[decoration={markings,mark=at position 1 with {\arrow[scale=1.5]{>}}},postaction={decorate},shorten >=0.5pt] (11) -- (13) node[anchor=south] {};
  \draw[decoration={markings,mark=at position 1 with {\arrow[scale=1.5]{>}}},postaction={decorate},shorten >=0.5pt] (6) -- (7) node[anchor=south] {};
  \draw[decoration={markings,mark=at position 1 with {\arrow[scale=1.5]{>}}},postaction={decorate},shorten >=0.5pt] (13) -- (14) node[anchor=south] {};
\end{tikzpicture}\]
In addition, it is easily seen that the functors $\Hom_{R}\left(L^{R}(\underline{a}),-\right):\mathcal{C}(R) \rightarrow \mathcal{C}(R)$ and $\Lambda^{\mathfrak{a}}(-):\mathcal{C}(R) \rightarrow \mathcal{C}(R)$ are additive way-out functors that commute with shift and preserve the exactness of degreewise split short exact sequences of $R$-complexes. Hence by Lemma \ref{3.3.4} (iv), we conclude that $\varsigma^{P}:\Hom_{R}\left(L^{R}(\underline{a}),P\right) \rightarrow \Lambda^{\mathfrak{a}}(P)$ is a quasi-isomorphism.

The second isomorphism is immediate since $\check{C}(\underline{a}) \simeq \check{C}_{\infty}(\underline{a})$ and ${\bf R}\Hom_{R}(-,X)$ is a functor on $\mathcal{D}(R)$.
\end{prf}

We note that if $\mathfrak{a}= (a_{1},...,a_{n})$ is an ideal of $R$ and $\underline{a}=a_{1},...,a_{n}$, then $\check{C}(\underline{a})$ as an element of $\mathcal{C}(R)$ depends on the generators $\underline{a}$. However, the proof of the next corollary shows that $\check{C}(\underline{a})$ as an element of $\mathcal{D}(R)$ is independent of the generators $\underline{a}$.

\begin{corollary} \label{3.3.6}
Let $\mathfrak{a}$ be an ideal of $R$. Then there are natural isomorphisms in $\mathcal{D}(R)$:
\begin{enumerate}
\item[(i)] ${\bf R}\Gamma_{\mathfrak{a}}(X) \simeq {\bf R}\Gamma_{\mathfrak{a}}(R)\otimes_{R}^{\bf L}X$.
\item[(ii)] ${\bf L}\Lambda^{\mathfrak{a}}(X) \simeq {\bf R}\Hom_{R}\left({\bf R}\Gamma_{\mathfrak{a}}(R),X\right)$.
\end{enumerate}
\end{corollary}

\begin{prf}
Suppose that $\mathfrak{a}= (a_{1},...,a_{n})$, and $\underline{a}=a_{1},...,a_{n}$. By Proposition \ref{3.3.5} (i), we have
$${\bf R}\Gamma_{\mathfrak{a}}(R) \simeq \check{C}(\underline{a})\otimes_{R}^{\bf L}R \simeq \check{C}(\underline{a}).$$
Now (i) and (ii) follow from Proposition \ref{3.3.5}.
\end{prf}

\section{Greenlees-May Duality Theorem}

Having the material developed in Sections 3.2 and 3.3 at our disposal, we are fully prepared to prove the celebrated Greenlees-May Duality Theorem. \index{Greenlees-May duality theorem}

\begin{theorem} \label{3.4.1}
Let $\mathfrak{a}$ be an ideal of $R$, and $X,Y \in \mathcal{D}(R)$. Then there is a natural isomorphism
$${\bf R}\Hom_{R}\left({\bf R}\Gamma_{\mathfrak{a}}(X),Y\right) \simeq {\bf R}\Hom_{R}\left(X,{\bf L}\Lambda^{\mathfrak{a}}(Y)\right)$$
in $\mathcal{D}(R)$.
\end{theorem}

\begin{prf}
Using Corollary \ref{3.3.6} and the Adjointness Isomorphism, we have
\begin{equation*}
\begin{split}
{\bf R}\Hom_{R}\left({\bf R}\Gamma_{\mathfrak{a}}(X),Y\right) & \simeq {\bf R}\Hom_{R}\left({\bf R}\Gamma_{\mathfrak{a}}(R)\otimes_{R}^{\bf L}X,Y\right) \\
 & \simeq {\bf R}\Hom_{R}\left(X,{\bf R}\Hom_{R}\left({\bf R}\Gamma_{\mathfrak{a}}(R),X\right)\right) \\
 & \simeq {\bf R}\Hom_{R}\left(X,{\bf L}\Lambda^{\mathfrak{a}}(Y)\right).
\end{split}
\end{equation*}
\end{prf}

\begin{corollary}
Let $\mathfrak{a}$ be an ideal of $R$, and $X,Y \in \mathcal{D}(R)$. Then there are natural isomorphisms:
\begin{equation*}
\begin{split}
{\bf L}\Lambda^{\mathfrak{a}}\left({\bf R}\Hom_{R}(X,Y)\right) & \simeq {\bf R}\Hom_{R}\left({\bf L}\Lambda^{\mathfrak{a}}(X),{\bf L}\Lambda^{\mathfrak{a}}(Y)\right) \\
 & \simeq {\bf R}\Hom_{R}\left(X,{\bf L}\Lambda^{\mathfrak{a}}(Y)\right) \\
 & \simeq {\bf R}\Hom_{R}\left({\bf R}\Gamma_{\mathfrak{a}}(X),{\bf L}\Lambda^{\mathfrak{a}}(Y)\right) \\
 & \simeq {\bf R}\Hom_{R}\left({\bf R}\Gamma_{\mathfrak{a}}(X),Y\right) \\
 & \simeq {\bf R}\Hom_{R}\left({\bf R}\Gamma_{\mathfrak{a}}(X),{\bf R}\Gamma_{\mathfrak{a}}(Y)\right).
\end{split}
\end{equation*}
\end{corollary}

\begin{prf}
By Corollary \ref{3.3.6}, Adjointness Isomorphism, and Theorem \ref{3.4.1}, we have
\begin{equation} \label{eq:1}
\begin{split}
{\bf L}\Lambda^{\mathfrak{a}}\left({\bf R}\Hom_{R}(X,Y)\right) & \simeq {\bf R}\Hom_{R}\left({\bf R}\Gamma_{\mathfrak{a}}(R),{\bf R}\Hom_{R}(X,Y)\right) \\
 & \simeq {\bf R}\Hom_{R}\left({\bf R}\Gamma_{\mathfrak{a}}(R)\otimes_{R}^{\bf L}X,Y\right) \\
 & \simeq {\bf R}\Hom_{R}\left({\bf R}\Gamma_{\mathfrak{a}}(X),Y\right) \\
 & \simeq {\bf R}\Hom_{R}\left(X,{\bf L}\Lambda^{\mathfrak{a}}(Y)\right).
\end{split}
\end{equation}
Further, by Theorem \ref{3.4.1}, \cite[Corollary on Page 6]{AJL}, and \cite[Proposition 3.2.2]{Li}, we have
\begin{equation} \label{eq:2}
\begin{split}
{\bf R}\Hom_{R}\left({\bf R}\Gamma_{\mathfrak{a}}(X),{\bf L}\Lambda^{\mathfrak{a}}(Y)\right) & \simeq {\bf R}\Hom_{R}\left({\bf R}\Gamma_{\mathfrak{a}}\left({\bf R}\Gamma_{\mathfrak{a}}(X)\right),Y\right) \\
 & \simeq {\bf R}\Hom_{R}\left({\bf R}\Gamma_{\mathfrak{a}}(X),Y\right) \\
 & \simeq {\bf R}\Hom_{R}\left({\bf R}\Gamma_{\mathfrak{a}}(X),{\bf R}\Gamma_{\mathfrak{a}}(Y)\right).
\end{split}
\end{equation}
Moreover, by Theorem \ref{3.4.1} and \cite[Corollary on Page 6]{AJL}, we have
\begin{equation} \label{eq:3}
\begin{split}
{\bf R}\Hom_{R}\left({\bf L}\Lambda^{\mathfrak{a}}(X),{\bf L}\Lambda^{\mathfrak{a}}(Y)\right) & \simeq {\bf R}\Hom_{R}\left({\bf R}\Gamma_{\mathfrak{a}}\left({\bf L}\Lambda^{\mathfrak{a}}(X)\right),Y\right) \\
 & \simeq {\bf R}\Hom_{R}\left({\bf R}\Gamma_{\mathfrak{a}}(X),Y\right).
\end{split}
\end{equation}
Combining the isomorphisms \eqref{eq:1}, \eqref{eq:2}, and \eqref{eq:3}, we get all the desired isomorphisms.
\end{prf}

Now we turn our attention to the Grothendieck's Local Duality, and demonstrate how to derive it from the Greenlees-May Duality.

We need the definition of a dualizing complex.

\begin{definition} \label{3.4.2}
A \textit{dualizing complex} \index{dualizing complex} for $R$ is an $R$-complex $D \in \mathcal{D}^{f}_{\square}(R)$ that satisfies the following conditions:
\begin{enumerate}
\item[(i)] The homothety morphism $\chi^{D}_{R}:R \rightarrow {\bf R}\Hom_{R}(D,D)$ is an isomorphism in $\mathcal{D}(R)$.
\item[(ii)] $\id_{R}(D)<\infty$.
\end{enumerate}
Moreover, if $R$ is local, then a dualizing complex $D$ is said to be \textit{normalized} \index{normalized dualizing complex} if $\sup(D)=\dim(R)$.
\end{definition}

It is clear that if $D$ is a dualizing complex for $R$, then so is $\Sigma^{s}D$ for every $s \in \mathbb{Z}$, which accounts for the non-uniqueness of dualizing complexes. Further, $\Sigma^{\dim(R)-\sup(D)}D$ is a normalized dualizing complex.

\begin{example} \label{3.4.3}
Let $(R,\mathfrak{m},k)$ be a local ring with a normalized dualizing complex $D$. Then ${\bf R}\Gamma_{\mathfrak{m}}(D) \simeq E_{R}(k)$. For a proof, refer to \cite[Proposition 6.1]{Ha2}.
\end{example}

The next theorem determines precisely when a ring enjoys a dualizing complex.

\begin{theorem} \label{3.4.4}
The the following assertions are equivalent:
\begin{enumerate}
\item[(i)] $R$ has a dualizing complex.
\item[(ii)] $R$ is a homomorphic image of a Gorenstein ring of finite Krull dimension.
\end{enumerate}
\end{theorem}

\begin{prf}
See \cite[Page 299]{Ha2} and \cite[Corollary 1.4]{Kw}.
\end{prf}

Now we prove the Local Duality Theorem for complexes. \index{local duality theorem for complexes}

\begin{theorem} \label{3.4.5}
Let $(R,\mathfrak{m})$ be a local ring with a dualizing complex $D$, and $X \in \mathcal{D}^{f}_{\square}(R)$. Then
$$H^{i}_{\mathfrak{m}}(X) \cong \Ext^{\dim(R)-i-\sup(D)}_{R}(X,D)^{\vee}$$
for every $i \in \mathbb{Z}$.
\end{theorem}

\begin{prf}
Clearly, we have
$$\Ext^{\dim(R)-i-\sup(D)}_{R}(X,D) \cong \Ext^{-i}_{R}\left(X,\Sigma^{\dim(R)-\sup(D)}D\right)$$
for every $i \in \mathbb{Z}$, and $\Sigma^{\dim(R)-\sup(D)}D$ is a normalized dualizing for $R$. Hence by replacing $D$ with $\Sigma^{\dim(R)-\sup(D)}D$, it suffices to assume that $D$ is a normalized dualizing complex and prove the isomorphism
$H^{i}_{\mathfrak{m}}(X) \cong \Ext^{-i}_{R}(X,D)^{\vee}$ for every $i \in \mathbb{Z}$.
By Theorem \ref{3.4.1}, we have
\begin{equation} \label{eq:3.4.5.1}
{\bf R}\Hom_{R}\left({\bf R}\Gamma_{\mathfrak{m}}(X),E_{R}(k)\right) \simeq {\bf R}\Hom_{R}\left(X,{\bf L}\Lambda^{\mathfrak{m}}\left(E_{R}(k)\right)\right).
\end{equation}
But since $E_{R}(k)$ is injective, it provides a semi-injective resolution of itself, so we have
\begin{equation} \label{eq:3.4.5.2}
{\bf R}\Hom_{R}\left({\bf R}\Gamma_{\mathfrak{m}}(X),E_{R}(k)\right) \simeq \Hom_{R}\left({\bf R}\Gamma_{\mathfrak{m}}(X),E_{R}(k)\right).
\end{equation}
Besides, by Example \ref{3.4.3}, \cite[Corollary on Page 6]{AJL}, and \cite[Proposition 2.7]{Fr}, we have
\begin{equation} \label{eq:3.4.5.3}
\begin{split}
{\bf L}\Lambda^{\mathfrak{m}}\left(E_{R}(k)\right) & \simeq {\bf L}\Lambda^{\mathfrak{m}}\left({\bf R}\Gamma_{\mathfrak{m}}(D)\right) \\
 & \simeq {\bf L}\Lambda^{\mathfrak{m}}(D) \\
 & \simeq D\otimes_{R}^{\bf L}\widehat{R}^{\mathfrak{m}} \\
 & \simeq D\otimes_{R}\widehat{R}^{\mathfrak{m}}.
\end{split}
\end{equation}
Combining \eqref{eq:3.4.5.1}, \eqref{eq:3.4.5.2}, and \eqref{eq:3.4.5.3}, we get
$$\Hom_{R}\left({\bf R}\Gamma_{\mathfrak{m}}(X),E_{R}(k)\right) \simeq {\bf R}\Hom_{R}\left(X,D\otimes_{R}\widehat{R}^{\mathfrak{m}}\right).$$
Taking Homology, we obtain
\begin{equation} \label{eq:3.4.5.4}
\begin{split}
\Hom_{R}\left(H^{i}_{\mathfrak{m}}(X),E_{R}(k)\right) & \cong \Hom_{R}\left(H_{-i}\left({\bf R}\Gamma_{\mathfrak{m}}(X)\right),E_{R}(k)\right) \\
 & \cong H_{i}\left(\Hom_{R}\left({\bf R}\Gamma_{\mathfrak{m}}(X),E_{R}(k)\right)\right) \\
 & \cong H_{i}\left({\bf R}\Hom_{R}\left(X,D\otimes_{R}\widehat{R}^{\mathfrak{m}}\right)\right) \\
 & \cong \Ext^{-i}_{R}\left(X,D\otimes_{R}\widehat{R}^{\mathfrak{m}}\right) \\
\end{split}
\end{equation}
for every $i \in \mathbb{Z}$.

Since $X \in \mathcal{D}^{f}_{\square}(R)$, we have $X\otimes_{R}\widehat{R}^{\mathfrak{m}} \in \mathcal{D}^{f}_{\square}\left(\widehat{R}^{\mathfrak{m}}\right)$, so $H^{i}_{\mathfrak{m}\widehat{R}^{\mathfrak{m}}}\left(X\otimes_{R}\widehat{R}^{\mathfrak{m}}\right)$ is an artinian $\widehat{R}^{\mathfrak{m}}$-module by \cite[Proposition 2.1]{HD}, and thus Matlis reflexive for every $i \in \mathbb{Z}$. Moreover, $D\otimes_{R}\widehat{R}^{\mathfrak{m}}$ is a normalized dualizing complex for $\widehat{R}^{\mathfrak{m}}$. Therefore, using the isomorphism \eqref{eq:3.4.5.4} over the $\mathfrak{m}$-adically complete ring $\widehat{R}^{\mathfrak{m}}$, we obtain
\begin{equation} \label{eq:3.4.5.5}
\begin{split}
H_{\mathfrak{m}}^{i}(X) & \cong H_{\mathfrak{m}}^{i}(X)\otimes_{R} \widehat{R}^{\mathfrak{m}}\\
 & \cong H^{i}_{\mathfrak{m}\widehat{R}^{\mathfrak{m}}}\left(X\otimes_{R} \widehat{R}^{\mathfrak{m}}\right) \\
 & \cong \Hom_{\widehat{R}^{\mathfrak{m}}}\left(\Hom_{\widehat{R}^{\mathfrak{m}}}\left(H^{i}_{\mathfrak{m}\widehat{R}^{\mathfrak{m}}}\left(X\otimes_{R} \widehat{R}^{\mathfrak{m}}\right),E_{\widehat{R}^{\mathfrak{m}}}(k)\right),E_{\widehat{R}^{\mathfrak{m}}}(k)\right) \\
 & \cong \Hom_{\widehat{R}^{\mathfrak{m}}}\left(\Ext^{-i}_{\widehat{R}^{\mathfrak{m}}}\left(X\otimes_{R}\widehat{R}^{\mathfrak{m}},D\otimes_{R}\widehat{R}^{\mathfrak{m}} \right),E_{\widehat{R}^{\mathfrak{m}}}(k)\right) \\
 & \cong \Hom_{\widehat{R}^{\mathfrak{m}}}\left(\Ext^{-i}_{R}(X,D)\otimes_{R}\widehat{R}^{\mathfrak{m}},E_{\widehat{R}^{\mathfrak{m}}}(k)\right) \\
\end{split}
\end{equation}
for every $i \in \mathbb{Z}$. However, ${\bf R}\Hom_{R}(X,D) \in \mathcal{D}^{f}_{\sqsubset}(R)$, so $\Ext^{-i}_{R}(X,D)$ is a finitely generated $R$-module for every $i \in \mathbb{Z}$. It follows that
\begin{equation} \label{eq:3.4.5.6}
\Hom_{\widehat{R}^{\mathfrak{m}}}\left(\Ext^{-i}_{R}(X,D)\otimes_{R}\widehat{R}^{\mathfrak{m}},E_{\widehat{R}^{\mathfrak{m}}}(k)\right) \cong \Hom_{R}\left(\Ext^{-i}_{R}(X,D),E_{R}(k)\right)
\end{equation}
for every $i \in \mathbb{Z}$.
Combining \eqref{eq:3.4.5.5} and \eqref{eq:3.4.5.6}, we obtain
$$H^{i}_{\mathfrak{m}}(X) \cong \Hom_{R}\left(\Ext^{-i}_{R}(X,D),E_{R}(k)\right)$$
for every $i \in \mathbb{Z}$ as desired.
\end{prf}

Our next goal is to obtain the Local Duality Theorem for modules. But first we need the definition of a dualizing module.

\begin{definition} \label{3.4.6}
Let $(R,\mathfrak{m})$ be a local ring. A \textit{dualizing module} \index{dualizing module} for $R$ is a finitely generated $R$-module $\omega$ that satisfies the following conditions:
\begin{enumerate}
\item[(i)] The homothety map $\chi^{\omega}_{R}:R \rightarrow \Hom_{R}(\omega,\omega)$, given by $\chi^{\omega}_{R}(a)=a1^{\omega}$ for every $a \in R$, is an isomorphism.
\item[(ii)] $\Ext^{i}_{R}(\omega,\omega)=0$ for every $i \geq 1$.
\item[(iii)] $\id_{R}(\omega)<\infty$.
\end{enumerate}
\end{definition}

The next theorem determines precisely when a ring enjoys a dualizing module.

\begin{theorem} \label{3.4.7}
Let $(R,\mathfrak{m})$ be a local ring. Then the following assertions are equivalent:
\begin{enumerate}
\item[(i)] $R$ has a dualizing module.
\item[(ii)] $R$ is a Cohen-Macaulay local ring which is a homomorphic image of a Gorenstein local ring.
\end{enumerate}
Moreover in this case, the dualizing module is unique up to isomorphism.
\end{theorem}

\begin{prf}
See \cite[Corollary 2.2.13]{Wa} and \cite[Theorem 3.3.6]{BH}.
\end{prf}

Since the dualizing module for $R$ is unique whenever it exists, we denote a choice of the dualizing module by $\omega_{R}$.

\begin{proposition} \label{3.4.8}
Let $(R,\mathfrak{m})$ be a Cohen-Macaulay local ring, and $\omega$ a finitely generated $R$-module. Then the following assertions are equivalent:
\begin{enumerate}
\item[(i)] $\omega$ is a dualizing module for $R$.
\item[(ii)] $\omega^{\vee} \cong H^{\dim(R)}_{\mathfrak{m}}(R)$.
\end{enumerate}
\end{proposition}

\begin{prf}
See \cite[Definition 12.1.2, Exercises 12.1.23 and 12.1.25, and Remark 12.1.26]{BS}, and \cite[Definition 3.3.1]{BH}.
\end{prf}

We can now derive the Local Duality Theorem for modules. \index{local duality theorem for modules}

\begin{theorem} \label{3.4.9}
Let $(R,\mathfrak{m})$ be a local ring with a dualizing module $\omega_{R}$, and $M$ a finitely generated $R$-module. Then
$$H^{i}_{\mathfrak{m}}(M) \cong \Ext^{\dim(R)-i}_{R}(M,\omega_{R})^{\vee}$$
for every $i \geq 0$.
\end{theorem}

\begin{prf}
By Theorem \ref{3.4.7}, $R$ is a Cohen-Macaulay local ring which is a homomorphic image of a Gorenstein local ring $S$. Since $S$ is local, we have $\dim(S)<\infty$. Hence Theorem \ref{3.4.4} implies that $R$ has a dualizing complex $D$. Since $R$ is Cohen-Macaulay, we have $H^{i}_{\mathfrak{m}}(R)=0$ for every $i \neq \dim(R)$. On the other hand, by Theorem \ref{3.4.5}, we have
\begin{equation} \label{eq:3.4.9.1}
\begin{split}
H^{i}_{\mathfrak{m}}(R) & \cong \Ext^{\dim(R)-i-\sup(D)}_{R}(R,D)^{\vee} \\
 & \cong H_{-\dim(R)+i+\sup(D)}\left({\bf R}\Hom_{R}(R,D)\right)^{\vee} \\
 & \cong H_{-\dim(R)+i+\sup(D)}(D)^{\vee}.
\end{split}
\end{equation}
It follows from the display \eqref{eq:3.4.9.1} that $H_{-\dim(R)+i+\sup(D)}(D)=0$ for every $i \neq \dim(R)$, i.e. $H_{i}(D)=0$ for every $i \neq \sup(D)$. Therefore, we have $D \simeq \Sigma^{\sup(D)}H_{\sup(D)}(D)$. In addition, letting $i=\dim(R)$ in the display \eqref{eq:3.4.9.1}, we get $H^{\dim(R)}_{\mathfrak{m}}(R) \cong H_{\sup(D)}(D)^{\vee}$, which implies that $\omega_{R} \cong H_{\sup(D)}(D)$ by Proposition \ref{3.4.8}. It follows that $D \simeq \Sigma^{\sup(D)}\omega_{R}$.

Now let $M$ be a finitely generated $R$-module. Then by Theorem \ref{3.4.5}, we have
\begin{equation*}
\begin{split}
H^{i}_{\mathfrak{m}}(M) & \cong \Ext^{\dim(R)-i-\sup(D)}_{R}(M,D)^{\vee} \\
 & \cong H_{-\dim(R)+i+\sup(D)}\left({\bf R}\Hom_{R}\left(M,D\right)\right)^{\vee} \\
 & \cong H_{-\dim(R)+i+\sup(D)}\left({\bf R}\Hom_{R}\left(M,\Sigma^{\sup(D)}\omega_{R}\right)\right)^{\vee} \\
 & \cong H_{-\dim(R)+i}\left({\bf R}\Hom_{R}\left(M,\omega_{R}\right)\right)^{\vee} \\
 & \cong \Ext^{\dim(R)-i}_{R}\left(M,\omega_{R}\right)^{\vee}. \\
\end{split}
\end{equation*}
\end{prf}

%%%%%%%%%%%%%%%%%%%%%%%%%%%%%%%%%%%%%%%%%%%%%%%%%%%%%%%%%%%%%%%%%%%%%%%%%%%%%%%%%%%%%%%%%%%%%%%%%%%%%%%%%%%%%%%%%%%%%%%%%%%%%%%%%%%%%%%%%%%%%%%%%%%%%%%%%%%%%%%%%%

\printindex

\end{document}